# FUZZY AND NEUTROSOPHIC ANALYSIS OF PERIYAR'S VIEWS ON UNTOUCHABILITY

**W. B. Vasantha Kandasamy**
**Florentin Smarandache**
**K. Kandasamy**

Translation of the speeches and writings
of Periyar from Tamil by
**Meena Kandasamy**

**November 2005**

# FUZZY AND NEUTROSOPHIC ANALYSIS OF PERIYAR'S VIEWS ON UNTOUCHABILITY

**W. B. Vasantha Kandasamy**
e-mail: **vasantha@iitm.ac.in**
web: **http://mat.iitm.ac.in/~wbv**

**Florentin Smarandache**
e-mail: **smarand@unm.edu**

**K. Kandasamy**
e-mail: **dr.k.kandasamy@gamil.com**

Translation of the speeches and writings
of Periyar from Tamil by
**Meena Kandasamy**

**November 2005**

*Dedicated to Periyar*

# CONTENTS

**Preface** 5

**Chapter One**
**BASIC NOTION OF FCMs, FRMs, NCMs AND NRMS**

**Chapter Two**
**UNTOUCHABILITY: PERIYAR'S VIEW AND PRESENT DAY SITUATION**
**A FUZZY AND NEUTROSOPHIC ANALYSIS**

# Preface

> "Day in and day out we take pride in claiming that India has a 5000 year old civilization. But the way Dalits and those suppressed are being treated by the people who wield power and authority speaks volumes for the degradations of our moral structure and civilized standards."
>
> ***Ex-President of India, the late K. R. Narayanan***
> The New Indian Express, Saturday, 12 Nov. 2005

K.R.Narayanan was a lauded hero and a distinguished victim of his Dalit background. Even in an international platform when he was on an official visit to Paris, the media headlines blazed, 'An Untouchable at Elysee'. He was visibly upset and it proved that a Dalit who rose up to such heights was never spared from the pangs of outcaste-ness and untouchability, which is based on birth. Thus, if the erstwhile first citizen of India faces such humiliation, what will be the plight of the last man who is a Dalit?

As one of the world's largest socio-economically oppressed, culturally subjugated and politically marginalized group of people, the 138 million Dalits in India suffer not only from the excesses of the traditional oppressor castes, but also from State Oppression— which includes, but is not limited to, authoritarianism, police brutality, economic embargo, criminalization of activists, electoral violence, repressive laws that aim to curb fundamental rights, and the non-implementation of laws that safeguard Dalit rights. The Dalits were considered untouchable for thousands of years by the Hindu society until the Constitution of India officially abolished the practice of untouchability in 1950.

The rigid rules of the caste system were codified in the *Manusmriti,* a Hindu religious text which decreed that Sudras and Dalits were to be distressed for subsistence. Max Weber, perhaps the most outstanding comparative sociologist of all time, clearly defined the social identity of the Hindu in terms of caste, and of caste in terms of ritual. "Caste, that is, the ritual rights and duties it gives and imposes, and the position of the Brahmins, is the fundamental institution of Hinduism. Before everything else, without caste there is no Hindu."

In this book we have analyzed the diverse manifestations of untouchability and the caste-system of Hinduism using Fuzzy and Neutrosophic theory.

In studying untouchability, we have chosen to view it through the eyes of the relentless crusader, Periyar E. V. Ramasamy (1879-1973). He was a social revolutionary who vociferously campaigned against untouchability and called for the annihilation of the caste system. For six decades, he powerfully influenced the course of politics in the South-Indian state of Tamil Nadu. His ideology and struggle to attain rationalism, caste-annihilation, self-respect and women liberation have made a permanent impact on the entire nation. Dr.K.Veeramani, who joined Periyar's movement at 11 years of age in 1944, is today the President of the Dravidar Kazhagam. He works vigorously to spread the revolutionary principles of his mentor and he has established numerous institutions as permanent memorials to Periyar.

The practices of untouchability and the caste system have a devastating effect on the nation's progress. That is why, in order to understand the multi-dimensional facets of untouchability and its consequences, we have given relevant excerpts from the translation of Periyar's writings and speeches.

This book is organized into four chapters.

In Chapter One we just introduce the basic Fuzzy and Neutrosophic tools used in the analysis of the social evil of Untouchability. Since the notion of caste is based on the mind, it is appropriate to use Fuzzy and Neutrosophic theory.

In Chapter Two we use the opinion of several experts to analyze the various aspects of untouchability. Here, we use the tools that have been described in Chapter One. Fuzzy Directed Graphs and Neutrosophic Directed Graphs of these Fuzzy and Neutrosophic models happen to be very dense. We have represented 16 such graphs in this book.

In Chapter Three we give a brief introduction about the life and struggle of Periyar. We have provided more than 220 pages of translations of his writings and speeches that dealt with the issues of caste and untouchability. In the final section of this chapter, we have provided a compilation of the various atrocities undergone by Dalits in present-day India.

The fourth and concluding chapter gives observations drawn from our mathematical results based on the Fuzzy and Neutrosophic analysis.

A list for further reading is also provided.

A study of this kind is being carried out for the first time. We have used mathematical tools to study this sociological problem because we wanted to perform a scientific study of this system of oppression.

In this book the terms Panchamas, Depressed Classes, Untouchables, Adi Dravidars, Paraiayar, Pallar, Chakkiliyar refers to the Dalits. We have retained the original words in the translations.

We feel that it was our lifetime duty to render this service to the great leader Periyar who devoted his entire life to working for the rights of the subjugated people.

W.B.VASANTHA KANDASAMY
FLORENTIN SMARANDACHE
K.KANDASAMY
26-11-2005

Chapter One

# BASIC NOTIONS OF FCMS, FRMS, NCMS AND NRMS

The main motivation of this chapter is to make the book a self contained one. In this chapter we just recall the basic concepts of the fuzzy tools and the neutrosophic tools used in this book for analyzing the vision of Periyar on untouchability. This chapter has six sections. In section one we just recall the definition of Fuzzy Cognitive Maps (FCMs) and illustrate them with examples. Section two gives the properties and some models of FCM. Section three describes the notion of Fuzzy Relational Maps (FRMs). In section four we introduce the concepts neutrosophy. Section five gives the definition of Neutrosophic Cognitive Maps (NCMs) and some of its properties. In the final section we give the definition of Neutrosophic Relational Maps (NRMs) and illustrate it by some models.

### 1.1 Definition of Fuzzy Cognitive Maps

In this section we recall the notion of Fuzzy Cognitive Maps (FCMs), which was introduced by Bart Kosko [26] in the year 1986. We also give several of its interrelated definitions. FCMs have a major role to play mainly when the data under consideration is an unsupervised one. Further this method is most simple and an effective one as it can analyse the data by directed graphs and connection matrices and give the hidden pattern of the system. As we are more interested to attract all types of reader we have not used very high mathematics.

**DEFINITION 1.1.1:** *An FCM is a directed graph with concepts like policies, events etc. as nodes and causalities*

*as edges. It represents causal relationship between concepts.*

***Example 1.1.1:*** In Tamil Nadu (a southern state in India) in the last decade several new engineering colleges have been approved and started. The resultant increase in the production of engineering graduates in these years is disproportionate with the need of engineering graduates. This has resulted in thousands of unemployed and underemployed graduate engineers. Using an expert's opinion we study the effect of such unemployed engineering graduates in the society. An expert spells out the five major concepts relating to the unemployed graduate engineers as

$E_1$ – Frustration
$E_2$ – Unemployment
$E_3$ – Increase of educated criminals
$E_4$ – Under employment
$E_5$ – Taking up drugs etc.

The directed graph where $E_1, \ldots, E_5$ are taken as the nodes and causalities as edges as given by an expert is illustrated in the following Figure 1.1.1:

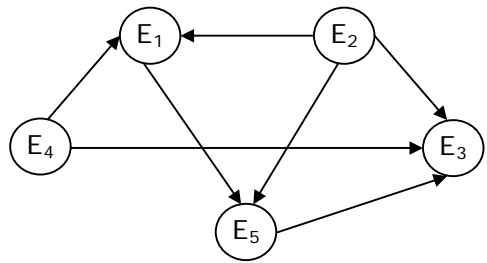

FIGURE: 1.1.1

According to this expert, increase in unemployment increases frustration. Increase in unemployment, increases the educated criminals. Frustration increases the graduates to take up to evils like alcohol drugs etc. Unemployment also leads to the increase in number of persons who take up to drugs, drinks etc. to forget their worries and unoccupied time. Under-employment forces them to do criminal acts like theft (leading to murder) for want of more money and so on. Thus one cannot actually get data for this but can use the expert's opinion for this unsupervised data to obtain some idea about the real plight of the situation. This is just an illustration to show how FCM is described by a directed graph. {If increase (or decrease) in one concept/node leads to increase (or decrease) in another concept/node, then we give the value 1. If there exists no relation between two concepts the value 0 is given. If increase (or decrease) in one concept decreases (or increases) another, then we give the value –1. Thus FCMs are described in this way.

**DEFINITION 1.1.2:** *When the nodes of the FCM are fuzzy sets then they are called as fuzzy nodes.*

**DEFINITION 1.1.3:** *FCMs with edge weights or causalities from the set {–1, 0, 1} are called simple FCMs.*

**DEFINITION 1.1.4:** *Consider the nodes / concepts $C_1, …, C_n$ of the FCM. Suppose the directed graph is drawn using edge weight $e_{ij} \in \{0, 1, –1\}$. The matrix E be defined by $E = (e_{ij})$ where $e_{ij}$ is the weight of the directed edge $C_i C_j$. E is called the adjacency matrix of the FCM, also known as the connection matrix of the FCM.*

It is important to note that all matrices associated with an FCM are always square matrices with diagonal entries as zeros.

**DEFINITION 1.1.5:** *Let $C_1, C_2, … , C_n$ be the nodes of an FCM. $A = (a_1, a_2, … , a_n)$ where $a_i \in \{0, 1\}$. A is called the instantaneous state vector and it denotes the ON - OFF position of the node at an instant.*
  *$a_i = 0$ if $a_i$ is OFF and $a_i = 1$ if $a_i$ is ON for $i = 1, 2, …, n$.*

**DEFINITION 1.1.6:** *Let $C_1, C_2, … , C_n$ be the nodes of an FCM. Let $\overrightarrow{C_1C_2}, \overrightarrow{C_2C_3}, \overrightarrow{C_3C_4}, … , \overrightarrow{C_iC_j}$ be the edges of the FCM ($i \neq j$). Then the edges form a directed cycle. An*

*FCM is said to be cyclic if it possesses a directed cycle. An FCM is said to be acyclic if it does not possess any directed cycle.*

**DEFINITION 1.1.7:** *An FCM with cycles is said to have a feedback.*

**DEFINITION 1.1.8:** *When there is a feedback in an FCM, i.e., when the causal relations flow through a cycle in a revolutionary way, the FCM is called a dynamical system.*

**DEFINITION 1.1.9:** *Let $\overrightarrow{C_1C_2}$, $\overrightarrow{C_2C_3}$, ... , $\overrightarrow{C_{n-1}C_n}$ be a cycle. When $C_i$ is switched on and if the causality flows through the edges of a cycle and if it again causes $C_i$, we say that the dynamical system goes round and round. This is true for any node $C_i$, for i = 1, 2, ... , n. The equilibrium state for this dynamical system is called the hidden pattern.*

**DEFINITION 1.1.10:** *If the equilibrium state of a dynamical system is a unique state vector, then it is called a fixed point.*

*Example 1.1.2:* Consider a FCM with $C_1$, $C_2$, ..., $C_n$ as nodes. For example let us start the dynamical system by switching on $C_1$. Let us assume that the FCM settles down with $C_1$ and $C_n$ on i.e. the state vector remains as (1, 0, 0, ..., 0, 1) this state vector (1, 0, 0, ..., 0, 1) is called the fixed point.

**DEFINITION 1.1.11:** *If the FCM settles down with a state vector repeating in the form $A_1 \rightarrow A_2 \rightarrow ... \rightarrow A_i \rightarrow A_1$ then this equilibrium is called a limit cycle.*

Methods of finding the hidden pattern are discussed in the section 1.2.

**DEFINITION 1.1.12:** *Finite number of FCMs can be combined together to produce the joint effect of all the FCMs. Let $E_1$, $E_2$, ..., $E_p$ be the adjacency matrices of the*

*FCMs with nodes $C_1, C_2, \ldots, C_n$ then the combined FCM is got by adding all the adjacency matrices $E_1, E_2, \ldots, E_p$.*

*We denote the combined FCM adjacency matrix by $E = E_1 + E_2 + \ldots + E_p$.*

*Notation*: Suppose $A = (a_1, \ldots, a_n)$ is a vector which is passed into a dynamical system E. Then $AE = (a'_1, \ldots, a'_n)$ after thresholding and updating the vector suppose we get $(b_1, \ldots, b_n)$ we denote that by $(a'_1, a'_2, \ldots, a'_n) \hookrightarrow (b_1, b_2, \ldots, b_n)$. Thus the symbol '$\hookrightarrow$' means the resultant vector has been thresholded and updated.

FCMs have several advantages as well as some disadvantages. The main advantage of this method is, it is simple. It functions on expert's opinion. When the data happens to be an unsupervised one the FCM comes handy. This is the only known fuzzy technique that gives the hidden pattern of the situation. As we have a very well known theory, which states that the strength of the data depends on, the number of experts' opinion we can use combined FCMs with several experts' opinions.

At the same time the disadvantage of the combined FCM is, when the weightages are 1 and –1 for the same $C_i$ $C_j$, given by 2 experts, we have the sum adding to zero thus at all times the connection matrices $E_1, \ldots, E_k$ may not be conformable for addition. Combined conflicting opinions tend to cancel out and assisted by the strong law of large numbers, a consensus emerges as the sample opinion approximates the underlying population opinion. This problem will be easily overcome if the FCM entries are only 0 and 1. We have just briefly recalled the definitions. For more about FCMs please refer Kosko [26, 27].

## 1.2 Fuzzy Cognitive Maps – Properties and Models

In this section we just give some basic properties of FCMs and give illustration of the model. Fuzzy cognitive maps (FCMs) are more applicable when the data in the first place is an unsupervised one. The FCMs work on the opinion of experts. FCMs model the world as a collection of classes and causal relations between classes. FCMs are fuzzy signed directed graphs with feedback. The directed edge $e_{ij}$ from causal concept $C_i$ to concept $C_j$ measures how much $C_i$ causes $C_j$. The time varying concept function $C_i(t)$ measures the non negative occurrence of some fuzzy event, perhaps the strength of a political sentiment, historical trend or military objective. FCMs are used to model several types of problems varying from gastric-appetite behavior, popular political developments etc. FCMs are also used to model in robotics like plant control. The edges $e_{ij}$ take values in the fuzzy causal interval $[-1, 1]$. $e_{ij} = 0$ indicates no causality, $e_{ij} > 0$ indicates causal increase, that is $C_j$ increases as $C_i$ increases (or $C_j$ decreases as $C_i$ decreases). $e_{ij} < 0$ indicates causal decrease or negative causality. $C_j$ decreases as $C_i$ increases (or $C_j$ increases as $C_i$ decreases). Simple FCMs have edge values in $\{-1, 0, 1\}$. Then if causality occurs, it occurs to a maximal positive or negative degree. Simple FCMs provide a quick first approximation to an expert stand or printed causal knowledge. We illustrate this by the following, which gives a simple FCM of a Socio-economic model. A Socio-economic model is constructed with Population, Crime, Economic condition, Poverty and Unemployment as nodes or concepts. Here the simple trivalent directed graph is given by the following Figure 1.2.1, which is the experts opinion.

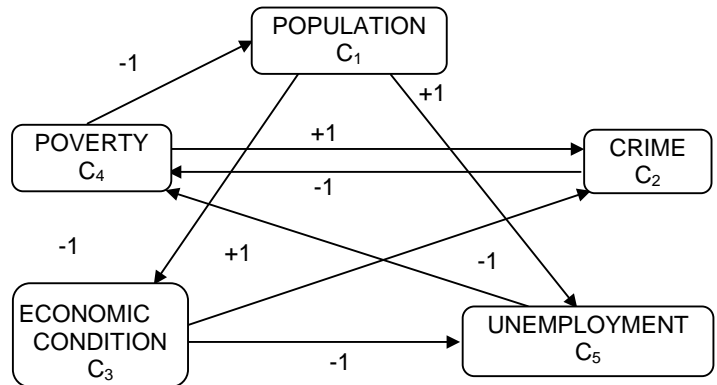

FIGURE: 1.2.1

Causal feedback loops abound in FCMs in thick tangles. Feedback precludes the graph-search techniques used in artificial-intelligence expert systems.

FCMs feedback allows experts to freely draw causal pictures of their problems and allows causal adaptation laws, infer causal links from simple data. FCM feedback forces us to abandon graph search, forward and especially backward chaining. Instead we view the FCM as a dynamical system and take its equilibrium behavior as a forward-evolved inference. Synchronous FCMs behave as Temporal Associative Memories (TAM). We can always, in case of a model, add two or more FCMs to produce a new FCM. The strong law of large numbers ensures in some sense that, knowledge reliability increases with expert sample size.

We reason with FCMs. We pass state vectors C repeatedly through the FCM connection matrix E, thresholding or non-linearly transforming the result after each pass. Independent of the FCMs size, it quickly settles down to a temporal associative memory limit cycle or fixed point which is the hidden pattern of the system for that state vector C. The limit cycle or fixed-point inference summarizes the joint effects of all the interacting fuzzy knowledge.

***Example 1.2.1:*** Consider the $5 \times 5$ causal connection matrix E that represents the socio economic model using FCM which is given in figure in Figure 1.2.1.

$$E = \begin{bmatrix} 0 & 0 & -1 & 0 & 1 \\ 0 & 0 & 0 & -1 & 0 \\ 0 & -1 & 0 & 0 & -1 \\ -1 & 1 & 0 & 0 & 0 \\ 0 & 0 & 0 & 1 & 0 \end{bmatrix}$$

Concepts/nodes can represent processes, events, values or policies. Consider the first node $C_1 = 1$. We hold or clamp $C_1$ on the temporal associative memories, recall process.

Threshold signal functions synchronously update each concept after each pass, through the connection matrix E. We start with the population $C_1 = (1\ 0\ 0\ 0\ 0)$. The arrow indicates the threshold operation,

$$
\begin{aligned}
C_1 E &= (0\ 0\ -1\ 0\ 1) &\hookrightarrow& (1\ 0\ 0\ 0\ 1) &= C_2 \\
C_2 E &= (0\ 0\ -1\ 1\ 1) &\hookrightarrow& (1\ 0\ 0\ 1\ 1) &= C_3 \\
C_3 E &= (-1\ 1\ -1\ 1\ 1) &\hookrightarrow& (1\ 1\ 0\ 1\ 1) &= C_4 \\
C_4 E &= (-1\ 1\ -1\ 0\ 1) &\hookrightarrow& (1\ 1\ 0\ 0\ 1) &= C_5 \\
C_5 E &= (0\ 0\ -1\ 0\ 1) &\hookrightarrow& (1\ 0\ 0\ 0\ 1) &= C_6 = C_2.
\end{aligned}
$$

So the increase in population results in the unemployment problem, which is a limit cycle. For more about FCM refer Kosko [26, 27] and for more about this socio economic model refer [74, 75].

This example illustrates the strengths and weaknesses of FCM analysis. FCM allows experts to represent factual and evaluative concepts in an interactive framework. Experts can quickly draw FCM pictures or respond to questionnaires. Experts can consent or dissent to the local causal structure and perhaps the global equilibrium. The FCM knowledge representation and inferencing structure reduces to simple vector-matrix operations, favours integrated circuit implementation and allows extension to neural statistical or dynamical system techniques. Yet an FCM equally encodes the experts' knowledge or ignorance, wisdom or prejudice. Worse, different experts differ in how they assign causal strengths to edges and in which concepts they deem causally relevant. The FCM seems merely to encode its designers' biases and may not even encode them accurately.

FCM combination provides a partial solution to this problem. We can additively superimpose each experts FCM in associative memory fashion, even though the FCM connection matrices $E_1, \ldots, E_K$ may not be conformable for addition. Combined conflicting opinions tend to cancel out and assisted with the strong law of large numbers, a consensus emerges as the sample opinion approximates the underlying population opinion. FCM combination allows knowledge researchers to construct FCMs with iterative interviews or questionnaire mailings.

The law of large numbers require that the random samples be independent, identically distributed random variables with finite variance. Independence models each experts individually. Identical distribution models a particular domain focus.

We combine arbitrary FCM connection matrices $F_1$, $F_2$, …, $F_K$ by adding augmented FCM matrices, $F_1$, …, $F_K$. Each augmented matrix $F_i$ has n-rows and n-columns n equals the total number of distinct concepts used by the experts. We permute the rows and columns of the augmented matrices to bring them into mutual coincidence. Then we add the $F_i$'s point wise to yield the combined FCM matrix F.

$$F = \sum_i F_i$$

We can then use F to construct the combined FCM directed graph.

Even if each expert gives trivalent description in $\{-1, 0, 1\}$, the combined (and normalized) FCM entry $f_{ij}$ tends to be in $\{-1, 1\}$. (Here $F = (f_{ij})$) The strong law of large numbers ensures that $f_{ij}$ provides a rational approximation to the underlying unknown population opinion of how much $C_i$ affects $C_j$. We can normalize $f_{ij}$ by the number K of experts. Experts tend to give trivalent evaluations more readily and more accurately than they give weighted evaluations. When transcribing interviews or documents, a knowledge engineer can more reliably determine an edge's sign than its magnitude.

Some experts may be more credible than others. We can weight each expert with non-negative credibility weight, weighing the augmented FCM matrix.

$$F = \Sigma\ w_i\ F_i.$$

The weights need not be in [0, 1]; the only condition is they should be non-negative. Different weights may produce different equilibrium, limit cycles or fixed points as hidden patterns. We can also weigh separately any submatrix of each experts augmented FCM matrix.

Augmented FCM matrices imply that every expert causally discusses every concept $C_1, \ldots, C_n$. If an expert does not include $C_j$ in his FCM model the expert implicitly say that $C_j$ is not causally relevant. So the $j^{th}$ row and the $j^{th}$ column of his augmented connection matrix contains only zeros.

The only drawback which we felt while adopting FCM to several of the models is that we do not have a means to say or express if the relation between two causal concepts $C_i$ and $C_j$ is an indeterminate. So we in section five will adopt in FCM the concept of indeterminacy and rename the Fuzzy Cognitive Maps as Neutrosophic Cognitive Maps as Neutrosophy enables one to accept the Truth, Falsehood and the Indeterminate. Such situations very often occur when we deal with unsupervised data that has more to do with feelings. Like political scenario, child labor, child's education, parent-children model, symptom-disease model, personality-medicine model in case of Homeopathy medicines, crime and punishment, in judicial problems (where evidences may be an indeterminate), problems of the aged, female infanticide problem and so on and so forth, which will be discussed in section 1.4 and 1.5. Now we proceed on to define a special type of FCM namely Fuzzy Relational Maps.

Here we mention some examples in which certain factors are indeterminate but where we have used FCM to find solutions. So that in section five we will use these for FCM model, Neutrosophic Cognitive Maps to study the systems.

## 1.3 Fuzzy Relational Maps

In this section, we introduce the notion of Fuzzy Relational Maps (FRMs); they are constructed analogous to FCMs described and discussed in the earlier sections. In FCMs we promote the correlations between causal associations among concurrently active units. But in FRMs we divide the very causal associations into two disjoint units, for example, the relation between a teacher and a student or relation between an employee or employer or a relation between doctor and patient and so on. Thus for us to define a FRM we need a domain space and a range space which are disjoint in the sense of concepts. We further assume no intermediate relation exists within the domain elements or node and the range spaces elements. The number of elements in the range space need not in general be equal to the number of elements in the domain space.

Thus throughout this section we assume the elements of the domain space are taken from the real vector space of dimension n and that of the range space are real vectors from the vector space of dimension m (m in general need not be equal to n). We denote by R the set of nodes $R_1,\ldots, R_m$ of the range space, where $R = \{(x_1,\ldots, x_m) \mid x_j = 0 \text{ or } 1\}$ for $j = 1, 2,\ldots, m$. If $x_i = 1$ it means that the node $R_i$ is in the ON state and if $x_i = 0$ it means that the node $R_i$ is in the OFF state. Similarly D denotes the nodes $D_1, D_2,\ldots, D_n$ of the domain space where $D = \{(x_1,\ldots, x_n) \mid x_j = 0 \text{ or } 1\}$ for $i = 1, 2,\ldots, n$. If $x_i = 1$ it means that the node $D_i$ is in the ON state and if $x_i = 0$ it means that the node $D_i$ is in the OFF state.

Now we proceed on to define a FRM.

**DEFINITION 1.3.1:** *A FRM is a directed graph or a map from D to R with concepts like policies or events etc, as nodes and causalities as edges. It represents causal relations between spaces D and R.*

*Let $D_i$ and $R_j$ denote the two nodes of an FRM. The directed edge from $D_i$ to $R_j$ denotes the causality of $D_i$ on $R_j$ called relations. Every edge in the FRM is weighted with a number in the set $\{0, \pm 1\}$. Let $e_{ij}$ be the weight of the edge $D_iR_j$, $e_{ij} \in \{0, \pm 1\}$. The weight of the edge $D_i R_j$ is positive, if increase in $D_i$ implies increase in $R_j$ or decrease in $D_i$ implies decrease in $R_j$ ie causality of $D_i$ on $R_j$ is 1. If $e_{ij} = 0$, then $D_i$ does not have any effect on $R_j$. We do not discuss the cases when increase in $D_i$ implies decrease in $R_j$ or decrease in $D_i$ implies increase in $R_j$.*

**DEFINITION 1.3.2:** *When the nodes of the FRM are fuzzy sets then they are called fuzzy nodes. FRMs with edge weights {0, ±1} are called simple FRMs.*

**DEFINITION 1.3.3:** *Let $D_1, ..., D_n$ be the nodes of the domain space D of an FRM and $R_1, ..., R_m$ be the nodes of the range space R of an FRM. Let the matrix E be defined as $E = (e_{ij})$ where $e_{ij}$ is the weight of the directed edge $D_iR_j$ (or $R_jD_i$), E is called the relational matrix of the FRM.*

*Note:* It is pertinent to mention here that unlike the FCMs the FRMs can be a rectangular matrix with rows corresponding to the domain space and columns corresponding to the range space. This is one of the marked difference between FRMs and FCMs.

**DEFINITION 1.3.4:** *Let $D_1, ..., D_n$ and $R_1, ..., R_m$ denote the nodes of the FRM. Let $A = (a_1,...,a_n)$, $a_i \in \{0, \pm 1\}$. A is called the instantaneous state vector of the domain space and it denotes the on-off position of the nodes at any instant. Similarly let $B = (b_1,..., b_m)$ $b_i \in \{0, \pm 1\}$. B is called instantaneous state vector of the range space and it denotes the on-off position of the nodes at any instant, $a_i = 0$ if $a_i$ is off and $a_i = 1$ if $a_i$ is on for i= 1, 2,..., n Similarly, $b_i = 0$ if $b_i$ is off and $b_i = 1$ if $b_i$ is on, for i= 1, 2,..., m.*

**DEFINITION 1.3.5:** *Let $D_1, ..., D_n$ and $R_1,..., R_m$ be the nodes of an FRM. Let $D_iR_j$ (or $R_j D_i$) be the edges of an FRM, j = 1, 2,..., m and i= 1, 2,..., n. Let the edges form a directed cycle. An FRM is said to be a cycle if it posses a directed cycle. An FRM is said to be acyclic if it does not posses any directed cycle.*

**DEFINITION 1.3.6:** *An FRM with cycles is said to be an FRM with feedback.*

**DEFINITION 1.3.7:** *When there is a feedback in the FRM, i.e. when the causal relations flow through a cycle in a*

*revolutionary manner, the FRM is called a dynamical system.*

**DEFINITION 1.3.8:** *Let FRM be as described in definition 1.3.5 be taken. Let $D_i R_j$ (or $R_j D_i$), $1 \leq j \leq m$, $1 \leq i \leq n$. When $R_i$ (or $D_j$) is switched on and if causality flows through edges of the cycle and if it again causes $R_i$ (or $D_j$), we say that the dynamical system goes round and round. This is true for any node $R_j$ (or $D_i$) for $1 \leq i \leq n$, (or $1 \leq j \leq m$). The equilibrium state of this dynamical system is called the hidden pattern.*

**DEFINITION 1.3.9:** *If the equilibrium state of a dynamical system is a unique state vector, then it is called a fixed point. Consider an FRM with $R_1, R_2,..., R_m$ and $D_1, D_2,..., D_n$ as nodes. For example, let us start the dynamical system by switching on $R_1$ (or $D_1$). Let us assume that the FRM settles down with $R_1$ and $R_m$ (or $D_1$ and $D_n$) on, i.e. the state vector remains as (1, 0, …, 0, 1) in R (or 1, 0, 0, … , 0, 1) in D), This state vector is called the fixed point.*

**DEFINITION 1.3.10:** *If the FRM settles down with a state vector repeating in the form*
$$A_1 \to A_2 \to A_3 \to ... \to A_i \to A_1$$
$$(or\ B_1 \to B_2 \to ... \to B_i \to B_1)$$
*then this equilibrium is called a limit cycle.*

Let $R_1, R_2,…, R_m$ and $D_1, D_2,…, D_n$ be the nodes of a FRM with feedback. Let E be the relational matrix. Let us find a hidden pattern when $D_1$ is switched on i.e. when an input is given as vector $A_1 = (1, 0, …, 0)$ in $D_1$, the data should pass through the relational matrix E. This is done by multiplying $A_1$ with the relational matrix E. Let $A_1E = (r_1, r_2,…, r_m)$, after thresholding and updating the resultant vector we get $A_1 E \in R$. Now let $B = A_1E$ we pass on B into $E^T$ and obtain $BE^T$. We update and threshold the vector $BE^T$ so that $BE^T \in D$. This procedure is repeated till we get a limit cycle or a fixed point.

**DEFINITION 1.3.11:** *Finite number of FRMs can be combined together to produce the joint effect of all the FRMs. Let $E_1, ..., E_p$ be the relational matrices of the FRMs with nodes $R_1, R_2, ..., R_m$ and $D_1, D_2, ..., D_n$, then the combined FRM is represented by the relational matrix $E = E_1 + ... + E_p$.*

## 1.4 An Introduction to Neutrosophy and some Neutrosophic algebraic structures

In this section we introduce the notion of neutrosophic logic created by Florentine Smarandache [55-59], which is an extension / combination of the fuzzy logic in which the indeterminacy is included. It has become very essential that the notion of neutrosophic logic play a vital role in several of the real world problems like law, medicine, industry, finance, IT, stocks and share etc. Use of neutrosophic notions will be illustrated/ applied in the later sections of this chapter. Fuzzy theory only measures the grade of membership or the non-existence of a membership in the revolutionary way but fuzzy theory has failed to attribute the concept when the relations between notions or nodes or concepts in problems are indeterminate. In fact one can say the inclusion of the concept of indeterminate situation with fuzzy concepts will form the neutrosophic logic.

As in this book the concept of only fuzzy cognitive maps are dealt which mainly deals with the relation / non-relation between two nodes or concepts but it fails to deal the relation between two conceptual nodes when the relation is an indeterminate one. Neutrosophic logic is the only tool known to us, which deals with the notions of indeterminacy, and here we give a brief description of it. Throughout this book *I* will denote indeterminacy of a relation / concept / node. For more about Neutrosophic logic please refer Smarandache [55-59].

**DEFINITION 1.4.1:** *In the neutrosophic logic every logical variable x is described by an ordered triple x = (T, I, F) where T is the degree of truth, F is the degree of false and I the level of indeterminacy.*

(A). To maintain consistency with the classical and fuzzy logics and with probability there is the special case where T + I + F = 1.
(B). But to refer to intuitionistic logic, which means incomplete information on a variable proposition or event one has T+ I+F <1.
(C). Analogically referring to Paraconsistent logic, which means contradictory sources of information about a same logical variable, proposition or event one has    T + I + F > 1.
Thus the advantage of using Neutrosophic logic is that this logic distinguishes between relative truth that is a truth is one or a few worlds only noted by 1 and absolute truth denoted by $1^+$. Likewise neutrosophic logic distinguishes between relative falsehood, noted by 0 and absolute falsehood noted by $^-0$.

It has several applications. A few such illustrations are given in [55-59] which is as follows:

*Example 1.4.1:* From a pool of refugees, waiting in a political refugee camp in Turkey to get the American visa, a% have the chance to be accepted – where a varies in the set A, r% to be rejected – where r varies in the set R, and p% to be in pending (not yet decided) – where p varies in P.
   Say, for example, that the chance of someone Popescu in the pool to emigrate to USA is (between) 40-60% (considering different criteria of emigration one gets different percentages, we have to take care of all of them), the chance of being rejected is 20-25% or 30-35%, and the chance of being in pending is 10% or 20% or 30%. Then the neutrosophic probability that Popescu emigrates to the Unites States is NP (Popescu) = ((40 – 60) (20 – 25) ∪ (30 – 35), {10, 20, 30}), closer to the life.

This is a better approach than the classical probability, where 40 P(Popescu) 60, because from the pending chance – which will be converted to acceptance or rejection – Popescu might get extra percentage in his will to emigrating and also the superior limit of the subsets sum

$$60 + 35 + 30 > 100$$

and in other cases one may have the inferior sum < 0, while in the classical fuzzy set theory the superior sum should be 100 and the inferior sum μ 0. In a similar way, we could say about the element Popescu that Popescu ((40-60), (20-25) ∪ (30-35), {10, 20, 30}) belongs to the set of accepted refugees.

*Example 1.4.2:* The probability that candidate C will win an election is say 25-30% true (percent of people voting for him), 35% false (percent of people voting against him), and 40% or 41% indeterminate (percent of people not coming to the ballot box, or giving a blank vote – not selecting any one or giving a negative vote cutting all candidate on the list). Dialectic and dualism don't work in this case anymore.

*Example 1.4.3:* Another example, the probability that tomorrow it will rain is say 50-54% true according to meteorologists who have investigated the past years weather, 30 or 34-35% false according to today's very sunny and droughty summer, and 10 or 20% undecided (indeterminate).

*Example 1.4.4:* The probability that Yankees will win tomorrow versus Cowboys is 60% true (according to their confrontation's history giving Yankees' satisfaction), 30-32% false (supposing Cowboys are actually up to the mark, while Yankees are declining), and 10 or 11 or 12% indeterminate (left to the hazard: sickness of players, referee's mistakes, atmospheric conditions during the game). These parameters act on players' psychology.

As in this book we use mainly the notion of neutrosophic logic with regard to the indeterminacy of any relation in cognitive maps we are restraining ourselves from dealing with several interesting concepts about neutrosophic logic. As FCMs deals with unsupervised data and the existence or non-existence of cognitive relation, we do not in this book just recall the notion of neutrosophic concepts.

In this book we assume all fields to be real fields of characteristic 0 all vector spaces are taken as real spaces over reals and we denote the indeterminacy by '$I$' as i will make a confusion as it denotes the imaginary value, viz $i^2 = -1$ that is $\sqrt{-1}$ = i. The indeterminacy $I$ is such that $I \cdot I = I^2 = I$.

**DEFINITION 1.4.2:** *Let K be the field of reals. We call the field generated by $K \cup I$ to be the neutrosophic field for it involves the indeterminacy factor in it. We define $I^2 = I$, $I + I = 2I$ i.e., $I + ... + I = nI$, and if $k \in K$ then $k.I = kI$, $0I = 0$. We denote the neutrosophic field by K(I) which is generated by $K \cup I$ that is $K(I) = \langle K \cup I \rangle$. ($\langle K \cup I \rangle$ denotes the field generated by K and I.*

*Example 1.4.5:* Let R be the field of reals. The neutrosophic field of reals is generated by $\langle R \cup I \rangle$ i.e. R(*I*) clearly R $\subset \langle R \cup I \rangle$.

*Example 1.4.6:* Let Q be the field of rationals. The neutrosophic field of rational is generated by $Q \cup I$ denoted by Q(*I*).

**DEFINITION 1.4.3:** *Let K(I) be a neutrosophic field we say K(I) is a prime neutrosophic field if K(I) has no proper subfield which is a neutrosophic field.*

*Example 1.4.7:* Q(*I*) is a prime neutrosophic field where as R(*I*) is not a prime neutrosophic field for Q(*I*) $\subset$ R (*I*).

It is very important to note that all neutrosophic fields used in this book are of characteristic zero. Likewise we can define neutrosophic subfield.

**DEFINITION 1.4.4:** *Let K(I) be a neutrosophic field, P ⊂ K(I) is a neutrosophic subfield of P if P itself is a neutrosophic field. K(I) will also be called as the extension neutrosophic field of the neutrosophic field P.*

**DEFINITION 1.4.5:** *Let $M_{nxm} = \{(a_{ij}) \,/\, a_{ij} \in K(I)\}$, where K (I), is a neutrosophic field. We call $M_{nxm}$ to be the neutrosophic matrix.*

*Example 1.4.8:* Let $Q(I) = \langle Q \cup I \rangle$ be the neutrosophic field.

$$M_{4x3} = \begin{pmatrix} 0 & 1 & I \\ -2 & 4I & 0 \\ 1 & -I & 2 \\ 3I & 1 & 0 \end{pmatrix}$$

is the neutrosophic matrix, with entries from rationals and the indeterminacy *I*.

However we just state, suppose in a legal issue the jury or the judge cannot always prove the evidence in a case, in several places we may not be able to derive any conclusions from the existing facts because of which we cannot make a conclusion that no relation exists or otherwise. But existing relation is an indeterminate. So in the case when the concept of indeterminacy exists the judgment ought to be very carefully analyzed be it a civil case or a criminal case. FCMs are deployed only when the existence or non-existence is dealt with but however in our Neutrosophic Cognitive Maps we will deal with the notion of indeterminacy of the evidence also. Thus legal side has lot of Neutrosophic (NCM) applications. Also we will show how NCMs can be used to study factors as varied as stock markets, medical diagnosis, etc.

## 1.5 Neutrosophic Cognitive Maps (NCM)

The notion of Fuzzy Cognitive Maps (FCMs) which are fuzzy signed directed graphs with feedback are discussed and described in section 1.1. The directed edge $e_{ij}$ from causal concept $C_i$ to concept $C_j$ measures how much $C_i$ causes $C_j$. The time varying concept function $C_i(t)$ measures the non negative occurrence of some fuzzy event, perhaps the strength of a political sentiment, historical trend or opinion about some topics like child labor or school dropouts etc. FCMs model the world as a collection of classes and causal relations between them.

The edge $e_{ij}$ takes values in the fuzzy causal interval [–1, 1] ($e_{ij} = 0$ indicates no causality, $e_{ij} > 0$ indicates causal increase; that $C_j$ increases as $C_i$ increases and $C_j$ decreases as $C_i$ decreases, $e_{ij} < 0$ indicates causal decrease or negative causality $C_j$ decreases as $C_i$ increases or $C_j$, increases as $C_i$ decreases. Simple FCMs have edge value in {-1, 0, 1}. Thus if causality occurs it occurs to maximal positive or negative degree.

It is important to note that $e_{ij}$ measures only absence or presence of influence of the node $C_i$ on $C_j$ but till now any researcher has not contemplated the indeterminacy of any relation between two nodes $C_i$ and $C_j$. When we deal with unsupervised data, there are situations when no relation can be determined between some two nodes. So in this section we try to introduce the indeterminacy in FCMs, and we choose to call this generalized structure as Neutrosophic Cognitive Maps (NCMs). In our view this will certainly give a more appropriate result and also caution the user about the risk of indeterminacy.

Now we proceed on to define the concepts about NCMs.

**DEFINITION 1.5.1:** *A Neutrosophic Cognitive Map (NCM) is a neutrosophic directed graph with concepts like policies, events etc. as nodes and causalities or indeterminates as edges. It represents the causal relationship between concepts.*

Let $C_1, C_2, \ldots, C_n$ denote n nodes, further we assume each node is a neutrosophic vector from neutrosophic vector space V. So a node $C_i$ will be represented by $(x_1, \ldots, x_n)$ where $x_k$'s are zero or one or $I$ ($I$ is the indeterminate introduced in Section 1.4 of chapter one) and $x_k = 1$ means that the node $C_k$ is in the on state and $x_k = 0$ means the node is in the off state and $x_k = I$ means the nodes state is an indeterminate at that time or in that situation.

Let $C_i$ and $C_j$ denote the two nodes of the NCM. The directed edge from $C_i$ to $C_j$ denotes the causality of $C_i$ on $C_j$ called connections. Every edge in the NCM is weighted with a number in the set $\{-1, 0, 1, I\}$. Let $e_{ij}$ be the weight of the directed edge $C_iC_j$, $e_{ij} \in \{-1, 0, 1, I\}$. $e_{ij} = 0$ if $C_i$ does not have any effect on $C_j$, $e_{ij} = 1$ if increase (or decrease) in $C_i$ causes increase (or decreases) in $C_j$, $e_{ij} = -1$ if increase (or decrease) in $C_i$ causes decrease (or increase) in $C_j$. $e_{ij} = I$ if the relation or effect of $C_i$ on $C_j$ is an indeterminate.

**DEFINITION 1.5.2:** *NCMs with edge weight from {-1, 0, 1, I} are called simple NCMs.*

**DEFINITION 1.5.3:** *Let $C_1, C_2, \ldots, C_n$ be nodes of a NCM. Let the neutrosophic matrix N(E) be defined as $N(E) = (e_{ij})$, where $e_{ij}$ is the weight of the directed edge $C_i C_j$, and $e_{ij} \in \{0, 1, -1, I\}$. N(E) is called the neutrosophic adjacency matrix of the NCM.*

**DEFINITION 1.5.4:** *Let $C_1, C_2, \ldots, C_n$ be the nodes of the NCM. Let $A = (a_1, a_2, \ldots, a_n)$ where $a_i \in \{0, 1, I\}$. A is called the instantaneous state neutrosophic vector and it denotes the on – off – indeterminate state position of the node at an instant*

    $a_i$ = *0 if $a_i$ is off (no effect)*
    $a_i$ = *1 if $a_i$ is on (has effect)*
    $a_i$ = *I if $a_i$ is indeterminate(effect cannot be determined)*
*for i = 1, 2,…, n.*

**DEFINITION 1.5.5:** *Let $C_1, C_2, ..., C_n$ be the nodes of the FCM. Let $\overrightarrow{C_1C_2}, \overrightarrow{C_2C_3}, \overrightarrow{C_3C_4}, ..., \overrightarrow{C_iC_j}$ be the edges of the NCM. Then the edges form a directed cycle. An NCM is said to be cyclic if it possesses a directed cyclic. An NCM is said to be acyclic if it does not possess any directed cycle.*

**DEFINITION 1.5.6:** *An NCM with cycles is said to have a feedback. When there is a feedback in the NCM i.e. when the causal relations flow through a cycle in a revolutionary manner the NCM is called a dynamical system.*

**DEFINITION 1.5.7:** *Let $\overrightarrow{C_1C_2}, \overrightarrow{C_2C_3}, \cdots, \overrightarrow{C_{n-1}C_n}$ be cycle, when $C_i$ is switched on and if the causality flows through the edges of a cycle and if it again causes $C_i$, we say that the dynamical system goes round and round. This is true for any node $C_i$, for $i = 1, 2,..., n$; the equilibrium state for this dynamical system is called the hidden pattern.*

**DEFINITION 1.5.8:** *If the equilibrium state of a dynamical system is a unique state vector, then it is called a fixed point. Consider the NCM with $C_1, C_2,..., C_n$ as nodes. For example let us start the dynamical system by switching on $C_1$. Let us assume that the NCM settles down with $C_1$ and $C_n$ on, i.e. the state vector remain as (1, 0,..., 1) this neutrosophic state vector (1,0,..., 0, 1) is called the fixed point.*

**DEFINITION 1.5.9:** *If the NCM settles with a neutrosophic state vector repeating in the form*
$$A_1 \to A_2 \to ... \to A_i \to A_1,$$
*then this equilibrium state is called a limit cycle of the NCM.*

Let $C_1, C_2,..., C_n$ be the nodes of an NCM, with feedback. Let E be the associated adjacency matrix. Let us find the hidden pattern when $C_1$ is switched on when an input is given as the vector $A_1 = (1, 0, 0,..., 0)$, the data should pass through the neutrosophic matrix N(E), this is done by multiplying $A_1$ by the matrix $N(E)$. Let $A_1N(E) = (a_1, a_2,..., a_n)$ with the threshold operation that is by replacing $a_i$ by 1 if $a_i \geq k$ and $a_i$ by 0 if $a_i < k$ ($k$ – a suitable positive integer) and $a_i$ by $I$ if $a_i$ is not a integer. We update the resulting concept, the concept $C_1$ is included in the updated vector by making the first coordinate as 1 in the resulting vector. Suppose $A_1N(E) \to A_2$ then consider $A_2N(E)$ and repeat the same procedure. This procedure is repeated till we get a limit cycle or a fixed point.

**DEFINITION 1.5.10:** *Finite number of NCMs can be combined together to produce the joint effect of all NCMs. If $N(E_1), N(E_2),..., N(E_p)$ be the neutrosophic adjacency matrices of a NCM with nodes $C_1, C_2,..., C_n$ then the combined NCM is got by adding all the neutrosophic adjacency matrices $N(E_1),..., N(E_p)$. We denote the combined NCMs adjacency neutrosophic matrix by $N(E) = N(E_1) + N(E_2) +...+ N(E_p)$.*

**NOTATION:** Let $(a_1, a_2, ..., a_n)$ and $(a'_1, a'_2, ..., a'_n)$ be two neutrosophic vectors. We say $(a_1, a_2, ..., a_n)$ is equivalent to $(a'_1, a'_2, ..., a'_n)$ denoted by $((a_1, a_2, ..., a_n) \sim (a'_1, a'_2, ..., a'_n)$ if $(a'_1, a'_2, ..., a'_n)$ is got after thresholding and updating the vector $(a_1, a_2, ..., a_n)$ after passing through the neutrosophic adjacency matrix $N(E)$.

The following are very important:

*Note 1:* The nodes $C_1, C_2, ..., C_n$ are nodes and not indeterminate nodes because they indicate the concepts which are well known. But the edges connecting $C_i$ and $C_j$ may be indeterminate i.e. an expert may not be in a position to say that $C_i$ has some causality on $C_j$ either will he be in a position to state that $C_i$ has no relation with $C_j$ in such cases the relation between $C_i$ and $C_j$ which is indeterminate is denoted by *I*.

*Note 2:* The nodes when sent will have only ones and zeros i.e. on and off states, but after the state vector passes through the neutrosophic adjacency matrix the resultant vector will have entries from {0, 1, I} i.e. they can be neutrosophic vectors.

The presence of I in any of the coordinate implies the expert cannot say the presence of that node i.e. on state of it after passing through N(E) nor can we say the absence of the node i.e. off state of it the effect on the node after passing through the dynamical system is indeterminate so only it is represented by *I*. Thus only in case of NCMs we can say the effect of any node on other nodes can also be indeterminates. Such possibilities and analysis is totally absent in the case of FCMs.

*Note 3:* In the neutrosophic matrix N(E), the presence of *I* in the $a_{ij}^{th}$ place shows, that the causality between the two nodes i.e. the effect of $C_i$ on $C_j$ is indeterminate. Such chances of being indeterminate is very possible in case of unsupervised data and that too in the study of FCMs which are derived from the directed graphs.

Thus only NCMs helps in such analysis.

Now we shall represent a few examples to show how in this set up NCMs is preferred to FCMs. At the outset before we proceed to give examples we make it clear that all unsupervised data need not have NCMs to be applied to it. Only data which have the relation between two nodes to be an indeterminate need to be modeled with NCMs if the data has no indeterminacy factor between any pair of nodes one need not go for NCMs; FCMs will do the best job.

### 1.6 Neutrosophic Relational Maps — Definition with Examples

When the nodes or concepts under study happens to be such that they can be divided into two disjoint classes and a study or analysis can be made using Fuzzy Relational Maps (FRMs) which was introduced and described in the earlier section. Here we define a new concept called Neutrosophic Relational Maps (NRMs), analyse and study them. We also give examples of them.

**DEFINITION 1.6.1:** *Let D be the domain space and R be the range space with $D_1,\ldots, D_n$ the conceptual nodes of the domain space D and $R_1,\ldots, R_m$ be the conceptual nodes of the range space R such that they form a disjoint class i.e. $D \cap R = \phi$. Suppose there is a FRM relating D and R and if at least a edge relating a $D_i\ R_j$ is an indeterminate then we call the FRMs as the Neutrosophic relational maps. i.e. NRMs.*

*Note:* In everyday occurrences we see that if we are studying a model built using an unsupervised data we need not always have some edge relating the nodes of a domain space and a range space or there does not exist any relation between two nodes, it can very well happen that for any two nodes one may not be always in a position to say that the existence or nonexistence of a relation, but we may say that the relation between two nodes is an indeterminate or cannot be decided.

Thus to the best of our knowledge indeterminacy models can be built using neutrosophy. One model already discussed is the Neutrosophic Cognitive Model. The other being the Neutrosophic Relational Maps model, which are a further generalization of Fuzzy Relational Maps model.

It is not essential when a study/ prediction/ investigation is made we are always in a position to find a complete answer. This is not always possible (sometimes or many a times) it is almost all models built using unsupervised data, we may have the factor of indeterminacy to play a role. Such study is possible only by using the Neutrosophic logic.

*Example 1.6.1:* Female infanticide (the practice of killing female children at birth or shortly thereafter) is prevalent in India from the early vedic times, as women were (and still are) considered as a property. As long as a woman is treated as a property/ object the practice of female infanticide will continue in India.

In India, social factors play a major role in female infanticide. Even when the government recognized the girl child as a critical issue for the country's development, India continues to have an adverse ratio of women to men. Other reasons being torture of the in-laws may also result in cruel death of a girl child. This is mainly due to the fact that men are considered superior to women. Also they take into account the fact that men are breadwinners for the family. Even if women work like men, parents think that her efforts is going to end once she is married and enters a new family.

Studies have consistently shown that girl babies in India are denied the same and equal food and medical care that the boy babies receive. Girl babies die more often than boy babies even though medical research has long ago established that girls are generally biologically stronger as new-borns than boys. The birth of a male child is a time for celebration, but the birth of female child is often viewed as a crisis. Thus the female infanticide cannot be attributed to single reason it is highly dependent on the feeling of individuals ranging from social stigma, monetary waste, social status etc.

Suppose we take the conceptual nodes for the unsupervised data relating to the study of female infanticide. We take the status of the people as the domain space D

$D_1$ – very rich
$D_2$ – rich
$D_3$ – upper middle class
$D_4$ – middle class
$D_5$ – lower middle class
$D_6$ – poor
$D_7$ – very poor.

The nodes of the range space R are taken as

$R_1$ – Number of female children - a problem
$R_2$ – Social stigma of having female children
$R_3$ – Torture by in-laws for having only female children
$R_4$ – Economic loss / burden due to female children

$R_5$ – Insecurity due to having only female children (They will marry and enter different homes thereby leaving their parents, so no one would be able to take care of them in later days.)

Keeping these as nodes of the range space and the domain space experts opinion were drawn which is given by the following Figure 1.6.1:

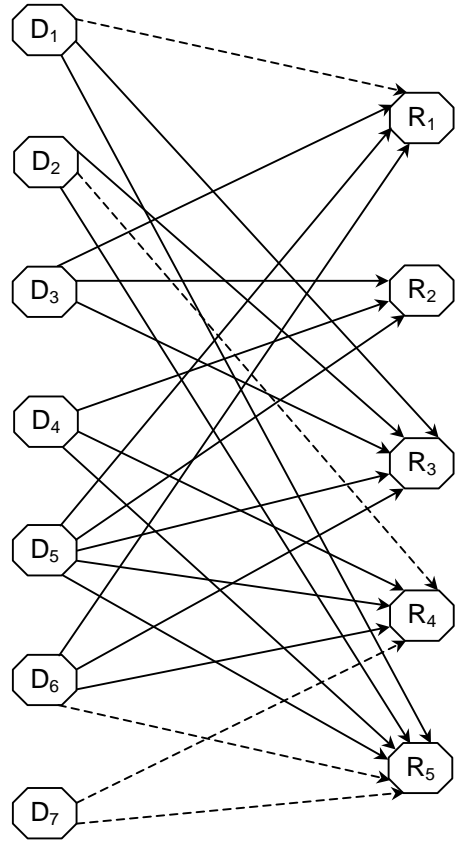

**FIGURE: 1.6.1**

Figure 1.6.1 is the neutrosophic directed graph of the NRM where the dotted lines show the indeterminacy relation.

The corresponding neutrosophic relational matrix $N(R)^T$ is given below:

$$N(R)^T = \begin{bmatrix} I & 0 & 1 & 0 & 1 & 1 & 0 \\ 0 & 0 & 1 & 1 & 1 & 0 & 0 \\ 1 & 1 & 1 & 1 & 1 & 1 & 0 \\ 0 & I & 0 & 0 & 1 & 1 & I \\ 1 & 1 & 0 & 1 & 1 & I & I \end{bmatrix} \text{ and}$$

$$N(R) = \begin{bmatrix} I & 0 & 1 & 0 & 1 \\ 0 & 0 & 1 & I & 1 \\ 1 & 1 & 1 & 0 & 0 \\ 0 & 1 & 1 & 0 & 1 \\ 1 & 1 & 1 & 1 & 1 \\ 1 & 0 & 1 & 1 & I \\ 0 & 0 & 0 & I & I \end{bmatrix}.$$

Suppose $A_1 = (0\ 1\ 0\ 0\ 0)$ is the instantaneous state vector under consideration i.e., social stigma of having female children. The effect of $A_1$ on the system $N(R)$ is

$$\begin{aligned}
A_1 N(R)^T &= (0\ 0\ 1\ 1\ 1\ 0\ 0) \hookrightarrow (0\ 0\ 1\ 1\ 1\ 0\ 0) = B_1 \\
B_1[N(R)] &= (2, 3, 3, 1, 2) \hookrightarrow (1\ 1\ 1\ 1\ 1) = A_2 \\
A_2[N(R)]^T &= (1+I, 1+I, 1, 1, 1, 1+I, I) \\
&\hookrightarrow (1\ 1\ 1\ 1\ 1\ 1\ I) = B_2 \\
B_2(N(R)) &= (I+3, 3, 5, 2I+2, 2I+4) \\
&\hookrightarrow (1\ 1\ 1\ 1\ 1) = A_3 = A_1.
\end{aligned}$$

Thus this state vector $A_1 = (0, 1, 0, 0, 0)$ gives a fixed point $(1\ 1\ 1\ 1\ 1)$ indicating if one thinks that having female children is a social stigma immaterial of their status they also feel that having number of female children is a problem, it is a economic loss / burden, they also under go torture or bad treatment by in-laws and ultimately it is a insecurity for having only female children, the latter two cases hold where applicable.

On the other hand we derive the following conclusions on the domain space when the range space state vector $A_1 = (0\ 1\ 0\ 0\ 0)$ is sent

$$A_1[N(R)]^T \hookrightarrow B_1$$
$$B_1[N(R)] \hookrightarrow A_2$$
$$A_2[N(R)]^T \hookrightarrow B_2$$
$$B_2[N(R)] \hookrightarrow A_3 = A_2 \text{ so}$$
$$A_2[N(R)]^T \hookrightarrow B_2 = (1, 1, 1, 1, 1, 1, I)$$

leading to a fixed point. When the state vector $A_1 = (0, 1, 0, 0, 0)$ is sent to study i.e. the social stigma node is on uniformly all people from all economic classes are awakened expect the very poor for the resultant vector happens to be a Neutrosophic vector hence one is not in a position to say what is the feeling of the very poor people and the "many female children are a social stigma" as that coordinate remains as an indeterminate one. This is typical of real-life scenarios, for the working classes hardly distinguish much when it comes to the gender of the child.

Several or any other instantaneous vector can be used and its effect on the Neutrosophical Dynamical System can be studied and analysed. This is left as an exercise for the reader.

Having seen an example and application or construction of the NRM model we will proceed on to describe the concepts of it in a more mathematical way.

**DESCRIPTION OF A NRM:**

Neutrosophic Cognitive Maps (NCMs) promote the causal relationships between concurrently active units or decides the absence of any relation between two units or the indeterminance of any relation between any two units. But in Neutrosophic Relational Maps (NRMs) we divide the very causal nodes into two disjoint units. Thus for the modeling of a NRM we need a domain space and a range space which are disjoint in the sense of concepts. We further assume no intermediate relations exist within the domain and the range spaces. The number of elements or nodes in the range space need not be equal to the number of elements or nodes in the domain space.

Throughout this section we assume the elements of a domain space are taken from the neutrosophic vector space of dimension n and that of the range space are from the neutrosophic vector space of dimension m. (m in general need not be equal to n). We denote by R the set of nodes $R_1, \ldots, R_m$ of the range space, where $R = \{(x_1, \ldots, x_m) \mid x_j = 0$ or 1 or $I$ for $j = 1, 2, \ldots, m\}$.

If $x_i = 1$ it means that node $R_i$ is in the on state and if $x_i = 0$ it means that the node $R_i$ is in the off state and if $x_i = I$ in the resultant vector it means the effect of the node $x_i$ is indeterminate or whether it will be off or on cannot be predicted by the neutrosophic dynamical system.

It is very important to note that when we send the state vectors they are always taken as the real state vectors for we know the node or the concept is in the on state or in the off state but when the state vector passes through the Neutrosophic dynamical system some other node may become indeterminate i.e. due to the presence of a node we may not be able to predict the presence or the absence of the other node i.e., it is indeterminate, denoted by the symbol I, thus the resultant vector can be a neutrosophic vector.

**DEFINITION 1.6.2:** *A Neutrosophic Relational Map (NRM) is a Neutrosophic directed graph or a map from D to R with concepts like policies or events etc. as nodes and causalities as edges. (Here by causalities we mean or include the indeterminate causalities also). It represents Neutrosophic Relations and Causal Relations between spaces D and R.*

*Let $D_i$ and $R_j$ denote the nodes of an NRM. The directed edge from $D_i$ to $R_j$ denotes the causality of $D_i$ on $R_j$ called relations. Every edge in the NRM is weighted with a number in the set $\{0, +1, –1, I\}$. Let $e_{ij}$ be the weight of the edge $D_i R_j$, $e_{ij} \in \{0, 1, –1, I\}$. The weight of the edge $D_i R_j$ is positive if increase in $D_i$ implies increase in $R_j$ or decrease in $D_i$ implies decrease in $R_j$ i.e. causality of $D_i$ on $R_j$ is 1. If $e_{ij} = –1$ then increase (or decrease) in $D_i$ implies decrease (or increase) in $R_j$. If $e_{ij} = 0$ then $D_i$ does not have any effect on* $R_j$. If $e_{ij} = I$ it implies we are not in a position to determine the effect of $D_i$ on $R_j$ i.e. the effect of $D_i$ on $R_j$ is an indeterminate so we denote it by $I$.

**DEFINITION 1.6.3:** *When the nodes of the NRM take edge values from {0, 1, –1, I} we say the NRMs are simple NRMs.*

**DEFINITION 1.6.4:** *Let $D_1$, …, $D_n$ be the nodes of the domain space D of an NRM and let $R_1$, $R_2$,…, $R_m$ be the nodes of the range space R of the same NRM. Let the matrix N(E) be defined as $N(E) = (e_{ij})$ where $e_{ij}$ is the weight of the directed edge $D_i R_j$ (or $R_j D_i$) and $e_{ij} \in$ {0, 1, –1, I}. N(E) is called the Neutrosophic Relational Matrix of the NRM.*

The following remark is important and interesting to find its mention in this book.

**Remark**: Unlike NCMs, NRMs can also be rectangular matrices with rows corresponding to the domain space and columns corresponding to the range space. This is one of the marked difference between NRMs and NCMs. Further the number of entries for a particular model which can be treated as disjoint sets when dealt as a NRM has very much less entries than when the same model is treated as a NCM.

Thus in many cases when the unsupervised data under study or consideration can be spilt as disjoint sets of nodes or concepts; certainly NRMs are a better tool than the NCMs.

**DEFINITION 1.6.5:** *Let $D_1$, …, $D_n$ and $R_1$,…, $R_m$ denote the nodes of a NRM. Let $A = (a_1,…, a_n)$, $a_i \in$ {0, 1, I} is called the Neutrosophic instantaneous state vector of the domain space and it denotes the on-off position or indeterminate state of the nodes at any instant. Similarly let $B = (b_1,…, b_n)$ $b_i \in$ {0, 1, I}, B is called instantaneous state vector of the range space and it denotes the on-off position of the nodes at any instant, $a_i = 0$ if $a_i$ is off and $a_i = 1$ if $a_i$ is on for i = 1, 2, …, n. Similarly, $b_i = 0$ if $b_i$ is off and $b_i = 1$ if $b_i$ is on for i = 1, 2,…, m. $a_i = I$ or $b_i = I$ gives the indeterminate* *state at that time or in that situation for i = 1, 2, …, m or i = 1, 2, …, n.*

**DEFINITION 1.6.6:** *Let $D_1,…, D_n$ and $R_1, R_2,…, R_m$ be the nodes of a NRM. Let $D_i R_j$ (or $R_j D_i$) be the edges of an NRM, j = 1, 2,…, m and i = 1, 2,…, n. The edges form a directed cycle. An NRM is said to be a cycle if it possess a directed cycle. An NRM is said to be acyclic if it does not possess any directed cycle.*

**DEFINITION 1.6.7:** *A NRM with cycles is said to be a NRM with feedback.*

**DEFINITION 1.6.8:** *When there is a feedback in the NRM i.e. when the causal relations flow through a cycle in a revolutionary manner the NRM is called a Neutrosophic dynamical system.*

**DEFINITION 1.6.9:** *Let $D_i R_j$ (or $R_j D_i$) $1 \leq j \leq m$, $1 \leq i \leq n$, when $R_j$ (or $D_i$) is switched on and if causality flows through edges of a cycle and if it again causes $R_j$ (or $D_i$) we say that the Neutrosophical dynamical system goes round and round. This is true for any node $R_j$ (or $D_i$) for $1 \leq j \leq m$ (or $1 \leq i \leq n$). The equilibrium state of this Neutrosophical dynamical system is called the Neutrosophic hidden pattern.*

**DEFINITION 1.6.10:** *If the equilibrium state of a Neutrosophical dynamical system is a unique Neutrosophic state vector, then it is called the fixed point. Consider an NRM with $R_1, R_2, …, R_m$ and $D_1, D_2,…, D_n$ as nodes. For example let us start the dynamical system by switching on $R_1$ (or $D_1$). Let us assume that the NRM settles down with $R_1$ and $R_m$ (or $D_1$ and $D_n$) on, or indeterminate or off, i.e. the Neutrosophic state vector remains as (1, 0, 0,…, 1) or (1, 0, 0,…I) (or (1, 0, 0,…1) or (1, 0, 0,…I) in D), this state vector is called the fixed point.*

**DEFINITION 1.6.11:** *If the NRM settles down with a state vector repeating in the form $A_1 \rightarrow A_2 \rightarrow A_3 \rightarrow … \rightarrow A_i \rightarrow A_1$*

*(or $B_1 \to B_2 \to \ldots \to B_i \to B_1$) then this equilibrium is called a limit cycle.*

Let $R_1, R_2, \ldots, R_m$ and $D_1, D_2, \ldots, D_n$ be the nodes of a NRM with feedback. Let $N(E)$ be the Neutrosophic Relational Matrix. Let us find the hidden pattern when $D_1$ is switched on i.e. when an input is given as a vector; $A_1 = (1, 0, \ldots, 0)$ in D; the data should pass through the relational matrix $N(E)$. This is done by multiplying $A_1$ with the Neutrosophic relational matrix $N(E)$. Let $A_1 N(E) = (r_1, r_2, \ldots, r_m)$ after thresholding and updating the resultant vector we get $A_1 E \in R$. Now let $B = A_1 E$ we pass on B into the system $(N(E))^T$ and obtain $B(N(E))^T$. We update and threshold the vector $B(N(E))^T$ so that $B(N(E))^T \in D$.

This procedure is repeated till we get a limit cycle or a fixed point.

**DEFINITION 1.6.12:** *Finite number of NRMs can be combined together to produce the joint effect of all NRMs. Let $N(E_1), N(E_2), \ldots, N(E_r)$ be the Neutrosophic relational matrices of the NRMs with nodes $R_1, \ldots, R_m$ and $D_1, \ldots, D_n$, then the combined NRM is represented by the neutrosophic relational matrix $N(E) = N(E_1) + N(E_2) + \ldots + N(E_r)$.*

For more about FCM, FRM, Neutrosophic Cognitive Maps and Neutrosophic Relational Maps, please refer Book [77].

Chapter Two

# UNTOUCHABILITY: PERIYAR'S VIEW AND PRESENT DAY SITUATION
*A Fuzzy and Neutrosophic Analysis*

In this chapter, we use experts' opinions, documented writings and speeches of Periyar and contemporary casteist atrocities to analyze the evils of untouchability and its consequences using Fuzzy Cognitive Maps (FCMs) and Neutrosophic Cognitive Maps (NCMs). For the first time, we use mathematical models to analyze untouchability. We specifically use Fuzzy and Neutrosophic models.

At the first stage, since the term untouchability is a vague and ill-defined concept because of its diverse manifestations, we use fuzzy and neutrosophic models. Our study would be very vast if we take into consideration the notion of untouchability from the time of the ancient lawgiver Manu (approx. 1 Common Era)[1]. So, we specifically study the social effects of untouchability during Periyar E. V. Ramasamy's lifetime (1879-1973).[2]

A relentless crusader for social justice, Periyar fought against the caste system for more than half a century. He rendered thousands of speeches condemning it and wrote hundreds of articles against it in his newspapers and magazines like *Kudiarasu, Viduthalai, Revolt,* etc. and his struggle had far-reaching impact in the Tamil society. Always a leader of mass organizations, he served as the State President of the Congress party, founded the Self-

---

[1] Wendy Doniger, *The Laws of Manu*, Penguin, 1991.
[2] For more biographical information refer p.103-113 of this book.

Respect Movement, was a leader of the Justice Party, later the Dravidar Kazhagam, he was undoubtedly a populist leader and a man of the masses. He fought against the caste system by organizing demonstrations and was arrested and incarcerated several times.

In this chapter we analyze his views and the present state of untouchability in India using the two technical tools: Fuzzy Cognitive Maps (FCMs) and Neutrosophic Cognitive Maps (NCMs).

This chapter is divided into six sections. In the first section, we analyze the evils of untouchability due to Hinduism and its codes. In section two, the discrimination faced by Dalits and Sudras in the field of education is studied. In section three, we analyze the social inequality experienced by Dalits and Sudras. In section four we analyze the discrimination faced by Dalits in the political field. Economic threat faced by Dalits due to untouchability is studied in section five.

All the above sections utilize the two tools FCMs and NCMs. In the final section we analyze the generalized social problem of untouchability using a different tool, Fuzzy Relational Maps (FRMs).

Experts' opinions are given verbatim. No change or modification has been made since it would make the data biased. Further, in most places, we provide an exact description of the concepts/nodes as specified by the expert who have their own views about untouchability.

In several models, the nodes were very intricately defined by the experts, so, the ON state of a single practice of untouchability made the resultant, hidden pattern of the dynamical system to give fixed point in which all the nodes (connected to practices of untouchability) were in the ON state. Most of the directed graphs obtained from the experts' opinions were highly dense.

We have refrained from giving our views since the sole purpose of this study is a mathematical analysis of untouchability.

## 2.1 Analysis of untouchability due to Hindu religion using FCMs and NCMs

We proceed on to categorize untouchability practiced due to Hindu religion. Nodes related to religious untouchability, as given by an expert who is a social activist, are as follows:

- $R_1$ - Religion
- $R_2$ - Superstition
- $R_3$ - Faith in particular religious sect
- $R_4$ - Discrimination in religion
- $R_5$ - No freedom of choice (caste is birth-based)
- $R_6$ - Untouchability
- $R_7$ - Caste system
- $R_8$ - Psychological oppression
- $R_9$ - Discrimination in social outlook
- $R_{10}$ - Practice of Varnashrama Dharma
- $R_{11}$ - Social identity
- $R_{12}$ - Social fear
- $R_{13}$ - Social binding
- $R_{14}$ - Social rituals
- $R_{15}$ - Solace

As described by the expert, these attributes are explained below.

**$R_1$ – Religion**

By the term 'religion' the expert refers to the social codes and spiritual activities one is supposed to follow by belonging to a particular caste.

This expert points out the difference between the native religion of the Tamil people and the casteist religion of the Aryans (which was based on the Manusmriti and Vedas), which gave sanction to untouchability, and hierarchy among the people. He says that there was no historical evidence for the existence of religion that could be labeled Hinduism.

He even adds that the social and spiritual practices followed by the indigenous people were distinctly different compared to the practices prescribed by the *Manu Smriti* (Laws of Manu). He also pointed out the vast differences between these two systems.

A woman was never considered a possession of her father, son or husband in the Tamil culture. She was free to choose her husband unlike in the Manuvadi culture. This is a contradiction in the case of women identity among the ancient Tamils and the Manu Smriti. This is strong evidence that the way of life of the Tamilians was in no way related to Manu Smriti.

Likewise there were several differences in rituals, gods, sacrifices, modes of worship etc.

This expert points out that the splendid Tamil text *Thirukkural* was written by a Paraiyar called Valluvar, whereas Aryan religion banned education to Sudras.

### $R_2$: Superstition

In every culture, superstition was a negative force that eclipsed the scientific development of human beings. The fear of snakes made them worship it with a superstition that it would consequently not harm them.

According to this expert, superstition has become so strong in this modern India that even monkeys are worshipped in the form of Hanuman, whose tiny statues litter the countryside and giant statue line the highways! This has no scientific rhyme or reason! Thus the $21^{st}$ century's wave of superstition is the worship of monkey. Hanuman statues up to 60 feet in height are worshiped in Tamil Nadu. Why does nobody question this stupidity?

Then, the cow is treated as sacred to Hindus and is considered a mother by all Hindus. So caste-Hindus don't consume its meat and consider its hide polluting. This has resulted in ban on cow slaughter that has affected the livelihood of Dalits. Several other irrational behaviors can be discussed.

For instance, the concepts of purity and pollution are merely religious taboos.

**R₃: Faith in a particular religious sect**

If anyone takes a statistics of the number of people converting to Christianity, Buddhism, Jainism or Islam from their native religion, Dalits would top that list. They convert to other religions because they find Hinduism to be very oppressive and they strive to escape from caste atrocities. Many of them feel that the change of faith has changed their lifestyle and socio-economic status. Besides, the vicious cycle of 'karma' can be broken by religious conversion.

Faith gives hope for betterment. Faith also provides them an improved identity and a new sense of belonging.

**R₄ – Discrimination in religion**

The Hindu religion follows the fourfold varna structure. At the head of the hierarchy come the *Brahmins* (priestly class), *Kshatriyas* (warriors), *Vaishyas* (merchants and artisans) and *Sudras* (slaves). Outside this caste-Hindu fold are the Untouchables, also called the *Panchamas*. These four varnas are in turn divided into thousands of castes and sub-castes.

According to classic Hinduism, the origin of the caste system is as old as the origin of the human race itself. The story of this cosmogonic myth regarding the creation of the physical elements of the universe from the Rig Veda goes as follows:

> "The gods made a Man [*Purusha*] who had a thousand heads, a thousand eyes, a thousand feet. When the gods spread the sacrifice, using the Man as the offering, spring was the clarified butter, summer the fuel, autumn the oblation.
> From that sacrifice in which everything was offered, the clarified butter was obtained, and they made it into those beasts who live in the air, in the forest and in villages.
> His mouth was the Brahmins, his arms were made into the Kshatriyas, his thighs became the

Vaishyas and his feet became the Shudras. The moon was born from his mind; the sun was born from his eye."

The pity is that Dalits don't even find mention in the creation. Even their shadow is considered polluting. Rigid mathematical computations decree that some castes are polluting at four feet, some at eight feet, and some as far as even sixty-four feet.

Discrimination against Dalits is dominant even today in the public and private spheres of life. Prior to Indian independence, Dalits where not even allowed into temples. Even now, in rural areas Dalits cannot walk on temple streets or enter caste-Hindu temples. When Dalits demand their rights to temple entry and temple property, it often results in the eruption of large-scale violence and bloodshed.

### $R_5$: No freedom of choice

Although the Constitution of India gives the freedom to preach and practice any religion of ones choice, in practice this is not so straight or simple. Dalits who convert to Christianity or Islam lose the benefits of reservation in education and employment because they are no longer treated as Dalits. So, it is a choiceless situation for them. If they leave Hinduism which is imposed on them, they lose the benefits of affirmative action. So, majority of the Dalits continue to remain in their same state, without converting to other religions.

### $R_6$: Untouchability

This concept is closely interconnected with $R_1$, religion. It is religion that has sanctioned untouchability and imposed it upon people. Only an untouchable is untouchable, but all the products of his labour (farming, leather work, etc.) are touchable. It proves that this concept is only to subjugate a section of people. The concept of untouchability has made some people into permanent slaves.

Not only where some people called untouchable, but a few of them were even made unseeables such as the *Nadar* and *Pudharai Vannars*.

The Brahmins created this idea of untouchability in order to protect themselves as a superior group. The very root of untouchability lies in the laws of Manu, which was a codification of existing caste practices. Untouchability is a by-product of the caste system.

### R$_7$: Caste system

One of the most peculiarly oppressive systems, all the castes and the thousands of sub-castes are religiously forbidden, to intermarry, or interdine, or engage in social activity with any but members of their own group. The caste system was clearly an ideological construct of the upper castes to maintain their sinister monopoly over cultural capital (knowledge and education), social capital (status and patriarchal domination), political capital (power) and material capital (wealth).

This expert quoted revolutionary Dr.Ambedkar: "Caste is not a division of labour, it is a division of labourers."

The concept of superiority and inferiority and the hierarchy of caste system has claimed millions of lives through riots, rampages and bloodshed. Needless to say that the major victims were only Dalits and Sudras.

### R$_8$: Psychological oppression

Because of the caste hierarchy, concepts of superiority and inferiority are ingrained in the minds of people. This psychologically suppresses the Dalits and the Sudras from reaching their true potential. Brahmins are considered 'intellectual' and they use their superiority to snub and deny education and opportunities to the lower castes and Dalits. Since the caste-system is birth-based, people have no chance of getting rid of caste or their inferior identity. So it oppresses them psychologically.

### $R_9$ – Social discrimination

(*According to this expert this concept is different from the concept of religious discrimination, i.e. node $R_4$*).

Dalits are denied their civil rights. They are treated as socially inferior people. For centuries Dalits have been living in separate, segregated settlements known as the *cheri*, which is situated outside the caste Hindu village (*oor*). In certain places, they are not allowed to walk on the streets wearing slippers. In backward areas where casteism is rampant, Dalits are not allowed to carry umbrellas, ride bicycles or wear new clothes. They are not allowed to draw water from common ponds, or graze their cattle on public grounds. In teashops in many villages, the two-tumbler system is rampant: one set of tumblers for the caste-Hindus and another for the Dalits. Even very old Dalit men are addressed in the singular and treated disrespectfully by small caste-Hindu boys. Even after death they are treated differently, because caste-Hindus and Dalits have separate burial grounds. Inter-caste marriages are strictly prohibited. In several cases, those who dare to love across caste lines have to face the dire consequence of death.

### $R_{10}$ – Varnashrama Dharma

It was the earlier, original form of the caste system. Only the four varna system developed into the thousands of *jatis* (castes) and sub-castes. Varnashrama Dharma decided ones position in life and the caste-based occupation they had to follow hereditarily. In a democratic nation, one must have the right to choose ones profession. But it is not so. Varnashrama Dharma robs people of the basic equality of opportunity. Dalits and Sudra students have to overcome all sorts of obstacles in order to get educated and employed.

All promotions, perks and increments etc. are denied to them, and they are labeled 'quota' people because they might have benefited from affirmative action.

Thus even if they come out of the clutches of the Varnashrama Dharma they are the worst victims of discrimination.

### $R_{11}$ – Social identity

Everybody requires a social identity, and sense of belonging to a larger power group. The Dalits and Sudras are unable to throw away their identity of being outcastes and 'lower' castes. The quest for an improved social identity and self-respect triggers religious conversions.

### $R_{12}$ – Social Fear

Out of social fear they identify themselves with a particular caste and lead a life of constant suffering. For instance the Sudras think that they will be treated as outcastes (a category lower to them) so they emphasize their caste status.

So it is the social fear (of being considered lower than what they are in the hierarchy) inherent in them makes them accept their caste. This is because the caste system has been imposed on them and they have begun to value it.

### $R_{13}$ – Social Bonding

The expert is of the opinion that they become socially bonded along caste lines because of social fear, fear of further exploitation, ill treatment, and harassment. For the sake of sheer existence, peace and to avoid further torture they accept the social bond.

### $R_{14}$ – Social Rituals

They are involuntarily forced to practice social rituals in accordance with religion mainly due to social identification,

social fear and social binding (nodes $R_{11}$, $R_{12}$ and $R_{13}$). The expert says that social celebrations and festivals are highly ritualistic and intended to develop the bonding within a community as well as mark their identity.

**$R_{15}$ – Solace**

Any human follows a religion and worships a god mainly for solace and peace of mind. People who suffer caste harassment and subjugation believe that god or religion can give them some respite.

The directed graph given by this expert is as follows:

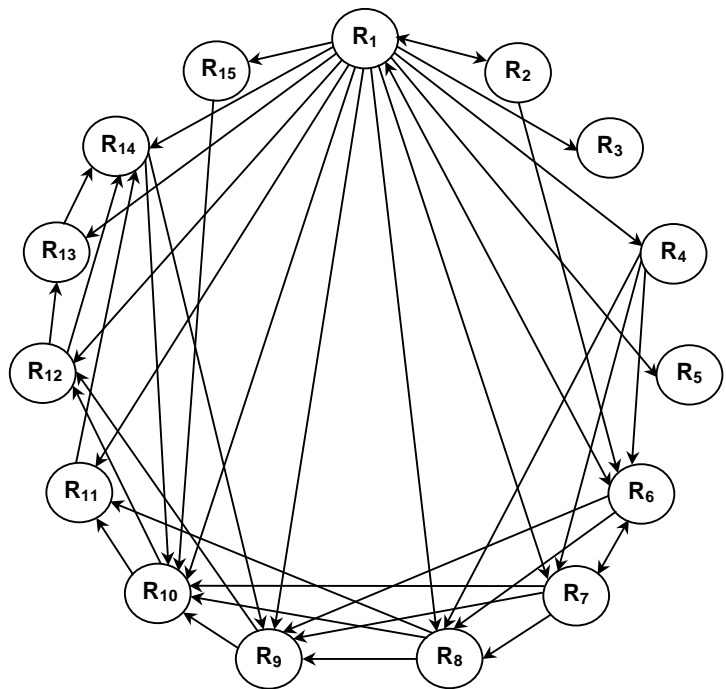

FIGURE 2.1.1

The related connection matrix R is as follows

$$\begin{array}{c} \phantom{R_1} \quad R_1 \; R_2 \; R_3 \; R_4 \; R_5 \; R_6 \; R_7 \; R_8 \; R_9 \; R_{10} \; R_{11} R_{12} R_{13} R_{14} R_{15} \\ \begin{array}{c} R_1 \\ R_2 \\ R_3 \\ R_4 \\ R_5 \\ R_6 \\ R_7 \\ R_8 \\ R_9 \\ R_{10} \\ R_{11} \\ R_{12} \\ R_{13} \\ R_{14} \\ R_{15} \end{array} \begin{bmatrix} 0 & 1 & 1 & 1 & 1 & 1 & 1 & 1 & 1 & 1 & 1 & 1 & 1 & 1 & 1 \\ 1 & 0 & 0 & 0 & 0 & 1 & 0 & 0 & 0 & 0 & 0 & 0 & 0 & 0 & 0 \\ 0 & 0 & 0 & 0 & 0 & 0 & 0 & 0 & 0 & 0 & 0 & 0 & 0 & 0 & 0 \\ 0 & 0 & 0 & 0 & 0 & 1 & 1 & 1 & 0 & 0 & 0 & 0 & 0 & 0 & 0 \\ 0 & 0 & 0 & 0 & 0 & 0 & 0 & 0 & 0 & 0 & 0 & 0 & 0 & 0 & 0 \\ 1 & 0 & 0 & 0 & 0 & 0 & 1 & 1 & 1 & 0 & 0 & 0 & 0 & 0 & 0 \\ 0 & 0 & 0 & 0 & 0 & 1 & 0 & 1 & 1 & 1 & 0 & 0 & 0 & 0 & 0 \\ 0 & 0 & 0 & 0 & 0 & 0 & 0 & 0 & 1 & 1 & 1 & 0 & 0 & 0 & 0 \\ 0 & 0 & 0 & 0 & 0 & 0 & 0 & 0 & 1 & 0 & 1 & 0 & 0 & 0 & 0 \\ 0 & 0 & 0 & 0 & 0 & 0 & 0 & 0 & 0 & 0 & 1 & 1 & 0 & 0 & 0 \\ 0 & 0 & 0 & 0 & 0 & 0 & 0 & 0 & 0 & 0 & 0 & 0 & 0 & 1 & 0 \\ 0 & 0 & 0 & 0 & 0 & 0 & 0 & 0 & 0 & 0 & 0 & 1 & 1 & 0 \\ 0 & 0 & 0 & 0 & 0 & 0 & 0 & 0 & 0 & 0 & 0 & 0 & 1 & 0 \\ 0 & 0 & 0 & 0 & 0 & 0 & 0 & 1 & 1 & 0 & 0 & 0 & 0 & 0 \\ 0 & 0 & 0 & 0 & 0 & 0 & 0 & 0 & 0 & 1 & 0 & 0 & 0 & 0 \end{bmatrix} \end{array}$$

R is a 15 × 15 matrix. The analysis only uses simple FCMs. If we pass on the state vector X = (1 0 0 0 0 … 0), where all the nodes are in the OFF state, except the attribute $R_1$ (religion) we see the effect of the state vector X on the dynamical system R

$$X R \; \hookrightarrow \; (1\;1\;1\;1\;1\;1\;1\;1\;1\;1\;1\;1\;1\;1\;1) \;=\; X_1$$
$$X_1 R \; \hookrightarrow \; (1\;1\;1\;1\;1\;1\;1\;1\;1\;1\;1\;1\;1\;1\;1) \;=\; X_1$$

The symbol '$\hookrightarrow$' denotes the resultant of XR after thresholding and updating the vector. Thus when the concept or the attribute 'religion' alone is in the ON state we see all other attributes come to the ON state. Thus the hidden pattern of the dynamical system states that the religion (here, Hindu Brahminical religion) induces on people the superstition, caste system, untouchability, psychological oppression, discrimination of all forms, Varnashrama Dharma, social fear and so on. Thus, as per this expert we can say what Periyar said was correct: Don't follow any religion, particularly Hinduism that is nothing but Varnashrama Dharma (caste system). That is why the Self-Respect Movement founded by Periyar fought against the oppressive interlinked concepts of religion and god.

Thus it is clear from the dynamical system, religion as given by Manu Smriti, Vedas, Gita and so on is the root cause of all discrimination, under-development and unrest prevailing in the nation. Dalits and Backward classes (Sudras) lead a life of insecurity, fear, segregation, and discrimination. Majority of them live below the poverty line, without even proper school education, no proper occupation, work over ten hours a day. The atrocities committed against the dalits shows the worst side of Hindu religion.

As suggested by the expert we analyze the state vector Y = (0 0 0 0 0 1 0 0 0 0 0 0 0 0) i.e., only the node untouchability is in the ON state and all other nodes are in the OFF state. We study the hidden pattern given by the vector Y on the dynamical system R.

$$YR \hookrightarrow (1\ 0\ 0\ 0\ 0\ 1\ 1\ 1\ 1\ 0\ 0\ 0\ 0\ 0\ 0) = Y_1$$
$$Y_1 R \hookrightarrow (1\ 1\ 1\ 1\ 1\ 1\ 1\ 1\ 1\ 1\ 1\ 1\ 1\ 1\ 1) = Y_2$$
$$Y_2 R = Y_2.$$

Thus when the node untouchability is in the ON state all other nodes becomes ON just like the earlier example. This shows that the practice of untouchability is the cause of all problems and humiliations.

Suppose we consider the node 'caste system' to be in the ON state and all other nodes remain in the OFF state we obtain the state vector Z = (0 0 0 0 0 0 1 0 0 0 0 0 0 0). The effect of Z on R gives

$$ZR \hookrightarrow (0\ 0\ 0\ 0\ 0\ 1\ 1\ 1\ 1\ 1\ 0\ 0\ 0\ 0\ 0) = Z_1 \text{ (say)}$$
$$Z_1 R \hookrightarrow (1\ 0\ 0\ 0\ 0\ 1\ 1\ 1\ 1\ 1\ 1\ 0\ 0\ 0\ 0\ 0) = Z_2 \text{ (say)}$$
$$Z_2 R \hookrightarrow (1\ 1\ 1\ 1\ 1\ 1\ 1\ 1\ 1\ 1\ 1\ 1\ 1\ 1\ 1) = Z_3 \text{ (say)}$$
$$Z_3 R \hookrightarrow Z_3$$

leading to a fixed point.

The hidden pattern of the state vector when caste system is in the ON state makes all other nodes come to the ON state. Thus the caste system is sufficient to ruin the nation.

On the other hand if the node 'faith in religion' is in the ON state i.e., we obtain the following fixed point

A    =    (0 0 1 0 0 0 0 0 0 0 0 0 0 0 0)
AR   ↪    (0 0 1 0 0 0 0 0 0 0 0 0 0 0 0)  =  A.

Thus 'faith in religion' does not have any impact on the social setup. Likewise we can study any node.

Take the state vector B = (0 0 0 0 0 0 0 1 0 0 0 0 0 0 0) i.e., the node 'psychological oppression', to be in the ON state and all other nodes are in the OFF state. Consider the effect of B on R.

BR    ↪    (0 0 0 0 0 0 0 1 1 1 1 0 0 0 0)   =  $B_1$ (say)
$B_1$R  ↪  (0 0 0 0 0 0 0 1 1 1 1 1 0 1 0)   =  $B_2$ (say)
$B_2$R  ↪  (0 0 0 0 0 0 0 1 1 1 1 1 1 1 0)   =  $B_3$ (say)
$B_3$R  ↪  (0 0 0 0 0 0 0 1 1 1 1 1 1 1 0)   = $B_4$(= $B_3$ say)

We see that 'psychological oppression' makes social fear and social binding and related factors to come to the ON state.

The same type of analysis can be carried out by taking one or several set of nodes in the ON state and analyze them.

We have used a C-program to analyze the various concepts related to this directed graph. For C-program refer [69, 76]. Detailed conclusions are given in the final chapter.

We analyze another expert's opinion for the same attribute viz. Untouchability and Hindu religion. First we consider the nine attributes given by this expert. He has clubbed several of the social factors

$P_1$  -  Religious cruelty
$P_2$  -  Untouchability
$P_3$  -  Caste system
$P_4$  -  Varnashrama Dharma

P_5 - Manu Dharma
P_6 - Samadharma (Equality)
P_7 - Atheism
P_8 - Social inequality
P_9 - Psychological oppression

These attributes have been described earlier in this section. The directed graph given by this expert is as follows:

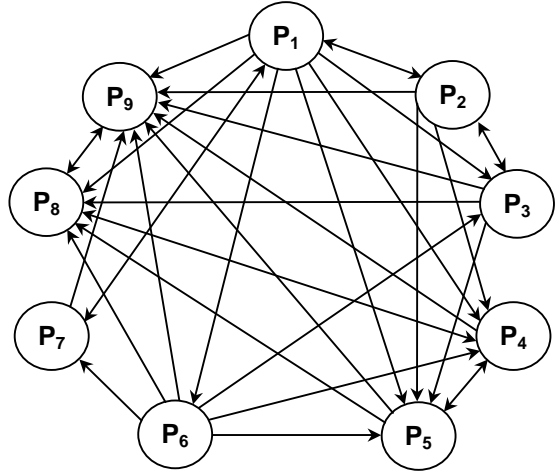

FIGURE 2.1.2

The related connection matrix of the directed graph, P is as follows:

$$P = \begin{array}{c} \\ P_1 \\ P_2 \\ P_3 \\ P_4 \\ P_5 \\ P_6 \\ P_7 \\ P_8 \\ P_9 \end{array} \begin{array}{c} \begin{array}{ccccccccc} P_1 & P_2 & P_3 & P_4 & P_5 & P_6 & P_7 & P_8 & P_9 \end{array} \\ \left[ \begin{array}{ccccccccc} 0 & 1 & 1 & 1 & 1 & 1 & 1 & 1 & 1 \\ 1 & 0 & 1 & 1 & 1 & 0 & 0 & 0 & 1 \\ 0 & 1 & 0 & 0 & 1 & 0 & 0 & 1 & 1 \\ 0 & 0 & 0 & 0 & 1 & 0 & 0 & 1 & 1 \\ 0 & 0 & 0 & 1 & 0 & 0 & 0 & 1 & 1 \\ 0 & 0 & -1 & 1 & -1 & 0 & 1 & 1 & -1 \\ 1 & 0 & 0 & 0 & 0 & 0 & 0 & 0 & -1 \\ 0 & 0 & 0 & 1 & 0 & 0 & 0 & 0 & 1 \\ 0 & 0 & 0 & 0 & 0 & 0 & 0 & 1 & 0 \end{array} \right] \end{array}$$

Let P be the connection matrix. According to this expert, Samadharma[3] increases with the annihilation (decrease) of the caste system.

We shall see the effect of the state vectors on the dynamical system P.

Let X = (1 0 0 0 0 0 0 0) be the state vector which depicts the ON state of the node 'religious cruelty' and all other nodes are in the OFF state. Effect of X on P is given by

$$XP \hookrightarrow (1\ 1\ 1\ 1\ 1\ 1\ 1\ 1) = X_1 \text{ (say).}$$
$$X_1 P \hookrightarrow (1\ 1\ 1\ 1\ 1\ 1\ 1\ 1) = X_2 \text{ (say).}$$

Clearly $X_1 = X_2$. Thus religious cruelty is the root cause of all factors like untouchability, caste system, Manu Dharma, Varnashrama Dharma, social inequality and psychological oppression.

Thus religious cruelty is the basis for all disturbances, distrust, discrimination, desperateness, divisions and depressions in society. It has forced them to seek means to protect themselves from further harassment by taking refuge in atheism. Naturally they defy religion so that the cruel clutches of religion are released and an atheist feels that at least he is liberated.

Likewise Samadharma also works as an antidote for Varnashrama Dharma and Manu Dharma. After a prolonged discussion about Hinduism this expert, a Periyarist changed the node religion to religious cruelty.

Consider the node Samadharma alone in the ON state and all other nodes in the OFF state.

$$Y = (0\ 0\ 0\ 0\ 0\ 1\ 0\ 0\ 0)$$
$$YP = (0\ 0\ -1\ 1\ -1\ 1\ 1\ 1\ -1)$$
$$\hookrightarrow (0\ 0\ 0\ 1\ 0\ 1\ 1\ 1\ 0) = Y_1 \text{ (say)}$$

---

[3] *According to Periyar, the concept of Samadharma or egalitarianism was postulated as an alternative to Manudharma and the Varnashrama Dharma. This word gained currency when he founded the Self Respect Movement.*

$$Y_1 P \hookrightarrow (1\ 0\ 0\ 1\ 1\ 1\ 1\ 1\ 1) = Y_2 \text{ (say)}.$$
$$Y_2 P \hookrightarrow (1\ 1\ 1\ 1\ 1\ 1\ 1\ 1\ 1) = Y_3 \text{ (say)}.$$
$$Y_3 P = Y_3.$$

This gives a fixed point, which converts all other nodes to the ON state. So the node Samadharma has a very positive impact. When it is in the ON state it destroys religious cruelty, untouchability, caste system and Varnashrama Dharma. Samadharma supports atheism, breaks social inequality and psychological oppression, thereby all states are activated.

Next we study the effect of social inequality on the system P.

Let us take the state vector

$$Z = (0\ 0\ 0\ 0\ 0\ 0\ 0\ 1\ 0)$$

which has only social inequality in the ON state and all other nodes in the OFF state. The effect of Z on the dynamical system P is given by

$$ZP \hookrightarrow (0\ 0\ 0\ 1\ 0\ 0\ 0\ 1\ 1) = Z_1 \text{ (say)}.$$
$$Z_1 P \hookrightarrow (0\ 0\ 0\ 1\ 1\ 0\ 0\ 1\ 1) = Z_2 \text{ (say)}$$
$$Z_2 P \hookrightarrow (0\ 0\ 0\ 1\ 1\ 0\ 0\ 1\ 1) = Z_2.$$

We arrive at the fixed point. The effect of social inequality on the system turns the nodes Varnashrama Dharma, Manu Dharma and psychological oppression to ON state.

According to this dynamical system, social inequality causes psychological oppression and these are caused due to Manu Dharma and Varnashrama Dharma.

Several other experts' opinions are taken for our study and it's modeled in the form of connection matrices. A thorough analysis of the inter-relation between attributes is given in the final chapter.

For these same 9 concepts we got the opinion of the same expert using Neutrosophic Cognitive Map (NCM), where he has the option to state if certain relationships given in the graph are indeterminate.

The neutrosophic directed graph is as follows:

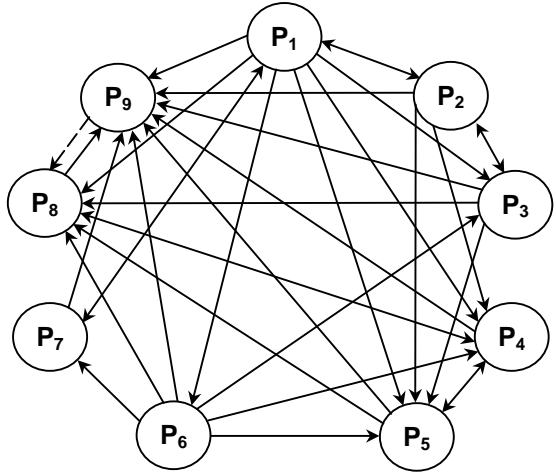

FIGURE 2.1.3

The related connection matrix is given by P'.

$$P' = \begin{array}{c} \\ P_1 \\ P_2 \\ P_3 \\ P_4 \\ P_5 \\ P_6 \\ P_7 \\ P_8 \\ P_9 \end{array} \begin{array}{c} \begin{array}{ccccccccc} P_1 & P_2 & P_3 & P_4 & P_5 & P_6 & P_7 & P_8 & P_9 \end{array} \\ \left[ \begin{array}{ccccccccc} 0 & 1 & 1 & 1 & 1 & 1 & 1 & 1 & 1 \\ 1 & 0 & 1 & 1 & 1 & 0 & 0 & 0 & 1 \\ 0 & 1 & 0 & 0 & 1 & 0 & 0 & 1 & 1 \\ 0 & 0 & 0 & 0 & 1 & 0 & 0 & 1 & 1 \\ 0 & 0 & 0 & 1 & 0 & 0 & 0 & 1 & 1 \\ 0 & 0 & -1 & 1 & -1 & 0 & 1 & 1 & -1 \\ 1 & 0 & 0 & 0 & 0 & 0 & 0 & 0 & -1 \\ 0 & 0 & 0 & 1 & 0 & 0 & 0 & 0 & 1 \\ 0 & 0 & 0 & 0 & 0 & 0 & 0 & I & 0 \end{array} \right] \end{array}$$

Consider the state vector $Z' = (0\ 0\ 0\ 0\ 0\ 0\ 0\ 1\ 0)$ on the neutrosophic matrix P'. We see the resultant remains the same. One can work with any other state vector and derive conclusions.

## 2.2 Analysis of discrimination faced by Dalits/ Sudras in the field of education as untouchables using FCMs and NCMs

In this section, we discuss the denial of education to Dalits and Backward classes. The literacy rate of Dalits is much lower compared to the general population. This expert says that the Brahmins and caste-Hindus systematically deny education to Dalits and Sudras right from school to professional levels. Dalits and Sudras are harassed, tormented, and driven away from their professions.

In the Manusmriti it was decreed that if Sudras were going to even hear the Vedas (the holy texts of the Brahmins) being uttered, molten lead had to be poured into their ears as a punishment. The expert quoted Jotirao Phule[4] who spoke of how the Dalit and Sudras were actively denied education by the Brahminical orthodoxy, deduces: "Without knowledge, intelligence was lost, without intelligence, morality was lost and without morality was lost all dynamism! Without dynamism, money was lost and without money the Sudras sank. All this misery was caused by the lack of knowledge."

Our expert also quoted Periyar who had said that even during the British and Congress rule, the Dalit and Sudra children were denied entry into certain schools. For instance, even today in villages Dalit children are segregated and made to sit apart in the classroom. In educational institutions run by Hindu religious leaders like the Sankaracharya, Brahmin and non-Brahmin college students are given different types of food and statewide agitations took place against it. [88]

It goes without saying how much discrimination and intellectual harassment they undergo in classrooms and college premises. So we are justified in mathematically analyzing the problem faced by Dalits and Sudras in education, which we have broadly classified as educational untouchability. Since the attributes chosen by different experts is different the notion of Combined Fuzzy Cognitive Maps (CFCMs) cannot be applied to this model.

---

[4] Phule (1827 -1890) the earliest pioneer of the non-Brahmin movement in India established the Satyashodhak Samaj (the society of seekers of the Truth).

Further, the expert opinion (unsupervised data) is full of sensitive, rationalist feeling not based on any science/ law. So we are justified in using fuzzy theory to study it. The attributes given by this expert as the causes of educational untouchability or educational discrimination is as follows.

- $E_1$ - Discrimination in education
- $E_2$ - Varnashrama Dharma
- $E_3$ - Manu Dharma
- $E_4$ - Samadharma
- $E_5$ - Fear of studies/ harassment by high caste teachers
- $E_6$ - Discouragement to read/study
- $E_7$ - Role of Brahmins and caste-Hindus
- $E_8$ - Untouchability
- $E_9$ - Economic condition
- $E_{10}$ - Social condition

The associated directed graph given by this expert.

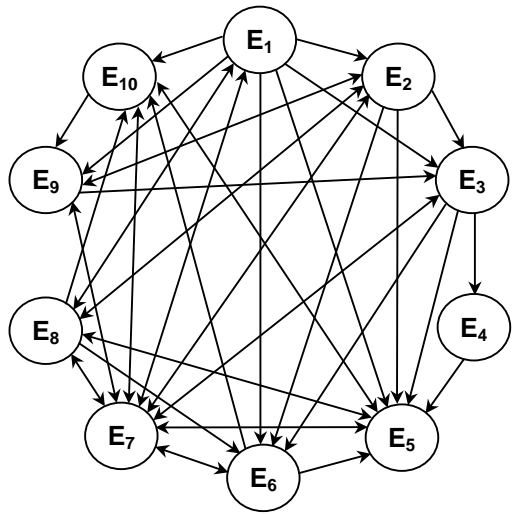

FIGURE 2.2.1

The related connection matrix for the above directed graph is given by E

$$E = \begin{array}{c} \\ E_1 \\ E_2 \\ E_3 \\ E_4 \\ E_5 \\ E_6 \\ E_7 \\ E_8 \\ E_9 \\ E_{10} \end{array} \begin{array}{c} E_1\ E_2\ E_3\ E_4\ E_5\ E_6\ E_7\ E_8\ E_9\ E_{10} \\ \begin{bmatrix} 0 & 1 & 1 & 0 & 1 & 1 & 1 & 1 & 1 & 1 \\ 0 & 0 & 1 & 0 & 1 & 1 & 1 & 1 & 1 & 0 \\ 0 & 0 & 0 & 1 & 1 & 1 & 1 & 0 & 0 & 0 \\ 0 & 0 & 0 & 0 & -1 & 0 & 0 & 0 & 0 & 0 \\ 0 & 0 & 0 & 0 & 0 & 0 & 1 & 1 & 0 & 1 \\ 0 & 0 & 0 & 0 & 1 & 0 & 1 & 0 & 0 & 1 \\ 1 & 1 & 1 & 0 & 1 & 1 & 0 & 1 & 1 & 1 \\ 1 & 1 & 0 & 0 & 1 & 1 & 1 & 0 & 0 & 1 \\ 0 & 1 & 1 & 0 & 0 & 0 & 1 & 0 & 0 & 0 \\ 0 & 0 & 0 & 0 & 1 & 0 & 1 & 0 & 1 & 0 \end{bmatrix} \end{array}$$

Now we take the vector X = (1 0 0 0 0 0 0 0 0 0), i.e., only discrimination in education is in the ON state and all other nodes are in the OFF state. The effect of X on the dynamical system E is given by

$$\begin{array}{lll} XE & \hookrightarrow\ (1\ 1\ 1\ 0\ 1\ 1\ 1\ 1\ 1\ 1) & =\ X_1\ (\text{say}) \\ X_1E & \hookrightarrow\ (1\ 1\ 1\ 1\ 1\ 1\ 1\ 1\ 1\ 1) & =\ X_2\ (\text{say}) \\ X_2E & \hookrightarrow\ (1\ 1\ 1\ 1\ 1\ 1\ 1\ 1\ 1\ 1) & =\ X_2 \end{array}$$

which is clearly a fixed point. We see from the hidden pattern of the dynamical system that discrimination in education is due to Manu Dharma, Varnashrama Dharma, untouchability, Fear of studies/ harassment by high caste teachers, role of Brahmins and caste Hindus, discouragement to read/study, economic and social condition. Thus the notion of Samadharma becomes ON to fight against the discrimination in education.

Next we take the state vector Y = (0 0 0 0 1 0 0 0 0 0), where the only node which is ON is the 'Fear of studies/ harassment by high caste teachers' and all other nodes are in the OFF state. We study the effect of Y on E

$$
\begin{aligned}
YE &\hookrightarrow (0\,0\,0\,0\,1\,0\,1\,1\,0\,1) &&= Y_1 \text{ say} \\
Y_1E &\hookrightarrow (1\,1\,1\,0\,1\,1\,1\,1\,1\,1) &&= Y_2 \text{ (say)}. \\
Y_2E &\hookrightarrow (1\,1\,1\,1\,1\,1\,1\,1\,1\,1) &&= Y_2.
\end{aligned}
$$

We see that when 'fear of studies/ harassment by high caste teachers' is in the ON state, it causes discouragement in education and so on.

Now we work with the node 'untouchability' alone in the ON state and all other nodes are in the OFF state. Let us denote the state vector by $Z = (0\,0\,0\,0\,0\,0\,0\,1\,0\,0)$.

The effect of Z on the dynamical system E is given by

$$
\begin{aligned}
ZE &\hookrightarrow (1\,1\,0\,0\,1\,1\,1\,1\,0\,1) &&= Z_1 \text{ (Say)} \\
Z_1E &\hookrightarrow (1\,1\,1\,0\,1\,1\,1\,1\,1\,1) &&= Z_2 \text{ (say)} \\
Z_2E &\hookrightarrow (1\,1\,1\,1\,1\,1\,1\,1\,1\,1) &&= Z_3 \text{ (say)} \\
Z_3E &\hookrightarrow (1\,1\,1\,1\,1\,1\,1\,1\,1\,1) &&= Z_3.
\end{aligned}
$$

Just the effect of untouchability is highly cruel in the context of discrimination in education. The ON state of Samadharma is the only positive attribute which has forced them to fight against both the imposition of untouchability and the discrimination in education. Since all the concepts/ nodes given by the experts are tightly interlinked, we see that the ON state of any of the node makes all the other nodes to ON state, leading to a fixed point. Moreover, it is important to note that in no place in the model is the hidden pattern a limit cycle. In the expert's opinion a recurrence of nodes as a limit cycle is impossible.

The reader may be under the assumption that the talk of discrimination in education is an exaggeration. However, this is not so. Even in the days of Periyar, in the Presidency College in Chennai, Mr. Namachivaya, a teacher of Tamil with several years experience could only earn a salary of Rs.80 per month. On the other hand, a Brahmin Kuppusamy Sastrigal, who was a Sanskrit teacher in the same college, and had lesser experience than the Tamil pundit would get Rs.400/- per month. That was one instance of the widespread discrimination faced by the Sudras and Dalits in the hands of the Brahmins. This is prevalent even today.

This expert was asked if any of the interrelations given in the directed graphs could be indeterminate. So we obtained the following neutrosophic directed graph from the expert.

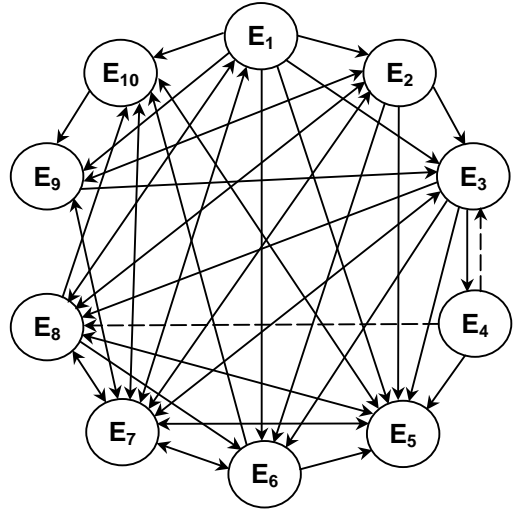

FIGURE 2.2.2

The related connection matrix of the neutrosophic directed graph is as follows:

$$E = \begin{array}{c} \\ E_1 \\ E_2 \\ E_3 \\ E_4 \\ E_5 \\ E_6 \\ E_7 \\ E_8 \\ E_9 \\ E_{10} \end{array} \begin{array}{c} E_1\ E_2\ E_3\ E_4\ E_5\ E_6\ E_7\ E_8\ E_9\ E_{10} \\ \begin{bmatrix} 0 & 1 & 1 & 0 & 1 & 1 & 1 & 1 & 1 & 1 \\ 0 & 0 & 1 & 0 & 1 & 1 & 1 & 1 & 1 & 0 \\ 0 & 0 & 0 & 1 & 1 & 1 & 1 & 0 & 0 & 0 \\ 0 & 0 & I & 0 & -1 & 0 & 0 & I & 0 & 0 \\ 0 & 0 & 0 & 0 & 0 & 0 & 1 & 1 & 0 & 1 \\ 0 & 0 & 0 & 0 & 1 & 0 & 1 & 0 & 0 & 1 \\ 1 & 1 & 1 & 0 & 1 & 1 & 0 & 1 & 1 & 1 \\ 1 & 1 & 0 & 0 & 1 & 1 & 1 & 0 & 0 & 1 \\ 0 & 1 & 1 & 0 & 0 & 0 & 1 & 0 & 0 & 0 \\ 0 & 0 & 0 & 0 & 1 & 0 & 1 & 0 & 1 & 0 \end{bmatrix} \end{array}$$

It can easily be seen that if the node Samadharma is in ON state, all other nodes are in the OFF state. We see using Neutrosophic Cognitive Maps several of the relations are indeterminate, it is in fact acceptable and in keeping with the ON state of the node $E_4$.

Next we study another expert's opinion. It was more diverse and had greater number of attributes dealing with untouchability practiced in education. The attributes given by this expert are:

1. Discrimination in Education
2. Untouchability
3. Role of Brahmins as heads of most institutions of importance and higher learning
4. Certain fields of education are unapproachable for Dalits and Sudras
5. Differences practiced in delivering education (discrimination)
6. Discrimination in school education
7. Discrimination in college education
8. Discrimination in law education
9. Discrimination in medical education
10. Discrimination in arts like music and dance
11. Discrimination in acting/ drama schools
12. Research (technical and scientific) devoid of any Dalit and Sudra representation
13. Research – Arts – place of Dalits and Sudras
14. Defence research
15. Colour discrimination in Tamil Nadu
16. Merit scholarship
17. Discrimination in food and accommodation in educational institutions
18. Educational marking system
19. Plight of Tamil teachers vis-à-vis Hindi/ Sanskrit teachers
20. State of physical education
21. Medium of instruction
22. Employment in educational institutions of national importance

23. Fear of learning due to discrimination and ill treatment
24. Partiality or discrimination in valuation of Dalits and Sudras answer scripts
25. Role of present day government

This expert is an educationalist who had been in the field of education/teaching and administration for over 5 decades. She is also a Periyarist and held discussions with us for many days on this topic of untouchability in education. During these discussions, astonishingly, she gave us over 30 attributes. We told her our inability to work with such a big number as 30. She hesitatingly reduced it to 25; we felt it impossible to argue with her or make compromise over our views. She offered to give reasons, so we had no other option other than to accept her views.

 Many of the attributes have already been described; hence we proceed on to describe the new nodes in a brief manner to avoid repetition.

 She says that Brahmins occupy most higher, powerful posts and hence oppress Dalits and Sudras. Brahmins grossly dominate certain fields, so it is unapproachable for the oppressed sections of society.

 She explains that 'differences are practiced in delivering education' in the following way: Education delivered in a primary school in a rural area (where majority of the students are from Dalit and backward communities) is far low in standard when compared to the education given in posh urban schools, where the children of elites and upper castes study. This is more emphasized in the case of elementary and secondary education. In rural schools, there are no proper teachers, no space, no infrastructure and above all, no motivation for education is provided to the Dalit and Sudra children in contrast with the city schools. Above all the medium of instruction for all these children is only Tamil. Unfortunately, these rural Dalit and Sudra children have to compete in public examinations with the well-educated Brahmin and upper caste and answer the same set of questions at the $10^{th}$ and $12^{th}$ standards.

At times, Dalit and Sudra students find it difficult to enter any professional course: engineering, medicine, law or architecture. At times, even if they qualify for a professional course by sheer merit and hard work, they are unable to pay the high fees. Also, in some cases, they student find it very difficult to cope with the English medium of instruction. Another sad thing to be noted is that these students are ill-treated and harassed not only by their teachers for their ignorance in English but also by their own classmates. The teachers are also indifferent to the problems faced by students from these marginalized backgrounds, and they don't try to give them a proper orientation.

According to our expert, this is also prevalent in every vocation and occupation such as law, medicine, acting and performing arts, science and technology, etc. She spoke about the cultural richness of these subjugated communities and how the Brahmins appropriated it.

In institutes of higher education, there is no reservation for Dalits and Sudras for doing research. Even granted there is some arbitrary reservation, the harassment they are made to undergo is so terrible that some of the students don't complete the course in the stipulated period; some do research for over several years together, become averse and quit. The reason behind this is absence of solid representation of Sudras and Dalits as teachers or research guides in these institutes.

In Tamil Nadu, the Dalit and Sudra students are also easily picked up and are made victims of discrimination, because majority of them are dark complexioned. So colour stands as an easy recognition of caste by the upper castes and Brahmins. So they are discriminated at the places of learning. This expert states that there is partiality in the valuation of answer scripts of Dalit and Sudra students.

Finally, the expert also pointed out that none of the Government policies of education had succeeded in reaching out to the lower strata of society.

Since this expert insisted on the 25 attributes and was against lessening it, working with a $25 \times 25$ matrix and a graph with 25 nodes was very complicated. So we took the links and worked with it using the pseudo code in C language [69, 76]. The data given by her was fed in the form of the program and several results were derived. These are ingrained in our chapter on conclusions.

As our expert quickly grasped the concept of NCM, she obliged us by working with NCMs. We noted the related matrix. When the state vector 'discrimination in school education' alone was in the ON state we found the resultant had four nodes viz. 10, 14, 16 and 19 were in the OFF state and the nodes 15, 20, 23 and 24 were in the indeterminate state. The rest of the nodes came to ON state.

### 2.3 Social inequality faced by Dalits and some of the most backward classes - an analysis using FCM and NCM

The third major attribute taken for our analysis is social untouchability. Everywhere in India, the Dalits are made to live in separate settlements called the *cheri* or *basti* on the outskirts of the village *oor* or *gram*. So, India itself is divided. This concept has been described earlier in p.48 of this book. Social inequality is also gruesomely displayed in caste riots and atrocities. According to a ten-year-old statistics from the Government of India: "On an average day, two Dalits are killed, three Dalit women are raped, two Dalit houses are burned and 50 Dalits are assaulted by caste-Hindus." This is just the tip of the iceberg. The actual numbers are indeed frightening. The nodes given by the first expert as important attributes are

- $T_1$ - Social inequality
- $T_2$ - Social discrimination
- $T_3$ - Absence of conscience and humaneness in caste-Hindus and Brahmins
- $T_4$ - Manu Dharma
- $T_5$ - Varnashrama Dharma
- $T_6$ - Brahmin/ caste-Hindu arrogance
- $T_7$ - Law is not favorable to the Dalits/ Sudras.
- $T_8$ - Dalits/ Sudras are not policy makers.

| | | |
|---|---|---|
| $T_9$ | - | Power is not in the hands of Dalits/ Sudras |
| $T_{10}$ | - | Fear of police atrocities |
| $T_{11}$ | - | Police join hands with power/ caste Hindus. |
| $T_{12}$ | - | Fear of losing their means of livelihood (employment as coolies in farms) |
| $T_{13}$ | - | Fear of further harassment (denial of water, provision, milk, etc. and social boycott) |
| $T_{14}$ | - | No education / education denied |
| $T_{15}$ | - | Fear of being arrested and put in prison on false charges by upper caste and Brahmins/ ruling party |
| $T_{16}$ | - | No support from the ruling party |
| $T_{17}$ | - | Poor economic condition |

This expert has given these 17 attributes as vital reasons for the social untouchability. The related directed graph is:

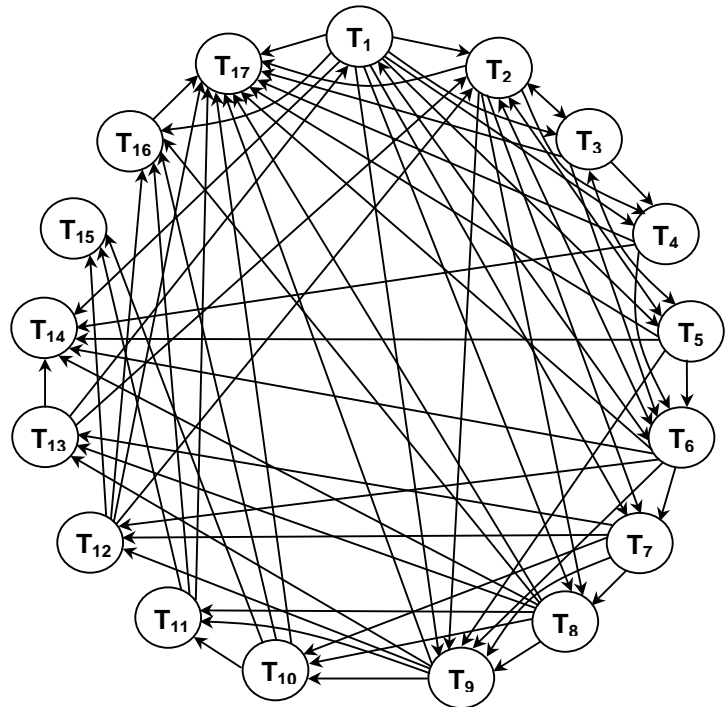

FIGURE 2.3.1

The related matrix is given by

$$\begin{array}{c} \\ T_1 \\ T_2 \\ T_3 \\ T_4 \\ T_5 \\ T_6 \\ T_7 \\ T_8 \\ T_9 \\ T_{10} \\ T_{11} \\ T_{12} \\ T_{13} \\ T_{14} \\ T_{15} \\ T_{16} \\ T_{17} \end{array} \begin{array}{c} T_1\ T_2\ T_3\ T_4\ T_5\ T_6\ T_7\ T_8\ T_9\ T_{10}\ T_{11}\ T_{12}\ T_{13}\ T_{14}\ T_{15}\ T_{16}\ T_{17} \\ \left[\begin{array}{ccccccccccccccccc} 0 & 1 & 1 & 1 & 1 & 1 & 1 & 1 & 1 & 0 & 0 & 0 & 1 & 0 & 1 & 1 \\ 0 & 0 & 1 & 1 & 1 & 1 & 1 & 1 & 1 & 0 & 0 & 0 & 0 & 0 & 0 & 1 \\ 0 & 1 & 0 & 1 & 1 & 1 & 0 & 0 & 0 & 0 & 0 & 0 & 0 & 0 & 0 & 1 \\ 0 & 1 & 0 & 0 & 0 & 1 & 0 & 0 & 0 & 0 & 0 & 0 & 1 & 0 & 0 & 1 \\ 0 & 1 & 1 & 0 & 0 & 1 & 0 & 0 & 1 & 0 & 0 & 0 & 1 & 0 & 0 & 1 \\ 1 & 1 & 0 & 0 & 0 & 0 & 1 & 0 & 1 & 0 & 0 & 1 & 0 & 1 & 0 & 0 & 1 \\ 0 & 0 & 0 & 0 & 0 & 0 & 0 & 1 & 1 & 1 & 0 & 1 & 1 & 0 & 0 & 0 & 0 \\ 0 & 0 & 0 & 0 & 0 & 0 & 0 & 1 & 1 & 1 & 0 & 1 & 1 & 0 & 1 & 1 \\ 0 & 0 & 0 & 0 & 0 & 0 & 0 & 0 & 1 & 1 & 1 & 1 & 0 & 0 & 0 & 1 \\ 0 & 0 & 0 & 0 & 0 & 0 & 0 & 0 & 0 & 1 & 0 & 0 & 0 & 1 & 1 & 1 \\ 0 & 0 & 0 & 0 & 0 & 0 & 0 & 0 & 0 & 0 & 0 & 0 & 0 & 1 & 1 & 1 \\ 0 & 1 & 0 & 0 & 0 & 0 & 0 & 0 & 0 & 0 & 0 & 0 & 0 & 1 & 1 & 1 \\ 1 & 1 & 0 & 0 & 0 & 0 & 0 & 0 & 0 & 0 & 0 & 0 & 1 & 0 & 0 & 0 \\ 0 & 0 & 0 & 0 & 0 & 0 & 0 & 0 & 0 & 0 & 0 & 0 & 0 & 0 & 0 & 0 \\ 0 & 0 & 0 & 0 & 0 & 0 & 0 & 0 & 0 & 0 & 0 & 0 & 0 & 0 & 0 & 0 \\ 0 & 0 & 0 & 0 & 0 & 0 & 0 & 0 & 0 & 0 & 0 & 0 & 0 & 0 & 0 & 1 \\ 0 & 0 & 0 & 0 & 0 & 0 & 0 & 0 & 0 & 0 & 0 & 0 & 0 & 0 & 0 & 0 \end{array}\right] \end{array}$$

Let S denote a $17 \times 17$ matrix, the effect of $X = (1\ 0\ 0\ 0\ \ldots\ 0)$ i.e., only the node 'social inequality' is in the ON state and all other nodes are OFF.

To find the effect of X on the dynamical system S, consider

$$\begin{array}{lll} XS & \hookrightarrow & (1\ 1\ 1\ 1\ 1\ 1\ 1\ 1\ 1\ 0\ 0\ 0\ 0\ 1\ 0\ 1\ 1) = X_1\ (\text{say}) \\ X_1 S & \hookrightarrow & (1\ 1\ 1\ 1\ 1\ 1\ 1\ 1\ 1\ 1\ 1\ 1\ 1\ 1\ 0\ 1\ 1) = X_2\ (\text{say}) \\ X_2 S & \hookrightarrow & (1\ 1\ 1\ 1\ 1\ 1\ 1\ 1\ 1\ 1\ 1\ 1\ 1\ 1\ 1\ 1\ 1) = X_3 = (X_2). \end{array}$$

Thus the hidden pattern results in a fixed point. This shows that if there is social inequality in the national set-up all other nodes become ON signifying the role of Manu Dharma, untouchability so on and so forth.

We take state vector Y = (0 0 0 0 0 1 0 0 0 0 0 0 0 0 0) i.e., only the node 'Brahmin or caste Hindu arrogance' is in the ON state and all other nodes are OFF, we see the effect of Y on the dynamical system S.

YS  $\hookrightarrow$  (1 1 0 0 0 1 1 0 1 0 0 1 0 1 0 0 1) = $Y_1$ (say)
$Y_1$S  $\hookrightarrow$  (1 1 1 1 1 1 1 1 1 1 1 1 1 1 1) = $Y_2$ (say)
$Y_2$S  $\hookrightarrow$  (1 1 1 1 1 1 1 1 1 1 1 1 1 1 1) = $Y_2$ .

This also gives a fixed point as hidden pattern. Thus caste Hindu/ Brahmin arrogance is sufficient to ruin the society entailing in social inequality, social discrimination and so on and so forth. Unless the oppressor castes and Brahmins correct themselves it would certainly ruin the nation. Power and money are concentrated in the hands of a few elites belonging to Brahmin and 'upper' castes. Over 70% of the people are Dalits and Sudras who suffer due to poverty, landlessness, and inferior social status.

Suppose the node 'poor economic conditions of the majority of Dalits and Sudras' and the node that they are 'not policy makers' is in the 'ON' state, we study the effect of it on the dynamical system i.e. on our society; even if all other nodes are in the OFF state.

Let P = (0 0 0 0 0 0 0 1 0 0 0 0 0 0 1). The effect of P on the dynamical system is:

PS   $\hookrightarrow$  (0 0 0 0 0 0 0 1 1 1 1 0 1 1 0 1 1)   =   $P_1$ say
$P_1$S  $\hookrightarrow$  (1 1 0 0 0 0 0 1 1 1 1 1 1 1 1 1 1)   =   $P_2$ say
$P_2$S  $\hookrightarrow$  (1 1 1 1 1 1 1 1 1 1 1 1 1 1 1)   =   $P_3$ say
$P_3$S  $\hookrightarrow$  (1 1 1 1 1 1 1 1 1 1 1 1 1 1 1)   =   $P_3$ .

Thus the hidden pattern is a fixed point, which makes all states ON in the resultant state vector. This shows that the attribute that they are not policy makers coupled with poor economic conditions makes all other states ON showing once again the importance of the fact that unless they improve economically and become policy makers it would be impossible to see any radical change in their lifestyle.

To be more precise, their social and economic conditions are spiraling downward with the advent of globalization and privatization. Several other aspects were analyzed using FCM and we arrived at conclusions given in the final chapter of this book.

We informed the expert that now we give him the option to introduce indeterminates in his associations or relations. Consequently, some relationships were marked as being indeterminate.

The following neutrosophic graph, represented by figure 2.3.2, was obtained:

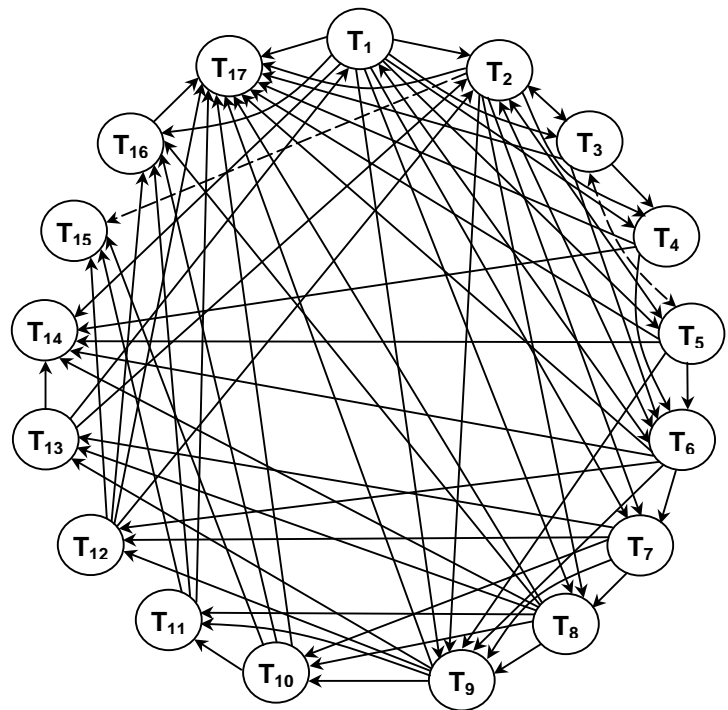

FIGURE 2.3.2

The related neutrosophic matrix is given in the next page.

$$\begin{array}{c}\phantom{T_1}\;T_1\;T_2\;T_3\;T_4\;T_5\;T_6\;T_7\;T_8\;T_9\;T_{10}\;T_{11}\;T_{12}\;T_{13}\;T_{14}\;T_{15}\;T_{16}\;T_{17}\end{array}$$

|   | $T_1$ | $T_2$ | $T_3$ | $T_4$ | $T_5$ | $T_6$ | $T_7$ | $T_8$ | $T_9$ | $T_{10}$ | $T_{11}$ | $T_{12}$ | $T_{13}$ | $T_{14}$ | $T_{15}$ | $T_{16}$ | $T_{17}$ |
|---|---|---|---|---|---|---|---|---|---|---|---|---|---|---|---|---|---|
| $T_1$ | 0 | 1 | 1 | 1 | 1 | 1 | 1 | 1 | 1 | 0 | 0 | 0 | 0 | 1 | 0 | 1 | 1 |
| $T_2$ | 0 | 0 | 1 | 1 | 1 | 1 | 1 | 1 | 1 | 0 | 0 | 0 | 0 | 0 | 0 | $I$ | 1 |
| $T_3$ | 0 | 1 | 0 | 1 | $I$ | 1 | 0 | 0 | 0 | 0 | 0 | 0 | 0 | 0 | 0 | 0 | 1 |
| $T_4$ | 0 | 1 | 0 | 0 | 0 | 1 | 0 | 0 | 0 | 0 | 0 | 0 | 0 | 1 | 0 | 0 | 1 |
| $T_5$ | 0 | 1 | $I$ | 0 | 0 | 1 | 0 | 0 | 1 | 0 | 0 | 0 | 0 | 1 | 0 | 0 | 1 |
| $T_6$ | 1 | 1 | 0 | 0 | 0 | 0 | 1 | 0 | 1 | 0 | 0 | 1 | 0 | 1 | 0 | 0 | 1 |
| $T_7$ | 0 | 0 | 0 | 0 | 0 | 0 | 0 | 1 | 1 | 1 | 0 | 1 | 1 | 0 | 0 | 0 | 0 |
| $T_8$ | 0 | 0 | 0 | 0 | 0 | 0 | 0 | 0 | 1 | 1 | 1 | 0 | 1 | 1 | 0 | 1 | 1 |
| $T_9$ | 0 | 0 | 0 | 0 | 0 | 0 | 0 | 0 | 0 | 1 | 1 | 1 | 1 | 0 | 0 | 0 | 1 |
| $T_{10}$ | 0 | 0 | 0 | 0 | 0 | 0 | 0 | 0 | 0 | 0 | 1 | 0 | 0 | 0 | 1 | 1 | 1 |
| $T_{11}$ | 0 | 0 | 0 | 0 | 0 | 0 | 0 | 0 | 0 | 0 | 0 | 0 | 0 | 0 | 1 | 1 | 1 |
| $T_{12}$ | 0 | 1 | 0 | 0 | 0 | 0 | 0 | 0 | 0 | 0 | 0 | 0 | 0 | 1 | 1 | 1 | 1 |
| $T_{13}$ | 1 | 1 | 0 | 0 | 0 | 0 | 0 | 0 | 0 | 0 | 0 | 0 | 0 | 1 | 0 | 0 | 1 |
| $T_{14}$ | 0 | 0 | 0 | 0 | 0 | 0 | 0 | 0 | 0 | 0 | 0 | 0 | 0 | 0 | 0 | 0 | 0 |
| $T_{15}$ | 0 | 0 | 0 | 0 | 0 | 0 | 0 | 0 | 0 | 0 | 0 | 0 | 0 | 0 | 0 | 0 | 0 |
| $T_{16}$ | 0 | 0 | 0 | 0 | 0 | 0 | 0 | 0 | 0 | 0 | 0 | 0 | 0 | 0 | 0 | 0 | 1 |
| $T_{17}$ | 0 | 0 | 0 | 0 | 0 | 0 | 0 | 0 | 0 | 0 | 0 | 0 | 0 | 0 | 0 | 0 | 0 |

The resultant shows $T_3$ and $T_6$ are indeterminates. When node $T_7$ alone is in ON state, in the resultant vector $T_{14}$ is in indeterminate state and all other nodes become ON.

Next, we get the second expert's opinion on the same topic. The nodes given by this expert are as follows:

$K_1$ - Manu Dharma
$K_2$ - Varnashrama Dharma
$K_3$ - Monopoly of power by Brahmins and caste Hindus
$K_4$ - Religious belief
$K_5$ - Social discrimination
$K_6$ - Caste discrimination
$K_7$ - Economic condition
$K_8$ - No access to power
$K_9$ - Government has no concern about Dalits /Sudras; has concern about privatization and globalization
$K_{10}$ - Social disrespect based on their caste

$K_{11}$ - Fear of losing livelihood / job
$K_{12}$ - Judiciary monopolized by Brahmins / caste Hindu
$K_{13}$ - Government reforms do not reach grassroots
$K_{14}$ - Landless labourers
$K_{15}$ - No facility for education / economic improvement

This expert feels that socio-economic status of Dalits and Sudras can never be improved until they become policy makers. Their status will remain the same; Brahmins and 'upper' caste Hindus and the government is indirectly the cause of such treachery.

This expert has given a very dense directed graph, which as follows

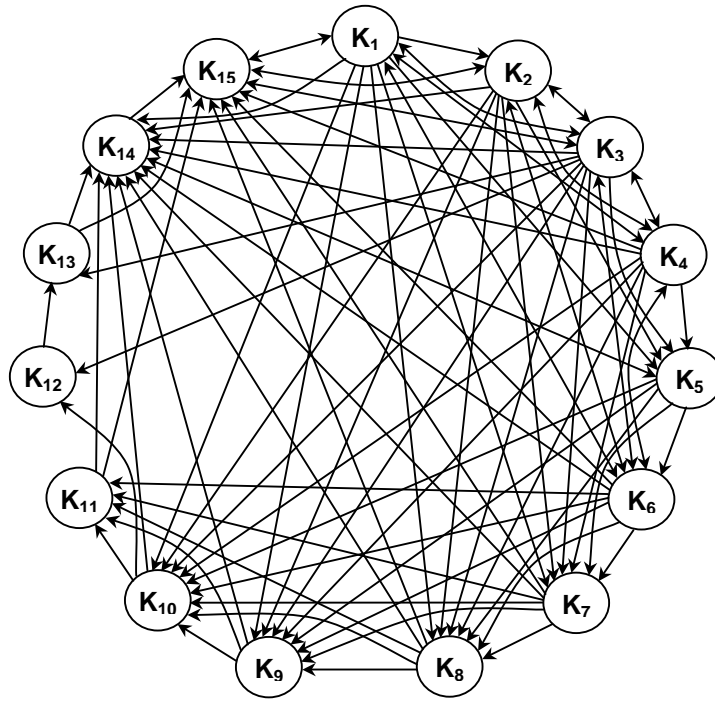

FIGURE 2.3.3

The related connection matrix is a $15 \times 15$ matrix which we denote by L.

$$\begin{array}{c} \phantom{K_{1}} \quad K_1 \ K_2 \ K_3 \ K_4 \ K_5 \ K_6 \ K_7 \ K_8 \ K_9 \ K_{10} \ K_{11} K_{12} K_{13} \ K_{14} K_{15} \\ \begin{array}{c} K_1 \\ K_2 \\ K_3 \\ K_4 \\ K_5 \\ K_6 \\ K_7 \\ K_8 \\ K_9 \\ K_{10} \\ K_{11} \\ K_{12} \\ K_{13} \\ K_{14} \\ K_{15} \end{array} \left[ \begin{array}{ccccccccccccccc} 0 & 1 & 1 & 1 & 1 & 1 & 1 & 1 & 1 & 1 & 0 & 0 & 0 & 1 & 1 \\ 0 & 0 & 1 & 1 & 1 & 1 & 1 & 1 & 1 & 1 & 0 & 0 & 0 & 1 & 1 \\ 1 & 1 & 0 & 1 & 1 & 1 & 1 & 1 & 1 & 1 & 0 & 1 & 1 & 1 & 1 \\ 1 & 1 & 1 & 0 & 1 & 1 & 1 & 1 & 1 & 1 & 0 & 0 & 0 & 1 & 1 \\ 0 & 0 & 0 & 0 & 0 & 1 & 1 & 1 & 1 & 1 & 0 & 0 & 0 & 1 & 0 \\ 1 & 1 & 1 & 1 & 0 & 0 & 1 & 1 & 1 & 1 & 0 & 0 & 1 & 1 & 1 \\ 0 & 0 & 0 & 0 & 0 & 0 & 0 & 1 & 1 & 1 & 0 & 0 & 1 & 1 \\ 0 & 0 & 0 & 0 & 0 & 0 & 0 & 0 & 1 & 1 & 1 & 0 & 0 & 1 & 1 \\ 0 & 0 & 0 & 0 & 0 & 0 & 0 & 0 & 1 & 1 & 0 & 0 & 1 & 0 \\ 0 & 0 & 0 & 0 & 0 & 0 & 0 & 0 & 0 & 1 & 1 & 0 & 1 & 0 \\ 0 & 0 & 0 & 0 & 0 & 0 & 0 & 0 & 0 & 0 & 0 & 0 & 1 & 1 \\ 0 & 0 & 0 & 0 & 0 & 0 & 0 & 0 & 0 & 0 & 0 & 1 & 0 & 0 \\ 0 & 0 & 0 & 0 & 0 & 0 & 0 & 0 & 0 & 0 & 0 & 0 & 1 & 1 \\ 0 & 0 & 0 & 0 & 0 & 0 & 0 & 0 & 0 & 0 & 0 & 0 & 0 & 1 \\ 1 & 1 & 1 & 0 & 0 & 0 & 0 & 0 & 0 & 0 & 0 & 0 & 0 & 0 \end{array} \right] \end{array}$$

We study the effect of any state vector on the dynamical system using L. Let X = (0 0 0 0 1 0 0 0 0 0 0 0 0 0 0) be a state vector in which the node 'social discrimination' is in ON state and all other nodes are in OFF state.

| | | | |
|---|---|---|---|
| XL | ↪ (0 0 0 0 1 1 1 1 1 1 0 0 0 1 0) = | $X_1$ (say) |
| $X_1$ L | ↪ (1 1 1 1 1 1 1 1 1 1 1 1 1 0 1 1) = | $X_2$ say |
| $X_2$ L | ↪ (1 1 1 1 1 1 1 1 1 1 1 1 1 1 1) = | $X_2$ (say) |
| $X_2$L | ↪ (1 1 1 1 1 1 1 1 1 1 1 1 1 1 1) = | $X_3$. |

The hidden pattern is a fixed point in which all attributes are in the ON state. It is clear that the attribute social discrimination has forced all the attributes to become ON. The main points are: government has no concern over their living conditions, no budgetary plans ever reach them and so on. Each factor represents negative effect of social discrimination practiced even today in full swing.

Let Y be (0 0 0 0 0 0 0 1 0 0 0 0 0 0 0) where only the node 'no access to power' is in the ON state and all nodes are in the OFF state. Therefore, the effect of Y on the dynamical system is given by

YL ↪ (0 0 0 0 0 0 0 1 1 1 1 0 0 1 1) = $Y_1$ (say)
$Y_1$L ↪ (1 1 1 0 0 0 0 1 1 1 1 1 0 1 1) = $Y_2$ (say)
$Y_2$L ↪ (1 1 1 1 1 1 1 1 1 1 1 1 1 1 1) = $Y_3$ (say)
$Y_3$L ↪ $Y_3$.

The resulting hidden pattern is a fixed point. If majority of the population viz. Sudras and Dalits have no access to power, all the nodes become ON. It is evident that the marginalized societies have been kept out of power sharing. Several such conclusions can be made on this study, which is listed in the final chapter.

Now when we permit the expert to give indeterminate relations he has given the following neutrosophic graph.

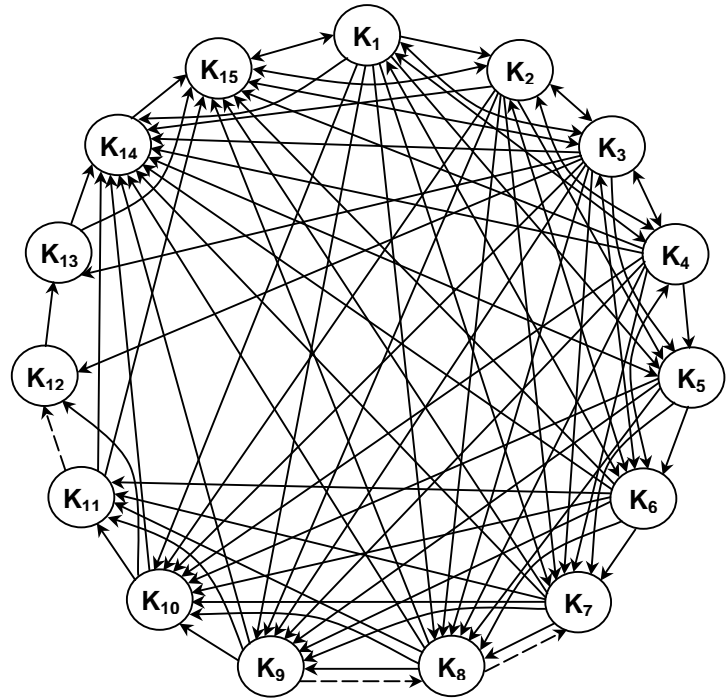

FIGURE 2.3.4

Using the neutrosophic graph we determine the following relational matrix N(E).

$$
\begin{array}{c}
\phantom{K_{1}}\;\;K_1\;K_2\;K_3\;K_4\;K_5\;K_6\;K_7\;K_8\;K_9\;K_{10}\;K_{11}\;K_{12}\;K_{13}\;K_{14}\;K_{15}\\
\begin{array}{c}K_1\\K_2\\K_3\\K_4\\K_5\\K_6\\K_7\\K_8\\K_9\\K_{10}\\K_{11}\\K_{12}\\K_{13}\\K_{14}\\K_{15}\end{array}
\left[\begin{array}{ccccccccccccccc}
0 & 1 & 1 & 1 & 1 & 1 & 1 & 1 & 1 & 0 & 0 & 0 & 0 & 1 & 1\\
0 & 0 & 1 & 1 & 1 & 1 & 1 & 1 & 1 & 1 & 0 & 0 & 0 & 1 & 1\\
1 & 1 & 0 & 1 & 1 & 1 & 1 & 1 & 1 & 1 & 0 & 1 & 1 & 1 & 1\\
1 & 1 & 1 & 0 & 1 & 1 & 1 & 1 & 1 & 1 & 0 & 0 & 0 & 1 & 1\\
0 & 0 & 0 & 0 & 0 & 1 & 1 & 1 & 1 & 1 & 0 & 0 & 0 & 1 & 0\\
1 & 1 & 1 & 1 & 0 & 0 & 1 & 1 & 1 & 1 & 1 & 0 & 0 & 1 & 1\\
0 & 0 & 0 & 0 & 0 & 0 & 0 & 1 & 1 & 1 & 1 & 0 & 0 & 1 & 1\\
0 & 0 & 0 & 0 & 0 & 0 & I & 0 & 1 & 1 & 1 & 0 & 0 & 1 & 1\\
0 & 0 & 0 & 0 & 0 & 0 & 0 & I & 0 & 1 & 1 & 0 & 0 & 1 & 0\\
0 & 0 & 0 & 0 & 0 & 0 & 0 & 0 & 0 & 0 & 1 & 1 & 0 & 1 & 0\\
0 & 0 & 0 & 0 & 0 & 0 & 0 & 0 & 0 & 0 & 0 & I & 0 & 1 & 1\\
0 & 0 & 0 & 0 & 0 & 0 & 0 & 0 & 0 & 0 & 0 & 0 & 1 & 0 & 0\\
0 & 0 & 0 & 0 & 0 & 0 & 0 & 0 & 0 & 0 & 0 & 0 & 0 & 1 & 1\\
0 & 0 & 0 & 0 & 0 & 0 & 0 & 0 & 0 & 0 & 0 & 0 & 0 & 0 & 1\\
1 & 1 & 1 & 0 & 0 & 0 & 0 & 0 & 0 & 0 & 0 & 0 & 0 & 0 & 0
\end{array}\right]
\end{array}
$$

Using this neutrosophic matrix N(E) the effect of several state vectors on N(E) is determined. We proceed on to the next section.

### 2.4 Problems faced by Dalits in the political arena due to discrimination – a FCM and NCM analysis

We take the next major attribute, political untouchability. Although reservation exists for Dalits at every tier of the political systems, right from membership in local village government (*panchayat*) to the parliament. Because it is a constitutional obligation, this reservation is mostly filled. However real assertive candidates fail to get elected and the dummy dalit candidates out up by the major political parties end up holding public office. Also, electoral violence against Dalit candidates is rampant and widespread. Those who manage to get elected are also either puppets in the hands of caste Hindu forces, or they are finished off. Worst instance of caste fanaticism against Dalit political assertion was witnessed in the murder of Melavalavu Murugesan and seven others on 30 June 1997, because they violated a caste Hindu order that they must not contest the post of Panchayat President. All over India, in several places elections in grassroots level bodies (*panchayat*) are not allowed to take place because these seats have been reserved for dalits (like in Pappapatti and Keeripatti). In fact, in the Tamil Nadu State cabinet, there is only one Dalit minister and he is in charge of dalit welfare. Dalits have never been allocated important portfolios, a clear evidence of discrimination because they are the largest vote bank: one in five Tamilians is a dalit.

We list the attributes given by the first expert. He says the political untouchability is closely connected with state sponsored atrocities i.e. "police cruelty" against Dalits due to caste.

$A_1$   -   State sponsored atrocities
$A_2$   -   Political murder
$A_3$   -   Political discrimination
$A_4$   -   Role of ruling party
$A_5$   -   Role of caste Hindus / Brahmins
$A_6$   -   Role of judiciary
$A_7$   -   Role of caste majority in that village
$A_8$   -   Role of untouchability
$A_9$   -   Influence of Caste Hindu power groups
$A_{10}$   -   Economic conditions of Dalits
$A_{11}$   -   Police discrimination
$A_{12}$   -   Insecurity to life and belongings of Dalits
$A_{13}$   -   Fear of holding political post due to Life threat
$A_{14}$   -   Torture by caste Hindus and Brahmins
$A_{15}$   -   No co-operation from government and police
$A_{16}$   -   Rule of Manu Dharma

It is pertinent to mention here that the very chair used by Dr. Ambedkar when he was the central minister was washed and purified with Ganga water. If this is the plight of the father of our Constitution one can imagine what will be the torture, harassments and discriminations invariably faced by poor uneducated Dalits.

Even granted a Dalit is elected, occupies an elected office he always faces a risk to his life, consequent of which he may be assassinated or brutally murdered on the streets. This is the deplorable plight of the Dalits who are properly elected. How can they make policies? Unless this situation is changed progress in the true sense cannot be achieved. So liberation has become very difficult for Dalits and Sudras in every sphere of life.

The directed graph given by the first expert is given below:

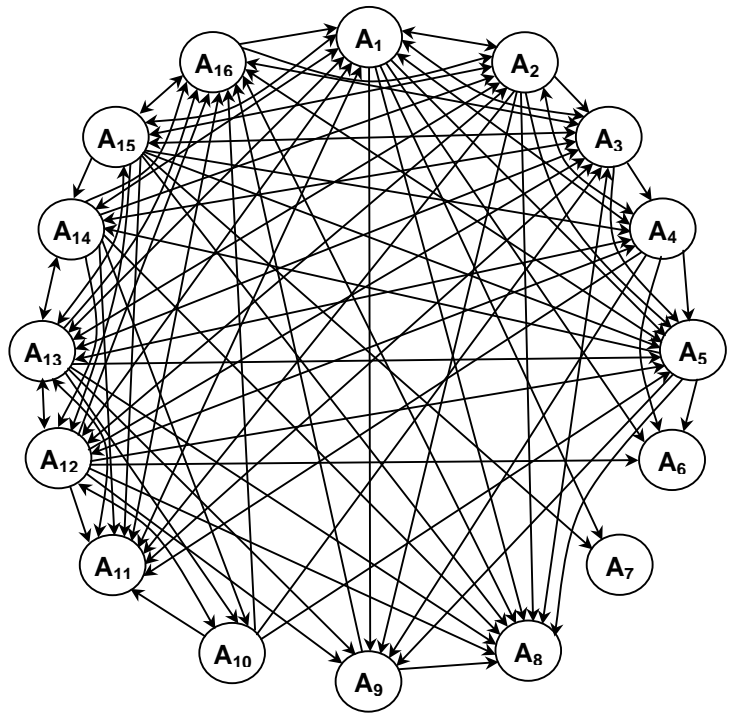

FIGURE 2.4.1

Now we give the related 16 × 16 matrix denoted by V.

$$
\begin{array}{c}
\phantom{A_{1}}\,A_{1}\,A_{2}\,A_{3}\,A_{4}\,A_{5}\,A_{6}\,A_{7}\,A_{8}\,A_{9}\,A_{10}\,A_{11}\,A_{12}\,A_{13}\,A_{14}\,A_{15}\,A_{16}\\
\begin{array}{c}A_{1}\\A_{2}\\A_{3}\\A_{4}\\A_{5}\\A_{6}\\A_{7}\\A_{8}\\A_{9}\\A_{10}\\A_{11}\\A_{12}\\A_{13}\\A_{14}\\A_{15}\\A_{16}\end{array}
\left[\begin{array}{cccccccccccccccc}
0 & 1 & 1 & 1 & 1 & 1 & 1 & 1 & 1 & 0 & 1 & 1 & 1 & 1 & 1 & 0\\
1 & 0 & 1 & 1 & 1 & 0 & 0 & 1 & 1 & 0 & 1 & 1 & 1 & 0 & 1 & 0\\
1 & 0 & 0 & 1 & 1 & 0 & 0 & 1 & 0 & 0 & 1 & 0 & 1 & 1 & 1 & 1\\
1 & 1 & 0 & 0 & 1 & 1 & 0 & 0 & 1 & 0 & 1 & 1 & 1 & 0 & 0 & 0\\
0 & 0 & 0 & 0 & 0 & 1 & 0 & 1 & 1 & 0 & 0 & 0 & 0 & 1 & 0 & 1\\
0 & 0 & 0 & 0 & 0 & 0 & 0 & 0 & 0 & 0 & 0 & 0 & 0 & 0 & 0 & 0\\
0 & 0 & 0 & 0 & 0 & 0 & 0 & 0 & 0 & 0 & 0 & 0 & 0 & 0 & 0 & 0\\
0 & 0 & 1 & 0 & 0 & 0 & 0 & 0 & 0 & 0 & 0 & 0 & 0 & 0 & 0 & 1\\
0 & 0 & 0 & 0 & 0 & 0 & 0 & 1 & 0 & 0 & 0 & 0 & 0 & 0 & 0 & 1\\
0 & 0 & 1 & 0 & 1 & 0 & 0 & 0 & 0 & 0 & 1 & 1 & 0 & 0 & 0 & 1\\
1 & 0 & 1 & 0 & 0 & 0 & 0 & 0 & 0 & 0 & 0 & 0 & 1 & 0 & 0 & 1\\
1 & 1 & 1 & 1 & 1 & 1 & 0 & 1 & 1 & 1 & 1 & 0 & 1 & 0 & 0 & 1\\
1 & 1 & 1 & 1 & 1 & 0 & 0 & 1 & 1 & 1 & 1 & 1 & 0 & 1 & 1 & 1\\
1 & 1 & 1 & 0 & 0 & 0 & 0 & 1 & 0 & 1 & 1 & 1 & 1 & 0 & 0 & 1\\
1 & 1 & 1 & 1 & 1 & 0 & 1 & 1 & 0 & 0 & 1 & 1 & 1 & 1 & 0 & 1\\
1 & 1 & 1 & 0 & 1 & 0 & 0 & 1 & 0 & 0 & 1 & 0 & 0 & 0 & 1 & 0
\end{array}\right]
\end{array}
$$

We consider the state vector $X = (1\ 0\ 0\ 0\ 0\ \ldots\ 0)$, where the node 'state-sponsored atrocities' is in the ON state and all other nodes are in the OFF state.

The effect of X on the dynamical system V is given by as follows

$XV \hookrightarrow (1\ 1\ 1\ 1\ 1\ 1\ 1\ 1\ 1\ 0\ 1\ 1\ 1\ 1\ 1\ 0) = X_1$ (say)
$X_1 V \hookrightarrow (1\ 1\ 1\ 1\ 1\ 1\ 1\ 1\ 1\ 1\ 1\ 1\ 1\ 1\ 1\ 1) = X_2$ (say)
$X_2 V \hookrightarrow (1\ 1\ 1\ 1\ 1\ 1\ 1\ 1\ 1\ 1\ 1\ 1\ 1\ 1\ 1\ 1) = X_3\ (= X_2)$.

The hidden pattern of the dynamical system is a fixed point and if 'state-sponsored atrocities' alone is in the ON state it makes all other nodes like political murder, political discrimination, role of ruling party, role of caste Hindus and Brahmins and so on to become ON. Hence 'state-sponsored atrocities' promotes political untouchability.

We take the state vector Y = (0 0 0 0 1 0 0 0 0 0 0 0 0 0 0 0) where the node 'rule of caste Hindus and Brahmins' alone is in the ON state. The effect of Y on the dynamical system V is

Y V ↪ (0 0 0 0 1 1 0 1 1 0 0 0 0 1 0 1) = $Y_1$ (say)
$Y_1$ V ↪ (1 1 1 0 1 1 0 1 1 1 1 1 1 1 1 1) = $Y_2$ (say)
$Y_2$ V ↪ (1 1 1 1 1 1 1 1 1 1 1 1 1 1 1 1) = $Y_3$ (say)
$Y_3$ V ↪ (1 1 1 1 1 1 1 1 1 1 1 1 1 1 1 1) = $Y_4$ = $Y_3$.

The hidden pattern results in a fixed point. We see that the effect of 'rule of caste Hindus and Brahmins' makes all nodes of the system under political untouchability come to ON state.

We obtained similar results when we studied with only the node $A_{11}$ in the ON state. This conclusion is supported based on documentary evidence. Violence against Dalits by the partisan state machinery, particularly the police are highly condemnable. Instigated by the ruling casteist forces, they police launch planned attacks and ransack Dalit settlements, foist false castes against Dalits, take the men into arbitrary custody, beat them, and torture the women and children. There are several such instances that have taken place in Tamil Nadu. For more, refer [65, 66].

When the expert was given the option to use NCM, i.e., the neutrosophic model, the 16 × 16 relational matrix was made into a neutrosophic matrix. When the node 'state-sponsored atrocities' was alone in the ON state and all other nodes were in the OFF state it was found that in the resultant vector, the nodes 2, 7 and 10 become indeterminate. Likewise several state vectors were analyzed using NCM.

Next we take another expert's opinion on political untouchability. The attributes/ nodes that the expert associated with the term "political untouchability" are:

$S_1$ - Political discrimination
$S_2$ - Political ill treatment
$S_3$ - Varnashrama Dharma
$S_4$ - Economic condition

$S_5$    -    Role of caste
$S_6$    -    Role of ruling party's caste votebank
$S_7$    -    Police atrocities on Dalits
$S_8$    -    Dalits prevented from voting
$S_9$    -    Dalits live under fear during election
$S_{10}$  -    Political violence

The expert explained $S_{10}$, the concept of political violence by speaking of the large scale electoral violence that takes place in the country. Several lives have been lost in such mindless, casteist violence. Dalit settlements are burnt during elections in order to terrify them and prevent them from casting their votes to a candidate of their choice. For instances, please refer [65, 66].

This expert differentiates the two concepts political discrimination and political ill treatment. For instance, if a Dalit who is elected as a Panchayat President faces a life threat and the police do not give him adequate security, it is an evidence of political discrimination. If an elected Dalit Panchayat President is not allowed to sit on a chair during meetings, and is instead asked to sit on the floor, or if he is abused by his caste name, then it is an evidence of political ill treatment.

Though we feel that political discrimination and ill treatment has very subtle difference, this expert wants to consider them as separate attributes. This expert views that the political discrimination is one of the causes for police atrocities on Dalits. If there was reasonable representation of Dalits in the ruling party certainly the state-sponsored police atrocities on these innocent Dalits can be prevented to a large extent.

The expert says if Dalits were economically powerful, certainly they will not become victims of such atrocities, and they would automatically find representation in power sharing and they would become politically empowered.

He says the role of caste is directly dependent on the ruling party's caste votebank. He is of the opinion that in India caste plays a major destructive role!

The directed graph given by the expert is as follows:

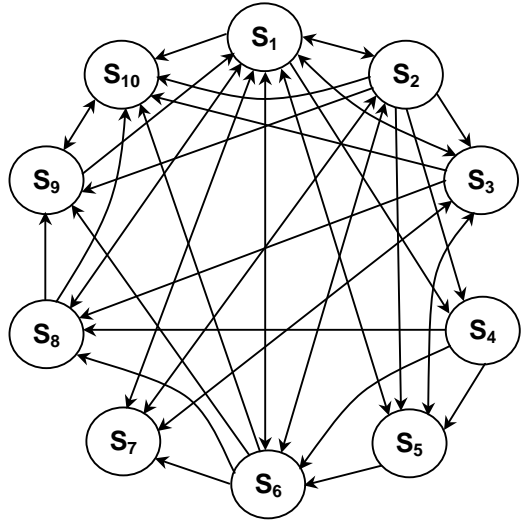

FIGURE 2.4.2

The connection matrix related to this directed graph is denoted by Q which is a $10 \times 10$ matrix.

$$Q = \begin{array}{c} \\ S_1 \\ S_2 \\ S_3 \\ S_4 \\ S_5 \\ S_6 \\ S_7 \\ S_8 \\ S_9 \\ S_{10} \end{array} \begin{array}{c} \begin{array}{cccccccccc} S_1 & S_2 & S_3 & S_4 & S_5 & S_6 & S_7 & S_8 & S_9 & S_{10} \end{array} \\ \left[ \begin{array}{cccccccccc} 0 & 1 & 1 & 1 & 1 & 1 & 1 & 1 & 0 & 1 \\ 1 & 0 & 1 & 1 & 1 & 1 & 1 & 0 & 1 & 1 \\ 1 & 0 & 0 & 0 & 1 & 0 & 1 & 1 & 0 & 1 \\ 0 & 0 & 0 & 0 & 1 & 1 & 0 & 1 & 0 & 0 \\ 1 & 0 & 1 & 0 & 0 & 1 & 0 & 0 & 0 & 0 \\ 1 & 1 & 0 & 0 & 0 & 0 & 1 & 1 & 1 & 1 \\ 0 & 0 & 0 & 0 & 0 & 0 & 0 & 1 & 1 & 1 \\ 1 & 0 & 0 & 0 & 0 & 0 & 0 & 0 & 1 & 1 \\ 1 & 0 & 0 & 0 & 0 & 0 & 0 & 0 & 0 & 1 \\ 0 & 0 & 0 & 0 & 0 & 0 & 0 & 0 & 1 & 0 \end{array} \right] \end{array}$$

We consider the state vector $X = (0\ 0\ 0\ 0\ 0\ 1\ 0\ 0\ 0\ 0)$ where the node 'role of ruling party's caste votebank' is in the ON state and all other attributes are in the OFF state. The effect of X on the dynamical system Q is given by

$X Q \hookrightarrow (1\ 1\ 0\ 0\ 0\ 1\ 1\ 1\ 1\ 1) = X_1$ (say)
$X_1 Q \hookrightarrow (1\ 1\ 1\ 1\ 1\ 1\ 1\ 1\ 1\ 1) = X_2$ (say)
$X_2 Q \hookrightarrow (1\ 1\ 1\ 1\ 1\ 1\ 1\ 1\ 1\ 1) = X_3 (= X_2$ say).

The hidden pattern is a fixed point in which all nodes are in the ON state. This shows that the role of ruling party's caste votebank plays a vital role in political untouchability.

We study the same dynamical system Q with another state vector $T = (0\ 0\ 0\ 0\ 0\ 0\ 0\ 0\ 1\ 0)$, where only the node 'Dalits live under fear during election time' is in the ON state and all other nodes are in the OFF state. The effect of Y on Q is given by

$Y Q \hookrightarrow (1\ 0\ 0\ 0\ 0\ 0\ 0\ 0\ 1\ 1) = Y_1$ (say)
$Y_1 Q \hookrightarrow (1\ 1\ 1\ 1\ 1\ 1\ 1\ 1\ 1\ 1) = Y_2$.

Clearly $Y_2 Q = Y_2$. All nodes are in the ON state in the resultant vector. This implies that when Dalits live under fear during elections, they are denied even the democratic right to vote for a candidate of their choice.

The expert says that invariably the ruling party rowdies, who take the situation near polling booth under their control, cast dalit votes. Thus political violence is the result of political untouchability. If the ruling party realizes this and keeps in mind that the nation is to be democratically ruled, certainly the concept of political discrimination and hence, political untouchability would be wiped off. The undue influence of religion, and particularly communalism in politics has been one of the factors that established caste as the prime player in politics.

This expert also suggested communal representation, as an alternative to the current existing system. He pointed out that this would limit the dominance by the various oppressor castes.

## 2.5 Study of Economic Status of Dalits due to untouchability using fuzzy and neutrosophic analysis

According to this expert the term economic untouchability means that Brahmins and caste-Hindus have prevented Dalits and Ssudras from gaining economic development because of practicing the religious taboo of untouchability. Since majority of Dalits live as landless agricultural labourers, they are economically downtrodden. They form the largest section of the unorganized sector of labourers. Their dismal representation in higher-income jobs is also reflective of their economic status. Dalits form the greatest share of those living below the poverty line. The expert also quoted Dr.B.R.Ambedkar to trace the history of economic subjugation of the Dalits.

Next, a major observation made by the expert was about the neo-economic policies of liberalization, globalization and privatization that has further robbed the Dalit people and left them much poorer.

He also described that bonded labour (rice-mills, brick-making and several small-scale industries etc.) still existed in several pockets of Tamil Nadu and almost all its victims were Dalits and most backward classes.

Economic boycott was a visible manifestation of economic untouchability. During caste riots in villages, the caste-Hindus refuse to hire Dalits to work on their lands, they refuse to sell provisions to Dalits, and thus utilize economic embargo as a tool of oppression.

Next, the expert began to criticize the concept of reservation for economically backward people? He says that it is clear-cut method to give reservation to the Brahmins. Over 95% of the poor people are only Dalits and Sudras. Hardly 0.001% of the Brahmin community could be living below the poverty line.

Thus according to this expert, economic untouchability means economic discrimination, economic backwardness, absence of economic self-sufficiency and so on. The expert feels that the major reason for economic untouchability lies in the four factors religious untouchability, political untouchability, social untouchability and educational untouchability.

The attributes given by this expert are:

- $F_1$ - Economic discrimination
- $F_2$ - Varnashrama Dharma
- $F_3$ - Caste system
- $F_4$ - Political power in the hands of caste Hindus
- $F_5$ - No proportion between labour and pay
- $F_6$ - State apathy towards their economic plight
- $F_7$ - Belief in destiny / karma
- $F_8$ - Lack of political power
- $F_9$ - Poor social status
- $F_{10}$ - No education
- $F_{11}$ - Policy makers are not Dalits/ Sudras
- $F_{12}$ - No planning for income generation activities

We now give the directed graph of this expert.

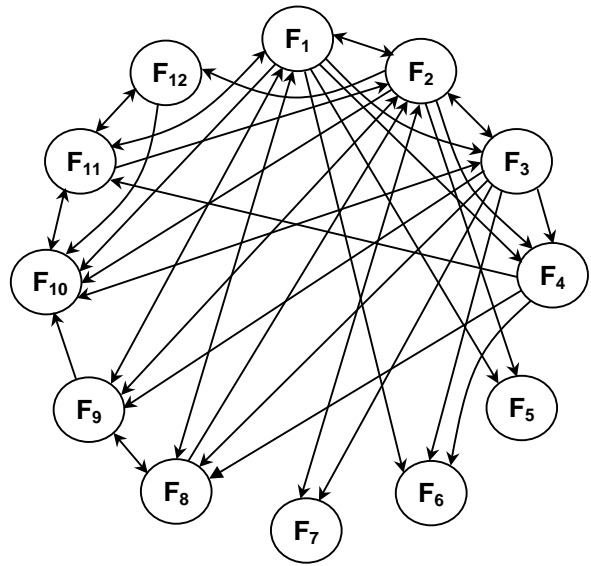

FIGURE 2.5.1

The related connection matrix of the directed graph is as follows:

$$\begin{array}{c} \\ F_1 \\ F_2 \\ F_3 \\ F_4 \\ F_5 \\ F_6 \\ F_7 \\ F_8 \\ F_9 \\ F_{10} \\ F_{11} \\ F_{12} \end{array} \begin{array}{c} F_1\ F_2\ F_3\ F_4\ F_5\ F_6\ F_7\ F_8\ F_9\ F_{10}\ F_{11}\ F_{12} \\ \begin{bmatrix} 0 & 1 & 1 & 1 & 1 & 1 & 0 & 1 & 1 & 1 & 1 & 0 \\ 1 & 0 & 1 & 1 & 1 & 0 & 1 & 0 & 1 & 1 & 0 & 1 \\ 0 & 1 & 0 & 1 & 0 & 1 & 1 & 1 & 1 & 1 & 0 & 0 \\ 0 & 0 & 0 & 0 & 0 & 1 & 0 & 1 & 0 & 0 & 1 & 0 \\ 0 & 0 & 0 & 0 & 0 & 0 & 0 & 0 & 0 & 0 & 0 & 0 \\ 0 & 0 & 0 & 0 & 0 & 0 & 0 & 0 & 0 & 0 & 0 & 0 \\ 0 & 1 & 0 & 0 & 0 & 0 & 0 & 0 & 0 & 0 & 0 & 0 \\ 1 & 1 & 0 & 0 & 0 & 0 & 0 & 0 & 1 & 0 & 0 & 0 \\ 1 & 1 & 0 & 0 & 0 & 0 & 0 & 1 & 0 & 1 & 0 & 0 \\ 0 & 0 & 1 & 0 & 0 & 0 & 0 & 0 & 0 & 0 & 1 & 0 \\ 1 & 1 & 0 & 0 & 0 & 0 & 0 & 0 & 0 & 1 & 0 & 1 \\ 0 & 0 & 0 & 0 & 0 & 0 & 0 & 0 & 0 & 1 & 1 & 0 \end{bmatrix} \end{array}$$

Let S denote the 12 × 12 matrix. The hidden pattern of any desired state vector is found using the dynamical system S.

We take the state vector X = (0 0 1 0 0 0 0 0 0 0 0 0) in which the node 'caste system' alone is in the ON state; and all other nodes are in OFF state. The effect of X on the system S is given by

XS  ↪  (0 1 1 1 0 1 1 1 1 1 0 0)  =  $X_1$ (say)
$X_1$S  ↪  (1 1 1 1 1 1 1 1 1 1 1 1)  =  $X_2$ (say)
$X_2$S  ↪  (1 1 1 1 1 1 1 1 1 1 1 1)  =  $X_2$.

$X_2$ is the fixed point of the hidden pattern. The effect of this vector X on the dynamical system is a fixed point in which all nodes come to the ON state. Hence we see that caste system is the major cause for economic untouchability. Unless the caste system is abolished; there can be no improvement in the economic conditions for economic discrimination will be at its peak.

Next we consider the state vector Y = (0 0 0 0 0 0 0 1 0 0 0 0) where the node 'lack of political power' alone in the ON state and all the other nodes are in the OFF state.

The effect of the state vector Y on the dynamical system S is given by

$$YS \hookrightarrow (1\ 1\ 0\ 0\ 0\ 0\ 0\ 1\ 1\ 0\ 0\ 0) = Y_1 \text{ (say)}$$
$$Y_1 S \hookrightarrow (1\ 1\ 1\ 1\ 1\ 1\ 1\ 1\ 1\ 1\ 1\ 1) = Y_2 \text{ (say)}$$
$$Y_2 S \hookrightarrow (1\ 1\ 1\ 1\ 1\ 1\ 1\ 1\ 1\ 1\ 1\ 1) = Y_2 \text{ (say)}.$$

The hidden pattern results in a fixed point. Thus we see that 'lack of political power' has not only made Dalits and Sudras economically poor and economic untouchables but also made them to be socially discriminated and suffer socio-economic problems.

Next, we take the state vector Z where the node 'No education' is alone in the ON state of the node and all the other nodes are in the OFF state.

The effect of Z = (0 0 0 0 0 0 0 0 0 1 0 0) on the dynamical system S is given by

$$ZS \hookrightarrow (1\ 1\ 1\ 0\ 0\ 0\ 0\ 0\ 0\ 1\ 1\ 0) = Z_1 \text{ (say)}$$
$$Z_1 S \hookrightarrow (1\ 1\ 1\ 1\ 1\ 1\ 1\ 1\ 1\ 1\ 1\ 1) = Z_2 \text{ (say)}$$
$$Z_2 S \hookrightarrow Z_2.$$

This results in a fixed point. The resultant vector shows that all nodes are in the on state when the node 'no education' alone is on state. It explains that 'no education' has the biggest impact on the economic condition. Dalits and Sudras lacking education are economically discriminated and become victims of casteism. Study of all other nodes in the ON state has been carried out and the conclusions are given in the final chapter.

He expressed his opinion using indeterminates for Neutrosophic Cognitive Maps, which is as follows. He feels the relation between the nodes 'belief in destiny/karma' and 'lack of political power' is an indeterminate; also the relation between the nodes 'belief in destiny/karma' and 'no education' is an indeterminate.

The directed neutrosophic graph is as follows:

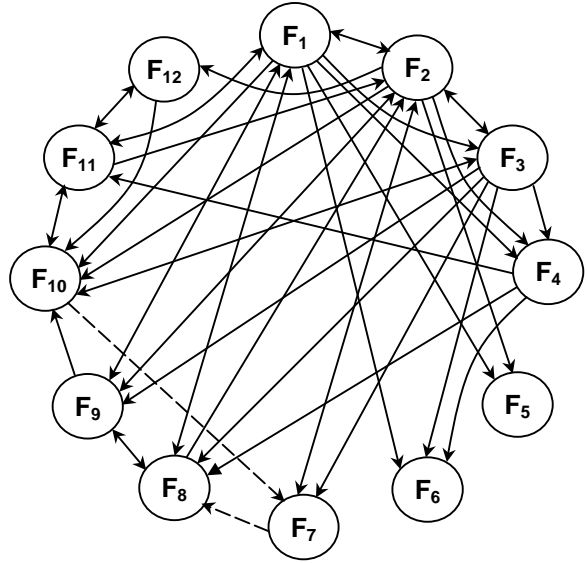

FIGURE 2.5.2

The related connection matrix S of the directed graph.

$$
\begin{array}{c}
\phantom{F_{1}}\begin{array}{cccccccccccc} F_1 & F_2 & F_3 & F_4 & F_5 & F_6 & F_7 & F_8 & F_9 & F_{10} & F_{11} & F_{12} \end{array} \\
\begin{array}{c} F_1 \\ F_2 \\ F_3 \\ F_4 \\ F_5 \\ F_6 \\ F_7 \\ F_8 \\ F_9 \\ F_{10} \\ F_{11} \\ F_{12} \end{array}
\begin{bmatrix}
0 & 1 & 1 & 1 & 1 & 1 & 0 & 1 & 1 & 1 & 1 & 0 \\
1 & 0 & 1 & 1 & 1 & 0 & 1 & 0 & 1 & 1 & 0 & 1 \\
0 & 1 & 0 & 1 & 0 & 1 & 1 & 1 & 1 & 1 & 0 & 0 \\
0 & 0 & 0 & 0 & 0 & 1 & 0 & 1 & 0 & 0 & 1 & 0 \\
0 & 0 & 0 & 0 & 0 & 0 & 0 & 0 & 0 & 0 & 0 & 0 \\
0 & 0 & 0 & 0 & 0 & 0 & 0 & 0 & 0 & 0 & 0 & 0 \\
0 & 1 & 0 & 0 & 0 & 0 & 0 & I & 0 & 0 & 0 & 0 \\
1 & 1 & 0 & 0 & 0 & 0 & 0 & 0 & 1 & 0 & 0 & 0 \\
1 & 1 & 0 & 0 & 0 & 0 & 0 & 1 & 0 & 1 & 0 & 0 \\
0 & 0 & 1 & 0 & 0 & 0 & I & 0 & 0 & 0 & 1 & 0 \\
1 & 1 & 0 & 0 & 0 & 0 & 0 & 0 & 0 & 1 & 0 & 1 \\
0 & 0 & 0 & 0 & 0 & 0 & 0 & 0 & 0 & 1 & 1 & 0
\end{bmatrix}
\end{array}
$$

Consider the state vector with node $F_7$ that is 'belief in destiny/karma' to be in the ON state and all other nodes are in the OFF state.

The effect of this on the neutrosophic dynamical system makes all nodes to come to the ON state in spite of the indeterminate factors in the NCM.

Now we proceed on to consider the attributes given by the second expert. This expert introduced several new nodes.

$H_1$ - Economic discrimination
$H_2$ - Psychological oppression
$H_3$ - Poor working condition with low pay
$H_4$ - Poor facility for education
$H_5$ - Caste system
$H_6$ - Manu Dharma
$H_7$ - Domination of Brahmins/ higher castes in political spheres so government plans to alleviate poverty are a farce
$H_8$ - No improvement in the economic conditions
$H_9$ - Poor economic conditions of Dalits and Sudras
$H_{10}$ - Constant fear of discrimination and ill-treatment
$H_{11}$ - Constant harassment by Brahmins/ higher castes
$H_{12}$ - No proper political representation
$H_{13}$ - No policy for their attaining equal social, political, economic status
$H_{14}$ - Financial assistances by government never reaches the poorest strata of the society
$H_{15}$ - Brahmins and higher caste should change their psychological hatred and give equal opportunities and strive to treat Dalits and Sudras with respect

The expert explained the node $H_{15}$ as follows: "Brahmins and caste Hindus should not exploit Dalits and Sudras by paying them meager amounts and extracting the maximum labour from them." He asked them to treat Dalits and Sudras with empathy. The practice of the caste system and untouchability has drained the economy of the nation. If caste is abolished, it could certainly make a drastic change in the economic, social and political conditions of Dalits and Sudras.

This expert has given the following directed graph using these attributes.

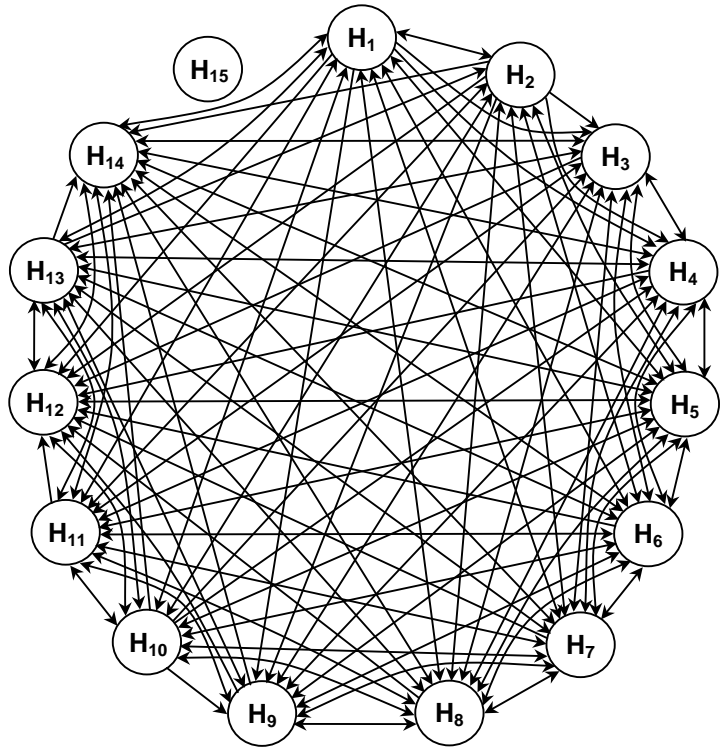

FIGURE: 2.5.3

The related 15 × 15 connection matrix W is given in the following page:

$$\begin{array}{c} \phantom{H_1} H_1\ H_2\ H_3\ H_4\ H_5\ H_6\ H_7\ H_8\ H_9\ H_{10}\ H_{11}\ H_{12}\ H_{13}\ H_{14}\ H_{15} \\ \begin{array}{c} H_1 \\ H_2 \\ H_3 \\ H_4 \\ H_5 \\ H_6 \\ H_7 \\ H_8 \\ H_9 \\ H_{10} \\ H_{11} \\ H_{12} \\ H_{13} \\ H_{14} \\ H_{15} \end{array} \left[\begin{array}{ccccccccccccccc} 0 & 1 & 1 & 1 & 1 & 1 & 1 & 1 & 1 & 1 & 1 & 1 & 1 & 0 \\ 1 & 0 & 1 & 1 & 1 & 1 & 1 & 1 & 1 & 1 & 1 & 1 & 0 & 1 & 0 \\ 0 & 0 & 0 & 1 & 1 & 1 & 1 & 1 & 1 & 0 & 1 & 1 & 1 & 1 & 0 \\ 0 & 0 & 1 & 0 & 1 & 1 & 1 & 1 & 1 & 0 & 1 & 1 & 1 & 1 & 0 \\ 1 & 1 & 1 & 1 & 0 & 1 & 1 & 1 & 1 & 1 & 1 & 1 & 1 & 1 & 0 \\ 1 & 1 & 1 & 1 & 1 & 0 & 1 & 1 & 1 & 1 & 1 & 1 & 1 & 1 & 0 \\ 1 & 1 & 1 & 1 & 1 & 1 & 0 & 1 & 1 & 1 & 1 & 1 & 1 & 1 & 0 \\ 1 & 1 & 1 & 1 & 1 & 1 & 1 & 0 & 1 & 1 & 1 & 1 & 1 & 1 & 0 \\ 0 & 0 & 1 & 1 & 1 & 1 & 1 & 1 & 0 & 0 & 1 & 1 & 1 & 1 & 0 \\ 1 & 1 & 1 & 1 & 1 & 1 & 1 & 1 & 1 & 0 & 1 & 1 & 1 & 1 & 0 \\ 1 & 1 & 1 & 1 & 1 & 1 & 1 & 1 & 1 & 0 & 0 & 1 & 1 & 1 & 0 \\ 0 & 0 & 0 & 0 & 1 & 0 & 1 & 1 & 1 & 1 & 1 & 0 & 1 & 1 & 0 \\ 1 & 1 & 1 & 1 & 1 & 1 & 1 & 0 & 1 & 1 & 1 & 0 & 1 & 0 \\ 1 & 0 & 1 & 1 & 1 & 1 & 1 & 1 & 0 & 0 & 0 & 1 & 0 & 0 & 0 \\ 0 & 0 & 0 & 0 & 0 & 0 & 0 & 0 & 0 & 0 & 0 & 0 & 0 & 0 & 0 \end{array}\right] \end{array}$$

Node $H_{15}$ has no relation with all other nodes. The expert says that once this is achieved it will put an end to economic discrimination; free Dalits and Sudras from oppression and naturally working conditions will be improved and consequently they will get education. The caste system will die a natural death. It will also have an effect on the poverty alleviation schemes of the Government. If such individual self-realization comes upon the caste-Hindus and Brahmins, there will be drastic change not only in the lives of downtrodden communities but also positively reflect on the nation's progress. Systems of casteist exploitation of labour are so strongly ingrained that caste-Hindus, including Brahmins, refuse to think otherwise.

We obtain results using the 15 × 15 connection matrix, W. Let us consider the state vector X = (0 0 0 0 0 0 0 0 0 1 0 0 0 0 0) where the node 'constant fear of discrimination and ill-treatment' is in the ON state and all other nodes are in the OFF state. The effect of X on the dynamical system W

$XW \hookrightarrow (1 1\ 1\ 1\ 1\ 1\ 1\ 1\ 1\ 1\ 1\ 1\ 1\ 0) = X_2$ say
$X_2W \hookrightarrow (1 1\ 1\ 1\ 1\ 1\ 1\ 1\ 1\ 1\ 1\ 1\ 1\ 0) = X_3$ say
$X_3 = X_2$ a fixed point.

The hidden pattern is a fixed point, which shows the close interconnection between the constant fear of discrimination undergone by Dalits and other concepts in this system.

If we assume that the node 'caste system' alone is in the ON state and all other nodes in the state vector Y were in the OFF state, the effect of Y on the dynamical system W is given by

$YW \hookrightarrow (1\ 1\ 1\ 1\ 1\ 1\ 1\ 1\ 1\ 1\ 1\ 1\ 1\ 0) = Y_1$ (say)
(where $Y = (0\ 0\ 0\ 0\ 1\ 0\ 0\ 0\ 0\ 0\ 0\ 0\ 0\ 0)$).

$Y_1W = Y_1$ (a fixed point). As long as caste system exists in India it is certain that economic untouchability will be practiced.

When we asked the expert to use neutrosophy he remarked the node $H_{15}$ which he has purposefully included to show its independency over other nodes can be interrelated with some of the nodes as an indeterminate. We see the effect also is an indeterminate of the node $H_{15}$ in most of the cases.

Now we proceed onto analyze the second expert whom we felt was little different and innovative even in specifying the nodes for economic untouchability which are as follows:

$S_1$ - Economic inequality
$S_2$ - Caste system and its relation to economic status (Varnashrama Dharma)
$S_3$ - Concentration of wealth solely in the hands of Brahmins/upper castes
$S_4$ - No political power
$S_5$ - Landlessness
$S_6$ - Lack of education
$S_7$ - Domination by 'upper' castes in all fields of power and prosperity
$S_8$ - Less pay and more hours of work

- $S_9$ - No awareness among Dalits and Sudras (their job is dependent on their subjugation to caste-Hindus/ Brahmins)
- $S_{10}$ - Exploitation by caste Hindus and Brahmins
- $S_{11}$ - Terrorized from even demanding their rights by 'upper' castes and Brahmins (as they would be jobless)
- $S_{12}$ - Undemocratic actions taken on Dalits on all issues by caste-Hindu/ Brahmin public and government
- $S_{13}$ - Untouchability
- $S_{14}$ - Legal codes of the Varnashrama Dharma

These attributes were defined as the main concepts related to the economic untouchability. A very dense directed graph is given by the expert.

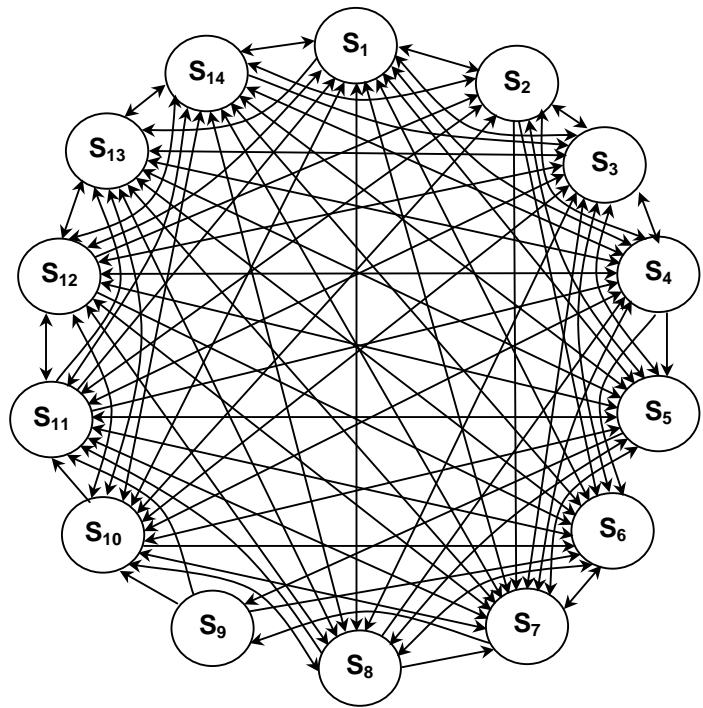

FIGURE: 2.1.5

The 14 × 14 relational connection matrix of the directed graph is given by matrix N.

$$N = \begin{array}{c} \\ S_1 \\ S_2 \\ S_3 \\ S_4 \\ S_5 \\ S_6 \\ S_7 \\ S_8 \\ S_9 \\ S_{10} \\ S_{11} \\ S_{12} \\ S_{13} \\ S_{14} \end{array} \begin{array}{c} S_1 \; S_2 \; S_3 \; S_4 \; S_5 \; S_6 \; S_7 \; S_8 \; S_9 \; S_{10} \; S_{11} \; S_{12} \; S_{13} \; S_{14} \\ \begin{bmatrix} 0 & 1 & 1 & 1 & 1 & 1 & 1 & 1 & 0 & 1 & 0 & 1 & 1 & 1 \\ 1 & 0 & 1 & 1 & 1 & 1 & 1 & 0 & 0 & 1 & 1 & 0 & 0 & 1 \\ 1 & 1 & 0 & 1 & 1 & 1 & 1 & 1 & 0 & 1 & 1 & 1 & 1 & 1 \\ 1 & 1 & 1 & 0 & 1 & 1 & 1 & 0 & 0 & 1 & 0 & 1 & 1 & 0 \\ 1 & 0 & 1 & 0 & 0 & 0 & 1 & 1 & 1 & 1 & 1 & 1 & 1 & 1 \\ 1 & 0 & 1 & 0 & 0 & 0 & 1 & 1 & 0 & 0 & 1 & 0 & 1 & 1 \\ 1 & 0 & 1 & 1 & 1 & 1 & 0 & 0 & 1 & 1 & 1 & 1 & 1 & 1 \\ 1 & 0 & 1 & 1 & 1 & 1 & 1 & 0 & 0 & 1 & 1 & 1 & 1 & 1 \\ 0 & 0 & 0 & 0 & 0 & 1 & 0 & 0 & 0 & 0 & 1 & 0 & 0 & 0 \\ 1 & 1 & 1 & 1 & 1 & 1 & 1 & 1 & 0 & 0 & 1 & 1 & 1 & 1 \\ 1 & 1 & 1 & 1 & 1 & 1 & 1 & 0 & 0 & 0 & 0 & 1 & 1 & 1 \\ 1 & 1 & 1 & 1 & 1 & 1 & 1 & 1 & 0 & 1 & 1 & 0 & 1 & 1 \\ 0 & 1 & 1 & 1 & 1 & 0 & 1 & 1 & 0 & 1 & 0 & 1 & 0 & 1 \\ 1 & 1 & 1 & 1 & 1 & 1 & 1 & 0 & 1 & 1 & 1 & 1 & 0 \end{bmatrix} \end{array}$$

We study the effect of the state vector X on the dynamical system N. The node 'Dalits are landlessness' is in the ON state, all other attributes are kept in the OFF state in X = (0 0 0 0 1 0 0 0 0 0 0 0 0 0).

The effect of the state vector X on the dynamical system N is given by

$$XN \hookrightarrow (1\;0\;1\;1\;1\;0\;1\;1\;1\;1\;1\;1\;1\;1) = X_1 \text{ (say)}$$
$$X_1N \hookrightarrow (1\;1\;1\;1\;1\;1\;1\;1\;1\;1\;1\;1\;1\;1) = X_2 \text{ (say)}$$
$$X_2N \hookrightarrow X_2.$$

This results in a fixed point. We see from the resultant state vector that when the node 'landlessness' is in the ON state, all other nodes like exploitation, less pay and more hours of work, gradually becoming indebted to the landlords due to poverty resulting in bonded labour and so on come to the ON state.

This clearly shows that if the Government distributes land to the Dalits, it will certainly empower them in all aspects. As per our expert, this will be a major step in putting an end to the economic exploitation faced by the Dalits.

Next, we take the node 'concentration of wealth solely in the hands of Brahmins/upper castes' to be in the ON state and all other nodes remain in the OFF state. We study the effect of the ON state of this node on the dynamical system. Let $Y = (0\ 0\ 1\ 0\ 0\ 0\ 0\ 0\ 0\ 0\ 0\ 0\ 0)$.

The effect of the state vector Y on the dynamical system N is given by

$$YN \hookrightarrow (1\ 1\ 1\ 1\ 1\ 1\ 1\ 0\ 1\ 1\ 1\ 1\ 1) = Y_1 \text{ (say)}$$
$$Y_1N \hookrightarrow (1\ 1\ 1\ 1\ 1\ 1\ 1\ 1\ 1\ 1\ 1\ 1\ 1) = Y_2 \text{ (say)}.$$

By thresholding and updating, we get the $Y_2N = Y_2$, resulting in a fixed point.

In the resultant vector all the nodes are in the ON state. Thus concentration of wealth solely in the hands of Brahmins/upper castes is sufficient to create and extend the economic untouchability. All these concepts are very closely linked and this is evident from the resultant vector.

Consider $T = (0\ 0\ 0\ 0\ 0\ 0\ 0\ 0\ 0\ 0\ 0\ 0\ 1)$ i.e., where the node, 'Code of Varnashrama Dharma' alone is in the ON state and all other nodes are in the OFF state. We now proceed on to study the effect of the state vector T on the dynamical system N.

$$TN \hookrightarrow (1\ 1\ 1\ 1\ 1\ 1\ 1\ 0\ 1\ 1\ 1\ 1\ 1) = T_1 \text{ (say)}$$
$$T_1N \hookrightarrow (1\ 1\ 1\ 1\ 1\ 1\ 1\ 1\ 1\ 1\ 1\ 1\ 1) = T_2 \text{ (say)}$$
$$T_2N = T_2 \text{ (a fixed point)}.$$

Thus we see that if Varnashrama Dharma alone is in the ON state, all other nodes come to the ON state. Thus it is evidently seen that Varnashrama Dharma alone is sufficient to ruin the lives of Dalits and the Sudras beyond any repair. It is unfortunate, but true that even today the sad state of our nation is due to the practice of Varnashrama Dharma as a result of which Dalits and the Sudras live a life of discrimination, poverty, constant fear of harassment and degradation. Several other conclusions based on this study is given in the last chapter.

This expert declined to work with NCM saying that the concepts were extremely obvious, and nothing was indeterminate in his opinion.

Now we proceed on to take a few of the attributes, common in the study of all the five aspects of untouchability viz., religious untouchability, political untouchability, social untouchability, educational untouchability and finally the economic untouchability. We use these special nodes and investigate the effect. Thus we take these five as major nodes and their common attributes of as the other nodes:

$N_1$ - Religious untouchability
$N_2$ - Political untouchability
$N_3$ - Educational untouchability
$N_4$ - Social untouchability
$N_5$ - Economic untouchability
$N_6$ - Role of the present day government
$N_7$ - Laws of Manu
$N_8$ - No power or decision-making vests in the hands of Dalits and Sudras
$N_9$ - Lack of education among Dalits and Sudras
$N_{10}$ - Poor economic condition of Dalits/Sudras
$N_{11}$ - Poverty / landlessness
$N_{12}$ - Not policy makers, even for them
$N_{13}$ - Samadharma
$N_{14}$ - Role of Brahmins and caste Hindus in power
$N_{15}$ - Social inequality/ discrimination
$N_{16}$ - Discrimination by caste, which is based on birth
$N_{17}$ - Brahmin and caste Hindu arrogance
$N_{18}$ - Varnashrama Dharma
$N_{19}$ - Inequality in all dealings with Dalits/Sudras

We had to lessen the number of attributes considerably in order to avoid the difficulty and cumbersome working out.

Let S denote the 19 × 19 matrix associated connection matrix of the directed graph given by the expert.

$$S = \begin{bmatrix}
0 & 1 & 1 & 1 & 1 & 1 & 1 & 1 & 1 & 0 & 1 & 1 & -1 & 1 & 1 & 1 & 1 & 1 & 1 \\
1 & 0 & 1 & 1 & 1 & 1 & 1 & 1 & 1 & 1 & 1 & 1 & -1 & 1 & 1 & 1 & 0 & 1 & 1 \\
1 & 1 & 0 & 1 & 1 & 1 & 1 & 1 & 0 & 1 & 1 & 1 & -1 & 1 & 1 & 1 & 1 & 1 & 1 \\
1 & 1 & 1 & 0 & 1 & 1 & 1 & 1 & 0 & 1 & 1 & 1 & -1 & 1 & 1 & 1 & 1 & 1 & 1 \\
1 & 1 & 1 & 1 & 0 & 1 & 1 & 1 & 1 & 1 & 1 & 1 & -1 & 1 & 1 & 1 & 1 & 1 & 1 \\
1 & 1 & 1 & 1 & 1 & 0 & 1 & 1 & 0 & 1 & 1 & 0 & -1 & 1 & 1 & 1 & 1 & 1 & 1 \\
1 & 1 & 1 & 1 & 1 & 1 & 0 & 1 & 1 & 1 & 1 & 1 & -1 & 1 & 1 & 1 & 1 & 1 & 1 \\
1 & 1 & 1 & 1 & 1 & 1 & 1 & 0 & 0 & 1 & 1 & 1 & -1 & 1 & 1 & 1 & 1 & 1 & 1 \\
1 & 1 & 1 & 1 & 1 & 1 & 1 & 1 & 0 & 1 & 1 & 1 & -1 & 1 & 1 & 1 & 1 & 1 & 1 \\
1 & 1 & 1 & 1 & 1 & 1 & 1 & 1 & 1 & 0 & 1 & 1 & -1 & 1 & 1 & 1 & 1 & 1 & 1 \\
1 & 1 & 1 & 1 & 1 & 1 & 1 & 1 & 1 & 0 & 1 & -1 & 1 & 1 & 1 & 1 & 1 & 1 \\
1 & 1 & 1 & 1 & 1 & 1 & 1 & 1 & 1 & 1 & 0 & 1 & -1 & 1 & 1 & 1 & 1 & 1 & 1 \\
-1 & -1 & -1 & -1 & -1 & -1 & -1 & -1 & -1 & -1 & -1 & 0 & 0 & -1 & -1 & -1 & -1 & -1 & -1 \\
1 & 1 & 1 & 1 & 1 & 1 & 1 & 1 & 1 & 1 & 1 & 1 & 0 & 1 & 1 & 1 & 1 & 1 \\
1 & 1 & 1 & 1 & 1 & 1 & 1 & 1 & 1 & 1 & 1 & 1 & 1 & 0 & 1 & 1 & 1 & 1 \\
1 & 1 & 1 & 1 & 1 & 1 & 1 & 1 & 1 & 1 & 1 & 1 & -1 & 1 & 0 & 1 & 1 & 1 \\
1 & 1 & 1 & 1 & 1 & 1 & 1 & 1 & 1 & 1 & 1 & 1 & 1 & 1 & 0 & 1 & 1 \\
1 & 1 & 1 & 1 & 1 & 1 & 1 & 1 & 1 & 1 & 1 & 1 & 1 & 1 & 1 & 0 & 1 \\
1 & 1 & 1 & 1 & 1 & 1 & 1 & 1 & 1 & 1 & 1 & 1 & 1 & 1 & 1 & 1 & 0
\end{bmatrix}$$

Let X = (0 0 0 0 0 0 1 0 0 0 0 0 0 0 0 0 0 0 0) be the state vector. Effect of X on the dynamical system S is given by

XS ↪ (1 1 1 1 1 1 1 1 1 1 1 1 0 1 1 1 1 1 1) = $X_1$ (say)
$X_1$S ↪ (1 1 1 1 1 1 1 1 1 1 1 1 0 1 1 1 1 1 1) = $X_2$ = $X_2$

is a fixed point. Thus Manu Dharma is the vital cause for all the social evils and once it is taken to be in the ON state of all the other nodes come to the ON state, except Samadharma (equality). Thus we see that Manu Dharma and Samadharma are directly opposed to each other. Manu Dharma is the root cause of all the 5 major aspects of untouchability.

Next we consider a state vector V = (1 0 0 0 0 0 0 0 0 0 0 0 0 0 0 0 0 0 0) where only the node 'religious untouchability' is in the ON state and all other nodes are in the OFF state. We find the effect of V on the dynamical system S

VS ↪ (1 1 1 1 1 1 1 1 1 0 1 1 0 1 1 1 1 1 1) = $V_1$ (say)
$V_1$S ↪ (1 1 1 1 1 1 1 1 1 1 1 1 1 0 1 1 1 1 1 1) = $V_2$ (say)
$V_2$S = $V_2$ is a fixed point.

X = (0 0 0 0 0 0 0 0 0 0 0 0 1 0 0 0 0 0 0) and XS ↪ X. Thus when the node $N_{13}$, i.e. Samadharma alone is in the ON state, the effect of this on the dynamical system gives the resultant vector to be the same. That is, all nodes remain in the OFF state, which shows that all social evils are in OFF state if Samadharma is alone ON.

Likewise we see that each of the nodes, which we have considered are very intricate and important because they totally affect all the other nodes. Conclusions based on the analysis of this dynamical system S are given in the fourth chapter.

In all these Fuzzy Cognitive Maps, which basically work with the aid of a connection matrix and the state vector, we have used C programme to determine the hidden pattern of the dynamical system. [69, 76]

## 2.6 Analysis of the social problem of untouchability using Fuzzy Relational Maps

We use another special type of mathematical tool called Fuzzy Relational Maps (FRMs). The main motivation for using this model in the place of FCM is that in several models we had to work with 25 × 25 matrix or a 19 × 19 matrix but if we use FRMs wherever viable we can certainly reduce the 25 × 25 matrix into a m × n matrix (m ≠ n) with m + n = 25. This will reduce the cumbersome working of the model. We use FRMs to study the Periyar's ideology and his social reform. The basic methods and movement by which Periyar fought against untouchability, denial of temple entry etc, are taken as the domain space of the FRM. The nodes of the domain space are given by

$P_1$ - Samadharma
$P_2$ - Self Respect Movement

- $P_3$ - Temple entry demonstration
- $P_4$ - Non-Brahmin Movement
- $P_5$ - Denial of existence of god
- $P_6$ - Demonstrations against untouchability
- $P_7$ - Opposition to denial of education for Dalits and Sudras (Laws of Manu)
- $P_8$ - Annihilation of caste system
- $P_9$ - Reformation of villages against superstitions/ rituals
- $P_{10}$ - Fight for communal representation

The range space nodes given by the expert are detailed below:

- $R_1$ - Fight against domination of Brahmins in fields of education and administration
- $R_2$ - Protest rule of Manu's Varnashrama Dharma
- $R_3$ - Stop the practice of untouchability
- $R_4$ - Abandon Varnashrama Dharma (discrimination by birth and occupation)
- $R_5$ - To fight against denial of education to Sudras and Dalits
- $R_6$ - Fight against denial of land and property to Dalits and Sudras
- $R_7$ - Fight against denial of temple entry and entry into sanctum sanctorum of the temple
- $R_8$ - Difference in food and place were food is served for Dalits and Sudra children in comparison with Brahmin children
- $R_9$ - Right to use all public roads by Dalits
- $R_{10}$ - Fight against denial of right to use common wells and lakes by Dalits
- $R_{11}$ - Rights denied to Dalits to wear new clothes, and carry umbrellas
- $R_{12}$ - Rights denied to Dalits to sit or use slippers in the presence of Brahmins or high castes even today in remote villages.

The first expert assigned the above attributes for the domain space and range space. The related directed graph given by the expert is given in the following figure:

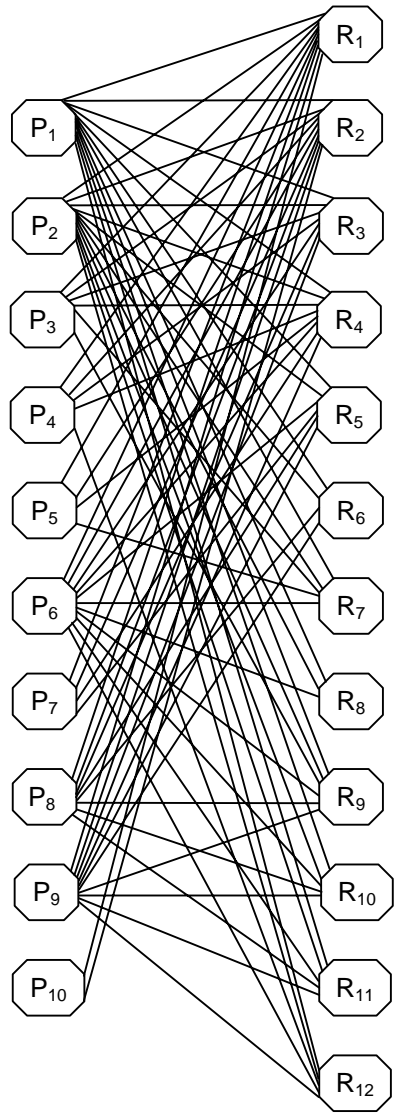

FIGURE: 2.6.1

The related connection matrix is given by M

$$
\begin{array}{c}
\phantom{P_{1}}\;\;R_1\;R_2\;R_3\;R_4\;R_5\;R_6\;R_7\;R_8\;R_9\;R_{10}\;R_{11}\;R_{12}\\
\begin{array}{c}P_1\\P_2\\P_3\\P_4\\P_5\\P_6\\P_7\\P_8\\P_9\\P_{10}\end{array}
\begin{bmatrix}
1 & 1 & 1 & 1 & 1 & 1 & 1 & 1 & 1 & 1 & 1 & 1\\
1 & 1 & 1 & 1 & 1 & 1 & 1 & 1 & 1 & 1 & 1 & 1\\
1 & 1 & 1 & 1 & 0 & 0 & 1 & 0 & 1 & 0 & 0 & 0\\
1 & 1 & 1 & 1 & 0 & 0 & 0 & 0 & 0 & 0 & 0 & 1\\
1 & 1 & 0 & 1 & 0 & 0 & 1 & 0 & 0 & 0 & 0 & 0\\
1 & 1 & 1 & 1 & 1 & 0 & 1 & 1 & 1 & 1 & 1 & 1\\
1 & 1 & 0 & 1 & 1 & 0 & 0 & 0 & 0 & 0 & 0 & 0\\
1 & 1 & 1 & 0 & 1 & 1 & 0 & 0 & 1 & 1 & 1 & 0\\
1 & 1 & 1 & 1 & 1 & 1 & 0 & 0 & 1 & 1 & 1 & 1\\
1 & 1 & 0 & 0 & 0 & 0 & 0 & 0 & 0 & 0 & 0 & 0
\end{bmatrix}
\end{array}
$$

Let M denote the 10 × 12 connection matrix. We study the effect of each state vector on the dynamical system.

Let X = (0 0 0 0 0 0 1 0 0 0) i.e., the only node which is in the ON state is the denial of education for all and all other nodes are in the OFF state. Now we obtain the effect of X on the dynamical system

$$
\begin{aligned}
XM &\hookrightarrow (1\;1\;0\;1\;1\;0\;0\;0\;0\;0\;0\;0) &=& \;Y\;(\text{say})\\
YM^T &\hookrightarrow (1\;1\;1\;1\;1\;1\;1\;1\;1\;1) &=& \;X_1\;(\text{say})\\
X_1 M &\hookrightarrow (1\;1\;1\;1\;1\;1\;1\;1\;1\;1\;1\;1) &=& \;Y_1\;\text{say}\\
Y_1 M^T &= X_2 = (X_1),\; X_1 M &=& \;Y_1.
\end{aligned}
$$

The hidden pattern is a binary pair that is a fixed point. {(1 1 1 1 1 1 1 1 1 1), (1 1 1 1 1 1 1 1 1 1 1 1)} is the fixed binary pair, which is such that all the nodes, both in the domain space and the range space are on.

Hence we see that according to Periyar the greatest social evil is the denial of access to education for everybody. This is evident from the fact that not only where they denied education in the ancient days according to the dictums of Manu Dharma, but even now they suffer because of denial of access to education. The education that Dalits and Sudras obtain is much lower in standard and nowhere in comparison with the education that upper-caste Brahmins receive. Schools run for the Dalits and Sudras lack proper infrastructure or well-qualified teachers. Similarly, the greatest sections of dropouts in school education are from the Dalit and Sudra communities. Majority of the Dalits and Sudras who live below the poverty line cannot think of educating their children even up to the 5$^{th}$ standard. When this attribute becomes ON it makes all attributes come to the ON state.

Now we see the effects of the attributes in the range space; suppose we take Y = (0 0 0 0 0 0 0 0 0 1 0 0 ) i.e. "right denied for the Dalits and Sudras to use common well and lakes" to be in the ON state and all other nodes in the range space are in the OFF state. The effect of Y on the dynamical system is given by

$YM^T$  ↪  (1 1 0 0 0 1 0 1 1 0)       =   X (say)
$XM$   ↪  (1 1 1 1 1 1 1 1 1 1 1 1)  =   $Y_1$ say
$Y_1M^T$ ↪  (1 1 1 1 1 1 1 1 1 1 )     =   Y.

This leads to the fixed point given by the pair {(1 1 1 1 1 1 1 1 1), (1 1 1 1 1 1 1 1 1 1 1 1)} which shows that the denial of the right to Dalits to use the common well or lake is an important discrimination which affects all other concepts.

Now suppose we consider the ON state of the node 'communal representation' and all other nodes are in the OFF state.

We study the effect of the state vector X = (0 0 0 0 0 0 0 0 1) on M.

$XM$   ↪  (1 1 0 0 0 0 0 0 0 0 0 0)   =   $Y_1$ (say)
$Y_1M^T$ ↪  (1 1 1 1 1 1 1 1 1 1)      =   $X_1$ (say)
$X_1M$  ↪  (1 1 1 1 1 1 1 1 1 1 1 1)  =   $Y_2$ (say)
$Y_2M^T$ =   $X_1$,

resulting in the fixed point, which is the binary pair {(1 1 1 1 1 1 1 1 1), (1 1 1 1 1 1 1 1 1 1 1 1)}.

'Fight for Communal Representation' for the Dalits and Sudras also makes all the nodes of the state vector ON.

In conclusion we can say that all the nodes given are so vital and the presence of anyone of them makes all other nodes to come to the ON state. In the place of FRMs one can use Neutrosophic Relational Maps (NRMs) and obtain results.

From our study and analysis of untouchability using fuzzy and neutrosophic logic we are surprised to see that almost all the graphs were very densely related with arrows, which obviously shows that every social evil is strongly bonded with other social evils.

Chapter Three

# PERIYAR'S BIOGRAPHY AND HIS VIEWS ON UNTOUCHABILITY

The first section of this chapter is a brief biography of Periyar. For more details about his life and works refer [2], [87] and [96]. The second section contains excerpts from the translated speeches and writings of Periyar that deal with this opinion on Untouchability. Some statistics of the atrocities against Dalits are given in the third and final section.

### 3.1 Life and Struggle of Periyar E. V. R.

> *"The enormous privileges given to Brahmins by the Vedas were sacrosanct only as long as they went unchallenged. The challenge rose in the Tamil country like a whirlwind, spearheaded by an iconoclast who questioned the Vedas and the gods as well. He took apparently extremist positions on some issues, but the fundamental nature of the social revolution he wrought was clear even to its victims. The political perspectives of Tamil Nadu and with it much of India, were altered with a seeming finality by Ramasamy Naicker. (1879-1973)"*
>
> - Makers of the Millennium, India 1000 to 2000, Express Publications, Madurai.

Periyar Erode Venkata Ramasamy (17.9.1879 - 24.12.1973) was a social revolutionary who influenced the politics of the south Indian state of Tamil Nadu for half a century. Periyar, literally meaning a great man was an honorific prefix bestowed on him by the people of Tamil Nadu. He was hailed by a UNESCO Award as 'the prophet of the new age, the Socrates of south east Asia, father of the social reform movement and arch enemy of ignorance, superstitions, meaningless customs and base manners.'

He was the second son of Venkata (Naicker) and Chinna Thayammal. They were a very affluent family; his father was a popular, rich merchant known for his religious philanthropy. Leaders like Motilal Nehru (father of Jawaharlal Nehru, independent India's first Prime Minister) and Mahatma Gandhi have stayed at Periyar's ancestral home in Erode. His school education stopped at the age of ten, and he subsequently entered his father's trade. Even as a small boy, he was naughty enough to question the contradictions and illusions that existed in the mythologies narrated by Brahmins. He married Nagammai in 1898. In 1904, he renounced his family and went to Varanasi, a sacred Hindu pilgrimage town on the banks of the Ganga. The insults and the differential treatment he faced, for not being a Brahmin, made him aware of the cruelty of caste. Giving up thoughts of renunciation, he returned to his hometown Erode, where his father gave him the trade rights. His business flourished and he became a well-known in Erode. In 1918, he formally entered public life by becoming the Chairman of Erode Municipality. The British Government made him an honorary magistrate and he held 29 honorary positions at that time in banks, libraries, war recruitment committee, merchant associations, school committees etc. His friends P.Varadharajulu Naidu and C. Rajagopalachari persuaded him to join the Indian National Congress led by Mohandas Gandhi. In 1919, he joined the Congress, and resigned all the public posts he held at that time. He gave up his lucrative wholesale dealership in grocery and agricultural products, and closed his newly begun spinning mill. His annual income had been Rs.20,000, a princely figure in those days where the price of eight grams of gold was not even Rs.10.

The Congress policy of Prohibition of Liquor was conceived in Periyar's house in Erode during Gandhi's stay there upon the request of Periyar's wife and sister. Gandhi subsequently announced that the Congress party men should undertake picketing in front of toddy shops throughout the country and urged the British Government to implement the policy of liquor prohibition. Periyar wholeheartedly undertook the constructive programme - spreading the use of Khadi, picketing toddy shops, and boycotting the shops selling foreign cloth and eradication of untouchability. He courted imprisonment for picketing toddy shops in Erode in 1921. When his wife as well as his sister joined the agitation, it gained momentum, and the administration was forced to come to a compromise. When Congress frontline leaders requested Mahatma Gandhi to stop that agitation, he seriously told them that the decision of stopping the agitation was not in his own hands, but with two women in Erode, implying Periyar's wife and sister.

In 1922, Periyar became the President of the Tamil Nadu Congress Committee. At its provincial conference held in Tirupur, he moved a resolution that urged people of all castes to be allowed to enter and worship in all the temples, as a measure to end birth-based discrimination. Citing the authority of Vedas and other Hindu scriptures, Brahmin members of the Committee opposed the resolution and stalled its passage. This reactionary stand of the members of upper castes provoked Periyar to declare that he would burn Manu Dharma Sastra, Ramayana etc. to show his disapproval to accept such scriptures to govern the life of people.

His radical outlook enabled him to even support the progressive policies of other parties such as the South Indian Liberation Federation (popularly, the Justice Party). Though he was a Congress leader, he welcomed the Justice Party's Act passed in the Madras State Legislative Council to form a Hindu Religious Endowment Board with a view to put an end to the age-old monopoly and exploitation of the Brahmins in the managements of Hindu temples. Likewise, in 1924, he appreciated the Justice Party's Government Order to implement the policy of Communal Representation in education and employment.

In the temple town of Vaikom in present day Kerala, the untouchable and 'low' caste communities were prohibited from walking on the temple streets. The Congress Party there began a peaceful agitation against it in 1924 and they invited Periyar to take up the leadership of the movement. He was sentenced twice to undergo imprisonment. The agitations went on for a year. Even Gandhi came to Vaikom. Periyar's intense involvement granted the basic civil right of using public roads to the 'low' castes and untouchables. His courageous struggle earned him the honour of being called the "Hero of Vaikom."

Out of the funds of the Tamil Nadu Congress Committee (TNCC), they started the National Training School at Cheranmadevi near Tirunelveli in Tamil Nadu, as an alternative to those run by the British Government. It was managed by a Brahmin, V.V.S. Iyer. In that school's hostel called *Gurukulam*, the Brahmin and non-Brahmin students were segregated. The Brahmin boys were treated in a better manner with regard to food, shelter and education. Angered over this, Periyar resigned from the post of Secretary of the TNCC. He later put an end to this discriminatory practice. Periyar was later elected the President of the TNCC.

Since 1920, he had attempted to move resolutions in the Congress conferences demanding Communal Representation for non-Brahmins in government jobs and education. But the Brahmins systematically defeated it: conferences in Tirunelveli (1920), Thanjavur (1921), Tirupur (1922), Salem (1923) and Thiruvannamalai (1924). When the same resolution was not passed in the Congress Conference at Kancheepuram in November 1925, Periyar quit the Congress party. Enraged, he thundered that his only future task was to destroy the Brahmin rule by all means. So, in December 1925, he organized a parallel conference of non-Brahmins and launched the Self Respect Movement. When he met the Congress leader M.K.Gandhi at Bangalore in 1927, he strongly argued that unless the poisonous caste-system called 'Varnashrama Dharma' was uprooted, the eradication of 'Untouchability' was not possible. With Periyar's active support, the minister S.Muthaiah (Mudaliar) implemented the Communal Reservation scheme of the Justice Party.

On 2 May 1925, he had started to publish a Tamil weekly magazine titled *Kudi Arasu* (Republic). In his second journalistic foray, he started the English magazine 'Revolt' on 7 November 1928.

In February 1929, the first provincial conference of Self-Respect Movement was organized in Chengalpattu. The same year, he introduced a new rationalist marriage system called the Self-Respect Marriage where rituals and Brahmin priests uttering mantras where forbidden. He encouraged inter-caste and widow remarriages. Marriage was stripped of its many rituals, it was secularized, and was a simple declaration between the life-partners. He insisted on a rational outlook to bring about intellectual emancipation. He sought to annihilate the hierarchical, birth-based caste structure and establish an egalitarian social structure.

The Self-Respect Movement also carried on a strong propaganda against the ridiculous superstitions practiced by the people. He strove to dispel their ignorance. He was against religion and the domination of Brahmins.

Because of his rationalist philosophy, Periyar has often been compared with Voltaire (1694-1778). Both had aroused their people to realize that all men are equal and it is the birthright of every individual to enjoy liberty, equality and fraternity. Both opposed religion virulently because the so-called men of religion invented myths and superstitions to keep the innocent and ignorant people in darkness and go on exploiting them. If Periyar's opposition of religion and superstition was an attack on the Brahmins, Voltaire too had condemned religious men in similar terms. In one of his articles, Voltaire wrote, "They (the religious men) inspired you with false beliefs and made you fanatics so that they might be your masters. They made you superstitious, not that you might fear god but that you might fear them."

The second provincial conference of the Self Respect Movement was held in Erode in May 1930. The same year he actively supported the Bill for the abolition of Devadasi system (young girls from a particular backward community were dedicated as dancing girls to temples).

His social work was multifaceted. Not only did he fight for equality and the rights of the oppressed people, but he also had a program for economic equality. Along with the veteran communist leader M.Singaravel, he organised industrial and agricultural labourers to stand against the exploitation of big capitalists and landlords. However, since the British took steps to ban the Communist Party of India and similar organizations, Periyar toned down his socialist activities. Because of his experience he knew that there were a lot of support to carry on the freedom struggle and organize the labourers, but not everybody would come forward to criticize the tyranny of Brahmins. In 1933, the British banned his magazine *Kudi Arasu*. He consequently started another magazine *Puratchi* (Revolution). In 1934, Jayaprakash Narayan, the reputed Socialist Leader asked Periyar to join the Socialist party. The same year, Periyar started the Tamil weekly *Pagutharivu* (Rationalism).

The Tamil weekly paper *Viduthalai* (Liberation) of the Justice Party, which was founded in 1935, started to be published by Periyar as a Tamil daily from 1.1.1937. When the Justice Party was voted out of office in 1937, Rajagopalachari headed the newly elected Congress Government. He announced that Hindi was going to be a subject of compulsory study in schools. Periyar vehemently opposed this imposition of a language and organized a campaign against it on a war footing. On 26 December 1935, at a conference convened in Tiruchi, Periyar declared that the only solution to defeat the dominance of Hindi over Tamil and the Dravidian race was to have a separate state: Tamil Nadu for the Tamils alone. More than 1200 people were imprisoned in 1938, of whom two youths Thalamuthu and Natarasan (a Dalit) lost their lives due to the beatings in prison. Periyar was also incarcerated. On 13 November 1938, while Periyar was in prison undergoing rigorous imprisonment, a conference of Tamil Women held in Chennai conferred on him the title of 'Periyar' (lit. Great Man).

On 29-12-1938, as Periyar continued to remain imprisoned, the Justice Party elected him its President in the Provincial Conference held in Chennai.

Periyar viewed the imposition of a language as a cultural subjugation. This lead to his demand for an independent Dravidar Nadu, a nation for the Dravidians. He believed in German scholar Max Mueller's (now controversial) Arya Invasion Theory that postulated the subjugation of Dravidian races (who Periyar identified with 'lower' castes) by an invading Aryan race (who Periyar identified with Brahmins and 'upper' castes). Periyar's anti-Brahmin views were based on the idea that all Brahmins are upper class Aryans and all non-Brahmins are lower-class Dravidians. He accused Brahmins of subjugating the Dravidian civilization. It must also be noted that his concept of Dravidians was not based on the purity of blood or on ethnicity. It was about an Aryan caste-ridden society being contrasted with an egalitarian Dravidian society based on the humane traditions of the Tamils.

The Congress Ministry in Madras and seven other Provinces resigned on 29 October 1939, following the outbreak of the Second World War. They protested against the British rulers involving India in the war without having consulted the High Command of the party. Consequently, the Governor and Governor General requested Periyar to come and form the ministry. The offer came once in 1940 and again in 1942. Even his friend C.Rajagopalachari personally requested Periyar to accept the offer and assured him of outside support. Periyar refused it on both the occasions. He reasoned that if he accepted power his aim of annihilating the caste system would receive a set back.

On his north-Indian tour, Periyar met Dr. B. R. Ambedkar in Mumbai. Both of them also met the Muslim League leader M. A. Jinnah on 8 January 1940. Periyar spoke of his commitment to create an independent State known as Dravidar Nadu. On 21 January 1940, the Madras Provincial Government, under the Governor's rule abolished the compulsory study of Hindi in schools. Congratulations poured in for Periyar's efforts, even Jinnah had sent a congratulatory telegram.

Because of Periyar's relentless struggle, the degrading practice of serving food separately for the Brahmins and the 'others' in restaurants at railway stations was abolished in March 1941.

In its provincial conference held in Salem on 27 August 1944, the Justice Party was officially renamed the Dravidar Kazhagam. Its members were required to drop the caste suffixes of their names. They were also asked to give up the posts, positions and titles conferred on them by the British rulers. A resolution was also passed that members of this movement should not contest elections, thereby transforming a political party into a mass socio-cultural organization.

The Dravidar Kazhagam flag was adopted in 1946. It was a black rectangular flag in the ratio 3:2, and there was red circle in the middle of the black background. Black represented the indignities to which the Dravidians were subjected by the Hindu religion, and the red stood for the liberation of the people. On 15 August 1947, Periyar said that the Tamils must observe Indian Independence Day as a day of mourning, because he felt that it was only a transfer of power from the British to the Brahmins and north-Indian merchants (Baniyas). The Congress Government of Madras Province banned the Black Shirt Volunteer Corps created by Periyar in March 1948. However, it only made the Dravidar Kazhagam more popular. As a result more than a lakh people participated in their conference held on 8 and 9 May 1948.

Periyar deeply regretted the assassination of Gandhi by a Brahmin. In a obituary, Periyar suggested that India should be named as "Gandhi Nadu (Nation)" in honour of Gandhi's martyrdom.

Because Hindi was reintroduced in the curriculum in June 1948, Periyar launched the second anti-Hindi agitation at Kumbakonam on 10 August 1948. Unable to bear the public pressure, the Government dropped the scheme.

Meanwhile, power-crazy elements in the Dravidar Kazhagam deserted it to form the Dravidar Munnetra Kazhagam in 1948. Periyar married Maniammai in 1949.

After the Indian Constitution was adopted on 26 January 1950, the Brahmins succeeded in striking down the Communal Government Order (GO). They argued that the policy of reservation in educational institutions for the disadvantaged communities violated the fundamental right to non-discrimination. The Madras High Court struck down the Communal GO saying that it was ultravires of the Constitution. The Supreme Court upheld this. Periyar toured all over Tamil Nadu. Rattled by his formidable opposition, the Constitution was amended for the first time by the Central Government. Constitution {First Amendment Act} was passed in 1951 adding the Clause 4 to the Article 15: "Nothing in this article or in clause (2) of Article 29 shall prevent the State from making any special provision for the advancement of any socially and educationally backward classes of citizens or for the Scheduled Castes and the Scheduled Tribes."

In 1952, Rajagopalachari introduced the hereditary education policy where children would study in schools in the forenoon and practice the occupation of their parents in the afternoon. Periyar was evidently against it. He viewed it as a perpetuation of the caste system and Brahminical domination. He vehemently opposed it.

His vociferous campaign against it was so powerful that Rajagopalachari had to quit the post of Chief Minister in March 1954. Consequently, Kamaraj came to power as the Chief Minister. He ensured that the casteist educational policy was withdrawn. Though Periyar had walked out of the Congress, he supported Kamaraj as an individual because of his commitment to people's issues. In fact Periyar labeled him the *Pachai Tamizhan* (Raw Tamilian).

On Vesak (Buddha's birthday), Periyar publicly broke idols of Pillaiyar in order to make people realize the idiocy of idol worship. He convened a Buddhist Conference at Erode. In 1954, he left for Burma to attend the World

Buddhist Conference where he once again met Dr. B.R. Ambedkar.

When the reorganization of states on linguistic basis took place on 1 November 1956, Periyar welcomed it.

Hotels in run by the Brahmins were called 'Brahmin Hotel' in the name boards, giving the impression that it was a superior caste. Periyar wanted to put an end to this and so launched a struggle. A symbolic protest was organized in Chennai on 5 May 1957 before one such Brahmin Hotel. Agitations continued daily and till 22 March 1958, 1010 volunteers had courted arrest. Periyar succeeded and the term Brahmin Hotel disappeared from nameboards all over the state.

More than 10,000 volunteers of the Dravidar Kazhagam burnt the provisions of the Constitution that helped safeguard the caste-system. More than 3000 of them were sentenced to undergo rigorous imprisonment. More than 15 people died in and out of jail due to incarceration. In fact, a special act was passed to convict them.

Periyar was sentenced on 14 December 1957 to undergo six months imprisonment in a false case where it was alleged that he had instigated his followers to attack Brahmins.

In February 1958, he undertook a tour of North India and addressed meetings at Kanpur, Lucknow and New Delhi where he propagated his ideology. In June 1960, Periyar asked people to burn the map of India to signify that the Union of India under the Central Government rule was a Brahmin rule that safeguarded casteism. Thousands of people were arrested for taking part in this demonstration.

Riding on the success of the anti-Hindi agitation that claimed hundreds of lives in the few weeks from 25 January 1965 to 15 February 1965, the Dravidar Munnetra Kazhagam won the state Legislative Assembly elections in 1967. C.N.Annadurai, who became the Chief Minister, went to his mentor Periyar and sought his goodwill and advice, apart from 'dedicating his ministry' to Periyar. Three of the achievements of his ministry where: the decision to rename the Madras state as Tamil Nadu; introduction the two-language formula of Tamil and English instead of the three-language formula of Tamil, English and Hindi; legalization of the Self-Respect marriage system.

In October 1967, Periyar accepted an invitation to visit Lucknow and address a conference of SCs, BCs, and Minorities.

In 1968, the Aryan epic Ramayana was burnt all over Tamil Nadu has a mark of protest against domination and oppression by the Aryans.

In 1970, Periyar started the fortnightly magazine *Unmai* (Truth). The then party General Secretary K.Veeramani released the first issue. In 1971, Periyar launched the English magazine Modern Rationalist.

In 1969, he had announced the agitation for entry into the sanctum sanctorum of temples. So far, only Brahmins had been priests and performed worship in Sanskrit. The Supreme Court Judgment on this topic was ambiguous and interpreted in favour of the Brahmins. Even in the Conference for Eradication of Social Indignity held in Chennai on 8, 9 December 1973, he exhorted people to fight for equal rights to enter sanctum sanctorum of the temples irrespective of their caste. On 19 December 1973, he delivered his final speech at Thayagaraya Nagar in Chennai. It was a clarion call, urging the people to fight for social equality. He took ill on the next day and breathed his last on 24 December 1973. Millions of Tamilians mourned his loss. It was the end of an epoch, a great man who was living history.

In his life period of 94 years, 3 months and 7 days he had totally traveled 8,20,000 miles for his propaganda work. (This distance is 33 times the diameter of earth and 3 times the distance between earth and moon). He had totally addressed about 10,700 meetings. The total duration of his lifetime speechings was 2,14,000. If it were allowed to run on a cassette it would run continuously for 2 years, 5 months and 11 days.

## 3.2 Translations of the Speeches and Writings of Periyar related to Untouchability

In this section we have provided the relevant, translated excerpts from the speeches and writings of Periyar that deal with untouchability.

**Temple Entry**

Few people feel uneasy with me, but I will not change my policies or principles. Even among the gods, they teach us differences. For example, everybody can touch the Hindu God Pillayar. Anyone can pray to him. In the Kasi temple in Uttar Pradesh, anyone can pray even in the sanctum sanctorum. But the cruel rules prevailing today in the temples in our place are: this man should not go inside, this man should not pray, only man of this community should pray, that too in a particular language and so on and so forth. For these social differences and ill treatments, we spend a lot of money and lose our self-respect.
[Hero of Vaikom, 1923]

**Conversion**

Who are the one crore Christians of our country? Who are the seven crores Muslims of our nation? Due to the atrocities by Brahmins and our religious leaders, they have been forced to take up other religions. They are our brothers. They are not from Jerusalem or Arabia.

How are we, the 24 crores people, respected? Muslims call us '*kufr*'. Christians call us 'ignorant of religion,' Whites call us 'coolies,' Britons call us 'uncultured' and Brahmins call us 'Sudra.' We do not feel disgraced by such things. In the name of religion, we are denied entry into certain streets even today. In the Vedas there are only four castes, but now there are over 4000 castes, these castes still strive to maintain their identity.
[Hero of Vaikom, 1923]

**Untouchability**

Things like untouchability, unseeability, not talking, not coming close has not left anyone in our country. Not only is it customary for a person to say that one below him [in the caste hierarchy] is untouchable and unseeable, and the same person is untouchable and unseeable to one above him, but all these people taken together are untouchable and unseeable, and unapproachable to the Europeans, which is the caste that rules us.

In this manner, to talk of eradicating untouchability is shameful. It is not only about the emancipation of the *Panchamas*, it is not only about removing their difficulties, but it is about removing the disgrace and atrocities inside each one of us: the notion of untouchability.

When this is said, they feel surprised: "Ah! Removal of untouchability! Can the Panchamas be allowed into the streets? Are they to be touched? Are they to be seen?"
[Kudiarasu, 21-6-1925]

**Untouchability**

When they say that some people should not walk in the streets and that some people should not be seen, we must think of what these people feel in their hearts and on what basis they say such things. People of other religions like Muslims, Christians and Parsis can walk in the streets. A pig, dog, cat, rat can walk in the street; these can be seen while eating or at other times. The Hindu Tamil Ian who was born and brought up in Tamil Nadu, and who has made Tamil Nadu his own for thousands of years, is told by another: 'You should not walk in the streets, don't come in front of me.' How can a human being digest this?

This is the philosophy of the Vaikom Satyagraha and the Gurukulam struggle. Remember that only based on this philosophy, Mahatma Gandhi carried out a Satyagraha in South Africa, and the Kenyan Boycott of Imperialism Day was celebrated.
[Kudiarasu, 5-7-1925]

**Communal Representation**

More than the communal representation for non-Brahmins, the communal representation for Untouchables is very important—we will say this even from the *gopuram* of the temples. They have not attained development in education, or employment, or in other streams of public life in proportion to their population in the society. As a result of this, these people who are one-fifth of the nation, forget the welfare of the nation; in search of the Government's grace, they go and fall into another religion and sprout as our enemies. The selfish do not worry about this. It is not just to say that it is betrayal of the nation. If the responsible public ignores it, I will say that it is a greater sin.

If this society had been given communal representation at least some 25 years ago, will we have this kind of difference of opinion, lack of unity, British tyranny and Brahmin cruelty in our country even today? The man who should not walk in the street, the man who should not be seen, the man who should not know his own religion, the man who should not see his own god: can such categories exist in India? Everyone with a social outlook must think about this. They need to implement this in the state conference and bring it into force because it is the duty of the patriots.

[Kudiarasu, 8-11-1925]

**Untouchability**

It is the foremost duty of non-Brahmins to abolish untouchability. Only the advancement of the Untouchables and the non-Brahmins is the advancement of the nation. Only the sorrow of the Untouchables is the sorrow of the non-Brahmins who form 97% of the population.

Only by the abolition of the untouchability, the non-Brahmins can be saved or redeemed. Only after the abolition of untouchability, our nation can achieve real *Swaraj* (Self-rule). So, all people who are worried about the abolition of untouchability and those who are called

Untouchables should come to Kancheepuram, conduct the conference and derive its advantages.

[Kudiarasu, 15-11-1925]

**Swaraj**

The Brahmins think that Swaraj is a means of asserting themselves as superior and it is a method of exploitation by establishing that the others are low. Does Swaraj mean using the common, uneducated people to conduct agitations and thus get better government posts and better salary for themselves?

The evidence that it is not Swaraj at present is clear from the fact that people live without self-respect, without equality; they keep stating that one is high and the other is low; majority of the people are termed as Untouchables, they are degraded worse than animals, they suffer all cruelties without any liberty and so on.

[Kudiarasu, 15-11-1925]

**Communal Representation**

Even after the British ruled for over 150 years, we have 7 crore people who remain Untouchables, they have not yet obtained or qualified for freedom. So, they are completely dependent on the government. The 24 crore non-Brahmin Hindus who have lost their equality in the British rule are termed as Backward Classes. They have to depend on the 3% Brahmins and the government for their living. What is the reason for all this? Is it due to the fact that they were not qualified to get these equalities by birth? Or, is it due to the fact that those in the higher positions cruelly ensured that these people remain undeveloped?

If the non-Brahmins and the Untouchables had been given communal representation at least 50 years ago what would be their position now? The posts held by B.N Sharma, and C.P.Ramasamy Iyer would have been held by [the Depressed Class leaders] M.C. Rajah and R. Veeraiyan.

If R.Veeraiyan had been in the place of C.P.Ramasamy Iyer, will they impose Section 144 [curfew/ embargo] in the streets of Palakkad? Would they [Untouchables] be punished for visiting temples? What is the reason for their not getting such posts? Is it due to lack of qualification? Is it not due to the absence of communal representation?

How did Mohammed Habibullah Sahib and Dr. Mohammed Usman get reputed positions with salaries of Rs.5,000 to Rs.6,000? How did they become qualified to get these posts? Is it not because of communal representation? Let the public think about it! Depressed Classes and Backward Classes cannot get equality unless they have communal representation.

Before the British came, were the Backward Classes, Backward Classes? The only reason that can be attributed is that the Brahmins pleased the British to such an extent that the Brahmins became forward class and obtained all the rights and comforts. However, the non-Brahmins who were not like them became Backward Classes. How can the Backward Classes improve without communal representation? Don't they have intelligence, strength, and education? Why have they not got their rights?

Indians are of several castes. They are divided into eight categories. If the government stops paying fat salaries it would not be difficult to give communal representation.
[Kudiarasu, 22-11-1925]

**Press**
The policy of the newspapers should be primarily for equality and improvement of the non-Brahmins and Untouchables. The direct administration of the newspapers should be in the hands of a few reliable persons who work only for the press and are not a part of any other public work or trade or government. The board that carries out this work should be renewed once in every five years.
[Kudiarasu, 22-11-1925]

**Communal representation**
For the sake of development and unity of the nation, it is resolved in this conference that in all the political posts, the Government must take into account the population of each of the three categories: the Brahmins, non-Brahmins and the Depressed classes and should follow communal representation in all posts. This conference will take attempts to pressurize the government to follow this.
[Kudiarasu, 6-12-1925]

**Communal Representation**
How do we maintain the Untouchables? Even after 200 years of British rule, can an untouchable stand in an election against a caste-Hindu and win? Does any Indian think of making the Depressed Classes, who are one-third of the population, stand in an election, support and elect them? Have anyone of the leaders shown this way? So, the only state of mind in which everyone lives is, "I will not give my support for your equality and if you wish to attain equality, we will destroy it." Is this human?

India will never get complete freedom until there are equal rights for every sect of people. If everyone has to get equal rights, the only medicine is communal representation. Trying to overrule and evade without providing communal representation is just like thinking that closing it without applying any medicine can cure a sore in the human body.

The people cannot realize true equality and freedom unless communal representation is given. In these days if anyone talks about the emancipation of his caste or about the development of the Depressed Classes he is termed a traitor. This situation not only shows the degraded state of the Backward Classes and the Depressed Classes, but also the supremacy of the high castes.
[Kudiarasu, 13-12-1925]

**Religious conversion**
We see that the government officials who were on the side of the Brahmins at the beginning of the Satyagraha are now against the Brahmins by holding the hands of the so-called Untouchables.

At that time, even we felt surprised about the power of the Satyagraha. We observe this kind of power because we patiently experienced the difficulties arising out of Satyagraha. If we had plunged into it with force, or in anger, or in malice, we could not have seen such powers at any point of time.

The objective of Satyagraha is not that we should walk in the streets in which the disgraceful dogs and pigs walk. It is just that in public life, there ought to be no difference between man and man. That philosophy does not end with walking on this street. So, it is the duty of humans to prove the rights in the street, in the temples also.

When Mahatma Gandhi and the Travancore Maharani met, she asked the Mahatma, "If we throw open the road now, you will immediately take efforts to enter the temple."

Mahatma Gandhi said, "Yes. That is my aim. But till I learn that the people who demand the rights to enter the temple have enough patience and forbearance and are prepared to make the necessary sacrifice, I shall not enter into that effort. Till then, I shall keep doing the work necessary for that. The government [officials] proved that it was only the Brahmins, and not they who were opposed to the Vaikom Satyagraha.

To go to another religion in order to attain human rights is a very disgraceful thing. If there is such a necessity, one can go to Christianity or Islam; I do not have any liking in going to the Arya Samaj.

Because by going to the Arya Samaj, not only should the meaningless sacred thread be donned, but the meaningless ritual, *sandhya vandanam* (an evening ritual performed by the Brahmins) must also be performed. Those who had donned the sacred thread and performed the sandhya vandanam at one point of time are today, the enemies of our freedom and reformation. If you think that you should also not come to that state, surely do not join that group.

[Kudiarasu, 6-12-1925]

**Hindu Mahasabha**

There are several branches of the Hindu Mahasabha. In all these sabhas and its branches, the Brahmins say that they are superior by birth and that the notion of untouchability that is mentioned in the Sastras should be maintained because they cannot disrespect the Vedas. T.R. Ramachandra Iyer, Chairman of the Chennai branch of Hindu Mahasabha is a monument who practices the concepts that Brahmins are of high caste, untouchability is ordained in the Sastras, and Untouchables should not walk in the streets and are not to be seen. The Kumbakonam chief is not any less since he has passed a resolution in his branch of the Hindu Mahasabha that untouchability should not be abolished because it is related with the Sastras.

Even Madan Mohan Malaviya who spoke about Untouchables with tears in his eyes has never accepted the fact that there is no difference between him and an untouchable. This is because he strongly believes in the Varnashrama dharma.

[Kudiarasu, 13-12-1925]

**Sacred Thread Identity**

Like how a new born child of the lower caste, is sparkling clean at birth, and as it grows, because of not bathing regularly, because of wearing dirty clothes, it is also called as a member of the untouchable caste, like that even though a Brahmin child is also born like ordinary children, after some time, when it dons a sacred thread it attains the status of being called upper caste.

But, if such a child is born in the homes of those who are not lower castes, even if it is not bathed properly and it wears dirty clothes, it does not attain the state of a low caste. We see that even today. Just like an unclean man who subconsciously feels that he is a little unclean, a man who wears the sacred thread feels that he is an upper caste. Not only that, but he also feels that others are lower caste.

How will the low status of a man improve because of wearing this sacred thread? Do all the non-Brahmins who wear the sacred thread, have parity with the Brahmins? Now, in our country the Singurao, Chettiar, Asariyar, Kuyavar, Vanniyar, Sowrashtrar, Thevangar, Komuti, Raju, Vaaniyar Pandaram, some sections of the Naickers, and the Valluvars (who are still called Panchamas): all these non-Brahmin brothers wear the sacred thread. Have they been accorded status like Brahmins?

Today, don't we see many cheats, prostitutes, lepers who are evidently Brahmins? What is the belief that this [sacred thread identity] will not become like that? In our country, even before this, in order to create equality, several great people like Sri Ramanujacharya put sacred threads to several Depressed Classes and through rituals, made them higher. Those who were strong in cunning and slyness, became Brahmins and are worse than the Brahmins before them; they have become the god of death for our society and our nation. Though they donned the sacred thread, those who were powerless remain as they were.

It is our opinion that there is no other work that could create greater blemish for our nation. What we understand is that for India's development, as important as it is for untouchability to be eradicated, it is more important for the sacred thread to be eradicated. In a sack of rice, if there is a measure of one *padi* (1/80 of a sack of rice) of stones; after spreading it on the floor, is it easy to pick out that one padi of stones, or is it easy to leave the stones and keep picking up rice?

[Kudiarasu, 27-12-1925]

**Self-Respect**

If we keep on one side of the balance (a measuring scale), the atrocities done to us by the government and on the other side the atrocities done to us by others [Brahmins], only the atrocities done to us by others would weigh more.

If a person with human form has no right to see (his) god and worship it, then how can he have any self-respect? What is the use of Swaraj in this social setup? What is the matter whatever the rule? The social setup has made a set of people remain without any self-respect. If these demons of society get Swaraj will it do any good to people leading a life without any self-respect? Seeking for, or attaining Swaraj is of no use when people like us do not have even the right to worship, see the god we believe in and walk in the streets; worse, we should not even be seen by them!

The foreigners (British) who ruled us did not ill-treat us like this; without their rule we would have continued to live in cruelty without even realizing the social injustice faced by us. We would not have had the Self Respect Movement.

The enemies of our Self Respect Movement are not the British, but only the Brahmins who are outsiders and are ruling us. Since the Brahmins are in favour of the government, the British are not bothered about the cruelties that the Brahmins subject us to. It is important for people to say that untouchability should be abolished, and that everybody should wear the *kadhar* (handspun cotton) dress for Self Respect. But the chief of the Self Respecter's duties is to fight against men who don't allow other men to walk in the same street, come near them, be seen by them, see god or worship him. What is the use of independence for these people?

Further, they cannot earn more because they simply labour for other people without any benefit. They cannot even lead a life of self-respect. A set of people is trying to use its cunning to exploit these people and continues to practice cruelties on them. In such conditions, the words 'Independence', 'Swaraj', 'rights' are words that will only give us disgrace and cruelties, so they are baseless. If anyone really wants to toil for the welfare of the nation, then, let him first work towards making every individual lead a life of self-respect.

[Kudiarasu, 24-1-1926].

**Caste atrocities**

It is of primary importance that the cruelties based on caste and the arrogance of high caste is abolished at first. Keeping all these faults intact, and yet voicing support for the rights of the South Africans and writing as if one has all concern for the Untouchables and shedding crocodile tears is only an atrocious act. It is not an act with any conscience. In these circumstances, do we need Swaraj or independence? Or, do we really need only advancement without self-respect? The common man should reflect upon it!

[Kudiarasu, 31-1-1926]

**Communal Representation**

Those who really want the nation to attain development will certainly think that it lies in the betterment of Depressed and Backward Classes. Nation here connotes every caste and community, not a particular powerful community. There are only 3 castes: Brahmins, Non-Brahmins and the Depressed Classes. They live with faithlessness and are unhappy among themselves. The Government need not talk about Andhras, Tamils or Kannadigas because everything falls within these 3 categories. Also problems would be settled if these three divisions are given communal representation. The Government need not feel that this will give raise to several subcastes.

Even among such subcastes, communal representation can be given to the subcastes. In Britain, there are only 3¼ crore people but they are governed by 700 members. No, the inclusion of all sections of people into power and politics is not difficult. In fact, only that can give the nation equality and real freedom.

[Kudiarasu, 14-2-1926]

**Communal Representation**

At present, Swaraj means only the rule of Brahmins. Even in the British rule, we see that people of specific castes cannot walk in certain streets, they should not use ponds for water and several other types of atrocities are practiced on them. If Swaraj, which is Brahmin Raj, is in force, they will not be afraid to torture these people. How to trust them?

Over seven crore people who are very calm are treated worse than animals by the Brahmins. Why should they stay like this? Can we change the calm people into animals and the dishonest cheats into gods by not giving communal representation in politics and employment? Through power and money, the cheats always use cunning means to pull down others and come up. Hence there is always a perennial problem of difference and discrimination at its height.

With the single fear that communal representation could come to practice, the Brahmins always tolerated and endured the government. If the government was not amiable to them, they tried to replace it. So, the government also feared the Brahmins and adjusted with them. If we rule our nation at any point of time it can be only through communal representation. This should be encrusted in diamonds.
[Kudiarasu, 18-4-1926]

**Hinduism and Caste**
Hinduism means keeping the signs of holy ash, or sandal etc., on the forehead and calling each caste as high or low.
[Kudiarasu, 30-5-1926]

**Origin of Caste**
Presently, Brahmins in this country, under the garb of patriotism have been preaching through their speeches, writings and through others that I am a sectarian and that I instigate sectarian clashes. You should think a little whether sectarian differences and enmity was created because of me or because of the Brahmins of our country!

If you will please get up and come with me—if we go around this village several boards which say, 'This is for Brahmins'; 'This is for Sudras'; 'Panchamas, Muslims and Christians will not be given food, snacks or water here'; 'Sudras should not draw water from here'; 'Sudras should not bathe here'; 'Sudras will not be admitted into this school'; 'Sudras should not read these topics'; 'Only Brahmins can go so far—Sudras should not go beyond this point'; 'Sudras should not reside in this street'; 'Panchamas should not walk in this street' are placed in every coffee shop that the Brahmins own, in every hall, in every tank and temple, rules are created, and people are divided; aversion and dislike are created; disgrace is created: Am I the one doing all this? Or, is it the Brahmins? Observe this.

The Buddhists and Jains eradicated the cruelty of Brahmins and said people are equal; love and brotherhood alone was God. The Brahmins could not accept it so they destroyed all these efforts.

At the time when the Brahmins said that they were upper-caste, those who were powerful refused to accept it and were about to use force, so the Brahmins cunningly said, "You shall be the Kshatriyas and rule, we will be your ministers and give you advice," and thus they cheated them and brought them within their influence.

When the wealthy and the influential asked, "How can you become a higher caste?" The Brahmins told them, "You shall be the Vaishyas, we also place the sacred thread on you, just like ourselves" and cheated them too. Later, the majority who were farmers and artisans were made to work for the three sections above them, including the Brahmins. When several among them did not accept this, the Brahmins said, "We have kept a section lower to you, you are their masters. You can treat them as you like." Thus they made the peaceful, guileless, non-aggressive people, to be called the Panchamas and exposed them to the Sudras, and thus cheated them. Finally, a voiceless section is suffering as Untouchables.

Today, you are considering only these people—who did such cruelty—as the head of your religion. If such cruelty has to be removed from us, we must gain dominance in religious and political matters. After gaining power, we must try to attain our self-respect. It is my humble and decisive opinion that whether it is for a society or for a nation, self-respect is more important than Swaraj.
[Kudiarasu, 15-8-1926]

**Biographical**
From the time I was seven years old, I knew about caste differences and untouchability. Though I was born in a very orthodox family, I criticized everyone and mingled with everyone without any distinction, so I was not allowed to get into the kitchen of my home. Except my father, all the family members used to wash the tumblers that I used. The persons who used to be jealous of my family's orthodoxy became peaceful or happy on seeing me.

[Kudiarasu, 22-8-1926]

**Self-Respect**

To fight against untouchability, one should have strength of mind, tolerance, and fearlessness since the majority of enemies of the self-respect movement are Brahmins of our own country and it is not an easy task to win over them. For this, we need unity and sacrifice. I can expect all these from you. But what more you need is self-confidence. Above all, the blind beliefs ingrained in you must be lost.

The first reason is your feeling that the Brahmins are superior by birth and your belief that there is someone inferior to you. This concept should first get out of your mind. In all your family rituals, the thought that the Brahmins are higher than you should be wiped away from your mind.

Please campaign against these. All thee people who believe that some others are lower than them will not come to help you. You need proper self-control. In every village form youth associations in which Brahmins are not allowed membership. So, do this. We need money for all this, because nothing can be done without money.

Since money or material is very important, you should even be prepared to beg to get it, but you should maintain the accounts properly. Enemies will try to falsely blame anybody who works powerfully, fearing them would not be of any use. We should bear all difficulties and loss. We should work hard and sacrifice anything and not mind the aspersions cast on us by our enemies or mind our names getting spoilt or defamed by our enemies but work without fear. This is the biggest form of sacrifice one can do.

Before getting from up the bed, everyone should think, "What should I do for self-respect?" If you have not done anything for self-respect consider that day to be wasted. This is the responsibility of youth. You should consider it a penance to give your life for self-respect.
[Kudiarasu, 12-9-1926]

**Congress and Untouchability**

What is the poison of untouchability in the Congress? Brahmins refuse to look at us and the Untouchables. You would have seen this in the context of the *Gurukulam* problem. When the Congress made a resolution that there is no difference about high/low in birth, Rajagopalachari, Rajendra Prasad and Sastri resigned from the party. What can one think of the Brahmins? Often, Srinivasa Iyengar said, "Don't mix up "politics" and "untouchability". Acharya said, "If a non-Brahmin child sees a Brahmin child eating, I would starve for a month" and ran away from the Congress Committee. From this, one can understand the untouchability prevailing in the Congress.
[Kudiarasu, 26-9-1926]

**Varnashrama Dharma**

More than 22 crores of our people have become Sudras and Untouchables by birth because of Varnashrama dharma and are treated worse than animals, insects, and worms. Our ancestors for over many thousands of years have fought against this untouchability and degraded status. To protect the Arya dharma, the Aryans, who are Brahmins, use their sabhas like *Varnashrama dharma Paripalana Sabha* and *Arya Dharma Paripalana Sabha* to maintain and establish that according to Varnashrama dharma they are higher castes and that Sudras are born to prostitutes, and Panchamas (Chandalas) should remain so (as Untouchables), and serve them without any reservations.

Brahmins are afraid that Sudras and Panchamas may achieve independence and equality through the Self Respect Movement which is against Varnashrama dharma. They want the Arya dharma to be protected and at the same time they do not want the Sudras and Panchamas to be wiped away. By their cunning and wickedness Brahmins create a lot of hostility and try to destroy all those movements, which work for abolition of Sudrahood and untouchability. They pass resolutions that they would never accept the abolition of Varnashrama dharma. All movements or sabhas of the Aryans are only to support and establish the concept of untouchability and Sudrahood. For them it is a sin to talk with us, see us or touch us.

They call our Self Respect Movement or any movement that tries to wipe away Sudrahood and untouchability as blemished movements. They say that it is heaven if a Sudra drinks the water in which a Brahmin has washed his feet. But, the Brahmins wrap their sacred thread around their ear when they talk with us or when they attend nature's call.
[Kudiarasu, 19-12-1926]

**Self-Respect Movement**

If true cadres, firm in ideology, campaign with integrity alongside honest leaders, certainly we would become independent in a year. For this, firstly no one should accept *kadhar*; secondly we should not accept highness or lowness attributed by birth. Anyone who has difference of opinion over these two ideologies is an enemy of the Self Respect Movement. These should be mentioned as a solemn pledge in everyone's membership form.
[Kudiarasu, 26-12-1926]

**Temple entry**

Why does the power of god run away when some people see god? Why do the temple and the god get polluted when some people visit the temple? Why does god die when few people touch him? Why do the gods, who have the same name and form, have different powers in different villages? The gods in the cities of Kasi, Jagannath and Pandaripur do not die if anyone touches them. People who go to these temples worship by themselves: they touch god directly, they sprinkle water on their heads, they shower flowers and they pray. But in our country, the same god with the same name dies when we touch it. The gods in the temples of Srirangam, Chidambaram, Perur, Bhavani, Kodumudi, Tiruchengode, Karur and such holy places do not die if the Nadars visit the temple and worship them. But the gods in the districts of Madurai, Tirunelveli, and Ramanathapuram alone die if the Nadars worship them. Why does the power and life of the gods and the power of the temples differ from village to village? And why does the Brahmin alone have the power to give a fresh lease of life to that god?

[Kudiarasu, 9-1-1927]

**Equality**

For equality, liberty and brotherhood to grow in the minds of people, feelings of birth-based difference must end. Communal clashes add fuel to the fire. This sort of ideas should develop among intelligent people.

[Kudiarasu, 1-5-1927]

**Congress and Untouchability**

The Congress, which should be open to all castes, is trapped in the hands of one or two Brahmins and has become a Brahmin Congress as reported by the magazine '*Tamil Nadu.*' When will it become common to all the castes? How many Christians are there in the Congress? How many Muslims are there in the Congress? How many Panchamas are there in the Congress? How many Hindus are there in the Congress? In such meetings are the majority from the Congress? Have they come away from the Congress? At any time, was the Congress in the hands of the non-Brahmins? At any time was it a non-Brahmin Congress?

[Kudiarasu, 5-6-1927]

**Hindu Mahasabha and other Hindu organizations**

Will the *Varnasharma dharma Paripalana Sabha* that is said to be common for all Hindus make Varadharajan,

Kalyanasundara Mudaliar, or me as its members? The few who run the institution say that it is the representation of the 24 crore Hindus of India and they teach Vedas, Sastras and the ceremonial recitation of the sacred texts.

Tradition and law also accept it. Great men like Jayavelu, Muthuranga Mudaliar, Adhi Narayana Chettiar, and O. Kandasamy Chettiar accept it as the representative, but does it mean great men like Varadarajulu, Veeraiyan, and M.C. Rajah accept it?

[Kudiarasu, 12-6-1927]

**Gandhi and Untouchability**

We need not meet the Mahatma [Gandhi] again and ask his opinion on this; we also state that we do not have even a little doubt about our opinion. Without keeping the dialogue with the Mahatma as a basis, that rebuttal was written keeping in mind the speech of the Mahatma in Mysore and what he wrote. But for the Mahatma, the Varnashrama dharma and caste differences are birth-based. We know this because of the Mahatma's speech and his writings.

The Mahatma calls the Aryan religion, which is the religion of the Brahmins, as the Hindu religion, and he calls himself a Hindu. For his propaganda, he often shows the Aryan myths—the puranas that is Ramayana, Bharatam, Bhagavata, etc.— as his texts and quotes from these. We dare to think that these slippery puranas will somehow push the Mahatma into swamp.

There is no chance for the Mahatma to know about the civilization of Tamil Nadu, and the ancient customs of Tamil Nadu. During the time he was in South Africa, through a few Tamilians such as Valliammai and Nagappan, he might have known about the bravery of the Tamilians. We can tell firmly that he would not have known of Tamil Nadu and the Tamil people anything else except through the Brahmins.

Who is there in the Mahatma's mutt to tell him of the nearly 2000-3000-year battle between Brahmin religion, which is the Aryan doctrine and the Tamil civilizations and customs? So, when we preach that Tamil Nadu does not require the Aryan doctrines as far as it is concerned we need to say that it does not require the Mahatma's doctrines.

We, who argue that the people should not be divided due to their profession, they cannot be divided because it is not right in the nature of the world and humanity, how can we accept it if it is said that 'caste exists in birth, occupation exists in birth, this was created by great men, it was created by God, it was created by sages?'

This is not the age of sages and Krishna etc. and that we are not in that age, and they are not here in this age—whether it is right or wrong or true or false, there is no proper evidence or need or possibility and if each and every rationalist human understands it by him/herself, how will such grandma tales have respect?

The essay 'The Mahatma and Varnashrama' was sent to the Mahatma also. We sent copies to his intimate followers and asked them to inform the Mahatma of the same, we have not done anything without the knowledge of the Mahatma. But his faithful followers teach him, "Do not bother about such issues, if you give it importance, it will give credential to their dissent; it will increase the integrity of the man who is dissenting and his magazine. So, throw it in the trash."

It is our opinion that the Mahatma does not have the time or power to transgress it. This issue is important for the freedom of our country and our self-respect. The philosophy of Varnashrama—that too the philosophy that Varna and caste are birth-based—is the 'Yama' (lord of death). Because we have come to this firm decision, we struggle to eradicate it. We write only with the sadness that for the nation's liberation and self-respect, it is the duty of the man to realize this and we need to condemn it.

[Kudiarasu, 28-8-1927]

**Religion**

I am not an agent for anybody or a slave to any religion; but I am subject to two philosophies: kindness and intellect. I

feel that it is my duty, concern and wish to talk about these. But I leave these duties to the reasoning of one and all. They spend several lakhs of rupees to build temples, which are hundreds of years old. What do we get out of these? Think of it!

You are a low caste who cannot offer prayers to these gods. By these temples it is established that we are lower castes than these Brahmins. Till the stones of the temple, the temple and its gods exist, our lowness and degradation will be permanent. It will continue for generations after me. We built temples; we perform Kumbhabishekams to establish ourselves as the low caste and the Brahmins as the high caste. We give lakhs of rupees for these temples in order for the Brahmins to use it. They continue to enjoy this money for many generations; as long as there are stones and grass in the Cauvery River!

[Kudiarasu, 11-9-1927]

**Caste Identity**

The question is whether community [caste] conferences taking place in our country are suitable for the nation's progress. Several people say that this is harmful. But, it is my opinion that there is no conference that is not a community conference (that is, conference without communalism).

Every community holds conferences to demand its rights. To browbeat another community from rising up, and to see that another community does not get its rights, several people hold many conferences.

But our Nadar Conference is not like that. Without harming others, this conference is taking place to secure the Nadar people's rights and welfare. Community conferences and communalism were created only because of the harassment by the upper castes. We don't organize conferences to create dislike. We organize this conference just to say, "No need of hate; instead say all are equal." The word 'unity' that is often repeated in our country, and the problems are merely outward garbs. True unity will be created only when people can help themselves. It is necessary to foster the feeling of unity among one another. To work towards that, every society should have the freedom. Every community must develop self-respect, equality and self-esteem. The nation shall attain progress only through that. Self-respect alone is the real nationalism. Swaraj depends only on self-respect. All are equal; there should be no feeling of upper or lower.

[Kudiarasu, 9-10-1927]

**Hindu Religion**

This religion which professes that everybody should not get educated; if the low castes get educated in spite of the religious codes their tongue should be cut off; if they hear slokas, molten lead should be poured into their ear; if they still get educated, their throat should be cut off and several other types of punishments are given to them. (These punishments are also described in the Laws of Manu). Except for the British who had come to our nation, we cannot imagine any form of education. Even one person out of thousand will not be educated due to our religion. So the suffering and disgrace we undergo in the name of religion is inexpressible.

[Kudiarasu, 23-10-1927]

**Hindu Religion**

Further, the idea that Brahmins are upper caste and others are of lower caste is being pumped into our bloodstream. Because of this, a small sect of people have tortured the majority of our people and made them lose their self-respect and live a life of poverty and difficulties. So some of them ask us to forget the differences with Brahmins, they ask the Brahmins and non-Brahmins to join together to find a remedy for such cruelties.

When we do not have the right to touch them, eat with them, pray with them and even walk in the streets in which they live[5] and enter the temples that are built for them. How can they beg us to forget the differences; you think about these! According to research, at one point of time or another, each and every caste were rulers. So, the Laws of Manu, which calls the rulers of this nation as Kshatriyas, is wrong. Likewise trade was in the hand of each and every caste, so saying that the Vaishyas were the people who looked after business is wrong!

Further, as of today if the Brahmin accepts all the divisions, he would not have the sole authority over everything; so now the Brahmins say that in the Kaliyuga there are only two divisions: Brahmins and Sudras. If we see who is a Sudra, we see that a Sudra is the son born to a Brahmin and his concubine, and his only duty is to serve the Brahmin. How can one's mind accept this? Further, he calls another group of people as Untouchables/ *Chandalas* and they call the Christians and Muslims as *Mlechas*.

I have no other option except to call people of all castes other than Brahmins as non-Brahmins. This problem of cruelty to one caste of people by another has existed for the past 5000 to 6000 years. This difference is against nature; it is non-existent in any other nation. Above all, it is against humanity. They have fought for unity and love among people. Let the Brahmins live with us; they had after all come to India begging. But what do they do, they talk of dharma and sastra and put an end to all our movements. We have not written in our places, 'this is for Brahmins' and 'that place is for Sudras' but they have written so even in the hotels, ponds and wells constructed by us. The Brahmins say, 'this is for Brahmins' and 'not for Sudras' and so on. The Sudras are treated worst than insects, dogs and pigs!

We are treated worse than these! Do we treat them badly? What have they done except talking philosophy? Those who lack self-respect will not know these difficulties. Who will understand the sorrow we undergo every day?

---

[5] *The Kudiarasu issue dated 24.4.1926 tells the story of a boy of the depressed classes who entered the Brahmin street and was beaten with a slipper by a Brahmin who felt that the street got polluted.*

Some tell us that if we abolish the differences among the lower castes then this will also go, but the differences among the lower castes are only imitation of the Brahmins and not done with any motivation, but only due to ignorance.

I am not a researcher. But the very term *Jati* (caste) is not a Tamil word. All Vedas and the Manu Dharma, which are the basis of *Jati*, are not in any way related to our culture. They do not even exist in our language. Our people have not written them. Somehow they have come on to our minds as a problem. The word 'Hindu' or 'Indu' is not there in any of our languages. Researchers say that it is the name of a river. In some places, they say that it means 'thief.'

To spoil our movements, our enemies have made some poisonous and silly comments.

Not even a single line about Hindu religion is given in any purana or history. Only during the British rule, the term/name Hindu came into being! It is found in railway platforms and hotels. The words, 'Hindu Women' and 'Hindu Men' are found on lavatories. I will suggest what you should do. We should uproot all the problems that are given to us in the name of Hindu Dharma. Never allow Brahmin youth to interfere in politics. Even in the rule of the British you torture us; when they go away, what will be the problems you will be giving to us?
[Kudiarasu, 30-10-1927]

**Nadar representation**

Because of Ramanathapuram King's death, the Presidency of the Ramathapuram District Board has fallen vacant today. We suggest to the Government that in the present circumstance, that too as far as Ramanathapuram District is concerned, than leaving that position for open election, it will be a more intelligent thing to nominate a suitable person.

Everybody knows that in some parts of Madurai, and in the Tirunelveli and Ramanathapuram districts, caste arrogance and caste atrocities are rampant. For instance, the cruelty of the Brahmins need not be described in order for a person to know about it. The cruelties committed by the Saiva Vellalas of the south who come in the next rung, make the whole non-Brahmin society ashamed. Caught between these two communities [Brahmins and Saiva Vellalas], the difficulties faced by the people of other communities are unthinkable.

Several cruelties that never take place anywhere in the Chennai Presidency is rampant in the above districts. Our Vaniya Chettiar brothers, who are Vaishyas, do not have the right to enter a Saiva temple in Tirchendur. Likewise, in several temples in Madurai and Rameshwaram, our Nadar brothers who are called Kshatriyas do not have the right of entry. But, they have the right to enter temples in places like Palani etc.

People are denied even the right to enter several streets in these two districts. Just because one or two Brahmins or Saiva Vellalas eat/drink with a few Nadar brothers, we cannot say that such problems have disappeared.

Still, these cruelties committed in the name of God, in the name of religion, in the name of Dharma Sastras, are going to be established in the name of Government's rule, in the name of laws and in the name of court judgments, how can we escape from this? This is our question.

If someone says that concern is absent among people of the community that has been subjected to cruelties, rather than using it as an excuse, in the name of common good, those concerned with the self-respect of the country must make it their first duty to do the needful for the country's welfare and the birth-rights of human society. Only with this intention, we implore that the post of President of the Ramanthapuram Zilla (District) Board must not be elected to; instead, it must be filled through nomination. If the vacancy is going to be filled through election, those who have the mindset and opportunity to eradicate these cruelties will not get that position. Only if the government is going to nominate, it can search and find a suitable person for the said post.

If the people who have been through tremendous suffering and the Depressed Classes have to attain equality, they have to be given power and position; else, it is not an easy task for reformation to take place through other means. It is our opinion that the communities, which we consider as being harassed by others, are generously given power, position and employment; automatically they will attain liberation and equality.

So, the Government must pay attention to these ideas, and bravely come forward to annihilate caste atrocities and caste arrogance and give necessary power, position, employment and support to communities that have been victimized for a long time, create a just influence and salvage them from being oppressed.

We also warn the government not to tell the excuse that they did not get a suitable candidate. Because when there are suitable persons in that community for the posts of Member of Legislative Assembly, Member of the District Board, Chairman of the District Board, and for other positions; and when there are people who by themselves manage large estates and businesses that have an annual income of Rs.50,000 and Rs.1,00,000, if the Government says that there is no one for this post alone, no one will believe the Government's integrity. So we believe that the Lord Governor and Minister of Local Governance will not lose this opportunity to prove their integrity.

[Kudiarasu, 16-9-1928]

**Temple Entry**

One of the biggest holy shrines in Orissa is the Jagannath temple. It attracts thousands of devotees every year. Here, the priests are only from the barber community. They offer food in heaps to god. That food is sold outside the temple as *prasad* and devotees buy and eat it. In this temple there is no uncleanliness; the barbers who act as priests, do the job of hair-cutting in their off-time. There is no caste difference here. In this temple, after one eats if there is leftover rice, it is put back into the same container or it is served to another.

Anyone can use his hand to take food from the vessel and eat it. Anyone can go to the temple to touch the statue of god and pray.

In this temple no one can perform the *vratha* (fast) or *tarpana* or *thithi* (observance of anniversary for dead ancestors). Unlike other temples, here we find the statues of Balarama, Krishna and Subadra. They are brothers and sister; Subadra, the sister of Balarama and Krishna stands in the middle and the statues are made of wood. Usually the statues in other temples are with their spouses. They say "*sarvam jaganathan*" when they enter the outskirts of the city, which means "all are equal and there is no difference between anyone."
[Kudiarasu, 30-9-1928]

### Indigenous gods
In South India, 90% of the people do not know about Savagism or Vaishnavism. They only know about gods like Maari, Karuppan, Madan, Madurai Veeran, and Pechi etc. The creation of Hindu religious texts *Devaram* and *Prabandam* was only to destroy and ruin the principles of Buddhism and to establish caste, religion, god and the epics and puranas in society.
[Kudiarasu, 7-10-1928]

### Communal Representation
It must be said, "books that teach higher and lower must not be read." If someone violates and reads, the books must be seized. All the heads of religious mutts, who have the feeling of upper and lower must be caught and jailed. If the public launches an agitation, we must send these religious heads on exile.

All the jewellery, chariots and lands of the Gods must be seized, it should be sold and used to educate the uneducated, and provide livelihood to those who lack it. A brave government that will dare to do such activities should be brought to this country.

We blame the Brahmins more because only the fire that they ignited has spread among ourselves. If we win the right of communal representation in employment by fighting with the upper castes, we must also give the required share to those who are the lower-castes.

What greater reform could be there than giving equality to all the communities in the policy of the all-party meeting?

[Kudiarasu, December 1928]

**Support of Varnashrama Dharma**

It is unfortunate to recall that the chief of the Theosophical Society Annie Besant said, "Adi Dravidas are unclean people and are Untouchables." Not only that, she recently addressed the university students in Varanasi and spoke on caste differences and Varnashrama dharma. Untouchability is supported because of Varnashrama dharma. So how can brotherhood and unity be inculcated or developed?

[Kudiarasu, 9-12-1928]

**Casteism among Non-Brahmins**

Brahmins discriminate us, the different castes among the non-Brahmins discriminate each other; if the upper castes wish to be equal to others they should not discriminate the low castes. You should practice equality. The differences among us are crueler than the differences that the Brahmins practice (against us). The fire they ignited has captured us. If we fight for communal representation and succeed in it we have given the due share to the low castes.

The Nattukkottai Chettiars who earn lot of money as profit spent it on the education of Brahmin by providing them food. Why don't they do the same service to the Adi Dravidar children? To satisfy their hunger, these small children go and work as coolies. If such support is extended will they not be educated and become ministers?

[Kudiarasu, 9-12-1928]

**Abolition of Untouchability**

For a long time, we have been fighting to prove that untouchability is not justified. We have been explaining that it is an atrocity. Yet, it doesn't look like this problem has been practically eradicated. Where compulsion existed, there it has been removed to an extent. We see that this untouchability takes firm root when we appeal to humaneness/ compassion to end this evil. We are not in anyway lesser than others to show that this evil of untouchability is baseless.

Only because of a courageous struggle, untouchability can be abolished. They say that "God is all powerful, he is impartial" and yet, how shameful it is when they say that "only God is responsible for the cruelty against the Untouchables." It is also being said that most probably he was the one who created this untouchability. If that is true, then we need to first annihilate that God and only then go ahead with the next task. If he doesn't know about these unjust happenings, we need to annihilate him even sooner. If it is not possible for him to remove this injustice or to control those who commit these atrocities, there is no need for him to inhabit any world. It is justifiable to annihilate him.

If there is any basis for saying that god or religion does not give way to the eradication of untouchability, then it should also be destroyed with fire no matter who said it or what it is. Without being firm in the task, and just verbally speaking about the eradication of untouchability, at any point of time our nation cannot attain development because of this deceit. If one gets a wound in the eye, do we not immediately use medicine to cure it? If we say that there must be no pain or burning sensation, then it is only a way for one to rot and die.

Even women are being segregated as an untouchable society. The eradication of untouchability is our duty. There is no other shameless act than following the superstitions blindly because someone said so, because it is written somewhere instead of researching and coming to a conclusion that as far as intelligence is concerned there is no basis for untouchability. Untouchability is not dependent on intelligence or evidence. But, it is only based on foolishness, arrogance, and cheating.

If one touches shit, it is enough to rinse off. But if they touch a man they need to bathe. Is it in the name of hygiene? Or is it on the basis of arrogance? If they believe the purana that Nandan was made into one of the Nayanmars and Pannan was made into the Tirupaanazhwar, why don't they let the brothers and grandchildren of Nandan into the temples and let them worship? How disgraceful is it to not let the people of the clan of a man considered as one of the 63 Nayanmars into the temples, and to misuse the names of the Nayanmars and Nandan and Azhwars in order to eat tamarind rice and pongal under the pretext of worshipping 63 stone idols?

Also, if it is being said that they stink when they come near, who is responsible? There is no place to bathe. They should not come to the roads; they should not come to the common wells and ponds. They don't have washermen, they don't have barbers. With so much of injustice will they smell fragrant instead of reeking?

Take the Sankaracharya and close him in a room for fifteen days and never let him bathe. Later, see whether he smells fragrant or if he is stinking. Even if God himself is not washed using water for ten days, it will start stinking.

It is also being said that they are drinking country liquor and arrack and eating meat. Who doesn't drink liquor? Who does not consume meat? Are they the ones who drink all the liquor that is produced? Those who drink and roll in the streets can enter the temples with liquor oozing out of their noses, people of other religions like Christians and Muslims can go halfway inside. But they who are called Hindus must alone stand outside the temple. What greater cruelty can be there in our religion?

Daily we see what the pigs and the hens eat on the roads. People who eat pigs and hens are not untouchable. But those who eat the cow, which eats green grass and cottonseeds, are untouchable. What is suitable for intelligence and justice must alone be accepted. Foolish, meaningless matters must never be accepted under the name of anything. Listen! In other countries, people have started being indifferent to the atrocities committed in the name of God, or in the name of religion. That is why they have attained advancement. I believe that you will not backtrack from eradicating untouchability that is against intelligence and justice.

[Kudiarasu, 17.2.1929

10.2.1929, Removal of Untouchability Conference]

### Self Respect Conference

The first Self Respecters conference was held in a grand manner on 7[th] to 11[th] of Feb in Chengulpet. It was talked about for 2-3 months and several resolutions were passed. As there was one law for man and another for women, one law for Brahmin and another law for non-Brahmin, one law for high caste and another for low castes and people suffer the impact, the Self Respect Movement wanted to change this and make equal laws on these inequalities.

    The caste differences like low or high is a cruelty that should be abolished. This has been accepted by all people in several fields. But several people do not accept those evidences which are against them. On the other hand, in order to make the cruelties long-lasting they hold conferences in the name of religion, campaign by sending messages to the government which clearly shows that our country will not easily abolish the cruelties of untouchability. Many do not understand the high-low basis of such cruelties. Likewise the Mahatma supports the divisions/differences of Varnashrama dharma. Reformers of the Self –Respect Movement passed a resolution in which the difference in caste system, and religious differences that are shown through symbols must be stopped.

    I want to abolish the practice of caste names, titles, symbols denoting particular sects, because in my opinion these have caused among the people an identification of caste, creed, and the real qualities or intellect of a person. This brought about a difference, which was baseless. This is

one of the reasons for the lack of unity among people. Such differences and symbols cannot be accepted by anyone who wishes for the unity and equality of one and all.

[Kudiarasu, 24-2-1929]

**Caste Identities**

Since I struggled for unity, someone who became jealous of me wrote an editorial about the abolition of caste titles and religious symbols that exhibit differences among people and questioned my participation in caste conferences. I advice that each caste can hold conferences, bring in resolutions to abolish the practices in their caste that distinguish them from other castes and spoils their unity. I said that only for such purposes one should hold caste conferences. The Vaishyas can hold the conference and decide in that conference that they should not wear the sacred thread and should not call themselves as Vaishyas. Likewise, for castes like Kshatriyas and Nadars.

In the Naidu caste conference, I requested them not to call themselves or identify themselves as Kshatriyas. Unless such conferences are held and resolutions are passed, how will people know about it? Only by going into each caste conference, the disease of the particular community can be cured. Now they hold Adi Dravidar conferences, which have brought about feelings for their development. Thus if caste difference is to be abolished, caste conferences should be held where resolutions are passed to put an end to the practices of caste which distinguish them for others. For example, in the southern states, if we had not held the conferences of non-Brahmins we would not have been able to know about the cruelties practiced on us by Brahminism. By frequently holding such caste conferences, we can know the problems and sufferings of the other castes so that we can unitedly help them come out of their problems. Except the selfish lazy caste that wants to live and eat on other's labour, one will not think of caste conferences as a danger. Religion, Vedas, Sastras and puranas, which support caste differences, should not be followed. Though this resolution cannot be flouted, those who make a living daily in the name of Vedas, Sastras and puranas cannot attain satisfaction by these resolutions so they would secretly make all their cunning propaganda.

The importance of the resolutions is not to follow the principles and sources of the Hindu religion. Since most importantly Hindu religion, based on Vedas, Sastras and puranas, insists on birth-based caste divisions of low and high, it is different from other important religions. It is not surprising if one thinks that this resolution blames the principles of Hindu religion. Hindu religion is not really a religion, and if one thinks that such a thing exists, the evidences are based on the selfishness of a few and due to our people's ignorance. The practice (of Hinduism) is only a confused feeling and for this confused feeling to be annihilated the Hindu religion and the principles which it is based on must be annihilated.

For example, no one has replied to the question whether there is any religion like Hinduism. In *Gnanasuriyan,* it is said that in the language of some nation, the degraded people called Hindu. That has been mainly used by Muslim religious leaders who say by the term 'Indian' means those "who do not belong to their country". English called us Swadeshi and in Arab and Urdu we were called as Hindu. As there was no name to refer to the Indian People, as there was no name for us in religion or based on our nation or on caste, the name Hindu was used to unite all of us and the principles of Aryans were named as the Hindu religion. Thus, this was forced on every Indian. For instance, there is no evidence about Hinduism some 400 or 500 years ago.

[Kudiarasu, 3-3-1929]

**Hinduism**

From the time the feeling of abolition of caste differences and untouchability came; the Brahmins conducted a conference on Varnashrama dharma and insisted on the practice of untouchability. If the Government tried to abolish such practices through law, the Brahmins made resolutions that the Government should not interfere in religious matters. Common man fears the religious codes and shuns any reformation. Having faith in religion, if one wishes to make any reformation, no one has so far succeeded even a little till date.

Thus, only using this experience we take up the issue of finding the quality of the Hindu religion and found that this garbage had only garbage and they were useless garbage. Is the principle of Saivaism useful for intellectual research and is it based on the feeling of kindness? No, it is not founded on any of these! There are 21 crores of Hindus. Can these 5 to 6 crores of Depressed Classes, who are termed as Untouchables, be accepted as Saivaites by the Saivaites themselves? They will not even accept them sitting by their side and eating food. How can the notion that Saivaite religion is a religion of kindness and is it fit for research can be accepted? I leave these discussions to the public; let them come to some conclusions.
[Kudiarasu, 7-4-1929]

**Saivaism**
Only after the research of Vedachalam and Subramania Pillai, these Saivaites who treated others worse than dogs have changed a bit. Thus there is some reformation in Saivaites only in the past 20 years. Till then, these Saivaites were of the opinion that it is a sin to eat with non-Saivaites and their Saivaism lived within the holy ash and the rudraksha seeds used by the Saivaites after they were given *dheeksha* by their gurus. These Saivaites did not even accept people of other religions as human.
[Kudiarasu, 7-4-1929]

**Varnashrama Dharma**
In those times, Kings were instructed to not only create differences in education and religion but were also forced to keep the castes separately and preach the idea of low and high castes and about Untouchables and Unseeables and associate each caste with a profession and so on.

Further, laws, punishments etc. were only based on it. Unless the religious laws in the Indian nation which preach low or high caste or caste differences by birth are refuted and abandoned by intellectuals, leaders, great men and above all the common man; our nation cannot get true independence; the nation's self-respect or the religion's respect cannot be attained.

Now the Brahmins join together, preach and confirm about the righteousness of the Varnashrama Dharma. They state that as per the Laws of Manu, Sudras are created to serve Brahmins and they try to establish this. They passed resolutions that separate temples must be built for the Untouchables. How can Brahmanism be abolished when such Varnashrama dharma is practiced? How can equality be there? How can unity prevail? How can people live with self-respect? These people only spoil India by such acts.

[Kudiarasu, 14-4-1929]

**Swaraj**

The word '*desiyam*' (nationalism) is to cheat the common man for their better living. The Brahmins coined this word in order to please the foreigners. Thus this term is used only to cheat the non-Brahmins!

According to them, the term nationalism only means cheating others without any conscience and leading a life of comfort and pomp at the cost of other's labour. It is not for the growth of education, growth of intellect, growth of research or for the improvement of technology or for equality or unity or for one's own effort or to realize the truth or for not cheating others or for anti-slavery or for anti-untouchability or for non-denial of common road/lake/well to anyone.

The term 'nationalism' was only to practice all the above for the deterioration of the common and poor man. Instead of trying to stop the practice of caste and religious differences, these nationalists try to strengthen and establish it. Who protested against the Devasthana bill? Who is the one against equality? Who is the one supporting prostitution done in the name of god in temples? Who is against equal representation to one and all? They are in fact traitors of the nation and not nationalists as they call themselves to be!
[Kudiarasu, 19-5-1929]

**Religious Conversion**

If people leave Hinduism they should follow some other religion. What religion should be recommended to people? Christianity was more like the Brahmin religion, so I recommended Islam to people. Because, in Christianity there are untouchable Christians, Nadar Christians, Vellala Christians and so on. So even in Christianity there is differences followed based on Varnashrama Dharma; so, in India, it is Hindu Christianity! If the Depressed Classes embraced Islam they would get social equality in a very short period.
[Kudiarasu, 19-5-1929]

**Non-Brahmins and Depressed Classes**

Over the past years, even though I have been invited some four or five times to this Zilla Adi Dravidar Conference, due to several reasons I was not able to attend it; so this time, I decided that certainly I had to reach the conference somehow. The leader of the welcoming committee lavished praises on me. I will say that so much of praise, other than making me feel shy, doesn't contain truth. In the matter of removal of untouchability — if there is anything, however small, that I have contributed, it will be something that has been done for our benefit, and it cannot be counted as being carried out for your benefit.

Because, the common philosophy in social life, of 'yourselves' and 'ourselves,' does not have even small distinctions. For example, like how you are Untouchables, in the same way, even we, who are a little higher class than you, remain Untouchables to one class — that is the Brahmin caste which claims itself to be born from God's face and talks of itself as the lords of the earth. Although we are allowed to go farther than you inside a temple, we have to stand behind the Brahmin. When you come into the temple, the temple and the God gets polluted; in the same way, if we go to a certain place (in the temple) the temple and the God gets polluted. Dining in front of you, dining with you, dining in your homes becomes an evil deed and a sin — in the same way, dining in front of us, dining with us, dining in our homes is considered a sin for the Brahmins.

Even where we are addressed by our caste names, we are being addressed in a manner more humiliating than you. They address you as Paraiyars and Pallars. Yet, the words Paraiyar and Pallar only signify your work/occupation and the place that you inhabit. Because of these names, Paraiyars and Pallars are independent and not fit of humiliation. Whereas, the name with which we are addressed, Sudran causes disgrace on our very birth, marks us a birth slave to one, as a son of a whore by birth and it carries with it only disgraceful facts. Paraiyan, denoting a person only recognizes him as the son of his lawful parents. But, Sudran means whoreson, prostitute-son, concubine-son, birth-slave, purchased slave and is full of several such derogatory meanings.

A section of people like you – the one being called Panchamas, its existence has no place in the Brahminical religion – that is the Hindu religion. The Hindu religion does not have any caste less than the Sudras. But, a group called the 'Chandala' arises out of the living patterns of the people of the four varnas. That is – children being born to a Brahmin male and a non-Brahmin female, Brahmin who did not read the Vedas, doesn't cultivate fire, doesn't perform the sandhya-vandanam, and who doesn't perform such related Brahminical tasks becomes a Chandala. Several Brahmin pundits, with erudition in discourses and expertise in the Vedas and Sastras, have confirmed this. So, if somebody says you are a Chandala they must be accepting that you are born of the liaison between a Brahmin male and a non-Brahmin female or vice versa. Otherwise, they must accept that the Brahmins who have slipped from the righteous path are you. If this is so, as on this day, 99¾% of the Brahmin society must only be Chandalas. So, in this situation, logically you cannot be faulted since you have a place and plenty of evidence.

But our condition is unthinkable, it is a great disgrace. That is why apart from saying that we are lower and more disgraced than you in societal life, I also say that the efforts undertaken to remove these disgraces are done more importantly for our section. Apart from this, if anybody has the concern/ anxiety – at least in small measure – for eradicating the current evils for us or for you, certainly it is impossible for him or her unless they are courageously prepared to destroy that which is supportive of this demeaning condition. This is because, although it is the scoundrel nature of man that has caused this humiliating state, for the purpose of establishing this as a religious tenet they attributed it to someone called the all-powerful God and it remains in force now.

If this situation has to be slightly changed, the aforesaid religion and God come and interfere. So, the religion that is said to be the basis for this condition and the God who is said to have created this religion must be opposed; unless we are ready to answer that. We are courageous and prepared to destroy this– it is impossible by any other path. Apart from this you also need the feeling of self-respect. Why are we low? Why must we toil for one person? This feeling must come. You must be considered as human beings just like anybody else.

If the people of the village torment you or treat you in a humiliating manner, you must stand up to oppose it. If you cannot do it, you must migrate to other towns. If you do not find a livelihood even there, you must shake off this cruel religion and go to a religion that has equality. If this is also not possible you must at least go abroad as coolies. If you do not have the guts to handle such sure methods, I will say that the burden thrust on you will not be easily eliminated. You cannot succeed in any work unless you are willing to suffer, to throw away the regulations and to sacrifice your life. Moreover, if you expect someone else will come and help you it is great foolishness.

You must get the courage to use yourselves. Some persons are discoursing to you that by getting yourselves educated, earning money, bathing, abstaining from drinking (alcohol) and by not eating meat your lowness will go. I will not accept these. If these are the reasons for your low status, then why do the others who are also having the same ill qualities, instead of reaching your low position remain as Brahmins? The only reason for your humiliation and low condition is absence of dignity and self-respect. Consider yourselves as human beings just like others; accordingly be courageous to behave like that. Acquire the bravery to bear the difficulties that arise. Soon you will be freed.

[Kudiarasu, 16-6-1929]

**Self-respect**

Moreover, you need the feeling of self-respect. Why did we become Untouchables? Why should we call someone '*samy*' (Lord)? Why should we work for someone? This feeling must come. You must consider that you are also humans like the others.

[South Arcot Adi Dravidar Conference, 16-6-1929] [6]

**Resolutions of the Chengulpet conference**
1. There are no differences in castes by birth.
2. Religion, Vedas, Sastras or puranas that teach caste differences should not be followed.
3. The 4 divisions, Brahmins, Kshatriyas, Vaishyas and Sudras; and Panchamas as given by the Varnashrama dharma should not be accepted.

---

[6] *Mr.Muthusamy, President of the Conference Reception Committee spoke: "Periyar, who has taken a vow to remove the horrible devil called caste, the cruel disease of untouchability, slavishness, poverty from this country, has renounced his body, property and energy for this and is working night and day and through campaigning against the Injustce of Manu and Sastra's doctrines and has established equality and brotherhood." (Adalarasan, Thanjai; Thanthai Periyarum Thazhthapattorum, Periyar Self-Respect Propaganda Institution, Chennai, 2nd Edition, 1992, p. 19)*

4. Untouchability among people must be abolished and common lakes, wells, schools, street, lodges and temples should be given complete and equal rights to be used by one and all.
5. As these cannot be attained by campaign, laws must be made to establish these.
6. As caste difference should not be practiced all symbols and caste names should not be used.

In case of Untouchables, they should be given land, food, clothing, and books, good and free education. The poromboke lands must be given to them freely. There should be no discrimination in hotels and restaurants.
[Kudiarasu, 25-8-1929]

**Religion and caste**
People are fully aware of the fact that religion is not based on birth and it is due to their mentality and propaganda. Further they know that any number of religions can be created. They also know only these people would go and join any newly created religion.

When we look at the Hindu religion and caste, we see that caste is more powerful than religion; for caste is created by birth and it is said that it is unchangeable! But religion is based on principles and it can be changed according to one's mind and heart. So, birth-based caste must be completely uprooted, all people must join together and be made to live with unity and brotherhood; this is also known to people. So what is the problem in uniting people by destroying the caste and Hindu religion?
[Kudiarasu, 3-11-1929]

**Hindu Mahasabha**
For the poor people who have migrated to Malaysia, South Africa etc. for their living, Brahminism has followed them there also. Even in those countries, Brahmins exploit the poor people in the name of god, religion, Sastras and Puranas. Whatever evil they have done to the poor in India, they continue to do it in other countries also.

In the Hindu Mahasabha, several castes even among the Hindus were not taken as its members. For example they did not accept the membership of Adi Dravidas, Maruthuvars and also Nadars for some time. But when they accepted Nadars as members they asked for lot of money just for the membership. Many even said that the Hindu Mahasabha would soon come to ruin and so they did not respect it.

[Kudiarasu, 9-2-1930]

**Abolishing Caste Discrimination**

This political decision assembly that makes rules and regulations in independent India should have amended/ banned the foundations that caused discrimination and division by birth as Brahmins, Sudras, Panchamas and Harijans.

They should have made a statement to this effect banning all differences based on birth. Can we say that their unconcern and indifference about such discrimination is due to the lack of concern about their low caste status? It can never be said so. But, Brahmins made them members of this political party only on the condition that they would not protest about the low caste disgrace or discrimination.

[Viduthalai, 7-8-1930]

**Simon Commission**

By the polices of the Simon Commission and the information given by Chakravarthy, the cruelties of caste differences, gender differences, rich-poor differences, landlord–farmer differences and above all cruelties practiced in the name of religion that have a strong foundation, have not been destroyed or thrown away. On the other hand, the policy of Simon Commission and Chakravarthy has only made the three categories of people: Brahmins, rich and educated to lead a more comfortable life. The poor farmer and the coolies who form 90% of the population have not benefited even a little!
[Kudiarasu, 23-11-1930]

**Abolishing Caste Discrimination**

If we analyze the caste system formed by the Brahmins, it is not based on rationalism or common sense and there is no dignity or self-respect for the Sudras.

According to the Brahmins, the Sudras are the people of the fourth caste who traditionally commit do mistakes or sins. What is said in the Dharmasastras is quoted here: "A bathed horse, an infatuated elephant, a lusty bull and a Sudra who is educated should not be kept near. Further, the Vedas say that women and Sudras should never perform prayers, penance, pilgrimage, *sanyas,* singing the praise of god.

Further it is said in the Vedas that a Sudra should not take bath before sunrise, should not pray, they should not read the Sastras. Likewise, thousands of rules and regulations are written on how a Sudra should lead his life only as a slave.

If we are really interested in the complete freedom or democracy of India—at first, even before the British quit India, i.e. in the rule of the British itself, we should ask them to make rules that will annihilate all caste differences. Instead, if we ask them to leave our land, it is just like drinking poison for our own deaths because 999 out of 1000 Indians are not interested in abolishing caste differences, on the other hand each caste strives to make itself superior, which is done by making some caste to be lower to them. If the rule is in the hand of casteist people, caste cruelties will never be banished!
[Kudiarasu, 30-11-1930]

**Cow Slaughter**

Many Naidus, Gounders, and Iyers say that they are of higher caste than the Paraiyars and the Chakkiliyars. They say that the lower birth is due to the deeds of previous birth.

Based on this, they deny the Paraiyars and Chakkiliyars their due freedom and liberty. All this is due to caste, religion and god. Sankaracharya says that killing a cow is wrong. He does not seek protection for all animals, but only for the cow. No one kills a cow that gives milk. They kill only old cows that don't give milk. Ninety-five percent of people in the world live on beef. If they want a law against cow slaughter and they go on a fast for this end; I feel that it looks like an act of madness. If the poor should eat meat, it can only be beef because it is cheap. Only this beef contains more proteins for the poor.

To say that the cow should not be killed because it feeds only on grass and vegetation, and that it does not eat insect etc. looks even more foolish.
[Kudiarasu, 11-1-1931]

**Hinduism and Caste**
From the time the feelings of God came into being in this nation, the feeling of low caste and high caste also came into being. No one can deny the caste discrimination and caste cruelties prescribed in the Laws of Manu.
[Kudiarasu, 11-1-1931]

**Hinduism and Caste**
To hear about the cruelty wreaked on the Adi Dravidar society by persons belonging to other castes makes me angry. But while thinking about who is responsible for this, I come to the conclusion that those who torment you are not the ones responsible; because they are doing it in the name of their belief, in the name of their religious feelings, in the name of their devotion to the basis of their religion, in the name of freedom of karma of the previous birth–purva punyam–fate, they are only exercising their right to do so and nothing else. In the same way, when subjected to brutality by others even you just think of it – that too only during that particular instant–as something atrocious and cruel.

You fail to deliberate upon issues like: 'what is the reason for that? Why does such an atrocious and cruel act take place? How to destroy this totally? What should be done?' If someone else points this out, you do not courageously undertake to accept it and put it into practice.

Pointlessly, during such an occasion when some of you individually face troubles–instead of paying attention to the basis of the problem–you shout one or two words and throw the blame on a particular individual or section, after the passage of four days this fades away and everything goes on as usual. Under the guise of speaking in your favour, even some 'persons of esteem' in public life talk a few words and then return to their jobs and live casually as always.

No matter what, when human freedom like walking on roads, drawing water from ponds, and touching fellow human-beings is being banned to a person in his own country, is it not funny to ask for freedom to govern from people of other countries?

You must get human freedom. If not, one must not remain in this country or religion or society. Other than this decision, if we think it is enough to somehow manage to live and continue in this manner it is shameless nature.

In these matters, it is cowardice to think of tolerance. How many days to endure this? How slow to go? Several people have argued for this over thousands of years. So many people have finished talking on dharma and justice. All that could only strengthen the prevailing condition. Like the proverb, quote Tamil original ('Biting the hide and biting the bellows, the dog became a hound'), nationalism demands political freedom and on the other hand, it takes steps to attain Swaraj based on Varnashrama dharma.

Even in this giddy state, accompanied by a yearning for awakening, meetings are held to demand 'Caste differences are needed. Paraiyar is needed. Manu dharma is needed'. The humiliation that has been imposed on us and on you must be eradicated. For that, you and your family must be ready to sacrifice your lives!

In today's world, several nations have decided that even the difference of rich-poor must not be there, and they are successfully working towards (putting an end to) it. At such a time, if you are having the principle 'Brahmin caste and Pariah caste; if one sees–touches– walks on the streets it is polluting', think of the differences in thought between the foreign countries and us!

These thoughts are arising in you because of your ignorance that is founded on your belief in religion and God. So, chase them away. Then you will get self-respect and social equality.

[Kudiarasu, 11-3-1931]

**Racism and Casteism**

The 'racial arrogance' that is the cause of sorrow of the Negroes in America is also the cause for untouchability in our country. This was a reason why the Aryans—who migrated to our country in ancient times—degraded the indigenous native tribal people as Asuras, Rakshashas, Mlechas, Panchamas, and Sudras. Untouchability, like the wooden steps of a ladder, has afflicted every person in this country. The Indian leaders who speak about destroying untouchability do not wish to destroy varnashrama dharma.

[Kudiarasu, 12-4-1931]

**Untouchability**

It is a fact said by all the people and political sociologists, and accepted by the public that the cruelties which are being committed in India due to untouchability—people not touching people— is a very big cruelty far worse than all the cruelties in the world and no other cruelty can be mentioned on par with it. But, it is evident to one's eyes that no person takes any effort to use this fact in any manner. Instead they simply speak empty rhetoric to cheat people and pass the time.

In a way, because of the agitations of the Self-respect movement and the extreme emotions that have arisen in the Untouchables a kind of flutter has been created among the public and it looks like there won't be much value for opposing this. So, the non-Brahmin leaders and cadres of the Congress, must independently come together at least at this juncture—and also involve people of other movements—to come forward and take up some efforts! In this context, those who are said to be Untouchables must suitably agitate or leave the Hindu religion that is the basis of untouchability, or take efforts to eradicate their disgraceful positions and cruelties. It is necessary that those who are called Untouchable in every place must jointly come to some decision, instead of expecting that extremist efforts would be undertaken by the English-educated or wealthy among them who expect jobs, appointments etc.
[Kudiarasu, 24-5-1931]

**Hinduism and Caste**

The evil they (Brahmins) do to non-Brahmins in the name of Hinduism is less than the evil they do to Muslims and Christians. But the cruelties they do to the Untouchables are greater than everything. The Hindu religion indirectly helps in the flourishing of both Christianity and Islam. So it is a profit for these two religions. So it is more important only for the non-Brahmins and Untouchables to destroy Hindu religion.
[Kudiarasu, 7-6-1931].

**Congress and Untouchability**

Of what use is the Swaraj that cannot bear to inflict anyone's feelings? It may claim itself to be independent. It is of good use for the Brahmins who can lead a comfortable life than before. But what is the plight of a Pariah? He cannot even have the right to walk in all streets. What is the plight of Sudra? He should do physical labour to satisfy the Brahmins. If a Pariah walks in the public streets, the feelings of Brahmins are hurt!

Once Jawaharlal said when he was talking with the Brahmins that no ones liberty would be affected or snatched away by Swaraj. But soon after he spoke in another gathering that the temples that does not allow the Untouchables would be razed. Gandhi also supports this.

Once I had gone to Dindugul with Srinivasa Iyengar on a Congress campaign and had to go to a Brahmin's house. I was given food in a different place outside the house and the leaf from which I ate in the morning was not cleared and in the same place, close to it, I was served supper also.

Likewise, once I had been to Tanjore with Venkatasamy for some Congress campaigning in Periyakulam. We alighted in a Brahmin lawyer's house. We were served breakfast in plantain leaves, outside the house, in an adjoining platform. These leaves were not removed and it was full of ants and flies. The same night, we were both served dinner near these leaves. I did not mind all this and worked for the Congress. So, if the Congress says that they don't practice caste differences, no one can believe it!
[Kudiarasu, 12-7-1931]

**Gandhi and Untouchability**
Under any circumstance, if we expect samadharma (equality) from Gandhi only we would become fools. Till the end of our life we cannot see the samadharma or unity in him.

He considers only the rich northerners and the southern Brahmins as humans. He mainly has their friendship. He thinks only of their problems as the international problems. All his plans are only to solve the problem of these two groups. But for the poor and low caste people to develop, he knows only two ways. One is spinning and the other is that the practice of untouchability is a sin. Even these two things would only become a philosophy in the presence of mill-owners and Brahmins.
[Kudiarasu, 26-7-1931]

**Religious Conversion**
When I advised the Untouchables to join Islam several of them got angry with me. I did not get angry with them.

I did not ask the Depressed Classes to join Islam to get heaven or mental satisfaction or to reach god, I mainly asked them to join Islam to destroy untouchability. It is just like doing *Satyagraha*. Even now I profess it.

It is difficult to make laws. Even if it is made it is difficult to be practiced. Because of this, one may get sorrow and defeat. But there is no difficulty for the Depressed Classes to convert to Islam. To follow this, there is no defeat or sorrow. What is the difficulty faced by others?

What is the problem whether one is an atheist or theist? Nobody suffers because of this! Even if he is a true Muslim or a false Muslim nobody is troubled by it! Even one need not have any mental change!

A person who wants to wipe off the disgrace and dishonour he suffers; if he converts to Islam at 5 p.m. he can walk in any street freely at 5.30 pm and he is free from the clutches of untouchability. What is wrong in becoming human?

Why should one not convert to Christianity or to Arya Samaj? I do not know what is said about Christianity in the Bible. But in practice they have in India especially, Paraiyar Christian, Nadar Christian etc.! Does the Muslim community have such differences? The Christians should not get angry with me.

Arya Samaj is hypocrisy. Because, in the Vaikom Satyagraha, the Arya Samaj Untouchables and the Pulayan Christians were not allowed to walk on the road near the temples. But the Muslims walked freely.

For instance, a Cheruman, an untouchable who converted to Islam walked along the road banned for Untouchables.

A Brahmin and a Nair came to the street and looked at him. A Muslim asked them, "Why do you sons-of-a-prostitute-pig look at him?" They bent their heads in shame and this untouchable walked on, smiling in that road.

[Kudiarasu, 2-8-1931]

**Caste System**

An idli-vendor Brahmin's son can become a High Court Judge. A horoscope-making Brahmin's son can become a minister. But think if a scavenger's son can become a Judge? That is the power of caste system.

[Kudiarasu, 6-9-1931]

**Caste and Class**

The caste system, which preaches low and high castes, acts as a supporting fort to the notion of rich and poor.

[Kudiarasu, 4-10-1931]

**Adi Dravidars and non-Brahmins**

We are uniting caste. Yes, we are trying. There is no doubt in that, but whether it will be done soon is a doubt. The human caste has to become one. We are boldly saying that those who prevent it are scoundrels and fools. We are being asked if we will give our daughters to Paraiyans.

This is a foolish question; I might say that it is a roguish question. This is because we are only going to make our daughters live with their loved one, and we are not people who will exercise the right of giving our daughters to men whom we like. The practice of considering women as objects and 'giving them to someone' must be eradicated; we are striving for it.

Friends! When we are speaking of Adi Dravidars, there is sense in the Brahmins getting mentally grieved. But, there is no meaning in the non-Brahmins getting upset. It is just plain foolishness and shamelessness. Because, in our society, apart from the Brahmin caste that constitutes 3% of the population, can the rest of the population have any classificatory title other than Sudra and Adi Dravidar? Please consider this, which is also based on your experience!

If you believe that without the label of Paraiyan being removed, the Sudra label will go away, then you are downright idiots.

Besides, if I talk at length, I would say that there is no evidence anywhere on who is a Paraiyan or Chakkiliyan and what their rights are! It is nothing but a show of strength and tactics. It is possible for the title of Paraiyan to disappear in a short while. The Sudra title given to you has many evidences–God, religion, Vedas, Sastras, puranas, history etc. Unless all these are completely destroyed you cannot expect the Sudra title on your head to get down. So, if anyone has a feeling of dignity, they wouldn't have asked us, 'You are uniting the castes!'

As a result, do realize that the speeches addressed regarding the benefit of Adi Dravidars and the efforts initiated for them are for the welfare of all non-Brahmins.
[Kudiarasu, 11-10-1931]

**Communal Representation**
What is wrong if any man asks for power sharing by saying that his religion is a minority, his caste is powerless and low; hence clearly ensures his share in the rule and sees to it that you do not suppress them by your rule? If anyone asks like this, they say labeled traitors to the nation? What can we say about these people? If one asks his share in his nation he is termed as traitor. All families that do not give the proper share to its members have come to ruin so if the share were denied or cheated out of any individual, religion or community, it is certain that such a nation would be ruined. At present in our nation, certainly one cannot refute the demonstrations carried against communal representation in such institutions.
[Kudiarasu, 8-11-1931]

**Abolition of Castes**
If we just abandon the Brahmins, then what about the Saivaites? Can they be there? They are also just like the Brahmins! So, all people who claim to be 'upper' caste are those who do not do any physical work but live on others' labour.

Castes should be abolished. A law should be made that everybody should toil physically. Without toil, there must be no place in this nation. Only then the Brahmin will get away and all the cruelties that are being practiced will go away. But keep in mind that there is another group which lives on the labour of others and which practices caste system. And it is this group that protected the Brahmins. They are the rich people, the capitalist, and the *mirasdhars* (landlords). They should also be reformed or annihilated.
[Kudiarasu, 5-3-1933]

### Congress and Untouchability

If the British rule were also based on Sanatana Dharma or on Manu Dharma, then certainly the notion of nationalism, law and non-cooperative movement would have flown away. For this reason, Gandhi was made Mahatma.

As the reputation and power of Gandhi lessened, he joined the anti-untouchability protest and temple-entry movement to renew his popularity. This caused some friction with Brahmins. That is why the Brahmin nationalists A. Rangasamy Iyengar, Satyamoorthy and K. Bashyam and others became silent. It has become very difficult to find Satyamoorthy's whereabouts. He became a full-time contractor to Raja Sir Annamalai Chettiar. So he was earning a lot. Rangasamy Iyengar became a full supporter of Sankaracharya so he earned money and fame!
[Kudiarasu, 19-3-1933]

### Christianity and Caste

Your priests know that you do not have the feeling of self-respect because of the madness of religion. They know that you will not do anything against their authority. For the sake of religion and God, you will bear any amount of sorrow and disgrace and will be speaking something merely orally—but any day you will not agree to transgress the regulation, or shake away or try to destroy the cause of your humiliation.

[Kudiarasu, 7-5-1933; Lalgudi Taluk Adi Dravidar Christians Conference, 23-4-1933]

**Gandhi and Untouchability**

The efforts of Gandhi were only to establish Brahminism and richness. That is the only thought that dominated his rich Brahmin disciples. The main reason that established such actions of Gandhi are: the discovery of spinning wheel for the poor, Bhagavad Gita for the Brahmins' existence, practice of caste system and varnashrama dharma in all types of industries and their representation in Swaraj. In the Round Table Conference, the acceptance of the Indian *samasthana* rulers (traditional kings) and zamindars established his intentions.

When asked why he had started the "salt satyagraha", he said, "If I don't start it, the nation would have faced a strong agitation towards establishing communism." Because of this, the rich lavishly supported and helped Gandhi.

Not only that, they said that they had abolished the low castes and started a Harijan Sevak Sangh by getting lakhs of rupees from rich industrialists like Birla. They purchased a few Untouchables. They did service only to the Hindu society by campaigning for puranas. They said that according to puranas, we were to work for others.

[Kudiarasu, 23-7-1933]

**Congress and Untouchability**

The Congress has not in anyway helped the poor, downtrodden people who toil like animals for over 14 to 18 hours in the rain or shine, forests, fields, hills, seas, workshops, industries and business centres. They don't have proper food, sanitation, medical aid, support for their children's education, proper clothes or proper house to stay; they live in broken huts, and die of poverty and starvation. There is no one to question this. What is the change that has taken place after Congress has come to power?

If the Indian National Congress served as the establishment for prostitution, no one can deny it.

Everybody should accept it because in the name of the poor and ignorant, the atrocities and the cruelties the Congress had done to them knows no bounds.

[Kudiarasu, 30-7-1933]

**Caste and Class**

In any nation, can the rich and the upper castes alone remain in a flourishing state, and the people who labour and are poor remain in a declining and degrading state? If any such institution or state or nation dies, the people would only be happy about it and they would by no means try to help it grow.

In fact they would bury it in a deep pit and it is the ardent duty of the intellectuals who seek equality to erect a headstone over it.

[Kudiarasu, 30-7-1933]

**Samadharma**

*Samadharma* (lit. all are equal) is above caste, religion, nation and varna and it is has destroyed the castes: Brahmin, Kshatriya, Vaishya, Sudra and Panchama. It has no difference between rich and poor, employee and employer, landlord and farmer. Samadharma annihilates all such notions. It has broken the mean notions of high and low caste, or rich and poor etc. based on karma. According to Samadharma, all are equal, everything is common, no difference exists and no caste or religion exists, in the world everybody is a comrade.

[Kudiarasu, 30-7-1933]

**Gandhi and Untouchability**

When Gandhi entered politics he chiefly spoke about social reforms: he said that we would not get Swaraj without the abolition of untouchability; we would not get Swaraj without the Hindu-Muslim unity, we would not get Swaraj without the abolition of alcholic drinks. He said that in every house spinning should be taken up as a compulsory work without which we would not get Swaraj. As days went by, he talked about the eradication of untouchability on one side but on the other, he supported varnashrama dharma.

He said the divisions of caste must exist; only the difference of high or low caste should not exist. He said that you couldn't become a Brahmin unless you are born a Brahmin. Then he said that if we get Swaraj, untouchability would go away by itself, automatically.

Lastly, because of the demonstrations by the Self-Respecters, even after getting their communal representation of reservation, the Untouchables have established a forum to fight for temple-entry. The demonstration started with fast and now it had gone to the assembly for making some laws. So it was decided that until some laws are passed they couldn't enter temples.

[Kudiarasu, 6-8-1933]

**Gandhi and Untouchability**
Gandhi's plan to drive away the devil of untouchability is only to make the religion pure and that has no relation with any feeling of sadness over untouchability. He ventured to make the Hindu religion pure because it should not come to a state of being ridiculed as a religion with flaws, since that would quickly destroy the Hindu religion. To save the Hindu religion, Gandhi was made Mahatma and by these acts, the Brahmins were able to collect lots of money to fill their stomach.

Gandhi had said over thousand times that the removal of untouchability is not a common dharma in the Hindu religion. Can the common man not follow that?

So, is the removal of untouchability only a religious campaign? Gandhi met people of different castes many times and has said that the caste differences are very essential and each caste should do only the profession of their forefathers. So, he had come to this Ashram and taught Untouchables the art of tanning cowhide.

I will bet that Gandhi has not till date spoken about the destruction of different castes but he always wanted the different castes viz. Brahmin, Kshtriya, Vaishya and Sudra; and the fifth caste of Chandalas, he wanted it merged with the Sudras.

When asked what is the place of Untouchables in the Hindu religion, Gandhi said that they should be combined with Sudras and the Sudras will continue to be Untouchables. The orthodox and people who follow Hindu dharma are better than Gandhi because the venom in them is well known even to fools and blinds; but the venom in Gandhi was sugarcoated so even an intellectual would eat it without doubt.

The venom of the orthodox will not kill the person who consumes it but the venom used by Gandhi was so sugarcoated that it will kill several generations. The campaign of Gandhi was more venomous than that of the followers of Manu.

[Puratchi, 10-12-1933]

**Cow worship**

Hindus means they will worship cows and they will show boundless respect to it. If anyone kills a cow or eats cow's meat the Hindus will be cruel to them and hate them. That is why they hate Muslims. That is why often there are Hindu-Muslim clashes. Among Hindus, the Paraiyars and Chakkiliyars are termed low castes because they eat beef. Based on these, several Sastras are written.

But these Hindus treat the bulls very cruelly. Bulls are used to pull carts, plough fields, pull water from wells and run the oil press. When the bulls are young, they castrate it in a very barbaric way. The amount of love and respect they show to the cow is not shown even to humans!

If the acts like the abolition of untouchability and temple entry to all is brought out as laws, they felt it will be a danger to orthodoxy according to Hindu scriptures. So, they are performing yagnas to escape that. This Hindu orthodoxy is very dangerous to non-Brahmins and it favours only the Brahmins. This is against everybody's self-respect.

That which was once considered to be civilization is now being ridiculed. For example, actions such as putting caste tags behind the name, wearing the sacred thread, using religious marks on the forehead are criticized for being contradictory to progressiveness.

[Puratchi, 10-6-1934]

**Press**

The main failure of the Round Table Conference, the disparity of people in political reform, the lack of faith between Hindus and Muslims, the differences between the Congress and Muslim league, the Untouchables seeking separate electorates are all due to religion and caste differences.

*Kudiarasu* was stopped because we wrote about Christianity and the magazine *Puratchi* was stopped because we wrote about Islam. Because we write about Hindu religion, the daily torture we undergo is innumerable. Officers from Brahmin community in the fields of police, railways, law and establishment trouble and torture us limitlessly.

[Pagutharivu, 9-9-1934]

**Non-Brahmin Representation**

The officials should supervise and ensure that non-Brahmin representation is duly given. For instance, if a Brahmin is a Deputy Collector, they should see to it that the Tahsildar is a non-Brahmin, likewise if a Brahmin is a District Judge than the sub-judge must be a non-Brahmin, it should be so in all power/ posts. If such a format were not followed, only the casteist emotions would dominate all spheres of life.

Now, there is some awareness among us for instance in the municipal elections. In Salem and Tirupur that was the base of Brahmins, not a single Brahmin could come to power. One can now see the status of the Congress and the Brahmins!

[Pagutharivu, 30-9-1934]

**Swaraj and Caste**

Even if a good General comes to the nation, the job of a scavenger will continue to be so: this was the dharma practiced by the nationalists and Mahatmas as *Swarajya dharma*.

British rule was better than Swaraj (Self-Rule). Because in the British Raj, a scavenger may become a minister; but in Swaraj that forces hereditary profession by no means will make a scavenger into a minister. That is why I want to be a traitor; I want to be an anti-Swarajist. I am not an enemy of the law that will destroy the difference between Brahmin and Pariah and profess equality among them and the Swaraj. I am not talking about this only in the Indian context but as the person who has analyzed it in the international context.

In every nation there is difference between poor-rich and capitalist-labourer but in no nation there is Brahmin-Pariah, low caste-high caste. These meaningless differences are created by caste culprits who demand Swaraj but are uninterested in destroying these differences. Until everybody is united under the banner of humane existence, I would like to be a traitor than accept the fraud of nationalism and Swaraj.

[Pagutharivu, 21-10-1934]

**Communal Representation**

Till date, the main goal of national heroes, nationalists and Congressmen was to be in good posts. They worked only based on that goal. But, when the non-Brahmins took the same mode of work they were called traitors. For, if there was communal representation, the Brahmins could enjoy only a certain amount of money. If the rule of communal representation is practiced then there will be no need for the Justice Party.

Whoever does any type of national service, even if it is Gandhi himself, without communal representation such service cannot be legal, proper or done with integrity.

Till the Untouchables and Backward classes continue to have poor or less representation in all posts, I would support them and will also be a traitor; because caste and equality in representation alone is the biggest thing for me. A group that labours is starving; conversely, another group is eating stomachful without doing any work. And the oppression by high castes and caste differences are so much that one is not able to tolerate them. How can one tolerate all these and remain a patriot?

[Pagutharivu, 9-9-1934]

**Hinduism and Caste Differences**

An untouchable Chokhamela, a Vaishnava devotee was not allowed into the temple. Vishnu himself led Chokhamela into the temple. Likewise in Saivaism, Nandan was taken into the temple. So the communal and caste differences in both Vaishnavism and Saivaism remain the same.

[Pagutharivu, 1935]

**Sectarianism**

We cannot consider sectarianism to be annihilated if the sectarian differences and disagreements such as Brahmin-Non-Brahmin, touchable-untouchable are solved. Who can contradict that sectarianism is rampant, both openly and as an internal consideration in today's politics?

[Kudiarasu, 26-5-1935]

**Communal Representation**

The Justice Party was created only to fight for communal representation. The revolutionaries Dr. Nair and S.P. Thiyagarayar, staunch Congressmen and nationalists were the main reason for creating the Justice Party. They wanted communal representation in education, science, technology and employment. Their main aim was to struggle for advancement of the Depressed and Backward classes. This party is a non-Brahmin Party.

[Kudiarasu, 23-6-1935]

**Congress and Untouchability**

In 1920, Gandhi collected a fund amounting to 1 crore rupees in the name of Swaraj Fund—this amount was spent only for Brahmin domination and their welfare. Once again in the year 1926-27, *Kadhar* (handspun cloth) Fund was collected amounting to 30 lakh rupees. This was also spent only for the betterment of Brahmins and they grew richer. Finally in 1934, Harijan fund was collected in the name of Untouchables. Around 20 to 30 lakh rupees were collected and it was also given to Brahmins to perpetrate their domination.

The village rehabilitation fund will also be only spent in the same way as the above three funds were spent. There can be no doubt about it. What was done for Untouchability was that they made it clear that no law should be brought against abolition of Untouchability. It ended with a slogan that people should eat hand-pounded rice and jaggery. Be very cautious about this collection of funds since it would make only the rich and the Brahmins economically powerful.

[Kudiarasu, 30-6-1935]

**Untouchability Eradication**

Unless religion is banished, untouchability will never be banished. If the Untouchables need some liberty it should be got only from the government i.e., law. So I advice them to get their freedom from the British rule since any other rule in the ancient days did not help them to come out of untouchability. Do not involve yourself in the demonstrations conducted by rich or other caste people. It is done only for their improvement.

[Kudiarasu, 28-7-1935]

**Hinduism and Untouchability**

If we see with a true consideration that the Depressed Classes must be liberated from the atrocities heaped on them by other people, it is the work of a revolution because the position of the Depressed Classes has been built on a great foundation. It has been built on such a strong foundation that the Depressed Classes, the people of the lower castes and those who are called Untouchables have attained the inferior quality in their birth itself. They were made to be born in such situation by God itself and that the actions of God and the rules of religion must not be changed by anyone and that it is not amenable to change.

If someone thinks the oppressed people can attain equality and the doctrine of untouchability can be removed from the society, because of mere words or propaganda, or by requesting the upper castes, I will only say that their life is a waste.

A few people of the Depressed classes and the Untouchables think that if they bathe and apply sacred ash, wear a *pattai* or *naamam* (marks on the forehead to indicate their religion) as per the customs, act like the orthodox, say that they do not drink or eat meat, roam about naming themselves 'samy'—their position will improve and that Untouchability will be eradicated. This is the madness of trying to cheat others and cheating themselves.

For a longtime we have seen many people act in this manner among the Depressed Classes, among the Adi Dravidars. There are several puranic and historic evidences for that.

We can say with courage that from thousands of years ago till now, in the matter of annihilation of Untouchability no work has been completed. The work of a few magicians could be used for selfishness; because of those garbs and piety no work took place. So, if God and religion and its basis—the Gita and Manu dharma Sastra are to be saved, think over of how the Sudra quality and the low-caste quality and untouchability can be changed?

Today, cent percent of those who are involved in the work of removal of untouchability and eradication of caste discrimination are idiots who believe in the Gita and the Manu dharma Sastras. For how many ever days they work, they will only be like fools who try to fill up water in pots that are broken. So, for the annihilation of untouchability or the annihilation of caste, first you have to annihilate your religion.

If you are not able to annihilate your religion, at least you must come out of the religion. Without your religion going, your untouchability or the pariah-hood will never be eradicated, this is a fact as firm as stone. If you want an example, those who were Untouchables in the human society could become touchable only when they shook off the religion that was imposed on them. So, do not think that religion can be protected and untouchability can also be eradicated and thus be cheated.

[Kudiarasu, 28-7-1935]

**Varnashrama Dharma**

How the kings ruled in the ancient days can be seen from the dramas staged by them in those days. In the scenes of the play, the king will ask his ministers, "Did you lavishly give the Brahmins subsidy, free distribution of food, Vedic schools? Are they happy without any mental strain? Do the other caste people serve them without any problem?" Thus the rule of kings was nothing but only the rule of varnashrama dharma.

[Kudiarasu, 8-9-1935]

**Patriotism**

In today's context, patriotism is based only on money and Brahmin campaign and not for the welfare of the commoner or to make amendments for some need or wrong in the social set up. The Brahmins have isolated themselves from society by the status that they are the only high caste and God created all other castes only to serve them.

Likewise the capitalists and rich think that they are different from the others and God has made them rich and it is their right to pay the workers what they feel like paying them.

The only difference between these two sets of people is that a Brahmin gets work done without paying them and the rich or capitalists gets work by paying them with what they like and both feel they are created by god. Thus these two groups form the nation's patriots.
[Kudiarasu, 29-9-1935]

**Poona Pact**

In the reforms offered by the Government, there were reforms with self-respect for you. That was spoiled by the patriotic reform called the Poona Pact. It is surprising for me that a few of you support the Poona Pact. The reason they give for that is very, very surprising. That is, you need to go and seek votes from the upper castes. The upper castes need to come to you and seek votes. It seems a mutual feeling will be created. This is like a mynah-catching stunt and there is no credibility or intelligence in this. If you have to enter an administrative office where a Brahmin is the authority you need a government order, and to walk in the public streets you need a law in the Penal code. Is it possible for you to go to walk into the *agraharam* (Brahmin settlement), cross the threshold of their homes, and ask for votes from a Brahmin who is lying on a couch and stroking his stomach.
[Presidential Address in Rasipuram Adi Dravidar Conference, 29-9-1935][7]

**Congress and Untouchability**

Our people have not realized the nature of Brahmin conspiracy against the annihilation of caste.

Gandhi works for the eradication of untouchability. But he doesn't change his opinion that the Depressed Classes must stitch slippers, they must work with hide. Our Congress enthusiasts worked against untouchability. They collected money. They used that money and caught comrade Sahajanandha and made him preach that Hindu religion and Hindu religious dharma must be saved and also helped him to run a magazine. They caught comrade M.C.Rajah and made him speak about the holiness of Hindu religion. Because such people were there, Dr.Ambedkar considered that it was not possible for him to destroy Hinduism, and that it was not possible to eradicate untouchability by remaining in Hinduism, he decided to quit Hinduism.

---

[7] *Rettamalai Srinivasan speaks on 6-5-1937: One has to remember that when Periyar was in Europe, he sent a telegram: Don't forget that more than the life of Gandhi, the life of 6 crores of the Depressed Classes is in your hands.*

The salary being given out of the Untouchability Eradication Fund is only being used for the Brahmins. […]

In any Congress meeting have they sat and taken food without being differentiated as Brahmins and non-Brahmins? Even in a recent Congress meeting, did they not feed the non-Brahmins by making them sit outside? In all the meetings, right from the Congress Conference, they have cooked only using Brahmins, have they allowed people of other castes to enter the kitchens?

Take the Self-Respect meetings or its conferences. Nadar, Naidu, Muslim, and those who are called Untouchables, all of them together cooked the meal. All of them together served the food and right from the 'Saivas' to all the so-called 'upper' castes sat together and eat.

At least from this you must understand the Congress's conspiracy of untouchability. Even if the self-respecters don't have faith, they praise some religions because they don't have untouchability. They support some political parties because it gives equal rights for the Untouchables. The Self-Respect movement tom-toms that anything that denies human rights must be destroyed because that is human nature.

[Kudiarasu, 19-1-1936]

**Traitors**

A few among you, who are Brahmins' coolies, live by obtaining wages from the Brahmins for saying that the British Government did not do anything for you and that the non-Brahmin agitation did not do any good for you. That group is not only amongst you, but is also amongst us. So, to an extent we need to give space for the life of people like that.

The reason for getting this right is because of the agitation of the Self Respect Movement. The Brahmins made a conspiracy of boycotting the Simon Commission since they didn't want their mistakes to be known in the Parliament in England. That did not work with the Self-Respecters.
[Tirunchengode Adi Dravidar Conference, 7-3-1936]

**Gandhi and Untouchability**

Instead of speaking about caste (because that will attract the wrath of Brahmins), Gandhi started saying that he will eradicate untouchability. For 10, 15 years he got publicity. He collected several lakh rupees. Till today, he could not say from his mouth that castes must be destroyed. Instead of that he says that it is his duty to save the caste (that is Varnashrama) system.

What is the reason? If his Mahatma title is to become permanent, he has to say that he will protect the castes. He has come to that state.
[Kudiarasu, 5-4-1936]

**Self Respect Movement**

Self Respect Movement was started for the people to lead a life of dignity and self-respect. There is nothing to laugh about it. It is a social organization to make people respect each other. Because of this social movement we face and suffer many problems.
[Kudiarasu, 5-4-1936]

**Gandhi and Untouchability**

In 1920, Gandhi made four plans for the betterment of human society. That is: removal of untouchability, Hindu-Muslim unity, prohibition of liquor, and *kadhar* (handspun cotton). What happened to these four?

The removal of untouchability has come to such a level that unless the Untouchables come of Hinduism, untouchability cannot be removed. When I said this in a Congress meeting in 1923, everyone detested me. When I said the same thing in 1928, 29 some people displayed patience, a few experienced difficulties. But today it has become one of the major problems in India.

If untouchability has to be eradicated, it cannot take place by remaining in Hinduism; it is not possible by the Congress; it is not possible by Gandhi. They need to leave the Hindu religion and go to another religion, that too the Islamic religion: this situation has been created. Keeping aside the concept of whether this is right or wrong, no one can refute that Gandhi's philosophy of untouchability eradication has created this situation.

[Kudiarasu, 5-4-1936]

**Round Table Conference**

What is our political agitation? What basis does it have other than 'which group should capture power?' If you want an example turn to the chapter on Gandhi and the Round Table Conference. What is there in that? The British Government called you (the Indians) to know what you wanted? We, who are the Indians, what did we go there and say? On behalf of the Hindus, Gandhi said, "our tradition must be protected."

Those who went on behalf of the non-Brahmins among the Hindus, said, "of the freedom you are going to offer, we need at least as per our percentage of the population."

Comrade Ambedkar who went on behalf of the Untouchables among the Hindus said, "You first distribute a certain quantity to our society, and then you can give any amount of freedom."

What was the answer of Gandhi to all this?

"All of them are traitors to the nation, they are all henchmen of the British Government; I am the only patriot. I am the sole leader of all the religions and all the sections of the Indian society. If you give the Swaraj that I demand, later I shall look into the matter of these people. He not only said this but also he took a vow, "If you give political representation to them without my consent, I shall commit suicide."

The Round Table became an eight-cornered table.

[Kudiarasu, 14-6-1936]

**Gandhi and Untouchability**

At first only the Justice Party worked for the Muslims, poor and the Untouchables. Only after that Gandhi started working like them.

After Gandhi became a 'Mahatma,' he said that without Hindu-Muslim unity, Swaraj could not be created. Likewise, he said to the Untouchables that Swaraj would come only when untouchability is eradicated, and that without untouchability being eradicated, Swaraj will not come and even if it comes then, it is not necessary.

Later, when he went to the Round Table Conference he said that without Swaraj, Hindu-Muslim unity would not come about. Jinnah got angry on hearing these words. He asked, "We need to know what is our position in the Swaraj Government?" because he knew that Gandhi's Ram Raj would only be a Hindu rule. Gandhi cannot cheat everyone, right? Seeing that they did not come to any decision, the British Government delivered justice to that society.

Next, to ruin the Untouchables state, he went on a fast. So, Ambedkar got afraid and he felt "we have committed many sins and are born in the untouchable society. What more sins will come to us?"

So he ran to Gandhi and signed the Poona Pact.

What do these people do for the removal of untouchability? They give them oil and soap and peanuts and rice and peppermint. Is this service for the Harijans? This costs one rupee a week. The person who distributes this material has a salary of Rs.36 per month.

Today, the difficulty of the poor cannot be really told. Nobody looks after the welfare of the poor. The torture of the poor people by the caste people is impossible to narrate.

We observed all this and felt that we can use the Congress, but they saw to it that a man with integrity does not come to the Congress. Anybody who is not a slave of the Brahmins does not have a place there.

But one question: If they win in the elections, let them say that they will give equal opportunity as per communal representation in all the positions. But what is possible when they threaten that this is betrayal to the nation?

The way you threaten us, go and threaten (Jinnah) Sahib, the Bengal Muslims: they will immediately break your teeth and give it in your hand.

[Kudiarasu, 18.10.1936]

**Social Equality**

If a riot breaks out, the upper caste people will supply the lower caste people with drinks and see that they become sacrifices. If their work is not achieved, they will be betrayers who will go under the feet of the enemies. Only these serve as the ancient heritage of India.

He, who thinks that India and the Indian people require benefit, must think as an Indian and understand the Indian situation and act accordingly. If he reads about foreign countries and remains a bookworm, it is a wasteful effort. When the state of our society is like the state of the foreign society, to use the foreign methods will be apt. Like the saying, "A blind man gazing a King's gaze" our Indian Paraiyar, Chakkili, Brahmin, Sudra speaking of economic equality, Marxism, Leninism in our villages is just plain masquerade and waste of time.

So the socialist youth of today please forget the idea of becoming rich for some days at least and involve themselves in the work of annihilating caste, creating social equality and social revolution. Because of unexpected reasons, if the situation is conducive we shall think about economy also.

[Kudiarasu, 19-11-1936]

**Untouchables**

Today, in this annual function, comrades U.P.A. Soundarapandian, N. Sivaraj, S.Gurusamy, T.N.Raman, Kunjitham, Vidwan Munusamy, Arockiasamy, Kalyanasundaram have spoken here. You must stay united. You must make a lot of women your members.

Several people who spoke here, talked about Self Respect and politics. Since it is seen that the majority of the members of this association are from the Depressed Classes, it is useful for you to support the Government. In politics, your group has been divided separately. Because there is a great difference between the objectives and status of your section and the other sections, you need to demand separate rights.

Though the Government gave you separate rights, the other upper caste people cheated you. The upper castes naturally have slyness and cunning. That is why they remain as upper caste. Because you don't have that, you are joined in the lower castes and you enjoy the benefits of that. Comrade Arockiasamy got angry when you were called downtrodden people. What is the point of getting angry? Obviously, are you not downtrodden? If that was not so, you would not need special rights? Look at how many obstacles you have in social life as per the law? You do not have the right of temple entry. After the Justice Party came, you have been given the right to enter the streets, ponds, and schools. Even now, in several places there is a ban on that.

If one has to speak justly, you alone are not downtrodden. Those who do not have temple entry rights in some places are not the only downtrodden. Even we who have rights to enter all the temples are also downtrodden. Even in the temples, there are several places where even we are not allowed. In hotels, we are also segregated and not allowed entry into some rooms.

Our sympathy is because we feel that if your disgrace is removed, our disgrace might also be removed. Do you not know of the practice of calling the server, under the pretext that your neighbour needs a food-item, and once he comes, you also get an extra helping from the server? That is how, if your problems and humiliation is removed, even our problems and humiliation will be automatically removed. That is why we keep talking about your grievances all the time.

How much you hate politics, how much you remain without joining political parties, that much your grievances can be solved: that is my opinion. If you do not give place to those who trample your heads, your difficulties will soon be taken care of. Otherwise, you will have to remain as a stepping-stone.

Though it is fifty years since the formation of the Congress, only after the Justice Party took part and began to ask for a share, in the social sector a great change could come about. Only after that, even in your state there have been many changes in these ten years. If this was not so, and you had to be behind the Congress, and if you had been behind the upper-caste and flattering them, think of what your position could be.

Think of how the Travancore temple's door was broken. The people there, including the Depressed Classes, Ezhavars, Nadars and others tried to break the religion and temple and god. Tens of thousands of people conducted conferences and resolved that Hindu religion is a sham, the temple is a sham, and that god is a sham.

Several people became Muslims and put on the Turkish cap. Several people grew their hair, their beards, got themselves a kirpan (knife) and became Sikhs. Some people converted into Christianity along with their family. Only after that the door of the temple was opened. All the upper castes in Tamil Nadu including the Brahmins congratulated the Travancore Maharajah. See where the secret of success lies. Likewise if you start demolishing Congress, religion, temples, and god you will automatically get all the rights without anyone's grace.

Without being so, nothing will happen because of your getting angry over being called downtrodden. The title of Paraiyars changed into Adi Dravidars, which changed into Harijans, which might change into some other name, but your problems and humiliation will not be removed.

Like the wood of an axe (which brings ruin to its own kindred), if a few are supportive of the Brahmins we must not be cheated. Because there are several such low, shameless people among ourselves we need to suffer by keeping such associations and establishments. If everyone has the feeling of self-respect, why should there be an association of the Depressed Classes or non-Brahmins?

Several among us are those who live by kissing the feet that kicks them. For the disgraceful occupation of such people, we need a remedy. That lies only in our non-cooperation and being obstacles.

Please remember that the Justice Party is a party of equality. Only after it was started, today the Paraiyan and the Brahmin are sitting equally in the same place. It is the speech of olden days to say that equality means that the tiger and the cow must drink water from the same place. But that is taking place today in the circus.

Because of that we don't say that it is a government of equality. But today the Paraiyan, Brahmin, Sastri, Sankaracharya and the Chakkili are all sitting in the same dais—How? Because of the whip? Or because of the fear of revolver? No, not at all. They themselves desire so, they spend ten thousand, and twenty thousand and wish to be seated so. The Brahmin remains standing and speaks to the Paraiyan saying O Lord, O Master. How was all this created? Before the creating of the Justice Party, in the Congress meetings, actions, policies, schemes was there one single word about the Depressed Classes: think and see? So I ask you, do you still have doubts whether the Justice Party is a party of equal justice? So, even now, even today, I only do the work of social equality.

Even I like to do the work of economic equality. Not only is the Congress a dead enemy of that, but it is also an establishment of treachery. Only if it is abolished it shall be comfortable to talk of economic equality. Still, without being against the law that is without becoming a victim of oppression by the Government, I am campaigning as much as possible for equality. I am going to continue doing so.
[Kudiarasu, 10-01-1937]

**Temple Entry**

In the time of Manu, or in the time of Raman, Krishnan, Harishchandran, Pandiyan, Naickar, could a 'Paraiyar' become a political minister? Or let it be any God. Even if the Paraiyan and Pulayan are going to come in front and pull out their tongue and die, in the Kerala country—that has tradition and Varnashrama madness like a monkey's grip, where the men and women think that they shall get a higher position if they are born to the Brahmin, where the people and the politics consider that it is a sin if they are seen by the Paraiyan, and if his shadow touches them, or if his words fall on their ears—could anybody have expected that it would have been possible for the Paraiyans to enter temples? This work has taken place only in our time.
[Kudiarasu, 18-7-1937]

**Caste and Education**

The Brahmin community in India is 100% literate. But why are 90% of other castes illiterate? 95% of the non-Brahmins are illiterates and 99½% of the Depressed Classes are illiterates. Is the government a cause for all these? What have we done to find the cause? Does the government give us education? Or does the government not have money to do it? Are people not interested in education? If we analyze these questions we can arrive at the conclusion as to why 90% of us remain uneducated! (The main reason for this is that good education is denied to these people, under the garb of caste).

If you look at the true cunningness in this—we can derive the following: it is mentioned in the law-giving Vedas and Sastras on which community alone should be given education. There were laws on who should teach and also who should be taught. So one can see that the true reason for the majority remaining uneducated is the Vedas.

The worst affected by these Sastras and Vedas are the Untouchables. They are making several plans to see that Untouchables remain uneducated. There is no school for them in several places. Even if there are some schools that give them education, the Untouchables cannot use public roads to reach the school since they should not use these roads.

Even if they try to get educated, they do not have the means or money to get education. They do not get enough food. To solve all these problems it has become basic to abolish the religion that alone stands as the stumbling block for development. Protecting the religion and saying that everyone should get education is an ignorance and hypocrisy. So we should abolish the rules given in the religion at first. To do so, we need courage, only when we have that courage we can achieve freedom and education.
[Kudiarasu, 22-8-1937]

**Brahmins and Untouchability**
Will the Brahmins be ready to give up any of the rights enjoyed by them? Likewise, will the Muslims be ready to divide any of the rights enjoyed by them? Will the Brahmins and the 'high' caste people agree to wipe off the disgrace and the cruelties suffered by the non-Brahmins and the Depressed Classes?
[Pagutharivu, 19-9-1937]

**Hinduism and Caste System**
The only reason for several divisions among the Indians is our religion. What is the reason for several castes: that too high-low, Brahmin-Paraiyar? Is not our religion the only cause for all this? Sastri, B.A., B.L., said that the low caste persons could touch the high caste persons only after death. It is impossible in this birth. The only basis is the Sastras and Vedas of the Hinduism. For, only Hindu religion gives the Brahmins all comforts, social status, money, freedom etc., so such a religion must be destroyed because it discriminates others as Untouchables. There is no difference in caste among animals as untouchable dog, Brahmin dog etc.
[Kudiarasu, 19-12-1937]

**Untouchability eradication**
Even now research is going on into knowing the caste of the child who will be born when a woman of one caste and a man of another caste unite. If someone like us says, "Why sir! People of other countries have progressed so much, why we are still indulging in this research? Is it correct?" They immediately get angry and say, "What do you know? How much have you studied? You come and advice us! The research of foreign countries has only two day's life."

From this, think of the progress of our country and the progress of other countries. In general, don't you agree that in our country, only seven out of hundred people have been educated? Secondly, don't you agree that our people don't have the means to eat, don't get a job or the wages, in tens of thousands they board the ships with their children and their pregnant women and go to lands that they haven't seen before or even heard of before and die there? Don't you accept that man considers another as low and shameful and castes are divided into upper and lower and are called Panchamas and Mlechas and Sudras, and in daily life they are belittled and cruelly treated in a manner that degrades their self-respect?

If that is so, you think and see whether is it possible for our country to progress by maintaining this state? If we think of the cruelty, irrationality and foolishness of untouchability, our mind does not give even a little place for us to forgive that or to be indifferent, or to say that 'it can be looked into tomorrow,' 'what is the need or urgency for it now,' or to postpone it.
[Pagutharivu, 1938, Issue 10, No 3]

**Untouchability eradication**
Though several cruelties that do not exist anywhere in the world and are unacceptable to humanity, justice and rationalism exist in the Indian nation—among them the issue that needs to be solved urgently, and in order to be respected by the peoples of the world is untouchability. To prove that Indians are not barbarians but are a humane and belong to a civilized society—and if they need to say that they are capable of protecting their own nation and establishing rule and governance like the majority of the peoples of the countries of the world — we can firmly state that there are two tasks that need to be completed at first. Among these, the first task is to abolish the practice of considering, as untouchable at birth, many groups with a population of several crores in the Indian soil. Secondly to abolish or refrain from giving them a treatment that is more dreadful than what is given to the animals that lack the sense of reason, and lower than given to the worms and insects that lack the sense of feeling.
[Pagutharivu 1938]

**Horror of Untouchability**
If we think of the horror and absurdity and brutality in the concept of untouchability, our mind does not allow us even a little to forgive it, or disregard it, or casually postpone it by saying, "what is the necessity or urgency, this can be dealt with later." Sometimes we feel that if one wants to make those people who consider others as untouchable and mete out cruelties, to realize the true difficulties of the Untouchables, the horrors experienced under the present British rule are not enough and we need a tyranny with a continuous military rule that is devoid of any little freedom or equality. Only then these people will become sensitive and achieve realization.

But because India is in the tight grasp of not only such brutality and lack of legitimacy, but also of stupidity, so no matter how much humiliation and cruelties take place it will be difficult to bring about a situation where such people realize the real sorrow. Yet, there is place for the belief that this state of affairs will change.
[Pagutharivu, 1938]

**Congress and Untouchability**

The Public Works Department Minister Honourable Yakub Hassan inaugurated the South Tanjore Zilla Congress Conference. In that conference, three Adi Dravidas were insulted for the sin of having taken part in a community inter-dining (*samabandhi bhojan*) and it has created a great flutter in south India. With the exception of the Congress papers, all the other newspapers carried condemnation statements. But the Congressmen who have girded up their loins to emancipate the 'Harijans' are maintaining a silence. No Congress paper has till now condemned this atrocity and written against it.

In the name of the Congress, the 'Harijan' members of the Chennai Assembly, or the 'Harijan' minister have not opened their mouth till now. It is being advertised that the 'Harijan' Minister Honourable Munusamy Pillai and the 'Harijan' Mayor comrade J.Sivashanmugam Pillai, are doing community dining with caste-Hindus. But the right of 'community dining'—is it only for the 'Harijan' minister and Mayor, is it not for the Harijan society? Is it just and fair when they accept the respect for themselves, a respect that their society does not have? Is it their opinion that whatever the status of their society, it is enough if they get position and money? Are they not shown the respect of community dining only because the 'Harijan' ministers and 'Harijan' mayors will clap hands for the Brahmins' farce and remain as henchmen?

When the injustice at the South Tanjore Congress Conference has caused mental turmoil to the entire Adi Dravidar community, it is not proper that the Honourable Munusamy Pillai and Mayor Sivashanmugam Pillai remain silent. If they are unable to search for a remedy for this humiliation, it is just and fair that they resign their positions. What are Honourable Munusamy Pillai and Mayor Sivashanmugam Pillai going to do?

[Kudiarasu, 16-1-1938]

**Brahmins and caste arrogance**

People other than Sudras and Untouchables will not know about the importance and greatness of this social movement.

In the olden days, the income tax collectors who were just paid Rs.10 to Rs.15 as monthly pay, who were invariably Brahmins used to address the income tax payer (naturally they should be very rich in comparison with the tax collectors) in singular respectless terms and also these rich men used to stand up with folded hands. Thus in these ways the arrogant Brahmins held respect. But now because of the Self Respect movement we can see the difference.

[Kudiarasu, 6-3-1938]

**Hinduism and caste**

In this country, there is only one caste that is not called untouchable – that is the Brahmin caste. Other than that people of all other castes are only considered to be Untouchables.

Europeans, Muslims, Indian Christians and those non-Brahmins who are called Hindus are only Untouchables to the Brahmins.

This is not being merely uttered from the mouth; we are made inferior by their saying that this is according to religious doctrines and Veda Sastra rules. If someone is said to be a Hindu, he will have a caste and he has to accept that. Still, not only the Hindu religion, but the Hindu Gods are also included in caste discrimination.

So, as long as a man wants to be in the Hindu crowd, he will come under caste discrimination, and as long as he is a non-Brahmin, he has to remain a lower caste. How many ever Gandhis arise, how much ever they are advertised as Mahatmas, caste discrimination and 'upper'-'lower' castes will not be annihilated any day. For this, till now several people have taken so much effort. Everyone who has tried so has been made into a slave of the Brahmins, and he who did not become a slave has been destroyed, they have not allowed anyone to rear his head.

[Kudiarasu, 13-3-1938]

**Practices of Untouchability**

For eating together with the upper castes, a man is tied up, beaten, his head is shaved, and mud and cow dung are mixed with water and poured on his head. Power is also only in the hands of high caste and it does not reach the low castes. Even today, the two-tumbler system is followed in hotels.

[Kudiarasu, 3-7-1938]

**Self-Respect Movement**

What has ripened the people's resolve to such an extent that they consider the throwing open of the Travancore temple, changing the gruesome state where once the Untouchables could not walk in the streets of that *samasthanam* (princely state), and yesterday, the Adi-Dravidars entering the Madurai Azhagar temple as very casual events?

Any ritual is not possible without the Brahmins, there is respect and honour only if the Brahmins come: such a state has changed. And today people think it is a shame to call the Brahmins for rituals because it is damaging to our self-respect. What brought about this change? Is this all? Who made the people who were Depressed Classes and considered disgraceful as a Mayor and Minister and gave them importance? Only the Self-Respect Movement started all this and countless such revolutionary activities in this country. Can anyone deny this? If anyone denies it, I pray that they may look at the propaganda we have been doing for the past 15 years and the resolutions we adopted in our conferences and attain clarity.

When we were indulging in so much of uproar, we did not have the necessary monetary strength, government support, or the help of intellectuals, the educated, the rich, the landowners, or the important people. All these work, we stood alone, without capturing politics, we were victims of oppression, and we accepted great difficulties. Today if all the Tamils, pundits, intellectuals, important people, rich, political leaders, several Tamils who are said to be protectors of religions are today going to be united, we can achieve more success than what we have so far achieved.

If we have unity, the capacity of living together without jealousy, discipline, selflessness, and no thirst for power, not only will we achieve success, it is certain that we will also become the guide for a great revolution in India.

[Kudiarasu, 16-07-1939]

### Justice Party

They preach that the Justice Party is an anti-Hindu movement. We should not even demonstrate. If we demonstrate, they say that it is legally wrong and they say we should be hanged! Intellectuals say that the prevailing laws are tyrannical.

[Viduthalai, 2-8-1938]

### Congress and Untouchability

I have often warned about Congress and Gandhi not because the British should leave India or not. We should have political power because the people who came here to make their living have practiced all types of cruelties against the natives and have succeeded in it. They have enslaved the natives and have overpowered them. Of the three (Aryans, Muslims and British) who ruled us, the Aryans were the cruelest to us. They have not only enslaved us but also made us barbarians by their religion and politically ruined us. Even if the British rule over us for 1000 years more or if the Muslims rule us for over 16,000 years still we do not want the Aryan rule, that is why the Congress and Gandhi who happen to be instrument and coolies of the Aryans must be totally abandoned from the root.

[Kudiarasu, 7-8-1938]

### Varnashrama Dharma

Gandhi asks the non-Brahmins to work without even thinking of saving a small silver coin, and follow their hereditary occupation on the basis of Varnashrama dharma. So he established the *Vardha* scheme. Based on this they banned alcohol and closed all the schools by stating that there was no revenue to run schools since there was no revenue due to prohibition (ban of liquor sales).
[Kudiarasu, 7-8-1938]

**Linguistic identity**

If a Tamilian considers himself to be of an unmixed/ pure Tamil race he should first come out of the religion which has no relation with him and which calls him as Sudra and Chandala. To be more precise, the low status of a Tamilian is because of the fact that he considers the Hindu religion as his religion and India as his nation and that labouring for other people and nation is the true service to nation. Indian Patriotism in Tamilian means serving people of other states and being a slave to people of other states.
[Kudiarasu, 28-10-1938]

**Casteism and Racism**

There is a proverb in the Bible:"First remove the beam from your own eye and then you can see clearly to remove the speck of chaff that is in your brother's eye. (Luke 6:42)" This proverb comes to our mind while the Congressmen speak of South Africa.

We ask why do they shed crocodile tears over what happens in distant South Africa while they do not worry about the people who are born in this country and are segregated, harassed and disgraced.

If the people of our country have to be treated with respect by people of other countries, the practice of degrading, oppressing and humiliating the people of our own country must be annihilated. We must gain the right to live with self-respect. Only when we destroy this, if we ask something to the foreigners they will give it weightage, only then we can have the right to ask.
[Kudiarasu, 04-02-1940]

**Working castes**

It is our opinion that the farm labourer, the viswakarma (carpenter), the weaver, the leather-worker are all the same. Not only have they been differentiated by virtue of caste, but they are also the necessary base and vital for the human society. They function in the same manner. Why should there be difference between them? Apart from the difference as per the religion and Sastras, what is the evident difference?

So, the emancipation of labourers shall lie in annihilation of caste differences that create the hierarchy of labour, and the hierarchy of labourers.

First, there should be no caste difference between the proletariats. There should be no higher or lower labour.

Everyone should unite. The religious Sastras that teach difference and degradation must be destroyed. Politics must be captured. Those who are called the higher castes and those who live lazily without doing physical labour must be considered as enemies and those concepts must be utterly destroyed. Every labourer class should think that the world exists because of them. Only because of them people are able to live. They must realize that they are the foundation for the shared human existence.

Because the working class people didn't realize it and they became henchmen of the Brahmins who love to lead a lazy life, even till today the labourers, or the labour methods have not attained any development. The non-Brahmin party that is the Justice Party, even when it got a little political dominance it carried out several efforts to abolish this caste difference. The people who were Backward and Depressed classes on account of their labour largely attained the benefits. Whoever the labourers are, whether it is caste-wise labor, or coolie labour, everyone who toils with their body are only non-Brahmins.

Rich men or upper caste or those who gained supremacy through the Laws of Manu, are not the authorities in the Justice Party. As far as possible, the Justice Party keeps the upper caste and the rich people without influence in the party.

If the labourers, the poor people, the lay people who form ninety percentage of the non-Brahmins fail to use this good opportunity and try to improve their status, it will mean that the Tamilians are only fit for being Sudras and Untouchables under the rule of Manu Dharma and that the time for a solution to their suffering and disgrace has not yet come.

Tamils! Intellectuals! Brave men! We fold our hands and welcome you to come and join this organization of the Tamils!

[Viduthalai, 16-2-1940]

**Aryans and Dravidians**

In this age, even today there are streets with hotels, idli shops, water pandals, gods, temples where Aryans and Dravidians are given different places. Dravidians should not enter the places where Aryans enter because if a Dravidian enters the place or area it becomes contaminated. So the Dravidians were punished. On seeing this, should not men with self-respect and human dignity think of self-respect and freedom than thinking of national politics? In the nation that was ruled by us, in which we lived, the cunning, deceitful, dishonest Aryans who came for their living call themselves as upper caste, and are well treated, whereas we are ill treated in our own country as low people. Till this stigma is removed do we need positions, posts or degrees?

Now it is our primary duty to tell our people about the social and religious discrimination between Aryans and Dravidians. For this we need a separate campaign, a separate establishment and a separate army. It will not become a crime if we realize our situation and if we struggle to remove our disgrace.

After the war is over, this should be our work. For this itself, we say, "The Dravidar Nadu is for the Dravidians." If the government fails to partition the country and create us a separate politics, we should be ready to die in that struggle.

It will only be a tyranny, no matter how the rule is, if someone comes to another person's land and dominates and makes the people of that land as Sudras and low castes.

Today, we do not have the name of Dravidians. We do not have the name of Tamils. We only have the names of Sudras, Untouchables, fourth caste, fifth caste etc. Is this right? Do you agree with this? Should this not be changed? Must we not work for it? If we are to consider it a shame that the British did not call our party's president, then how great a shame it is that the 3.5 crore people of our society are called Sudras, low caste, untouchable and kept apart by a begging group (Brahmins)? Is the Minister's post, the post of member of the Viceroy's Council, the post of Advisor, Governor etc. more important or more urgent or necessary?

Let the state of the old people, the shameless money crazy, power-crazy be anyway. You, who are the youth, come forward with firmness, courage and resolve. In ten years we can carry out a great social revolution and attain honour and liberation.

[Kudiarasu, 28-7-1940]

**Congress and Casteism**

In this century, even in this 20$^{th}$ century, we carry the title of being a disgraceful caste. Temples, tanks, hotels and such places, separate place is allotted for Brahmins and for the Sudras—so we are made into Sudras. This should be changed. We need to struggle for this. Non-Brahmin people in Congress, without even a little shame, join the chorus of the Brahmins and shout, "Annihilate Imperialism." Congress non-Brahmins who do not have the capacity or shame to abolish advertisement boards displaying Brahmins with durba grass and showing us as lower caste, are they going to abolish imperialism?

That too, when comrades Rajagopalachari, Satyamoorthy and Jawaharlal worship the temple of imperialism and circumambulate it from left and right, what a farce it is when these Congress non-Brahmins say that imperialism should be annihilated. If they have even a little shame, let them burn the durba grass, let them pull the advertisement boards and fling it off. Instead of doing so, don't they not shamelessly ingratiate with the Brahmins? Moreover, Comrade Lingam who spoke about the municipality elections that are to take place this month in Erode asked whether the Municipality would not leave way. We do not require that. I am willing to leave and give way. Let any Congress leader get up and say, let them accept that at least as far as this place is concerned in the coffee shops and Brahmin hotels, the boards that say 'Brahmins — Non-Brahmins' will be removed and flung off. Tomorrow itself I shall ask our party nominee to withdraw and today I shall make the ten Justice Party councilors in the Municipal Council to submit their resignation.

Instead of that, in the manner of saying, "Your property is mine. My property is also mine," they feel that we should not bother about our being a low caste. We should also accept that they are higher caste. If this is said to be freedom should we keep nodding like temple-cows. For me, somehow, as I think over this again and again it gives me anger. So, it is our duty that we need to win over any difficulty and sacrifice our all.

[Kudiarasu, 11-8-1940]

**Adi Dravidars**

In our country along with the difference between the Dravidian people as Dravidar and Adi Dravidar, the Adi Dravidar society has continued to be a very large number. In the Dravidian nation, how are the Dravidians untouchable for the Aryans who came from outside? This situation is a very shameful situation for the Dravidian society and it protects the Aryans by dividing and treating the Dravidians as untouchable people.

So, the name Adi Dravidar itself must be changed. Both must be called Dravidians or Tamilians. All social differences and discrimination between Dravidians and Adi Dravidars must be eradicated. It is my desire that they must become one. In both these cases, in my opinion this is the objective of the Justice Party, so I say these in the name of that party. It must be made a special duty of our party to ensure that in matters like employment and education for the Adi Dravidars special privileges must be shown so that they attain equality with us.
[Kudiarasu, 25-8-1940]

**Social Reform**
We have so far not got any social reforms through political reforms. Any small social justice we have obtained means it is only through the demonstration and never by law.
[Kudiarasu, 25-8-1940]

**Brahmin domination in employment**
The statistics collected by the government in the year 1939 about the pay and the number of posts held by Brahmins and non-Brahmins, is given by the following table:

| S. No. | Pay | No. of Brahmins | No. of Non-Brahmins |
|---|---|---|---|
| 1. | Gazetted officers and administrative officers pay from Rs.300/- to Rs.5000 | 609 | 398 |
| 2. | Posts in the pay scale greater than 100 | 3667 | 2492 |
| 3. | Pay above Rs.35 like *Gumasta* | 9183 | 8042 |
| 4. | Pay less than Rs.35 peon, coolies (running errands) | 1513 | 33662 |

So from the statistics it is evident that the Brahmins who are just 3% of the population occupy more number of the better positions i.e., 3% holding 609 posts and 97% holding only 398 posts. When it comes to coolies etc. the 3% represent 1513 posts whereas 97% hold 33662 posts. Even in these conditions we (non-Brahmins) are portrayed as hunters for jobs and they are described as nationalists and martyrs of the nation who renounce everything for the sake of nation.
[Kudiarasu, 25-8-1940]

**Common dining: Brahmins and Non-Brahmins**

In the village of Thiruvaiyaru in the Tanjore district, there is a Sanskrit College, which functions by taking grants from the Tanjore King. That college which is under the control of the Tanjore Zilla Board taught only Sanskrit. When the late Dravidar Mani Sir Paneerselvam was the President of the Tanjore Zilla Board, he ordered that Tamil should also be taught there. The *agragaram* (Brahmin settlement) fumed, but he did not get afraid.

On December 2, Tanjore Zilla Board President Comrade Nadimuthu Pillai passed a resolution in the Zilla Board meeting that there should be no caste difference in the college hostel and that both the Brahmins and the non-Brahmins must eat together. That resolution comes into effect from 27 January.

In this hostel, there are 70 non-Brahmin students and 45 Brahmin students.

After this common dining resolution came into effect, these 45 Brahmin students refused to eat in the hostel in order to condemn the resolution of the Zilla Board. They went to the agraharam and ate. The agraharam is agitated and roaring seeing this resolution, which says that Brahmins and non-Brahmins must eat together.

What is the idea? It means that even today, Brahmins are higher caste than the Tamils and that they are a separate caste. It only shows this arrogant opinion. It only shows their *varnashrama* mindset that if they sit together and eat with the non-Brahmins, they will get polluted, and they will come to harm. If they are together with the non-Brahmins and they eat, will the food turn into poison? Will their intestines blacken? What will be harmed?

What lesson are you going to teach for the Brahmin arrogance that refuses to eat with the Tamils? Does your shame not arise! Do you not feel ashamed and sad? When would the Brahmin caste--that refuses to accept that you are also men--stop it's bragging? When will you attain self-respect? Think and see. Resolve.

[Viduthalai, 29-01-1941]

**Press**

While speaking at the meeting in Ranade Hall, Mahaganam Sastri said,

"A great danger has come to the country now. Somehow I am feeling greatly afraid, I am feeling greatly worried. I feel that those who need to save me have betrayed me."

This has been published in the Swadeshamitran dated 27, page 4, 3rd column. But Sastri has not explained in his speech of what danger has come about that has caused him to worry.

Today the truth is that the British are in danger, in worry and in major loss: we can consider that itself as the danger that has come upon us. But Sastri or the people of his caste do not understand what danger there is. In railway hotels, discrimination and differences have been removed. The agitation from now onwards that such differences must be removed has become stronger. This is a danger for Sastri and his caste.

The agitation "Dravidar Nadu is for the Dravidians" has been gaining strength. Of the Indian party leaders, two of them that is Dr. Jinnah and Dr. Ambedkar have supported this. This is the second danger. Apart from these two, it is not understood what dangers the Brahmins have. To escape from this danger, Sastri charts the course: Those in jail instead of worrying about Gandhi will from now onwards not do Satyagraha, they will give an oath to the government, come out and accept ministerial positions and govern."

So Mahaganam Sastri and other Brahmins, consider it as a danger because Brahmin domination will not be there.

Nothing else. But we the Tamilians should ensure that the Brahmin domination must not continue. This is the important work that has to be carried out and we should be prepared to pay any price for the same.

If they think 4 Brahmins can assemble in a small room called the Ranade Hall, give information to the newspapers on what they spoke and send items, 4 columns, 5 columns, 7 columns will be published. They can boss over the people in the name of Swaraj and the nationalism.

In the Ranade Hall meeting, there were 4 Brahmins. Sastri presided over the meeting; comrade T.T.Krishnamachari was the speaker. One comrade Natesa Iyer delivered the welcome address. Another comrade K.S.Ramasamy Sastri was the audience. In the newspapers this came as a seven-column news item. They think that 'Indian' means only the four of them. They think that their statements are the Vedas. Our question is will this falsehood go on even now?

[Viduthalai, 01-05-1941]

**Press**

In Tamil Nadu, the Tamilian is considered as the fourth and fifth caste; he is considered as the fifth caste; he is considered untouchable, unapproachable; and unseeable and is disgraced. If we have to say it in brief, a man is not considered as a man, he is considered lower than a dog, a pig and even shit.

Whatever a man has or not, if he has the feeling of shame it is enough. A man with shame will not like a life of slavery, he will not accept another person humiliating him, he will not bear someone sitting on his head and not letting him get up. That is why we say that self-respect is the life of the Tamilian.

Think a little of how the Tamilian's shame is today. In the coffee clubs he is segregated as a low caste. In the temples he is made to stand behind the upper castes, because he is considered lower caste. The same way in the hostels that are being run for students, in the halls that have been built for charity, in the dinners being held by Brahmins espousing nationalism, Tamilian is considered as a low caste and he is segregated.

Even if we do not want to fight with them, they will not leave us. They will start a fight on their own accord. We are not cowards who will give up the fight.

In short, for our struggle, we have a few supporters and several enemies. In this situation a newspaper is a powerful weapon for our fight; the daily paper is only *Viduthalai*.

Our struggle for social rights is our revolution. We have not resolved to carry out this revolution by bearing arms or by the cowardliness called Ahimsa. Through creating emotion among the Tamilians, through showing the present disgraceful state of the Tamils, we want to unite them. We want to unite them and remove their social humiliation through a revolution. Only to create the emotion of revolution, our Viduthalai will work tirelessly. From today, Viduthalai is being published from Chennai. All our comrades should support it.

[Viduthalai, 01-05-1941]

**Caste-Annihilation**

He who does not strive to annihilate the caste system is not a Tamilian. One has to say that there is no Tamil blood in his body.

The superstitions of the Tamils must be eradicated. There is no limit to the superstitions of the Tamils. The Tamilian buys ticket for *moksha* (salvation) from the Brahmin. The Tamilian marries his god and lets the god sleep. The Tamilian takes the god to the home of prostitutes and keeps. For the sake of moksha, the Tamilian mixes urine and cow-dung and drinks it. He squanders away his earned wealth for all this. The Tamilian considers himself as a low caste. Must all this not go? Only for this the Tamilian movement, the Self-respect movement, the Justice movement are there. Only their efforts have caused such reform marriages and agitations for the removal of social disgrace.

The British don't consider us as a low-caste by birth. They don't obstruct giving social rights to us. But they allow the Aryans to dominate and exploit us. The British can look into and remove our social disgrace through laws, but we are not united enough to compel the British. Only because of their unity, the Aryans are treating us so badly. They threaten the government also. Not only that. A lot of Tamils are slaves of Aryans, coolies, in religion, society and politics and they betray us to the Aryans. If this situation comes to an end, and all the Tamils unite, within the blink of an eye, our lowness can be removed.

Though we took several efforts because nothing could be achieved, we are today agitating that "Tamil Nadu is for Tamils alone."

[Viduthalai, 10-05-1941]

**Current Affairs**

Comrades! I do not understand what to speak under the title of "Current Affairs"—one can speak about anything because of the title.

The war is going on. The Congressmen are acting out the drama of Satyagraha. Jinnah is demanding Pakistan. We are asking for Dravidar Nadu. The Hindu Mahasabha people are causing Hindu-Muslim riots. The Congressmen say that somehow the prohibition (on liquor) must go. We say that we need Muslim-Dravidian unity.

The government is collecting troops and money for the war.

Jinnah and we say that we should work together with Adi-Dravidars. Srinivasa Sastri and Sivasamy Iyer say that for a Brahmin and non-Brahmin to sit together and eat is against the law.

The Aryans went to court and quashed the order that the Aryans and Tamilians can pray together at the Tiruvannamalai temple.

The affluent, major landowners, and zamindars among our Tamil people do not have any feeling of shame. They flirt with our enemies. They even dare to complain to the government and to betray us. Our Governor is a good man, he will do justice, but it looks like they will spoil him. I don't go and see him unless he calls me. We need to join with the Muslims and Adi Dravidars and not bother about anybody else and do major agitations.

We should consider the partitioning of our country as our Swaraj and complete independence.

The leadership position, handicrafts, business, industry, and comfortable jobs have gone to our enemies. We are headless corpses. We will see that the support of Jinnah at this juncture is a victory.

The Act for the removal of social humiliation introduced by Comrade M. C. Rajah which was passed unanimously has been useless. The temple entry act of Comrade Rajagopalachari is also going to become likewise. So, all of you should support our agitation for the removal of social ills. Even if the need arises for us to go to jail in thousands, we must bravely stand at the forefront. Now, going to jail is the real and intense agitation for the government and for the lay people. Those who are afraid of that are equivalent to corpses. So be ready to go to jail.
[Viduthalai, 15-05-1941]

**Social Reform**
I do not have faith in reform. The little faith that I had has gone away. I believe only in revolution. Something must be done, only after that I should die. Death is close at hand. I have laboured not just a little. Only I know that. I am alone capable of respecting it. Only with me there is good strength. So, only I have the credentials to tell a decision about myself. Not those who are passive onlookers who attain good if it comes and shrug off any connection if evils come. Those who criticize and advice me, let them pay a little attention to their status and rights.

Only the worry that I have to soon do something that strikes me as the right thing to do, makes me decide some path and I do not have any ill-will towards anybody. How can one think that any work can be done with the people who feel worried and complain about the hatred towards

*Periya Puranam* and *Ramayana*? In the year 1922 itself, I said that these ought to be burnt. In 1928, in several presidencies outside and in our state, it has been burnt. In 1935 it has been burnt by the Adi Dravidars. Those who are protesting now are only the Aryan slaves. In this situation, they need to quit, or I must leave them and quit.

Why do I need a relationship that is devoid of goodwill and honesty? I need to tell things soon because my disease matures. You will ask me, "Why did you not join the Muslim society?"

If I go alone what can I do there? Even there I need at least ten thousand people, only then I can be with rights. You will say that mullahs and superstitions exist even there.

If we go with a sizeable crowd, even our country will become a true Muslim state like the Muslim countries of Egypt and Turkey. There is no necessity for superstitions in the rules of the Muslim religion. Even if it is there, I know about Egypt and Turkey. I say that it is enough if Egypt, Turkey and Persia remain Muslim. There is no need for garbs there. There is the essentiality only for the heart. The government officials will behave with integrity then. There will be no Brahmin menace. They will tremble or run away to Hindustan (north India). We will not have any type of disgrace or problem. We can be a large egalitarian society. We can be independent Muslims. We will also get the support of Pakistan in the north. Exploitation will be eradicated. I say this only for the time when nothing is possible. I say this because nobody should say that we do not have any plan. We must not be so.

Youth must think over this. Comrades Ambedkar and several personalities involved in public service have also said the same. Still I am not saying this—as some enemies and poisonous people allege—because Islam and the Islamic God will easily forgive the sins and make heaven nearer as compared to other religions.

Man is a creature who has to live together. That has been spoilt by Hinduism, gods of Hinduism, the great people and big shots of Hinduism, its leaders and they have made us into slaves for the Brahmins, themselves and those who are favorable to them. To change that, quit Hinduism to become people who are fit for human social life. I say this only for that reason. If you dare to do this, the Brahmin will step down and come. Even the Government will certainly become credible. What more do you want?
[Viduthalai, 1-5-1943]

**Abolition of Caste and Communal Representation**
The plans of Justice Party differ from the Congress mainly in two ways. In Congress, all old puranas and caste differences should be protected. Each caste should protect its ethics and have the right to practice its codes and customs; it also guarantees them that it will protect such communal and ritualistic rights. Suppose we accept this plan, it implies that we have accepted our Sudrahood and the untouchable caste. Also it means that we accept varnashrama dharma, Pariah, Brahmin and low/ high status. Today's Congress demonstrations are only to protect Brahminism and Panchamahood and not for any public good.

If in the Indian republic, there was one line that states 'there will be no caste like Brahmins and Pariahs and there would be no castes, so no difference,' then I will not talk about communal representation or any such reservation. In that case, the Congress can find a place to protest, or to question my plans. On the contrary, Congress is protecting and establishing the old dharma as well as guaranteeing their right to practice differences in castes! If that is the case, what are our rights? And when we question these differences why do you call us traitors?

The Brahmins, after completely protecting their high casteness, go to the low castes and ask them to accept their low status and their heinousness. Are we such big fools to accept this from the Brahmins who call us as lowborn for the sake of their living! The crores of non-Brahmins should first take up the mission of annihilating caste and treat everyone with equality.

What is Swaraj other than this? Every one of you should think about this. Think of all the non-Brahmins who do not accept the Congress policies. They first accept the Brahmins' supremacy by which they directly accept their Sudrahood and thereby accept they are born of prostitutes and thereby they accept it as their ardent duty to serve the Brahmins. That is why when I ask for communal representation and annihilation of castes, by no means is it a contradiction to their policy.

Caste difference can be annihilated only by law and not by common man's acceptance. Based on caste, one class of people do not do any work but at the same time they exploit others and live on the labour of the majority, which labours very hard and goes without even a square meal a day. So to give communal representation and annihilate caste legally is not a difficult task.

If Congress abolishes castes by its initiatives or by law, why are we going to demand communal representation? This basic idea would go by itself. As long as caste is alive, there is no wrong in asking for communal representation in reservation in all employment and educational institutions. Is it wrong to question the Congress that is trying to save caste and caste differences? We can achieve complete equality by giving communal representation in all posts, which will certainly annihilate the disgrace of untouchability, the ruin due to alcoholism and lack of education. If we live in unity, we can certainly make reformations and thereby lead a proper life. The Congress says we are slaves to rich! How can we accept this? Only those who get money from the Brahmins for their own living say this. So we need not bother about all these.
[Dravidar Nadu, 30-5-1943]

**Dravidians and Hindus**
The resolution passed in the Tiruvarur conference in 1940, said that it is wrong to call ourselves as 'Indians', 'Hindus', or 'people of Bharath' and instead we should call ourselves 'Dravidar', 'Rationalist' and 'people of the Dravidian country.' As we do not have such feelings we are made to disgracefully suffer all the complaints and humiliation that hinders our progress. If we say we are 'Hindus', then according to the Varnashrama Dharma of the Hindu religion, Brahmin, Kshatriya, Vaishya, Sudra and Panchamas are the 5 divisions and we are in the $4^{th}$ and the $5^{th}$ divisions so we become the degraded lowest caste. These concepts exist according to law, religion, Sastras and the existing rules and regulations. This made us accept it wholeheartedly with tolerance. So we call ourselves only as Dravidians and not as Hindus.

[Kudiarasu, 29-9-1943]

**Religious discrimination**

In the name of religion and in the name of law, among the Hindus, there is the fourth caste, the fifth caste (Sudra, Panchama)—that is the Brahmins call us Sudra and 'low' people—and we receive only the corresponding legal rights, social rights and treatment.

[Kudiarasu, 6-11-1943]

**Casteist disgrace**

Though several among us consider them to be low castes, degraded castes and treat them so and likewise give judgments in the courts, they think, "if our respect goes, let it go, if we get wealth and position it is enough," and they try to forget and hide their shame. If someone reminds them, or exposes it, they only get angry or try to take revenge, but they do not have a correct idea of their disgrace.

[Kudiarasu, 6-11-1943]

**Casteism**

Only due to caste some people are rich and some remain poor.

[Kudiarasu, 25-3-1944]

**Caste domination**

The nation is deteriorating and degenerating with several social problems because of caste. Each caste feels that persons of their caste should come to position, so it dominates other castes. Thus there is problem in every village.

[Kudiarasu, 11-11-1944]

**Partition of Dravidar Nadu**

When we speak of the partition of Dravidar Nadu we have to make a note of

1. Madras Presidency alone shall be Dravidar Nadu.
2. If people of other states have to step into the Dravidar Nadu they have to get visa.
3. Only with the permission and after paying required duty, they can come to Dravidar Nadu to buy/ sell things.
4. The increase or decrease of the area of the country is dependent on our comforts and wishes.
5. Foreigners doing business or handicrafts will depend only on our needs.

Till we get complete independence the security will remain as it is. Muslims, Adi Dravidas, Christians, Buddhists are all Dravidas. Their religion and their religious activity will continue to be the same, that is, as per their wishes.

[Kudiarasu, 2-12-1944]

**Quit Hinduism**

If a foreigner rules our nation, for his benefit he controls us by divide and rule policy. Likewise, the native people of our nation are permanently controlled by the foreigners i.e. the outsider Aryans by using tools like god, religion, Vedas and Sastras that has divided us into classes.

　　Because of this, the majority of us have become backward class people and Untouchables and they rule us. So we have become base people having no right even to call ourselves as human. A few of us are not aware of it. Shed away the concept of god, religion, Puranas, Sastras and epics. I came here with all difficulties just to say this.

If you all do not follow this, even after a thousand years, whatever conferences or campaigns you hold, or demonstrations you carry out, whatever be the political independence you get, whatever be your economic development, whatever be your degrees and positions, the heinous disgrace will not be wiped in your society. This is certain! This is certain! All of you are going to do the same mistakes your forefathers have done. So, you would do the same mistakes in your lifetime also. You can only do mistakes like them but not rectify any of it.

So, all gods and religions that stands as evidence for our disgrace must be ruthlessly cut away. It is utter foolishness to think of correcting or reforming Hinduism or its related concepts like religion, god, or Sastras.

The only intelligent thing to do is to come away from Hindu religion to save ourselves from any form of disgrace. If we have to utter any other word for Hindu religion, it is nothing but Aryan or Brahminism. If you have any doubts regarding this, kindly see the dictionary written by intellectuals. Don't become victims of the traitors. What happened to the efforts of Buddha, Sankara and Ramanuja? Only to annihilate Buddha, heroes like Rama, Krishna and puranas like Ramayana and Gita were created. When Ambedkar visited Chennai, I elaborated to him on what I had said in a massive conference in 1923 that unless one burns the Ramayana, untouchability couldn't be abolished. Accordingly there are large groups that burn the Ramayana and also large support for it.

Over 15,000 people have left Hinduism and have come to the Self Respect Movement. They have thrown away their Sanskritised names and have taken a Tamil name. They have left all symbols, which identifies them as Hindus.

Our people want a separate state mainly because they want to lead a life without disgrace. They should not be the victims of Aryan supremacy and socio-economic ill treatment. The leaders of the Depressed Classes, Dr.Ambedkar and Sivaraj have also spoken in this manner.

Quit Hinduism and call yourself to be of any religion you like. By accepting the fact that we are born in Hinduism, we accept Varnashrama dharma. To wipe away the disgrace we have to abandon Aryan supremacy by abandoning god, religion, Sastras and puranas.

[Kudiarasu, 13-1-1945]

**Elections**

I will say something of importance to that society. Those who are the leaders of the Depressed Classes must not hold any office, and they must not have affinity, desire or necessity to hold office. Apart from leaders, the others may get position, titles, power, salary, and rewards, in the name of their society. But the leaders must not look that side. Ninety percent of the time, the Government gives position, title, job, salary, reward, to the leaders mainly to punish and oppress that section of people, make them amiable, make them enemies of their enemies, and not for any other reason. At the least, those who get such positions must leave the leadership to someone else and look after their work. If the same person takes both the leadership position and the position that gives authority, title and salary, it is difficult for him and his society to have respect in the eyes of the Government and the public.

So, this is a crucial time where we have to fight by uprooting and not respecting the laws and policies. At this time, the Legislative Assembly, the nature of the Legislative Assembly, is not possible for us, and even if it is possible it is not suitable for us. I say this on the basis of 20 years of continuous experience. So, at least to an extent, the people who are the Depressed Classes should leave title, position, legislative assembly, etc. to those who believe it or have the need for it. They need to stand out and get ready to struggle. I request the true and strong-willed comrades.

This is my answer to the comrades who asked me about elections to the legislative assembly.

[Kudiarasu, 2-2-1946]

**Caste System**

It must be mentioned that though the British have come to this country and ruled for 200 years, even today thousands of people continue to be Backward Classes and Depressed Classes, divided into thousands of groups and are differentiated as high and low, and on the basis of this there is struggle between them.

[Kudiarasu, 30-3-1946]

**Partition**

Whether the British leave truly or falsely, if Pakistan is created and Dravidar Nadu is attached to it, somehow the Aryan atrocity will be removed from the Dravidian land.

In this 20$^{th}$ century, the Dravidians—who are the indigenous people of this country, who are the descendants of kings who ruled not only this land, but several other countries—face a disgrace and barrier to humane existence and cruelty of slavery that even the Muslims who are called *Mlechas* and base by the Aryans don't face; and if this has been the practice for thousands of years, even after the rule of the British it has existed for two hundred years, and if Swaraj is being demanded in order to preserve this, not only will the Dravidian hate such a kind of Swaraj, why will he not ask for the rule of Pakistan, or rule of Afghanistan, or any other rule that is non-Aryan? This is our question.

[Kudiarasu, 13-4-1946]

**Gandhi and Untouchability**

The Aryans have purposefully fabricated several false stories without any honesty, they state these as the history of the nation in schools, so these have become ingrained in the brains of the people in such a way that they think it is true that the Brahmins are the highest honorable caste and the sons of this soil are degraded people and Untouchables. How can a rationalist with self-respect accept all these?

Festivals are celebrated according to the puranic stories, likewise, if political power from puranas is imposed, when will there be development?

Even we have become slaves of the Congressmen i.e. we live like corpses. Even today, nothing else but leadership in politics remains primary for Gandhi. He claims himself to be the leader of Hindus, Muslims, Christians, Sikhs and Untouchables. On this basis he is treated respectfully. But in politics can we accept bhajans and prayers to Rama if he (Gandhi) is a leader of one and all?

[Viduthalai, 12-6-1946]

**Dravidar Nadu**

If we had the feeling of our nation (Dravidar Nadu), feeling of independence and self-respect, can we have Nehru Park, Gandhi Chowk and Tilak Square here? Can we call this complete independence? By looking at this, does not your mind become tensed? Does your stomach not burn?

In South Africa, there is the cruelty of the whites and in Dravidar Nadu there are Aryan atrocities. Dravidians in Dravidar Nadu are tortured by Aryans. Dravidians in Dravidar Nadu are called by Aryans as Untouchables, Chandalas, sons of prostitutes, and slaves. Does any Dravidian become mentally tensed because of this? Are Aryans in a 'guru' status to the Dravidians? Is the Aryan a holy man, a master, boss, judge, collector, magistrate and family guru to a Sudra? What evidence do we need to state that we are shameless society and that our country is a barbaric slave country?

[Kudiarasu, 6-7-1946]

**Hinduism and Untouchability**

Gandhi often says that untouchability is the curse on Hinduism. He says this because of the hidden reason that Hinduism has to be protected. If Gandhi was a true person or if he was intelligent, he would have understood that if the Hindu religion was not there untouchability would not be there and that if untouchability was not there, Hinduism would not be there. Hinduism exists in order to protect Untouchability, and Untouchability is not there to protect Hinduism, or as a curse on Hinduism. This is because we know that anywhere in the world where there is no untouchability, Hinduism is not there, and where there is no Hindu religion, there is no untouchability.

No one can show a specific doctrine anywhere that is supposed to be Hinduism.

But to create Untouchability and caste difference at birth, to protect it and to degrade a few people, and to conjoin in one society the people and the degraded people, these two aspects are properly and carefully stressed and this conspiracy that was created through experience is given the name Hinduism. And there is no place or truth or evidence to say that someone created it at some point of time as a religion. So, it is impossible to have a Hindu religion without untouchability or untouchability-eradicated Hindu religion.

If untouchability is abolished there will be no need for anyone to save Hinduism. If untouchability is eradicated, not only will the Hindus have to join one of the religions of Islam, Christianity, Brahmo Samaj, Jeevakarunya Samaj, Universal Brotherhood, and there will be no place or necessity for them to say that they are Hindu.

Some may say "there are several Gods in Hinduism, unlike in other religions. So, at least for several Gods to exist there will be something called Hinduism." These several gods were also created only on the basis of untouchability and caste difference and there is no one God who is without caste difference. After caste distinction and untouchability are annihilated, a lot of its foundation will be lost so the concept of polytheism will not have place and there will not be a religion separately for that purpose.

So, because of caste discrimination, the Hindu religion is the dream weapon of a malicious lazy group that oppresses others and makes them into slaves and grows fat on others' labour.

[Kudiarasu, 24-8-1946]

**Congress and Untouchability**

Congressmen used to blame all the previous regimes saying that they were not concerned about the eradication of untouchability. They boasted, "If only we capture the ministerial cabinet, in one second we will eradicate untouchability through laws."

We have never considered that merely due to getting the temple entry right, our indigenous people will get education or improved living standards. Even now we do not think so. Because, if see whether all the people of the different sections of society who can enter into the temple today have developed in the field of education and economy, we can understand that temple entry is not an important problem today. Moreover, we see from experience that because of temple entry, only the intelligence, money and time of people are getting ruined.

Since the temple is the abode of religion, god, codes and tradition that are (said to be) obstacles to the rights that have to be given to the indigenous people, we need to support the temple entry struggle.

[Viduthalai, 12-9-1946]

**Abolition of Untouchability**

It is not possible to eradicate untouchability unless a law that makes untouchability a criminal offence is passed.

Today, the Legislative Assemblies have become Congressized. The ministers are only Congressmen. So, if their whole-hearted feeling that untouchability has to be eradicated, why are they still silent? At least from this, will our Depressed Class comrades understand that 'Harijan Seva' is a mere deception?

[Viduthalai, 12-9-1946]

**Naming/ Terming**

I was shocked to see in the records that they mentioned the social status as Christian, Brahmin, and Muslim etc. but for all other it was mentioned as non-Brahmins. We non-Brahmins are the sons of this soil and how are we termed. I suggested that we should be called as Dravidians and the Brahmins must be called as non-Dravidians or Aryans.

[Viduthalai, 18-9-1946]

**Congress and Caste**

The Dravidian comrades in the Congress! Why do you stand as an obstacle for us? We do not need the British; we do not need them for anything. We, who ruled the country, serve as peons, butlers and constables after the British came. But the section of people who were begging, today serve as High Court Judge, Advocate General, Diwan, Minister, Sankaracharya and Bhagwans.

The Aryans know that once the British quit, only you, the Dravidian comrades of the Congress, will cut the tuft and the sacred thread that are symbol of the Brahmin's advancement and our lowliness.

The British are supportive to the Brahmins. This is a pact between the white Aryan and the yellow Aryan. What is going to be produced from now is the same pact. All the rest are gimmicks. I said this two years ago. Ambedkar and Jinnah say this now. Swaraj is only an effort to protect today's social structure. An example for this is that when we said that there should be no difference in food in Tiruvaiyaru, Mahaganam Sastri, P.S.Sivasamy Iyer, T.R.Venkatarama Sastri and others said to the British, "You are violating our agreement." This is what we need to break.

In a country with a great deal of caste and religious differences, what is wrong if the Muslims, a group that has only one doctrine, decide to separate themselves, that too in places where they are greater than sixty percentage? Why this obstacle? Will the country split? Will the country explode? If there is opposition to this, it is because of the revolt created by the British and the Aryans. Do the British not know of the strife and vulgarities here? Will the British not laugh when they see that there is a caste called the Brahmin, and a caste called the Pariah and that a caste cannot enter hotels?

[Kudiarasu, 9-10-1946]

**Sudras and Untouchables**

As said in the Ramayana, even if a pair of slippers rule (the nation) we are not worried. But the Sudra title that we, and the working class people of our race have in the Sastras, in the laws, and in practice must be annihilated. The Brahmin title of the lazy cheats must be abolished. Is there any Congress doctrine or policy that caste differences and the stigma of birth must be annihilated from Sastras, religion and Gods? Has any leader said this at any time? If you want send a telegram to Gandhi now. See if he agrees. Right from Raman's time, the upper and lower castes such as Brahmins, Paraiyans and Sudras existed. The Ramayana says that only after the Sudra who prayed to God was killed and the dead Brahmin child came back to life. Right from Harishchandra's time, the Paraiyan caste exists. So, only to remove such humiliation, the Dravidar Kazhagam functions as a true, incomparable establishment.

[Kudiarasu, 12-10-1946]

**Dravidians**

If the Madras Presidency has to function with complete independence, anyone who has to enter the Madras Presidency—whether it is the Birlas, Tatas, Mahatmas, Nehrus, British—they need to get a visa. No north Indian should be our leader, President or Mahatma.

Man should live with equality. In the human race, there should be no Paraiyar, Sudra, Chakkiliyar or Brahmin or Untouchables. Like other nations of the world, we should live like a separate independent nation. Among our 5 crore people anyone can be a minister, a Mahatma. Don't we have qualified people?

Here the Dravidians are the fourth caste. A group of people who are not Dravidians are the first caste. 90% of the population who are Dravidians is considered as Sudras,

Panchamas and degraded people by the laws and the Sastras.

Should the Dravidians not live like humans? Should the Dravidar Nadu not become a nation without degraded people? In the present days there is no value for throne or crown. Today's need is only democratic rule. Whoever is an honest moral humanitarian can wield power.

The Indian subcontinent is 3000 miles long and 2000 miles broad. From ancient days, people have developed a lot of differences here. One does not eat in another's house, does not perform inter-caste marriages. Several languages, several castes, several upper and lower differences exist. They don't touch each other. If they happen to touch, they perform purification in various ways as mentioned in the Vedas.

There is no nation in the world with such differences. Can a big subcontinent with so many differences be a single nation; if it is to be under a single rule how can we accept it? If Dravidar Nadu is to be partitioned according to my request, what discomfort does any Indian undergo?
[Kudiarasu, 19-10-1946]

**Social Equality**
If we take the social status of lay people in western countries, there is no social difference based on birth. All are equal. Among the public, there is no feeling of high or low by birth. Everyone has every social right. Nobody has any kind of discomfort or disgrace because of his birth. The cruelty of untouchability does not exist there. Whether it is the temple, or the hotels, or any public place everyone has equal rights, no difference is entertained. On the basis of birth, education, and status, priesthood or religious headship is not created, and this practice of exploitation has not been formed (there). Anyone can attain any education, any job and any position.

Because the people remain as one society, and because they have equal rights, their unity and discipline and common aim has grown and those countries' have attained development because of their society and they are a completely independent society.

This country's society has been divided in the name of birth. So, here the people have been divided into several groups. Here, to respect anybody the first task is to consider his birth. In society, status has been ascribed based on birth; several discomforts have been created for man on the basis of his birth. Whatever a man's integrity, he is either exalted or made lowly because of his birth.
[Kudiarasu, 2.11.1946]

**Birth-based discrimination**
In India just by seeing a person—by his dress, symbols he puts on his forehead and language—one can say to which religion he belongs! To exhibit the external and internal religious feelings is considered to be good. There are several gods and their history is abundant. Because of these, society has been divided and unity is lost; they oppress and dominate others based on the lowness and highness of caste. The nation has been ruined only due to religion. This will only increase the difference and the notion of low caste.

While talking about political parties one can only talk more on caste divisions and not the principles. Apart from this, there is a separate political assembly for Adi Dravidars, which has many crores of members, likewise there is a special but a different assembly for the Brahmins. Sufficient money is not spent for the education, social equality of Untouchables, etc. Likewise, for the Christians. Though the Varnashrama dharma-ists are very small in number, yet they have a Mahasabha: the Congress that caters to all their needs.

The Brahmin lives in luxury and spends at least 70 times on himself compared to what an untouchable spends on himself. If we consider the social status of a common man in western countries there is no discrimination by birth, all are equal. There is no low or high in respecting a common man. They enjoy equal rights.

There is no discomfort or discrimination due to birth. There is no untouchability. The cruelty that one should not come near another does not exist. Whether it is the temples or streets or shops, all are equal. One is not exalted in education, or status or religion because of one's birth. In general, one does not exploit the other. Anyone can achieve education or any post. That is why these nations have become developed ones. But in our country social status is based on birth. That is why we are divided into several castes. So the first duty of anyone in India is to find his birth! Based on his birth, he is respected or disrespected!

Because of these, men face severe discomfort and discrimination. Based on birth, untouchability is practiced! Some cannot be touched, some cannot even be seen, these Untouchables cannot use common roads, lakes, temple, hotels, etc. Even in classrooms they are discriminated. Even jobs are based on their birth. A class of people are given menial jobs; like scavenging is allotted to the Depressed Classes. The nation has no common objective. All common objectives are for the comfort of one class of people, the Brahmins. If untouchability exists how can there be unity and love among people?

So one not only respects other religions, but also thinks very low of his religion and considers it as his enemy. So, only Hinduism is the cause for all differences, the lack of unity and common objective that has hindered all developments and made the nation remain as a backward one.

[Kudiarasu, 9-11-1946]

**Need for unity**

If we still continue to remain divided it will only create misery for us. If we are divided because of political and religious differences, and keep fighting and hitting one another, it is very very disgraceful to our race. Apart from the conspiracy of enemies, and shameless, selfish, mean-mindedness what else is the reason for the splits and fights among us?

Look at the Brahmins: do they have such fights among themselves? Has a Brahmin ever beaten another? Has he troubled another Brahmin? No matter how selfish, a Brahmin's selfishness is used for the welfare of his society.

This is a very crucial time. If we don't unite at this time, if we don't protect ourselves, we will come to be in an even lower position.

If all the Muslims are killed and if it is ensured that not one of them remains in the Dravidar Nadu—Come to a conclusion that this has been done! What next? What will happen next? Think and see!

You and I and Kamaraj and Muthurangam and Bhaktavatsalam, the Pandarams, the Rajah Sir, Maharaja Sir, Sir Ramasamy Mudaliar, Sir Shanmugam, Kalyana Sundaram, are they not Sudras? Are they not traditional, Hindu Sudras? Comrades Ambedkar, Sivashanmugam, Munusamy Pillai, Goormaiya, Sivaraj are all Paraiyars, Chakkilis, Panchamas. Are they not people of the last caste? If these are unchanged, why the need for any Swaraj, or carnage, or confusion, or loot? For what? To make the Brahmin into the Brahman, and the Brahman into a *Bhudeva* (Lord of the earth)? Don't we need intelligence? Don't we need shame? Are we really low people? If we listen to the Brahmins and dance accordingly, poke the eyes, cut the legs, set ablaze and destroy our own people, do we need carnage, nationalism and Swaraj?

Comrades! After we get Swaraj, if Comrade Avinashilingam Chettiar says that we (Dravidians/ Sudras) do not have the status or the talent, after we get complete Swaraj what will the Brahmins not say? Today, the Brahmin is gleefully butchering the Muslims. Tomorrow they are going to call the Dravidians butchering the Dravidians as "October Sacrifice." I will say that this will certainly happen.

[Viduthalai, 23-11-1946]

**Annihilation of Caste**
When the people got the emotion of political freedom, the necessity, of feeling that the British government was totally not needed, was created. Likewise, because the social feeling was created social self-respect was created and people hated the feeling of low caste, and in anger they said that they didn't need this religion. This emotion comes to anyone who has the feeling of self-respect. There is no use if anyone feels sorry for it.

If each and every person is a Panchama, he will understand how important it is for him to become a Christian or a Muslim. When someone is getting thrashed for being a Paraiyan or a Chakkili, he scavenges and remains an untouchable, and someone else comes and lectures to him not to become a Muslim or a Christian; this is devoid of status or integrity.

Paraiyan, Chakkiliyan, Pallan, Pulayan, Cheruman, Kuravan, Naavidhan, they are named and called humiliatingly. Only today are they being called Untouchables? Right from the time of the Ramayana and Mahabharata, right from the time of the Vedas, Sastras and Puranas, several thousand years have passed. Who will give the compensation/ retribution for that? Which government will give the compensation for that? Which party will do that? Don't we know the information about the great leaders of yesteryears like Vijayaraghavachari, Pundit Malaviya? They only campaigned for the Varnashrama Dharma. But the Brahmins are doing propaganda that they are nationalists. Our people who are crazy are also singing along for the sake of money, or due to foolishness.

There is no use in getting angry with me. I will bet that in today's situation, there is no other way to eradicate the state of low castes than joining Islam or Christianity.

Gandhi says that untouchability is not there, but not a single day does he say that upper-caste lower-caste is not there. But instead he says, "I will establish the Ram Rajya which is based on Varnashrama Dharma. I live only to establish it, for that reason alone we need Swaraj." Can any Congressman deny this?

Is there any policy of annihilation of caste in the Congress? Untouchability is different; caste is different. Without understanding this there are several Sudra-Panchama comrades in the Congress who are cheated and get angry with us. I pity them. As long as Hinduism is there, the differences of Brahmin, Sudra and Panchama will not go. Caste difference will never go until *Gita*, *Ramayana*, *Manu dharma* and so on are trampled. That is why Gandhi and the Congress cautiously do not talk about caste.

He who does not like to become a Muslim or Christian and likes to be a Sudra or Panchama can be silent just like that! I am not worried about it. But why should he prevent others? Why must he blame others? This is my question.

As far as Swaraj is concerned, the Congressmen must keep working till they chase out the British. Why must they save the Hindu religion, protect the title of Sudra and Panchama, and prevent people from going to other religions? From this, is it not known that the Congress party is a religious establishment; it is an establishment that protects Hindu religion, caste and varnas? What mistake is there in Jinnah and Dr.Ambedkar calling the Congress a Hindu establishment?

When Dr.Ambedkar said, "I am not a Hindu; the people who accept me as their leader are not Hindus. I have the right to embrace any religion"—why should the Congressmen oppose it? That too, why should a Congressman who accepts himself as a Sudra and who considers Ambedkar and his society as Paraiyars, Chakkilis, Panchamas oppose it?

The Dravidar Kazhagam or the Scheduled Caste Federation has taken the annihilation of caste as their plan of Swaraj, so they need to talk and struggle for its sake. What else? If the Dravidian Congress comrades have concern, let them release a statement.

That is, let comrades Avinashilingam, Bhaktavatsalam and Kamaraj release a statement that "In the Swaraj obtained by the Congress, there will be no Sudra, no Panchama, there will be no religious basis for this. There will be only one God and one caste." Those with shame and rage can go and demand this in writing from Gandhi, Rajagopalachari and Prakasam. Instead of that why does the Congress get angry and assault us—the Dravidian and Scheduled castes who say that we are not Sudras, we are not Panchamas.

No one needs to worry about the Scheduled castes embracing Islam or Christianity. Other than that, they have no other option. No one has shown another way. Their organization is there to eradicate their Panchama-ness and not for (contesting seats to the) Legislative Assembly or (acquiring) the post of Minister. We have to realize that only a few will get it.

If a Hindu-Muslim riot takes place in this nation, this will be the next work. That is the people will have to convert their religion. The riot takes place only to hasten that process. Because of the Hindu-Muslim riots, the Scheduled Caste Minister, MLA, Buddhist, Jain, Muslim, Christian everyone will become one religion. This is confirmed. Except the Scheduled Caste Ministers, Secretary, MLA everyone else will join Islam. This is definite.

*Note: After Periyar finished speaking, some 100 Scheduled Caste comrades suddenly came towards the stage beating leather drums of various kinds. It caused an agitation in the meeting. Periyar called for a calm. After that, in front of Periyar they made a heap of those instruments. One of them got on the dais and said,"We are called Paraiyars only because we beat the Parai. From now, we are not going to beat the Parai. Certainly we are not going to beat the Parai for the house of death, or for birth or for God. We have now taken a resolve. So, we are burning that instrument." Having said so, they poured four or five liters of kerosene on it and set it on fire. Everyone applauded seeing it flaming away. They shouted, "We will convert after we get the permission of Dr.Ambedkar" as they departed.*
[Viduthalai, 23-11-1946]

**Representation for Depressed Classes**
Dr.Ambedkar got the permission that Adi-Dravidars must be given 12.5% positions and employment. On receiving this order, the Chennai Government said, "we are not getting qualified people" and they put that order in the dustbin. If the Chennai ministers ask for the proportion, they say that it is against nationalism.

If they see people who beat them, immediately they will give everything. If they see slaves, Sudras, Panchamas, they will trample them. This is Indian nationalism. This is also the Swaraj of the Presidencies.

[Kudiarasu, 30-11-1946]

**Conversion to Islam**

Let anyone suggest some other remedy for the atrocities of these Brahmins.

In the Dravidian nation, a Muslim can lead a life of his own comforts and not according to the wishes of a mullah. I know the rules followed by 10 Islamic nations, the rules only cater to the majority. It is not like the Hindu religion, where $3\frac{1}{8}$% of the Brahmins say to the non-Brahmins, who form 93% of the population "to stand at a distance," "you should not study," "you are Sudra," "you are a son of a prostitute," "you cannot have equality," "you cannot have a job," " you cannot have communal representation."

Islam preaches brotherhood, peace and kindness among one another. It is a path for unity. This cannot be refuted. I do not support or advocate Islam. What I say is the perfect and only truth. I do not have love or affection or faith or any other form of relationship with Islam. But in Hinduism the cruelties of the Brahmins are so atrocious, cruel, and unbearable like the poison of an ugly cobra. Islam is the only medicine to save Hindus from the poison of the cobra. To correctly know Islam, please tour through Turkey, Egypt and England.

"Islam is the way to remove disgrace" is the only mantra that can make you lead life like a human. I did not advocate this today or yesterday. I have been saying this for the past 28 years i.e., from 1919 or 1920 till date (1947). I have not lost popularity. The main reason for this is the fact that the only medicine to remove social disgrace is Islam.

All medicines that have been used to cure this disease of disgrace have only faced a failure.

If the Scheduled Castes are given communal representation in employment and education then these Brahmins can sit luxuriously in the rest of the places. What will be the plight of the Sudras? What is their future? What is the representation of Sudras in politics?

There are representatives of the Scheduled castes and Brahmins. There are representatives for the Sankaracharyas. Where is the representation for the agitated Sudras? Think with self-respect! Not only post and positions, but also think about the ideals! We should live without disgrace. The only medicine for this is Islam.

Go home, think and discuss with all great people and come and inform about your decision. Thank you for listening to me without throwing stones and making problems.

100 salaams.
[Ina Izhivu Ozhiya Islam Nanmarunthu! 1947]

**I am not a Hindu**
Why should there be stumbling blocks of "credentials and merit" for the Dravidian alone? Credentials and merit is not needed for the Muslims, it is not needed for the British, it is not needed for the Anglo-Indians, and it is not needed for the Adi-Dravidas who threaten to leave to another religion. But the Dravidian needs to be qualified and meritorious in order to even get enrolled in school—if the Swarajist Brahmin Prime Minister uses the Dravidian Education Minister in his Government to bring about such a policy, you can notice how forcefully the Laws of Manu and treachery are in practice.

Though the Adi-Dravidas and the Muslims have been very badly treated, they have won a prize. Though it is not the first prize it is a good prize. Their position has improved. Now they will get the correct communal representation, they will certainly get even more. Now, facilities to hear their recommendations and grievances have been formed.

I will say that unless we become either the Panchamas (that is the Adi Dravidas) or Muslims we shall not have a life of respect, wealth, or progress.

Comrades! Yesterday in this same place, in the midst of about 20,000 people I stressed this in front of our leader Vedachalam who presided over the meeting. I said the same thing yesterday morning, when I presided over the Dr. Ambedkar Students Hostel celebrations that took place at the municipality office. Some people said that because of this, I am earning a bad name among the Dravidians. Even today morning some people said this. A few people remain with me and say this and carry on a poisonous campaign. I am not worried about this.

I am not even a little worried that I should earn a good name from such Dravidians who do not have any place for emotions and do not bother about respect. My life is not based on such a fake reputation. Serve the Dravidians, in that I do not have any botheration other than the anxiety and expectation that the disgrace of the Dravidians must be removed. I am not worried about my life; I do not have the obstinacy to achieve something else. Today the disgrace the Dravidian faces is that he has to be a Hindu who has to be a Sudra.

Dr. Ambedkar, the sudden fortune that has come to the Adi Dravidars said, "I am not a Hindu, I am not a Panchama, I am not connected with any of the discrimination of this religion." The temples were opened, the minister asks, "Give the list, I shall give employment." Sardar Patel says, "I shall do good to you more than the law permits. What do you need? Ask", Gandhi says, "I am a Adi Dravida, a Bhangi"—all this happens because of the powerful mantra "I am not a Hindu." I said the same thing five years earlier to Dr.Ambedkar in 1925. He, who said it five years after me has succeeded. Yet, they are going to say, "I am not a Hindu" –at least orally—and they are going to get all the rights.

The benefits that the Muslims have received and will receive are due to this. That is, because of saying "We are not Hindus, we are not part of Indian society." I suggested Dravidar Nadu in 1937, Jinnah Sahib suggested Pakistan in 1940. Today Pakistan lies under his feet. Jinnah hits back by saying, "Scrub nicely, wash and clean and come." The Dravidian people live disgracefully—as Sudras—as the $4^{th}$, $5^{th}$ caste—as the children of prostitutes, servants and slaves; they are denied education, wealth, dignity, wages for work—they are treated worse than animals—they have lost their rights, their humaneness, they are forced into such a humiliating life—apart from being a Hindu, what else is the reason for this?

Though Dr.Ambedkar has won now, if he says that he is a Hindu after this enthusiasm is over, he shall again be included in the list of Panchamas.

[Ina Izhivu Ozhiya Islam Nanmarunthu, 1947]

**Caste Identity**

The Union Government decided and passed a resolution never to mention subcastes in voters lists, religious institutions, school records, court records but only write as 'Hindu', 'Muslim', 'Sikhs' or 'Christians'. When I heard about this I whole-heartedly welcomed this resolution. This was done mainly to protect Hinduism, since the Untouchables were converting themselves to Muslims and Sikhs. They broke their heads only to save, and to make this nation remain, a Hindu nation. The Congress never thought of annihilating caste.

They were least bothered about caste and the differences it made in society. They did not even dream of annihilating caste. What is the difference if a man is enslaved to another due to wealth or due to caste? They never think about these problems. To abolish caste, the government should make laws to punish those who put their caste names behind their name.

Everyone who gets married should perform only inter-caste marriage. Laws and rules must punish those who marry within their own caste.

Wearing the holy thread, and putting religious marks on the forehead must be legally abolished and those who do it should be punished. Thus, if caste is annihilated, only the Brahmins will be the worst affected. They entered our country seeking pasture for their cattle; they have become the highest caste and the sons of the soil have become low caste people. What is the reason for all these? The author of "*Politics of British*," Prof. T.K. says, "A revolutionary will never become a Pope and Pope will never become a revolutionary."

Thus a Brahmin can never become a revolutionist. Just like the Pope who is all-powerful the Brahmins are much more powerful. They will state, "Can a donkey become a horse? Or, do all the 5 fingers in any hand look alike? Likewise there are different castes." Till a last Brahmin lives in this nation, they will practice the divide and rule policy. As Buddha and Gurunanak have said the Vedas and Sastras are complete falsehood. We, the Dravidar Kazhagam people, alone have the courage to say so.
[Viduthalai, 10-1-1947]

**Education**
Today's education is fully based on Varnashrama Dharma: one man should labour and one should live lazily. So one group should be educated and the other group should be fools with education denied; so the very education system is like Varnashrama Dharma.
[Viduthalai, 25-1-1947]

**Hinduism and Caste**
Only because of Hindu religion there are several caste clashes. Otherwise there is no reason for religious fights. So, taking everything into account it is better that the non-Brahmins set themselves free from the prison of Hindu religion.

The gods and temples exist mainly to protect caste and the concept of high and low based on birth.

We are not allowed to enter the temples or pray to gods. So why should we tolerate such discrimination and disgrace when it speaks so much on our honour? In this problem, I will not discuss the question of existence of god, for it is an individual's wish. But temple entry and disgrace associated with it affects our honour.

One who goes to the temple is considered shameless because only after accepting that he is a low caste or an untouchable, he enters the temple or prays from outside.

If it is a disgrace or dishonour one should cease to live; but on the contrary, man accepts all disgraces and dishonours and accepts the religion and god. How can it be sensible? Devotees should think about it. The Brahmins do not bother about anything except the label that they are the highest caste.

[Kudiarasu, 22-2-1947]

**Depressed Classes**
In your welcome address, you mentioned about the struggle in Vaikom and about the work you have done for the welfare of the Depressed Classes. Though it is a historical event, I or my comrades or my wife Nagammai are alone not the reason. I cannot forget that several women took part in the Vaikom struggle. We did not get that victory so easily. The Vaikom struggle was held to get the rights of walking in the streets, drawing water from common wells and for temple entry for our Depressed Classes comrades. The Travancore Queen of that day used several oppressive measures but finally she had to accept everything except temple entry.

After some days we struggled for that also. Even at that time the queen did not accept. Finally when 400, 500 people became Muslims she got to her senses and made an order for temple entry. Its echo is what is called temple entry today. The trend of today's temple entry is surprising. Still, the philosophy of our agitation is certainly winning.

So, in such a situation, I have heard that the Depressed Classes do not have the right to enter the Brahmin streets in your village, Thirumangalam. In this year 1947, when the government has made laws to abolish untouchability, if the Aryans of this village are so firm and inhumane, to what can we compare their caste arrogance?

The agitations by the Self-Respect and Dravidian Movement have completely eradicated the barriers of caste differences in hotels. It must be the duty of the people who are bothered about the welfare of the nation to eradicate the cruelties practiced in such villages on the Depressed Classes by (those who are called) Aryans and upper caste.

If such atrocities take place anywhere, please bring it to the attention of the Dravidian Movement. We tell all this in a conference and bring it to the politician's notice and make them pass orders.

If only the indifferent attitude of the politicians is going to be the answer, we will attain human rights by transgressing even their 144 ban.

I say in brief, whoever enslaves the Depressed Classes, the Dravidian Movement will fight against it. Panchayat members in every village must take concern and give us support in this task.

[Viduthalai, 5-7-1947]

**Depressed Classes and the Dravidian Movement**

The Depressed Class comrade S.S.Mani, who gave funds to me on behalf of the Dravidar Kazhagam of this village, said that I have worked more for the Depressed Classes. He also said that all the Depressed Classes members of this village have resolved to work for the Dravidian Movement.

I, or the Dravidian Movement, have not worked separately for the Depressed Classes. If you want, say that the schemes of the Dravidian Movement try to secure the rights of those who are the Depressed Classes. I cannot lie to you that I am working separately for the Depressed Classes. I do not have such an idea. There is no necessity for it to be there. Why do I say so? You need not consider it in another way. We cannot accept the difference of Adi

Dravidar and Dravida. Our scheme is that all are Dravidians.

The Dravidians are not without upper caste people who have caste arrogance. From now onwards, I know how to rectify them. So, it is my intense feeling that we, who are of one race, should not have divisions among us even at the level of words.

I will explain in detail. The important doctrine of the Dravidian movement is that there should be no Paraiyan, Brahmin, upper caste, lower caste, and that the existence of Sudra and Panchama (identities) must be totally destroyed. All are of one race, all are people of the same society: this doctrine must be implemented. So there is no need to unnecessarily boast that the Dravidar Kazhagam or I did this for the castes among us.

The Aryans and their Congress establishment can say so, cheat and feel proud about it. Because they belong to another race, they are of another country. They enslave us in the name of God, they ensure that we do not have a living in our country, because they live here like termites. As if it is doing charity, that group says, "we did this and that for the Adi Dravidars," and they try to hide the conspiracy of domination. Why should our Dravidian Movement also get down to such a low level? So comrades we do not need differences of Adi Dravidar and Dravidar within us. Make them attain the sentiment. Do not say Paraiyan, do not say Sudra.

Also, my dear friend and intellectual Dr.Ambedkar is the leader of the Scheduled Caste Federation. With great evidence he said that the Depressed Classes are not Hindus. He is one of the important people who reduced the Aryan arrogance. Great intellectual. I had great expectations of his cooperation to destroy the cruelties of this Aryanism. That Aryanism whose main root he shook—the same Aryanism establishment—the Congress establishment—of the Hindus has today controlled Ambedkar. Even he has a relationship with the Congress. If it has to be mentioned, my mind is saddened on seeing him speak the dry north-Indian concept that India should not be partitioned. I fear that in some more days he will oppose the partition of Dravidar Nadu.

Still, it looks like in his situation in the north-Indian connection, may be that is how he needs to be. I am not coming forward to say something wrong or complain about it. If those who call themselves the leaders of the Depressed Classes in our country, unnecessarily blame the Dravidian movement—whose only aim is the liberation of the Depressed Classes—its policies, and me; can we blame Dr.Ambedkar who is attached to the north?

I wish to say one thing whether they blame or don't blame the Dravidian Movement or me. Comrades! I give the guarantee that till its last breath, the Dravidar Kazhagam shall work for annihilating the low caste titles of Pallan and Paraiyan and for their development. Till now, I have not said that you should not join the Scheduled Caste Federation. You attain the benefits coming out of that. The Depressed Class comrades have the right to enjoy the benefits of the Dravidar Kazhagam's labour, whether they join it or not.

Finally, one thing is certain. During the struggle for the partition of the Dravidar Nadu people like Dr.Ambedkar might even oppose it, but I have the definite belief that the indigenous people will remain with me and the Dravidian movement and work. The indigenous people and the Dravidian movement (or the Self-Respect movement) are like nail and flesh and nobody can separate them. So, majority of the benefit of Dravidar Nadu is only for the indigenous people. If not today, at least in the future, those who oppose it will understand.
[Viduthalai, 8-7-1947]

**Hinduism and Untouchability**
We often write about the pathetic state of those who are called Untouchables in the "Hindu religion." Today in our country, the people of the Depressed Classes having realized their disgraceful state have started to agitate against being treated as Untouchables, unsociables and being oppressed by the "Indian" people, because they want to rescue themselves.

Recently in Poona, a Satyagraha took place for temple entry. Its vigor has not yet decreased. It has aroused the feelings of people in that region. In north India, in several places like Varanasi, increased signs of rebellion by the Depressed Classes can be observed. In south India, particularly in our Tamil Nadu, within a few years of the spread of Self-Respect Movement, the public have been concerned about the issues of the Depressed Classes. They have also got the courage that can make them escape the disgraceful situation.

Those who are called the Saivas and those who are called the Vaishnavas have realized well that they cannot escape by merely showing the stories of Nandanar and Thirupaanalwar. Moreover, Adi Dravidars and the other Depressed Classes are also calling for conferences of their own communities and working for the development of their society. We also learn that day by day the agitations or emotions related to temple entry are increasing in places like Chidambaram, Erode, Villupuram and Tiruchi in Tamil Nadu.

We also come to know that in the Travancore *Samasthanam* (Kingdom), lakhs of 'Ezhavars' and 'Pulayans' who are ostracized as Untouchables are indulging very vehemently in demonstrations related to temple entry. Though it might be another matter that it is not possible for their agitations to be successful. Those who closely observe the philosophy of their efforts cannot deny that the time to determine the end of Hindu religion is drawing near.

We have often expressed our opinion about temple entry. It is not that if the 'Untouchables' enter the temple they will immediately get salvation or the grace of God. But 'Hinduism' says that our God resides in the temple. We also see people saying that God will die if the Untouchables enter the temple, or get close to God's statue. We will say that this temple entry movement has been started to prove to the Brahmins, the Saivaites, and the Vaishnavaites, that it is better for these Gods of fragile life to die and to be destroyed than for them to exist.

More than the other public places, the temples are given importance. We think that we need to express at this time the truth behind this reason.

We accept that the Depressed Classes attempt to enter the temple not because of the devotion they have over God. Since they have enough Gods and brokers and priests in their own community to enslave them and to prevent their rationalism from being used, they have no necessity for our Gods.

The depressed classes and those who are called Sudras follow different ceremonies of marriage. Those who are 'Sudras' are vehement that they should make Brahmins conduct their marriages. This is very shameful.

In the name of law, Sastras, gods, they are ostracizing us as Paraiyan, Panchama, and Sudra; should those who have intelligence not boycott them? We must realize that because we are Hindus we have all these ills. Only when we become conscious that we are Dravidians these ills run away from us. Muslims do not have the capacity of leniency. Like the Hindu, even the Christian priest is unconcerned about the existence of Christian Paraiyan, Christian Brahmin and Christian Mudaliar. But among the Muslims, there is no Muslim Brahmin, Muslim Paraiyan. Like that we should also be one (uniform) Dravidians.

We must not keep Hindu names. Hindu garbs and religious identities are against us. Has any of the Alwars, Nayanmars, avatars of the Gods, and the Mahatma at least said that the disgraceful titles of Sudra, prostitute-son, Paraiyan, Panchaman must be abolished? Only the Self-respect movement and the Dravidar Kazhagam are doing that work.

If the Brahmin touches shit, he just washes his hand. If he touches a Panchama, he says that a bath is essential.

That madman thinks that if the temple is open, this madman will believe. Such cheating will not take place any longer. Because of temple entry, the Panchama and the Sudra might have become one. But, the Brahmin has not yet let go. We still don't have share in the money he gets through the hundi (temple money/ charity). Unable to bear the atrocities several people of the Depressed classes are becoming Muslims. These Brahmin newspapers are blacking out the news.

We need to realize by ourselves. If anyone asks, "Why this black shirt?" then we need to answer, "I am of a disgraceful caste. I have the stinking job of working with leather. The Brahmin has the job of police superintendent. It seems I am a Paraiyan. The lazy Brahmin who does not work is a higher caste. I feel sad for this. I get angry. In order to remind myself of, and in order to retrieve myself from, this disgraceful position I am wearing this black shirt."

I cannot remain peaceful if someone just says that Periyar should live long. Our disgrace must be completely removed. Only that day I will be peaceful. In order for our disgrace to be eradicated, we need to give our lives for a revolution.

[Viduthalai, 23-7-1947]

**Dravidar Kazhagam**
Among the duties of this Kazhagam, the important one is to remove the disgraceful position of Sudra and Panchama that the Dravidian people have on account of their birth and to give them an equal, high status. Otherwise, the Dravidar Kazhagam is not opposed to other organizations.

The socialists say that we need economic equality. We say that it is not enough and that we also need social equality.

I ask the Dravidian comrades in the Congress to pay good attention. As long as the British were there, you kept saying that first the British should go only then you will look into the internal affairs of this country. Now, the British have gone. You can announce, "Now, the disgrace has gone away." You can announce, "There is no Brahmin, no Sudra, no Panchama." They can tom-tom, "Freedom has come. From now, no difference, no divisions, all are equal." Why have you not done so? Why was it not possible to do so? Think a bit.

Who is responsible for this disgraceful situation: for only 15% of our people are educated, and 100% of the Brahmins are literate? Who is responsible for our being the fourth and fifth caste? It is only 200 years since the British came. It is only 700 years since the Muslims came. But our disgrace has been there for 2000, 3000 years. So, the British and the Muslims are not responsible for our disgraceful situation. That is why even after the British left, the disgrace did not go. Is this not clear?

Then who is responsible? The Aryans came 5000 years ago in search of pasture to graze their cattle; they came to beg like nomads. These Aryans started Hinduism and divided the people into four castes and effortlessly made others work for them. Are not the Aryans, who enslaved people in the name of God and religion, responsible for this?

They made a law that "Jobs must be given on the basis of merit and qualification. Depressed Classes do not have the right of communal representation"—the freedom was used only for this. There are no primary schools for the children of our indigenous people to study. In many places, colleges are being built for the education of the children of the upper castes and the rich.

They get scholarships of one lakh and ten lakhs and they study and they go to work. Our people have to graze cows in the village, they need to draw the *pankha*, or they have to wear pavadai, take care of children and simply put salaams. We need to work under them; this is the educational system of the Congress.

If we kill each other, the Brahmin counts how many of us will pay for the *thivasam* rituals. He will only laugh that the Sudras are killing each other. What else will he do? Even the words, "Oh it is pathetic!" does not exist in his dictionary.

If we observe, is it not because of our propaganda that a few of the Dravidian comrades in the Congress enjoy positions and are in charge of posts? So comrades! The disgrace on account of birth must go away from society. Our party must work for that. The domination by the north-Indians in politics must be destroyed. It is our political scheme the Dravidian must have a separate government.

Till both of these are completely fulfilled, the Dravidar Kazhagam shall continue to work fearlessly. For that we will plead and get the support of anybody and work without any desire for position. So comrades! If you think what I have been saying so long is right, join the Dravidar Kazhagam. As a symbol of that, wear the black shirt. Fly a black flag in your homes. Work for the Dravidian people.
[Viduthalai, 20-8-1947]

**Brahmin domination**
We need to give rice and pulses to be offered to God. But when God eats, we should not see. We build temples. The consecration ceremony is done at our cost. But if we touch it, God will die. To do the purificatory rites, we need to give money to the Brahmins. So, he himself should do the prayers. He will take all the money that comes to his plate. With that his son should study B.A. Is this just?

Not only do we struggle to give, but also we are disgracefully called Sudras, Panchamas and we starve. The Brahmin who eats without toiling and calls himself a high caste is allotting seats for our salvation. Is this justified? His wife will wear an 18-cubit sari and will walk happily like a soldier. Our wives wear a four-cubit sari tightly that leaves them half covered, half naked and they bend and walk. Is this fair? With the money we give, he gambles or drinks or educates his daughter. Without money for education, our son grazes the buffalo. With the milk, curd, butter, we fill the Brahmin's stomach. Is this fair? Have we not come to this low state by not using our own intelligence and having believed everything that the Brahmin said was Hinduism. Only to permanently establish this they have made nationalism and Swaraj.

So, you must awake! There must be no place for the Brahmin in this country. There is no place for the Sudra. There is no place for the Panchama. Like other countries there will be space only for people, only for people who live in accordance with discipline/ order, justice and intelligence. We must create this situation. We are working only for that.

[Viduthalai, 22-8-1947]

**Dravidar Kazhagam**

The Dravidar Kazhagam in not a political party; it is a party related to society. From the beginning, I have said that politics in only a part of the social setup. Others say that only politics is important and that if power is captured social reforms can be made in one line. We will understand what is more important by seeing the plans that the present politicians have made in order to prevent our social development. This is because the public has not understood that the right of social equality is important. Our important work is to create humanism.

If we create unity in society and obtain the support of people, whoever is in power, we can treat them as we like. In this situation if we have power in our hands, it will be sufficient just to tackle the opposition and remain in power. However hard they work, will people of integrity, honesty, and with feelings of society get position in politics? Dr.Ambedkar himself was unable to secure a position despite his feelings for his community. Will it be possible for others? Who benefits in the present independence, except the gang with its separate ideology? The traitors of other parties will benefit.

Comrades, think over this! Was our lowness removed because of the freedom movement? Was the lowness/ disgrace removed because the Dravidians captured power? Because it has been secured firmly with nails, more than gaining power, or betraying our ideology and falling in the enemy's feet to gain power, we need to make the public accessible to us, pressurize those in power and make them work for us—we need to work towards this. We need stronger, consistent propaganda for that. More than the belief I have in power, I have greater belief in pressurizing those in power. Changing the mind of the public is essential and it is possible through propaganda.

Let us take Muslims and Christians. They are considered as *Mlechas* according to the Hindu religion. That is, they are lower than us. It seems they are lower than the Sudras and the Panchamas. What did they do that we didn't do? What do they eat that we don't eat? Are they not Dravidians who converted their religion, being unable to bear the title of Sudras and Panchamas? They were made into Mlechas only because they were Dravidians. So comrades, Muslims and Christians are only of our race. They did not come from Arabia or Gujarat. They are not also low. The Aryans who said to us that they (Muslims and Christians) are low and made us hate them are the *Mlechas*. They have written a dictionary that Aryans are Mlechas (foreigners). We have started work that no one has started before. We have touched a place which any great man, any avatar of God, any Alwar, any Nayanmar, any Acharya, any fighter who went to the gallows was afraid to touch.

Though it was also considered to be essential till now, we started this work of social reform, of destruction of varnashrama dharma, in the name of self-respect, in the name of feeling of community. After several difficulties and worries, only now the seed of self-respect that we planted has started to take root. When it is just a plant, instead of everyone breaking a twig for a toothpick and letting it be spoilt, everyone who passes by it should pour a jug of water, only then it will become a tree and give fruit. Reform will take root.

I pray that comrades who doubt if we love independence must think a little. 25, 30 years ago, I also rolled a little in the slush of independence, the slush of freedom in the belief that if freedom comes, disgraces will easily go away.

I immediately came out when I observed that the methods to establish varnashrama were expertly handled,

and I started the Self Respect Movement, a social reform movement. Children and youth must understand this well. I will repeatedly say that all this speech of independence is only to make way for Aryan imperialism to be rampant.

This is not to say that our work will attain success in our own period. This is what we need to do till we live, then leave the rest of the work to the following generations. Only if we work like that we will not get tired. In order to have children, a few women circumambulate the peepul tree and then touch their lower abdomen (expecting children). That is our state too.

What great change have we created among the people? How many people have we made to think, "I am not a Hindu, I am not a Sudra"? So, I request our comrades that we need to work without the desire for positions, without getting tired, with honesty, commitment and firmness of mind for the cause we have taken up.

[Viduthalai, 05-09-1947]

**Manusmriti and Caste System**

There are several evidences to prove that Dravidians are called Sudras in the name of Hinduism. In the 10th Chapter of *Manu Smriti*, under the title *'Sakar Jati'* it is said that those who do not follow the dharma of caste are known as Dravidians. It is mentioned that the children born to a Sudra man and a Brahmin women are Chandalas. A child born to a Brahmin and a Sudra is the degraded fisherman caste.

Further it is mentioned therein that those born through not following the caste dharma in the Aryan nation are the Chakkiliyars (who work with the skin and hide of animals) and Paraiyars who wear the clothes of the dead and eat the remaining food. Further, in the 44th sloka in the 10th chapter it is mentioned that the rulers of the Dravidian nation are Sudras. The language they speak is called a low language and they are termed as low caste people.

Just like Dravidians, the term Andhra is also mentioned. An Andhra is one who hunts animals in the forest and sells them (10th canto, 48th sloka). So, the Dravidians and Andhras are low caste people who are Untouchables. This is clearly mentioned in the Laws of Manu.

Further these people should live unseen, away from the city, under the trees or near the graveyards and tend dogs and donkeys and should not own cattle. This is mentioned in the 50th Sloka.

These people should dress themselves only in the clothes worn by the corpses, eat in broken vessels, should not use any form of metal vessels, should wear jewels only made of iron and brass. For their livelihood they should always roam. While auspicious functions are taking place they should not be seen and one should not talk with them. They should be forced to marry within their castes. They should not be given anything directly. The remaining food must be put in a broken vessel. They should not come into the village in the nights. Even if they are dressed well they should be considered only as degraded people.

Gandhi said in his speech in Tirupur in Tamil Nadu, what is the way for liberation for a Sudra female child? She must marry a Brahmin man, and her daughter must marry a Brahmin likewise, if the same procedure is carried out for seven generations then in the seventh generation he/she becomes a person of Brahmin community.

A Sudra by doing the work of a Brahmin does not become a Brahmin; likewise a Brahmin by performing any menial work of the Sudras never becomes a Sudra. A Brahmin never becomes a Sudra and a Sudra never becomes a Brahmin (Canto 10 Sloka 73)

A Brahmin should never do the work of ploughing. He can do that job by using some others. Because he thinks that by ploughing, the iron end inflicts the mother earth. (Canto 10 Sloka 84)

If a low caste man does the work of a higher caste man, then the king should chase him out of the nation after taking all his belongings (Canto 10 Sloka 96).

The Sudra has no right to perform rituals. (Canto 10 Sloka 126).

A Sudra should not earn more even if he is talented, because his property would cause pain to a Brahmin. (Canto 10 Sloka 129).

A Brahmin has the right to plunder the properties earned by a Sudra (Canto 11 Sloka 13).

The term Asura means only the Sudras, this has evidence in the 20th Sloka of 11th canto.

It is mentioned that it is proper and virtuous to usurp the property of persons who do not perform yagnas. If the Laws of Manu becomes the Hindu dharma and Dravidians become Hindus what will be the plight of the Dravidians?

[Kudiarasu, 20-9-1947]

**Depressed Classes**

No one will consider a society to be a developed society if there are Depressed Classes in it. If the majority in a society are Depressed Classes, it is really a great harm to that society.

If I or the members of the Dravidar Kazhagam are being criticized today by Aryanism it is not because we are blaming God. It is not because we question, "Why does God need temples, prayers and prostitutes?" But because we ask, "Why are we Sudras? Why are our mothers Sudra women? Why are our comrades Chandalas? Why are you alone Brahmins?"

It is because we question the varnashrama dharma. What can we say if you also criticize us for being atheists without understanding this? We can only say that the Brahmins have created real Sudras among us to such an extent and made us into a shameless society.

How much ever belief we have in God and religion, if someone has shame, will he dedicate his daughter, or his sister to the God and leave her for the service of the village? Has any Brahmin let his daughter to be become a prostitute and dance the *Sadir* in front of God? Has any Brahmin made you sit on the palanquin and carried you? Has any Brahmin used you as a priest and performed a ritual? When he has disgraced you so much, when he still continues to disgrace you, if you have any shame will you make your intelligence subservient to him and see him take you for a ride? This is what we ask.

Go to the temple if you like. We are not asking you not to go there. In the height of ecstatic bhakti, dance in joy. We are not preventing you from that. But why do you bow your heads to that dirty Brahmin? Why do you become a slave to his language? This is what we ask you. Is this atheism?

If the sages and mahatmas of those days need to be rid, then, are the sages and mahatmas of today any better? Would they have said one word that this varnashrama dharma must be destroyed? This Rajagopalachari? Leave them; would this Munusamy Pillai or Sivashanmugam Pillai ask why there must be Panchamas? If they had spoken one word against the Varnashrama Dharma they could not have reached to the position where they are today. What about today? Even if they say one word in opposition, the same day the Aryan gang will start blaming them! What happened to those who opposed them?

I will say what happened to those who supported them, listen.

Do you know of Kabilar who opposed Aryanism? How many of you know Thiruvalluvar who said that all beings are the same in birth? How many of you have the read Kabilar's *ahaval* or Thiruvalluvar's *kural*? Do you know of Ovvai who said that there are only two castes (women and men) and thus opposed the four-fold caste system?

On the contrary, even the cowherd in your home knows about Raman or Krishnan who supported Aryanism and highly respected the beggar Brahmins. Your prayer book is the Gita of the thief Krishna, the Gita of the immoral Krishna, the Gita of the Krishna who said that "I only created the Sudra, there are four castes." He is your god of the Kali Yuga. According to you, not only were the graceful Kabilan and the godly bard Thiruvalluvan and the artist Ovvai, born in the womb of a Paraiyan woman; you don't even know their songs.

It seems Thiruvalluvar, Kabilar, Ovvai, were born to a Sudra woman called Aadhi and a Brahmin called Bhagwan. Tholkappiyan whom you praise was born to the Aryan. When they are working so hard to see that all the greatness comes only to the people of their community, and even on seeing that if you continue to be their tails, how can others consider you to have some shame? Comrades if you have shame, please cooperate with us in order to gain human status. Otherwise, at least stay apart. Don't become the henchmen of Aryanism: Vibhishanas and Hanumans.

The uncivilized barbaric Negroes don't have the difference of Brahmin and Sudras. The Eskimos, who eat people raw and live amidst the dark snow, do not have these differences of high and low.

Dear friend who calls himself the son of Bharat Mata! The warrior who hails Bharat Mata! How can children of four-five castes have been born to your Bharat Mata? If it is so, do the children of Bharat Mata have four fathers? How can it be of one caste? Will you accept this? Why all this cheating: that one mother must have children of four castes? Why do you shamelessly try to boss over us?

What labels you Sudras? Only Hindu religion has made you into Sudras. Only because you call yourself a Hindu, you need to comply with the varnashrama dharma. You are made a Sudra only because you accept the gods of Hinduism. Brother, only because you go to worship these gods, you need to stand outside the prayer room, bow and keep beating your cheeks.

You think of the foolish, sly story! It seems Brahma created you! It seems the Brahmin was created from the face of the Brahma. It seems the Kshatriya came out of the shoulders of Brahma and the Vaishya was born out of Brahma's hips. Only because you were born out of his legs you are the Sudra. When Brahma has a wife, why should the husband do this work? Is it possible? Is Brahma male or female?

Even if that is ignored, why have they created it in such a way that one can give birth in several places? It seems that Krishna says in the Gita that it was he who created the four castes. If all the four were born of Brahma, then why does the Brahmin alone be called the son of Brahma, the race of Brahman? Because they were born in the lower organs, are they lower-caste? The jackfruit tree bears fruit even in its roots. It also bears fruit in the top. Is there difference in quality and taste? Is there no limit to foolishness and atrociousness? We who build temples and give away wealth, why should we be made into those born from the feet of god? The Brahmin who has never known labour, has never given even one coin for the god, cheats the God and also fattens himself with what we give this God; why should he be made the one who came out of the face of God? Can a god who did this or allowed this to be done, be our god? Can a god who says that his life will go if he touch him, be our God?

Today my friend of the Vannar caste, and friend of the Chakkili caste can send a telegram and talk on the phone. They can also today board an airplane and fly over your heads and the heads of your gods and the cupolas of your temples at the speed of thousand miles per hour. When they have been given the right to enjoy all this, why should the dirt-laden Brahmin alone have the right to touch this immobile stone God?

Only our people were the origin for all these specialties—for bravery, courage, energy, intelligence, manliness, affection and beauty. At what point of time were we, their descendants, made into Sudras? Who made us Sudras? Aryans, or the Gods?

Only after the Aryans migrated to this land could they have bullied us and called us children of prostitutes. Only after that we became Sudras. If we were Sudras before that, would our Gods have called us children of prostitutes? If it had been so, would we have accepted it as such? Three thousand years ago, from somewhere in middle Asia, the Aryan tribals were a begging crowd that came here in search of pasture to graze their cattle. For having hated their sacrifices and soma *bhan* (ambrosial liquor) they called our ancestors as children of prostitutes and Sudras.

In this age of rationalism and science, because we accept it without shame and remain their slaves, they are still ruling over us. You also remain dumb and not question them. Till now who thought about this disgrace? Because I have recruited several people and I do the work of removing this disgrace, some people have realized this disgrace and have joined me.

Seeing this, the Aryan gang feels that its influence is coming down and it begins to debase me. My friend, should you join with them and start shouting? For him (the Aryan), there is a loss of income, his superiority goes away and in that sadness he shrivels his mouth and blabbers. If we say his faults, why should you, a Dravidian, get angry? Why do you worry? If you are worried, go and tell him what is fair and make him mend his ways. Instead of doing that, why do you take my life?

Do I struggle for myself? I work for you and your family and your descendants too, that the Sudra title should go away. If I blame the Hindu religion that speaks of caste difference and Sastras and puranas that serve as the basis and the Gods who accept these, why do you get angry, brother (the Dravidian comrade of the Congress)? If I feel ashamed because they call me a Sudra, I list out the faults in them. If you do not have shame, if they say you are the son of a prostitute, and if you accept it, you can stand apart in a corner.

If you say that one should not talk about religion, does it mean that god is the private property of your home? The god who is omnipresent, will he not be present for me? Can he allow a group to dishonestly exploit him for the sake of their livelihood? Why do you jump at me if I speak against those who loot us in the name of God? Why should we accept a god who says that we are Sudras? By saying that we are all Hindus, why must the Brahmin make me alone donate rice and pulses? Brother, you answer that first.

I will say that it is barbaric to believe in such a God. Hinduism is a great deception. The Vedas and Sastras of Hinduism is all number one cheating. The puranas is a slush in which a jealous work is going to creep. Salvation and hell are all conjurors' tricks. Sin and virtue are for those without integrity. All this cannot approach a self-respecting man of integrity. Even if they approach, they will be exposed and it cannot destroy him. Can you contradict this? Who first condemned the devadasi system? Only the Self-Respecters. Brother, do you know how much of opposition was created in those days? Instead of being afraid of all that opposition, we fearlessly agitated. That is why this group of believers has itself passed a law that the devadasi system must be abolished. Like this how much we started—today, how they are being fulfilled? Think over this. The Hindu religion and the Varnashrama Dharma are responsible for your being Sudras. Only the Self-respecters touched the basis of it, and said that the religion of Varnashrama itself must be abolished. Your Alwars, Nayanmars, Sankara or Ramanuja or the heads of the mutts or the pundits did not say this. All of them said that a reformed Hindu religion was needed.

If the Brahmins in our country are not sectarian, why do they have a separate Hindu Mahasabha, and apart from that a Brahmin Protection Association and a Hindi Prachar (Sanskrit Prachar) Sabha?

[Viduthalai, 23-12-1947]

**Swaraj and caste**
Does Swaraj imply that a Brahmin should rule a Paraiyar? If a cat rules a rat will it be Swaraj? If a Zamindar rules a farmer can it be Swaraj? If a proprietor rules a labourer can this be called Swaraj?

"I am a Sudra (a prostitute's son) according to the Aryan's dictionary. My friend is an untouchable (he is a Panchama and Chandala). My Muslim friend is *Mlecha*, another a Brahmin." Keeping all these people in mind we talk about democracy without any shame or honour. If you have some honour or intellect, first declare that all are equal in this nation and there is no place for distinction according to caste. Give the voting right to everybody and then declare democracy.

"Swaraj" which stealthily came from behind the screens at midnight by seeking an auspicious day has not entered the scabbard like a sword till today. The robbery, murder, rape that started that day has not stopped. Only enmity is growing rapidly and there are no institutions that reduce these problems. Why should a Hindu and a Muslim exist? Why should an owner and a labourer exist? Did anyone question why a Brahmin and a Paraiyar, a superior and an inferior should exist? If they do not have a difference of feeling will there ever be any quarrel and dispute?

Freedom fighters! Have you ever made the public realize about humanity? Have you changed your policy as per the interest of the people depending on the times and existing circumstances and on the people's natural feelings of usual worldly desires? As long as such differences exist among people, humanity will continue to seethe. Realize this and act accordingly. To achieve permanent peace, begin by removing all the differences and do the necessary.

[Viduthalai, 19-1-1948]

**Labourers**
One should realize that according to the Hindu religion, all Dravidians are Sudras only because of their occupation. These labourers are considered as low, degraded people only due to Hinduism and Varnashrama Dharma. It is plain ignorance to think that one is a low caste person because he does menial jobs.

In no other nation is a particular job assigned to a group of people, and another type of job to a different group of people. So in other countries there is no stigma associated with occupation or labour. So no occupation has fame or shame in other countries.

The profession of making shoes and slippers is assigned to one sect of people and in particular it is assigned to the Panchamas i.e., Untouchables or outcastes. Only in our country, this profession is considered as a low one but by no means is that profession a low profession. On the contrary the profession of *archaka* (temple priest) that is carried out mainly by Brahmins is considered by the Hindu religion as a position of fame and it is assigned only to a high caste Brahmin.

Thus, one cannot be a high or low caste because of occupation. This realization should come to the labourers. Only labourers should rebel, revolt and fight against the situation. Why should they physically labour in order to live, and based on that be made as low caste? These days labourers only fight for better pay not realizing that better pay (or increase in pay) leads to increase in the cost of commodities.

Further these labourers are mad about religion, God, past sins, blessings and destiny. So, only the Brahmin leads a life of laziness without toil and rolls over with good food and a lot of money. On the contrary, the poor labourer believes in destiny, labours the whole day saying it is "god's will" and is not able to eat even a square meal a day. The Brahmin lives in arrogance saying that god has blessed him and lives on the labour of others. The labourers toil without cunningness and lead a life of poverty. Thus these selfish and cunning Brahmins have made (religious) laws to lead a life of comfort at the expense of other's labour. Do these have any rhyme or reason? That is why I question, "Can gods will or law, exist for Brahmins leading a lazy life and non-Brahmins toiling without any comfort?" So god should first ensure that the person who labours the whole day is well paid and the lazy Brahmin can be fed meagerly with the remainders.

Today's nationalist comrades and communist comrades sing the national anthem "to destroy the world that cannot provide food for an individual."

The Brahmin (Subramania) Bharati showed this way. How are they going to destroy the world? Is it by cutting the ocean into the land? Or by putting the sun into their oven? If they are going to destroy the world will not the labourer also get destroyed? What a foolish thing!

[Viduthalai, 20-1-1948]

**Sudras and Brahmins**

The class of people who were Sudras about 2000 years ago continues to remain Sudras even now. Those who were Brahmins 2000 years ago, continue to be Brahmins even now. Our colours, habits, love and hate and our behaviour have changed but our caste that is birth-based alone remains the same. Is this proper? Is it all right? Now, we cannot accept the disgrace. Of what use is a religion that does not give place for our respect! When we assemble like this, our mind and concern should be only about our self-respect. How will disgrace go?

[Viduthalai, 29-2-1948]

**Swaraj**

I was very clear that in the Swaraj, Panchamas and Sudras should be only one sect of people. Only for that I wore the black shirt.

[Viduthalai, 11-3-1948]

**Independence**

We thought that after we get independence we would become free from the stigma of birth. But after independence we have lost the freedom of even talking about the stigma of birth. For saying so I am going to be imprisoned!

What is the meaning of this? So lowness can exist! Sudra can exist! But we should not talk about the troubles created by it. What does this mean? Western countries have made plans! How much are we going to become developed in 50 years? We assess how we were in 1900, 1800, 1700 and 1600? What was our status then? What should we do? We carry out these sorts of plans! We go back to the days of Emperor Ashoka! What is our relation with Ashoka?

When other nations are planning for national development, our nation is planning regressively. We are living in a pathetic condition. They are playing on the ignorance of the common man! Should we not begin to think? In the name of ancient civilization they are dragging us to the olden days. Should not the intelligent people think about this? Should not we struggle and toil to stop this? Should we not think carefully and worry about the cunning and selfish motives by which we are treated and ruled? Be united and toil hard to wipe out the stigma of low birth, which has no scientific reason or rhyme. Live like a single family. Learn unity from the Brahmins. The large crowds that always gather to listen to my talks overwhelm me every time.

[Viduthalai, 11-3-1948]

**Aryans and Dravidians**

The Aryans call only the non-Aryans as Dravidians and say that they are fit only for degraded work and that they are the demons and *rakshasas* (devils) described in the Ramayana. Even the monkeys in Ramayana are nothing but the Sudras. (Jawaharlal) Nehru mentions this in a letter to his daughter. That is, the *lambadis* (nomads) who came here seeking green pastures fought with us, enslaved us and when we agreed to work for them, they termed us Sudras and Panchamas. Those who were against them were termed as Demons, rakshasas and monkeys!

[Kudiarasu, 8-5-1948]

**Brahmin domination**

The domination and the cruelties practiced by the Aryans on the natives are more than the problem faced by the natives under the rule of the British or the Muslims or even in times of war!

These Brahmins who came as beggars into India made the hooliganism grow. Should we also cultivate our people with this hooliganism? We have more responsibilities because it is our nation and 97% are our people. That is why we have to preach and teach to our people about morality, self-control and honesty. As the Brahmins have come to our nation for begging they have no feelings for our people or for our nation. So you should not become angry with them. Only the god who created the Brahmins has created me. Without his permission, my mission is to annihilate the Brahmins and the Sudrahood. I campaign only for this. How can I alone act against the rules of god? You go and ask him. Ask God who has created the Brahmins? Why did he create a Ramasamy to annihilate them?

Nearly 1½ lakh people have gathered here. Will any Brahmin magazine write a line about this? Suppose a small commotion or riot takes place, immediately all the Brahmin magazines will report it saying that there was a big caste clash!

With so much restraint they hide the information about us. Suppose 3¾ Brahmins assemble in the Rasiga Ranjana Sabha immediately all Brahmin magazines will write a page or more. Can anyone question these atrocities?

Doing all these atrocities you (Brahmins) claim, "The hand that holds the grass will hold the skies." What an arrogance all of you have? Is it proper? When will you come out in the open? How long shall we wait for you to come out in the open? The sooner you come, the better for us. Our disrespect will be abolished sooner. We are 97% and you are 3% (of the population). Even if 3 of us are destroyed we have 94% who can lead a life without stigma. If all the 3% of you are destroyed there will be only 0% of you.

If the grass-holders hold the sky, what will the beholder be? If a finger becomes a grinding stone what will the grinding stone become? What will be your plight? Think about it!

Everyday only Brahmins give problems to Dravidians but Dravidians never give any problem to Brahmins—everybody should realize this. Brahmins are plantain leaves and Dravidians are plants with thorn. Whatever the fight is, only the plantain leaves will be destroyed. If you try to clash with us, you will be eradicated from the root. The Brahmins should realize this. Brahmin comrades! Don't think we are fools and you are all intellectuals! Don't be destroyed due to arrogance! You say that you treat everyone equally. Then, why do you keep tufts? How long can you rule us by keeping some of us as your stooges? If they slowly come to our side what will be your plight? The black shirt that we wear is not for war; it is a symbol of disgrace. We feel disgraced due to disrespect, we feel sad. We have decided to wipe it off. That is what this black shirt indicates.
[Viduthalai, 19-5-1948]

### Nationalism
The word 'Desiyam' is a fraud since it is of Sanskrit origin. It has introduced the Varnashrama Dharma in our great traditional, honorable life. The Dravidians were made Sudras and the ancient people of our country were made Panchamas; all these are conspiracies of the Brahmins.
[Viduthalai, 20-7-1948]

### Anti-Hindi Agitation
I feel that my demonstration should be ingrained in gold letters.

Our demonstration is justifiable. People from all parties should join this demonstration because it is for our self-respect. This demonstration was against imposition of Hindi language on us.
[Viduthalai, 20-7-1948]

### Self-Respect Movement
All demonstrations conducted by me and the Self-Respect Movement never faced defeat. We succeeded in entering temples, taking water from common wells and so on.
[Viduthalai, 22-8-1948]

### Dravidians
All people, except the Aryans, are Dravidians. This includes Muslims and Christians also.
[Viduthalai, 5-10-1948]

### Religious Untouchability
If a Brahmin claims himself to be a non-Brahmin, we have to first question, "Where is the holy thread for a non-Brahmin?" and tell him to first cut it. According to Arya dharma even the Muslims and Christians are considered

Untouchables. Break the stone gods that say there are Untouchables and use these stones to lay roads.

[Viduthalai, 5-10-1948]

**Marriage**

Marriage rituals never existed among Dravidians. We see that the marriage rituals for the *dhobi* (washermen) and the barber are different. This itself is evidence that marriage rituals did not exist for Dravidians. Only the Aryans imposed all the marriage rituals on us.

Dear Dravidian comrades! If you want to be moral and have human instinct you have to abandon love for Sastras, puranas, rituals, religion and gods because only they and their laws disrespect you and they were made by the Aryans to enslave you. Should we not leave all these useless activities and start to work for reform? No one forces us to perform marriage in this form. If you want to improve yourself think and reason out with your intellect; if you want to raise yourself, listen to me. Never blindly accept any act passed on by your ancestors.

[Viduthalai, 24-10-1948]

**Congress and Caste**

The leader of the Depressed Classes Dr. Ambedkar said that the Congress is like a burning house and that soon it was going to be destroyed to ashes; hence he advised his friends not to join it. He said this when he was holding position due to the Congress' favour. He said this mainly with concern for people's development. He condemned the backwardness and tyrannical rule of the Congress and praised the Dravidar Kazhagam which alone for acting against the oppression of Congress.

[Viduthalai, 2-11-1948]

**British**

In Britain, a man feels ashamed if he is not educated. In Madras city, only the British run hospitals for leprosy patients, orphanages, hospitals, and colleges and take care of children abandoned at birth. They give education to the Untouchables who were not even permitted to walk in the streets or come near others. The Untouchables are treated equally, and above all, with human dignity.

[Kudiarasu, 28-5-1949]

**Religious discrimination**

No god would be so dishonest as to treat the toiling class as a low caste and the lazy caste that lives on other's labour as an upper caste.

[Kudiarasu, 18-6-1949]

**Aryan domination**

How did the Aryans make us into slaves? Not by winning any war, but by their cunningness. They made myths and epics into art and stories and imposed it on us.

   Slowly, we accepted their philosophy, god and dharma. Once they succeeded in that, they made the Laws of Manu through which we were made a low caste. Once they made us accept the Laws of Manu, they used the same law to make us into a low caste people. When we realized this and protest against it now, they impose their Aryan language on us. By these acts they are trying to ruin us.

[Viduthalai, 16-1-1950]

**Aryan domination**

The Aryans became successful by imposing their culture on us. By destroying our culture they began to rule us. Because we forgot our culture and accepted the Aryan culture now we have been made into low people, sons of prostitutes, Sudras and Panchamas.

[Viduthalai, 27-1-1950]

**Caste discrimination**

In the Gurukulam, the Brahmin students were fed with '*uppuma*' and the non-Brahmin students were fed with old rice. For these Brahmin boys they teach Vedas and for our children they teach Devaram and Thiruvasagam. There were different water pots for Brahmin and non-Brahmin boys. This difference in treatment existed even in those days.

[Viduthalai, 30-3-1950]

**Communal Representation**

Caste that gives distinction or low status should be abolished. Till people get equal opportunities in employment, education and politics, backward class people must be given communal representation.

This nation has 3% Brahmins, 16% Depressed Classes and 72% Non-Brahmins. Posts should be given only in this ratio. But on the contrary, the 3% people say that they are alone qualified and talented, so they usurp all the highly paid posts in every department. How can anyone accept the statements of these 3% Brahmins who occupy more than 97% of such posts?

[Viduthalai, 9-4-1950]

**Untouchability**

Like doctors, the healers tried to abolish untouchability by preaching equality in all walks of life.

But we question: What is untouchability? How did it come?

We see that it has come only due to religion. It is due to Sastras. How did Sastras come? It is the rule of god. So the religion that supports Sastras and the Sastras made by god should be abolished if untouchability is to be abolished.

[Viduthalai, 16-4-1950]

**Sudra**

If one goes to Bombay and calls somebody a Sudra he will not get angry, or if one goes to Punjab and calls somebody a Sudra he will not get angry, but if you go to a Dravidian country and call somebody a Sudra he will feel insulted. Thus only one out of nine people feel bad to be called Sudra. In a republic country, why should the Brahmin, the Sudra, and the Panchama exist?

Should the government not give the police a pair of scissors to cut the holy thread, whoever is wearing it?
[Viduthalai, 12-11-1950]

**Punishment for caste discrimination**
It is written 'Brahmins Club' in bold letters and so on. Whoever practices caste in this nation, whether Brahmin or Sudra or Panchama, should be imprisoned for a year.
[Viduthalai, 12-11-1950]

**British and Brahmins**
The British conspire with the Brahmins. This is the existing relationship between a white Aryan and a yellow Aryan. Only such agreements are made here. All other talks are strategy. I said this two years ago. Jinnah and Ambedkar say it now.

Swaraj is to protect the social setup. For instance, in Thiruvaiyaru when I said there should be no difference in food; Sastri, B.S.Sivasami Iyer and T.R.Venkatarama Sastri questioned the British by saying that they were over riding the agreement. We should break this. Whoever stands as a hindrance to the rule and peace of the nation by having the name of his community must be put behind bars for six months.

Will not the British laugh at the existence of the Brahmin caste and the Paraiyar caste and the fact that one caste cannot enter into hotels?
[Viduthalai, 14-2-1951]

**Untouchability**
An intelligent man has intelligence not in all the fields but only in few fields. He will not use his mind in certain fields, because applying ones mind in such fields is a sin and at times, dangerous! So, only out of fear, he does not apply his mind in certain fields and has never used it in some fields.

That is why we have remained Sudras, Panchamas, degraded people, low people, backward class people, Untouchables and so on. We do not research in these directions. Even if we research in these areas, we are punished. Because of the fear that we would be punished, we do not research in these areas. So, can we say people are brainless? They have brains; they use it in certain fields and do not use it in other fields at times due to their fear.
[Viduthalai, 27-3-1951]

**Hinduism and Caste System**
If the Varnashrama Dharma is to be annihilated, Hinduism must be annihilated. As long as Hinduism exists, one should accept the caste system. Caste can be destroyed only by the destruction of Hinduism. Gods and Sastras that are the basis of Hinduism must be destroyed. Maintaining one and destroying another cannot be done. If anyone says so he is a fool.
[Viduthalai, 30-5-1951]

**Struggle for Rights**
If we have to get our rights, it means that we have to snatch our rights from our enemies. It only means the destruction of their supremacy. Only then we can get our rights. This is not an easy task. It is not only a great task, but also several people have to die and certainly there will be bloodshed. Apart from the support of caste, religion, god, rule and Sastras, the power and support that our enemies have is very great. Our government also would only support them. Our enemies have more unity. We do not have any strength. Only we should have mental courage.
[Viduthalai, 30-5-1951]

**Constitution of India**
They have protected Sudrahood in the decision of making a political party. They formed a committee to write the Constitution of India. In the four-member committee all were from the orthodox Brahmin community. They were T.T. Krishnamachari, Sir Alladi Krishna Iyer, K.M Munshi and Tirumalai Rao. Then they added three more, viz. Dr. Ambedkar, one Muslim and one Christian. The representation from non-Brahmins i.e., Dravidians was nil, though they form nearly three-fourths of the population. On the contrary, four out of the seven members are Brahmins, and they represent less than 3% of the population.

Dr. Ambedkar fought for the Adi Dravidars. They informed him that they would give anything for his people without protest but on the only condition that he should not talk about others (Dravidians, that is non-Brahmins). Accordingly, Dr.Ambedkar got reservation as per representation of the population for his people.

Thus the Adi Dravidars obtained what they wanted. But they say it is wrong for us to ask for reservation as per our percentage in population, that is communal representation. Whoever asks for such reservation is termed sectarian.
[Viduthalai, 22-9-1951]

**Swaraj and Untouchability**
I was a very rich businessman. Even as a businessman I did a lot of social work. I joined the Congress in 1920. Gandhi said that he would not support Swaraj. Even if we demand Swaraj it will be only after the social reformation of abolishing untouchability and caste differences. I was not of the opinion that one should be in politics to serve people.
[Viduthalai, 14-10-1951]

**Caste Annihilation**
In any other nation in the world do people put holy ash on their forehead or any other mark of religion on their forehead? Such types of devotees, like the Nayanmars and Alwars, have not been found in other nations. As the caste system was very strong, all those who tried to annihilate caste were annihilated.
[Viduthalai, 21-2-1952]

**Brahminism**
Brahminism has enslaved us by religion, rituals, god and customs and not by force or guns or arms.
[Viduthalai, 22-2-1952]

**Inter-caste marriage**

Caste system cannot and will not be annihilated by inter-caste marriage. A new caste will be created by this intercaste marriage. The only means to annihilate caste is through laws and the Constitution. It is difficult to annihilate caste while retaining the Sastras. Only these non-Brahmins have caste by which they degrade others. In the case of Brahmins, all the divisions like Rao, Iyer, and Iyengar are only treated as Brahmins. Inter-caste marriages have not changed the differences in caste.

[Viduthalai, 27-2-1952]

**Annihilation of Caste**

Brahmin is a high caste that can directly send telegrams to god. A very dirty Brahmin, a leprosy-infected Brahmin or a dishonest Brahmin is deemed high-caste. He should be respected by other castes.

Likewise if a very intelligent, talented man is born in another caste he is deemed a Paraiyar or Sudra or prostitute's son.

He is expected to serve the Brahmins with folded hands. Such wonderful situation cannot occur in any other nation! 'Sudra' means serve the Brahmins endlessly. Everyone must do physical labour except the Brahmin. The situation will change only by banishing caste, god, religion and Sastras.

[Viduthalai, 2-7-1952]

**Politics**

People should do only one thing: they should hate dishonest traitors, honourless and characterless men who come up as representatives of the nation. Because people do not punish or protest against dishonest traitors and characterless people, they fearlessly gain recognition in society. Politics has gone to its lowest state and no religion or caste is backward in this activity. In society, honest and good politicians should come to power and these dishonest traitors must be driven away.
[Viduthalai, 7-7-1952]

### Aryans and Dravidians
In Dravidian country, especially among the Tamil-speaking people there is a strong demarcation of caste and cultural differences between the Aryans and Dravidians. This would be helpful at sometime or other. So my struggles will not be wasted.
[Viduthalai, 15-7-1952]

### Caste and Labour
In any country, the workers union will stop with the difference of rich and poor but in our country it does not stop with rich and poor but it protects the castes of Brahmin and Sudra and Paraiyar. This difference should be destroyed.
[Viduthalai, 7-10-1952]

### Names and Caste
You would have observed that varnashrama dharma is followed in naming us. All nice and beautiful names are assigned only for the upper castes. The names of gods and religious names are kept for the upper caste and these names should not be used for the Sudras and Panchamas. It is a sin if they keep these names. For instance, names like Raman, Krishnan, Saraswathi, Parvathi are the names of upper caste people and names like Karuppan, Mookan, Veeran, Katteri, Karuppayi are assigned for the low caste people.
[Viduthalai, 24-3-1953]

### Betrayers
Even today, they create Hanumans and Vibhishanas among us and they tell them to hit us. We need to tackle this. In our struggle, the Brahmins will not come like the Devas (celestial beings) who came according to the Puranas. They will only instigate and make the traitors like Vibhishana and Hanuman to fight against us. Our problem is: If we have to oppose and attack, we will only be attacking our own people. Should we fight among ourselves? That is why I control my people. Those who are coming to fight with us today will change at some point of time. For how long can the Brahmins feed them and give them wages and make them fight against us?
[Kudiarasu, 4-8-1953]

**Politics/ Personal**
Comrades! Only if we take efforts to solve these differences in my lifetime, it can be solved. I say this truly. Only if we take some efforts during my lifetime we can achieve success without carnage or bloodshed. If I die, what might happen afterwards? Why? After the person is not there, only the strength of the text remains. So, there shall come a situation where might becomes right. There are no responsible people, no one is concerned and takes effort. There is no one other than those who cause carnages and show their presence, and those who gain profit and publicity through this.
[Viduthalai , 8-9-1953]

**Labourers**
A labourer should not be rationalistic because if he were rationalistic it would affect the company. They deal with the problem only at the pay level and they will not allow a labourer to think why should he not be an owner. Suppose a scavenger's son gets educated, who will do the scavenging job? So they want these practices to be steadily carried out. In India we should demolish the caste through agitations!
[Viduthalai, 3-10-1953]

**Abolition of Brahmins**
In India there were 562 kings. After Independence they put down their crown and got around Rs. 50,000 to Rs. 2,00,000 and left us. A day came when there were no kings in history. We have changed to such situation. Zamindars have run away! Next, rich people are there. A day will come when we would ask why do you need 1000 or 2000 acres of land? He too will be abandoned. When such changes take place why need for a beggar Brahmin!

[Viduthalai, 8-10-1953]

**Educational policy**

As soon as Comrade Rajaji introduced the hereditary educational policy, only I said at first, "This educational policy is a casteist educational policy; this has to be opposed and abolished." The people supported me. Several people of other parties, why, even several people from Congress opposed this scheme.

What is this scheme? What do you think is the basis for this scheme? Is this educational scheme not a reconstruction and protection of varnashrama?

According to this educational scheme, every child should study for sometime every day, and the rest of the time they should practice their hereditary caste occupation. Who physically labours and slaves in the name of caste? Only we, who are called Sudras.

If we ask: "What? Should we keep doing the caste occupation while Brahmins alone get positions, employment, authority and go higher and higher? Must it be so? Is this justified?" Then, they will immediately say just one word, "You are sectarian."

I wish to say, "We don't think that physical labour is disgraceful, but why should we alone do that work?" Everybody can do it according to their percentages (in the population). Instead of that, if it is said that we, the Dravidian people alone must do such physical labour; and that Brahmins need not do any of it, they can lead a high life just remaining seated in one place and without any dirt touching their nails, what justice is there in this?

The British who were 6000 miles away from our country ruled us. Though he controlled us and ruled us, he treated us as human beings. He ate from our people. He did not say that he would not touch us!

Only those who are called Paraiyars by the Brahmins were the butlers and cooks in British homes. So the outsider Brahmin says he should not touch us. Our touch is polluting. He says he is born in the face; so he is the first caste, you are born from the leg so you are the lowest caste. He sees to it that our people do not come up in any field; he controls, rules and enslaves us. To do so he has taken hold of the law, rule, government, religion and god. How can we have this type of foreigners?

The work we do is not for our benefit. It is not for our selfish, self-centered profit or comfort. To wipe out this degraded state of all Dravidians who are Sudras, Paraiyars and Chakkiliyars; if we wipe out Sudrahood, will not all Dravidians become better? If we had this feeling and control, will Rajaji and the Brahmins have come to this level?

How much of opposition to Rajaji's minister post? Though the situation is such that his post can topple any day. If he can manage all that and still have a tight grip, what is it for? Only to save his race, only for the protection of the Brahmin society. As long as he is in power, he wants to fill up Brahmins in all the places. He wants to repress the Sudras, and with that motive alone he remains in power. Among us we do not have this mind, this feeling, this racial fanaticism! He fights for the welfare of the Brahmin race; our people do not come forward to struggle like that!
[Viduthalai, 17/18-11-1953]

**Educational policy**
Rajaji says that people in the villages do not need education. He has said that the student in the village should cut hair, he should wash clothes, and he should make pots and pans. The village school is only for three hours duration, the rest of the time our children should graze donkeys, this is called the New Primary Education Scheme. We organized a conference in Erode only to oppose this educational scheme.

Can we keep looking at law and democracy in such a situation? The work will be successful only through a great revolution. Though 100, 1000 of us will have to be sacrificed, we need to be prepared.

That is why even I gave a three months notice; and said that if they want to shoot, let them shoot. I am going to start a struggle. This educational scheme has to be destroyed by all means.

[Viduthalai, 26-2-1954]

**Educational policy**

What should we do now? We should be ready for the agitation. If a man has shame, he has to go to the jail at least once. That is the test of integrity for a man now.

If our elders backtrack by chance, only the youngsters must certainly do this. Otherwise, the next generation after us will only have to graze cattle. This is the only information I have to say about the hereditary educational system.

Today, do those who call themselves social reformers strive to abolish caste? Will the Brahmin agree to end caste? That is why I say that any one except us does not have this scheme. We say that for the annihilation of caste, the entire basis of caste must also be annihilated. No one will say that the basis for these differences such as religion, Sastras, God must be annihilated. No one will accept this. So, it is a difficult problem. Because that has to be said substantially to the people, we are revealing to people the doctrines of Buddha.

Buddha said this 2000, 3000 years ago. Even the government accepts Buddha's philosophy. They have embossed in the national flag the chakra of Emperor Ashoka, who walked the path of Buddha. Does it not mean that they have accepted his philosophy? As though this is not enough, it is also being said that Buddha is the tenth avatar of Vishnu. If this Ramasamy (Periyar) says something, it is wrong; if the Dravidar Kazhagam says it is wrong. But the Buddha himself said so, and he said it candidly.

[Viduthalai, 7-3-1954]

**Buddhism**
We, the working people, have been made into low caste; we are starving. We don't have clothes to wear; we don't have a place to live. But the Brahmin who doesn't work has all kinds of comforts and honour.

Only you built these temples. Only you gave the money. When it is so, can we just leave the god who says that we are low caste and untouchable? "Ungrateful god, only I built the temples and the tanks for you. Only I spend money on you." Should you not ask why you are lower caste or why you should not be touched? You believe what the Brahmins say, that is, if you ask such questions the god will get angry.

It was resolved in the Buddhism Propaganda Conference at Erode that such a god needs to be abolished. The World Buddha Sangha's President Mallala Sekara presided over that conference. Mr.Rajbhoj, who is the Secretary of the All India Depressed Classes Association and a Member of Parliament, inaugurated the conference. He spoke, "The policies of the Buddha are your policies, the Self-Respect policies."
[Viduthalai, 14-3-1954]

**Annihilation of Caste**
If our people realize the truth that they remain poor and enslaved only because of caste, all the castes will be harmed. Since the enjoyment and dominance of Brahmins will be destroyed due to this, the poet Subramania Bharati carefully sings in favour of the upper castes. "*Aayiram undu ingu saadhi, enil anniyar vandhu pugalenna needhi* (There are thousands of castes here/ Why should the outsiders come and tell justice?") that is, "There will be thousands of castes in our nation. For that, an outsider must not come and question. What is his work here?" he writes. If someone says that the raft has a hole but only the river's water must not enter it, can we say that it is an intelligent thing?

Except us, who toils for the annihilation of all this? The reform that others suggest, the British themselves destroyed those differences and went away. He (the British) gave to the Depressed Classes the post of President in Panchayat Boards, District Boards, Municipalities. He made our friend Sivaraj a Member of the Legislative Assembly. He wrote in his law that there should be a fine for saying 'Paraiyan': that exists even now.

[Viduthalai, 16-3-1954]

**Education**

"Only the Brahmins should read, the Sudras should not read; the Brahmins should not teach the Sudras." This is what the Manu Dharma says.

Rajaji did this only because of the fear that the position of his society will disappear. If the *Vanaan* washes clothes, if the *Paraiyan* beats the drums, if the *Chakkili* stitches shoes, if the *Ambattan* shaves, only then they will get the feeling that they are a low caste. If they also get educated and come ahead, the highness of the upper castes will disappear. So, Rajaji interfered in the basics by introducing a 3-hour occupational education.

How did Rajaji who did not contest an election become a Minister? There is no place to remove these atrocious people through law. Who prepared this Constitution? Did only four Brahmins not prepare it? Dr.Ambedkar who belonged to the Depressed Classes was there. He had thought that if his caste gets 5% representation it is enough. But when he was given 15% he got exhilarated and he placed his signature.

Later, I said that I would burn that Constitution. He said, "Now, only I will burn that Constitution." Why? Because he understood its truth, its strength.

The servant of this Rajaji, C. Subramanian was the one who raised his hand to say that Sudras don't need communal representation. Those who had gone on the merit of our votes were sitting there like the five Pandavas during Draupadi's disrobing. All this takes place in this country. If it had been another country, four or five (death) anniversaries would have been observed for such people. When it is so, how can we fight legally?

When I said that the idol of Vinayaka should be broken, they talked of beating my effigy with slippers. I don't bother even a little about that. It is not being beaten for my personal matter. Is it being beaten because of the work I did for public good? I only considered it as a reward. Even if I am myself beaten I will get longer life and years.

[Viduthalai, 17-4-1954]

**Ramayana**

Ramayana was a conspiracy story created by Brahmins only to save Varnashrama dharma. Only to establish their doctrines, they say that Vishnu himself took the avatar of Rama.

Rama rules after having brought Sita from Ravanan. At that time a young Brahmin child dies. Its father goes to Rama and says, "Because you or your populace are indulging in activities against the law, my son has died today; a Brahmin has to live for 100 years, if he dies prematurely it means that the King of that land is acting in violation of the Sastras."

When Rama goes around the country, in a forest he sees a Sudra called Sambhukan praying to God in order to attain *moksha* (salvation). When Rama realizes that he is a Sudra, he tells him: "If a Sudra does penance to god it is violative of Sastras. If he wants to attain moksha, he should be a slave of a Brahmin and follow his commands, and if the Brahmin is satisfied, the Sudra may attain moksha." So Rama said, "You are violating the Laws of Manu" and with one arrow he severed his head. It seems that as the Sudra Sambhukan died, the dead Brahmin child came to life. In this manner, the story of Rama was fabricated to save caste and religion.

[Viduthalai, 11-10-1954]

**Ramayana**

How does the Ramayana end? Because Sambhukan, a Sudra prayed to god a Brahmin boy died. Rama found this Sudra and killed him, so this dead Brahmin boy immediately came to life. What is the idea they want to drive in the minds of the public?

A Sudra should not worship god directly! He should worship god only through a Brahmin.

[Call for revolution, 1954]

**Brahmin domination**

Today, even a Brahmin who begs has coffee twice a day but think of a farmer who toils all day in the sun and does not know what is coffee; he drinks only gruel with salt. On the other hand, any Brahmin working in a hotel or temple educates his children well and they become collector, judge, minister or President of India. But we toil and remain in the lowest position.

Is the ex-Chief Minister Rajagopalachari son of a Lord? No, he is just a purohit's (Brahmin priest's) son. He studied with the aid of scholarship. So also is the grandfather of Nehru, our ex-Prime Minister. President Rajendra Prasad's father was a tribal Brahmin. So, they have all risen to this position only because of their caste. Once our people become a little rich, they never feel that they belong to low caste. Have we ever thought over this? Never!

We want money somehow, and others should say that we are rich. We will spend money in constructing temples, performing *Kumbhabishekams* (consecration), and give 100, 200, 300, …, 1000 acres of land to the Brahmins. We once led a very high life. Only now, that too after the coming of the Brahmins, have we become a low caste. In any country, has a worker who works for the nation ever become a low caste?

These Brahmins who came for green pastures for their cattle made us build temples and make idols of god. Finally they said that if we touch these gods, the gods will die or become polluted by our touch. The Hindu pundits and the rich accepted all these, consequent of this we have become Untouchables and ruined ones. There are animals but does there exist a Brahmin donkey or an untouchable donkey?

[Viduthalai, 6-5-1954]

**Division of Labour**

Why are you a labourer? You are a labourer because you are of the labourer's caste,. Otherwise, does any Brahmin like Sharma, Sastri, and Rao etc. ever do your work? Today the labourers pull carts. They do laundry. They do manual scavenging. They work as peons. All of them are Sudras! If we look at why they are in these jobs, the reason is caste. Why does he sweep the streets? What reason can be attributed? If we ask 'Why does any Brahmin not sweep the streets? Why does any Brahmin woman not sell buttermilk?' he will say, "This is a backward, regressive force."

You must think well. You need to question your minds.

For several ages, there has been 'Paraiyan' and 'Paraichi.' Did anyone say that this must be eradicated? Only the man who does not bother about the infamy he will get will dare to say this. If it comes from the mouth of another, he needs to attain the necessary benefit.

Even now, will those in public life say that caste should be annihilated? Ambedkar made the law. In that he was not allowed to write that caste should be annihilated. They told him, "Write what you want for your caste alone." They allowed him to write 15% for his own caste. Why did they allow? Because, if it is given to that society, there is no one among them. Later, he revolted a little. They have beaten him such that he can't rear his head. Even such an Ambedkar himself did not say so. Who else can say?

Why could the Buddha not annihilate caste? He got cheated by talking of *Ahimsa* (non-violence). The Jain was also similarly cheated. They were in the opinion of not even killing an insect.

Today, everyone is our enemy. Why? Because I say that caste should be annihilated. Though they know that it is good for everybody if this takes place, they don't extend support. Instead they oppose it. Now we have got into a major struggle. That is why we oppose. This is the only job for four High Court Judges. These Ministers also have the same job. The Dravidar Kazhagam has written in order to compensate their preaching, their newspapers and books etc.

Yesterday, they have dismissed a constable. When we ask what the charge was, they said that they had removed him on the suspicion that he participated in activities of the Dravidar Kazhagam. Any teacher, any professor, any clerk who joins us is immediately dismissed. This has to be made known to the people. The Brahmin society itself must be destroyed. Sudras must not remain lazy without worrying about these people who created so much fraud. If you cannot, at least leave your children with me. With them, I will finish off these Brahmins.

[Viduthalai, 7-6-1954]

**Annihilation of Caste**

Sun Yat Sen was born in China and caused the reformation. Kemal Pasha was born in Turkey and made changes. But in Tamil Nadu, Siddhars, Valluvar and Buddha have been born and they wanted to abolish caste yet no changes have come. All the workers are our people but still they continue to be low caste. In the Malaya state, an Ezhava should not walk in the street but now due to changes he can walk in all the streets. But untouchability still exists in Tamil Nadu.

What happened to Gandhi? He became a Mahatma but was a victim of three bullets. The Laws of Manu say that a Sudra should not read, but Gandhi said everyone should read. What happened? He was shot dead by a Brahmin. Gandhi said that Hindu religion and Muslim religion were one and the same. So he was shot dead.

Protest against change will always exist, so don't worry about that. Each and every one of you should strive for equality and annihilate the difference between low and high castes. All of you remain united as Tamilians. Help those who struggle. The basis of my service is only annihilation of caste.

[Viduthalai, 23-12-1954]

**Division of Labour**

A section of people who fatten up without working are called upper caste—Brahmins. Those who keep on toiling, those who starve, those who go from village to village and lose their life in the midst of all this are low caste—Sudras—Panchamas—the fourth caste—the fifth caste; what a contradictory activity! Should we not abolish this?

\*\*

When there are representatives in the Legislative Assembly who say that it is pathetic for a man pulling a rickshaw, why is there no representative who says that manual scavenging is disgraceful, humiliating and unhygienic?

\*\*

It is a disgrace to the human society, when there exists the occupation of manual scavenging in a country. This method is not necessary. All over the country there should be a modern scientific method of automatic waste-disposal.

\*\*

The Dravidians should not go to those temples where Brahmins function as priests and where there is compulsion that only Brahmins must perform prayer rituals. Even if they go, they should not have the prayers or rituals performed by Brahmins.

[Viduthalai, 19-1-1955]

**Annihilation of Caste**

The world's virtuous, the greatest of Mahatmas, has he said one word that caste should be annihilated? He only said that untouchability must not be there. Has he said a word that caste should be annihilated? Moreover, he only said, "I struggle for the protection of Hindu religion. I came for the establishment of the Varnashrama dharma. It is my primary aim to create Ram Rajya." He blatantly spoke of whatever was necessary for the protection of caste and he also said that he was working for the same.

   Only because he said so, he could live this long and get the title of Mahatma. The Brahmins assassinated him only because he, having realized that Hinduism is itself a fraud, said in spite of himself that caste problems should be annihilated.

In the same manner, the Brahmins ensured that there was no evidence about every great person who came to annihilate caste. Buddha and Valluvar came like this and the Brahmins victimized them severely. Till now, Valluvar's *Kural* was not known widely. Only because of our ruckus, people have started reading it. Though Buddha was born in this country and preached here, his ideology is present in other countries but not even to a little extent here. What is the reason? In the course of time Buddha's philosophy could not spread here without any barriers.

If I survive even now, it is a different matter. There is a plan, and a betrayers' group, to annihilate me. There will be conspirators waiting as to when I will get caught so they can put an end to me. This is because the work that I do causes a lot of trouble to the Brahmins. Without any means of living they have to flee the country. That is why I know that there will be such an opposition from them.

Today if I leave all my policies and say, "Religion should be there. The greatness of the Hindu religion created by God himself is really immense! The Brahmins are men of avatar! They are next to God! Only yesterday night, Ambal appeared in my dream and said these truths" it is enough. Immediately, I could become a Mahatma. They will say that Mahatma E.V.Ramasamy is the tenth avatar of Vishnu. My portrait will hang in the entrance of the house of every Brahmin. Daily they will worship it with garlands, incense and offer it flowers and fruits. Seeing me say so, the government will also say that the devil has left us and it will plan a major fraud. They will think, "Only one man was a great enemy, from now we are free" and they will enjoy.

I say this because there has been a conspiracy against such a great man in the world like Gandhi who lived in front of us. At such a time, even though I feel that even I have a risk, day-by-day I am more intense because I don't care about such threats. By saying my truths to the people, even they have understood the atrocities of the Brahmins and the dishonesty in the Sastras, puranas and they have reached a stage of feeling Self-respect. From now, they will only strive to remove their disgrace. I feel that I do not have the necessity to work like I worked earlier. Till my death, I will dare to say the basic mantras and the policies of the Kazhagam. That is why, I research fabricated stories like Ramayana that were intended to save religion and establish the low state of Dravidians, and expose their obscene, uncivilized stuff first through the *Kudiarasu* and through other published books. Now I go from village to village, get several Ramayana books and point out the errors in them.
[Viduthalai, 24-3-1955]

**Rationalism**
Like the other countries, even we have to be intelligent people and be one caste, the human caste. We need to give space for rationalism and have the human quality of following anything as per rationalism—that is how our doctrines have been made. Only on the basis of that doctrine, we say that Brahmins have to be driven away from our country, puranas must be burnt; and that the idols of Gods must be broken and using that streets must be laid.

Temples must be demolished and razed to the ground. Their limitless wealth must be spent on improvement of human society. For all this to take place, we carry out agitation and propaganda and publish our opinion through newspapers. So, it must be the duty of Dravidians to encourage and support these.
[Viduthalai, 15-4-1955]

**Non-Brahmins**
The necessity for us to take part in such activities is primarily because we want our low status and disgraceful caste titles like Sudra, Panchama to be removed. We were the first ones to plan that such low caste titles should be abolished.

The feeling of Brahmins and non-Brahmins started in the year 1894-95. Those who realized the cruelty of the Brahmins strived to condemn it. They had given it in writing to the Government.

[Viduthalai, 20-4-1955]

**Politics**

Nehru condemned caste on his recent visit here. Till now they did not have so much of concern in condemning the caste system, we do not understand what is the sudden reason for their talking about it now. The annihilation of caste atrocities and caste itself—We are not able to say how far this will be true.

As the next stage, the general election is going to come. We feel that only with this in mind, Nehru speaks these ensnaring words. Majority of the votes is only with the Depressed Classes. They have realized the truth that if they have the confidence of upper castes and Brahmins, that alone will not fetch them votes. So, to trap the Depressed Classes in their nets, they speak so.

[Viduthalai, 28-10-1955]

**Buddha**

Did anyone ask if low caste and high caste existed historically? Only Buddha questioned. He was a king's son. He questioned several things. He asked, "Why is he an old man? Why was that man a servant? Why is he blind?" Buddha asked, "Why is he a low caste man?" They replied that god created it. Then he asked, "Where is that god who created the low caste?"

Then they spoke about the *athma* (soul). He asked, "What is it? Was the athma seen?"

They chased away such persons who worked for the nation! Buddha Jayanthi is celebrated in the month of May in the days 24,25,26 and 27. They requested me to join the celebrations and I accepted the request made on behalf of their association.

Also, I request Dravidar Kazhagam comrades to join the celebrations and campaign about the advice and laws given by the Buddha.

The most important of Buddha's principles are

1. Use your intellect to analyze everything clearly.
2. If your intellect says it is correct then accept it.

Don't believe the terms like God, Athma, Devas, Heaven, Hell, Brahmin, Sudra and Panchama that you cannot understand. They are only imaginary terms. Use your common sense to analyze every object. Do not follow anything because God said so, Vedas says so or Manu says so and trust them. What your mind says believe that alone.
[Viduthalai, 19-4-1956]

**Aryan domination**

Only the Aryans made the Dravidians barbaric by denying them education and all sorts of development and made them a low people.

But in these 35 years of the Self Respect Movement, people who are in castes lower than the Sudras have become Chief Minister, one Sudra has become a High Court Judge, one Sudra has become Inspector General of Police, one Sudra has become an Education Officer. This is actually against the Laws of Manu. In Tamil Nadu such things have happened against the cunningness of these Brahmins.

Dr. Ambedkar boldly said in the Lok Sabha, "I don't believe in god or spirit" and completely came away from the Aryans.
[Thozhan, 20-5-1956]

**Hinduism/ Vedas**

The person who follows the Vedas doesn't know anything about Hindu religion. Only due to the mantras and Sastras, we have Brahmin and Paraiyar, high caste – low caste, countless gods, lakhs of temples and trash.
[Viduthalai, 3-6-1956]

**Caste discriminatory practices**

Who said that the Paraiyan, Chakkili, etc. should enter into the temple? Only this Ramasamy said so at first. That too, I said that the Gods must be broken and roads must be laid with the debris. Only after we said that we would leave this religion and become Muslims did the others open their mouths. Even at that time, several people said that worms would infest my tongue. Before we carried out the self-respect propaganda and reform, there was segregation in teashops, they were not allowed to board buses. You will not know this because 30 years ago all of you were children. At that time the Depressed Classes could not walk on the streets. They have to walk only over the sewers. They were not allowed into theater. Notices that 'Panchamas and those with contagious diseases are not allowed' would be displayed. Only the Brahmins cooked in hotels. Brahmins would eat separately, and we would eat separately. Does this take place now? Who changed all this? Which party worked for this? Which religious head said this?

If those who scavenge don't do so, will it not stink? If those who sweep don't do so, garbage will collect here and there. If those who wash clothes don't do so, we will have to wash it ourselves. If the farmers do not cultivate, then we will not get food. If the weavers do not weave, there will be shortage of clothes. But if the Brahmins are not there, what work would not take place?

In other countries, does work go on because of Brahmins? Brahmins don't have any occupation. What work does he do apart from spinning a yarn in sleep (lazing around)? It is a good thing the British put numbers, otherwise for every milestone and furlong stone our people themselves will stand there, call a Brahmin to ask him, "Can we go ahead only if we give money?"

So, several people among us should struggle to the extent of losing our lives. Then Brahmins will come forward to say, "We will also work along with you. Our women will also come to transplant crops, they will also work."

A man, who does not do anything for God, flies in the sky. Only our people, 2000 of them will pull a wooden chariot for a distance of a furlong. These people make a riot as if mad men were fed liquor. So, we must reform our people. We should see that our children are not superstitious. We should ourselves ask why they smear ash. He must think and research and know the answer. Only then intelligence will develop among our people.

[Viduthalai, 18-7-1956]

**Social Reform**

In all my demonstrations for social reformation I have spent my body, money, mind and spirit and I am waiting to spend this life. Instead of fighting, people get refuge from Sankaracharyas.

[Viduthalai, 30-7-1956]

**Morality in Hinduism**

The human qualities of one man cheating another, upper caste-lower caste, heaven-hell, ritual and further the evidences in our religion about homosexual relationships, bestiality, illicit sex are read and enjoyed, lived and heard by the devotes so they have no moral values or fear for moral values or even dislike or hatred for immoral and sinful acts.

Our caste system has made the low caste into respectless fools. High-caste men are culprits who cheat others mercilessly with hearts of stone.

Unlike Christianity or Islam, love of brotherhood or affection for others, or kindness, or respect for humans are never found in this religion. In this religion, the dharma is that one man should treat another as an untouchable, cheat him, exploit him and live on others' labour etc. This is being preached and practiced. How can we have good conduct or moral standard?

[Viduthalai, 3-8-1956]

**Partition**

If caste is not annihilated and this country does not become ours we need to take swords and petrol in our hands. See, in Ahmedabad slippers have become bouquets. Nehru begs, "Other countries will tease us seeing the happenings here. Please wait for five years. I feel ashamed to change the law which I made." Is it because he is afraid of the law, or is it because he is afraid of the violence?

How did the Muslims get Pakistan? Only with the sword! In Navakali did not 10,000, 12,000 people die. Take even today! Pakistan might have taken two beatings from India, but it is not that India has recaptured Pakistan. To cheat us, they kept talking *ahimsa* (non-violence) all the time.

If we take a sword or a stick, will a Brahmin remain standing? That is why he spoke of non-violence! A man who does not keep a knife with him is a coward! If I have to say more, I am going to say that even if you do not wear a black shirt you should keep a knife with you. Only then we will not have traitors among us.

If we make anyone influential, he turns into a traitor. If we take a bachelor to the ceremony of seeing a girl as a prospective bride, will he scrutinize the bride for himself, or for us?

In 1938, when the British were there, before independence came, I said that we need a separate country. The Congressmen say to me that what I say is correct. When there is so much of growing support if we lust for a legislative assembly seat, a parliament seat, and permit licence what does it mean? If the country goes like this, the Sudra will have to become the son of a prostitute.

[Viduthalai, 2-9-1956]

**Caste discrimination**

Let us assume that Nehru is going abroad on a tour. Someone sees Nehru and asks him, "It seems the people who scavenge in your country are of one caste. Those who sweep the streets are of one caste. It seems there is a caste called the Brahmins. If they touch others, it seems they would get polluted. Is this true?" What will he say then? Can he deny this? There is no other option for him than to say, "Yes, it is there. I am struggling to annihilate it."

We have to find out from where caste came. If a few Paraiyars apply sacred ash on their foreheads, or the Mark of Vishnu will the Paraiyan quality go away? They think that their disgraceful aspect has gone if four-five people become ministers and members of the Legislative Assembly. That only means the others are low people. After independence they called themselves Paraiyar, Chakkiliyar, Adi Dravidar. Because it was against the Brahmins, they called themselves Harijans. Calling themselves Depressed Classes, Adi Dravidars, Harijans is only a mere change of name. What else is it?

If a chaste woman calls another woman, "Chaste woman, come here!" what does it mean? It only means that she herself accepts that she is not a chaste woman. Even in 1956, if you call a *pappan* as a Brahmin, it only means that you are accepting yourself to be a Sudra, a son of prostitute, a low caste.

If this small thing cannot be changed, when and where do we annihilate caste? Who will listen if we ask them not to say their caste names?

If we write to the government, "He has put a board in his hotel saying it is a Brahmin hotel and he wears a holy thread," they write back saying that is his right. If that is so, to establish our right what can we do apart from erasing the board and cutting the holy thread and the tuft?
[Viduthalai, 3-9-1956]

**Development of the Depressed Classes**
In this country, they have allocated two or three crore rupees for the development of the Depressed Classes. To whom did they give it? They have placed it in front of the Brahmin who calls himself upper caste and on seeing us says, "Don't touch. Stand apart." If the Government had true concern about the advancement of the backward people shouldn't they have given the money to me, since I have spent a lifetime working for the annihilation of caste? Otherwise they should have given it to the Depressed and Backward Classes and told them to improve their communities. Even now some people say Brahmins are better than our people! If they say these words even in 1956 what does it mean?

Should I not see a society where caste has been annihilated? You are young so you can see it later. I am in a state of (passing away) today or tomorrow. If we say, "Brahmins quit" must we ask "To where?" Is it intelligence? The Brahmin can ask, "I have the right to be in this country." There is a meaning in his question. Can we Dravidians ask that question? What have you achieved by standing hidden for so long?

If caste has to be abolished, the Sastras must be burnt. The gods must be broken down. Likewise, if we want freedom, we must burn the national flag. Pakistan was created only by killing 1000, 2000 people. When Jinnah said that Pakistan was necessary, did one Muslim oppose it? If one Muslim had opposed it, will any Muslim ever attend any life or death ritual in that man's house? You have seen it with your eyes. What does it mean if you speak like this? If I say, "Take the knife," you criticize.

[Viduthalai, 4-9-1956]

**Caste discrimination**
They said that after Swaraj comes, they would make the Paraiyan sit in the center of the house and would make him eat in the common dining. What has the Government done? Even now, the Sudra title has not changed! All this trouble is only because there are 3¾% Brahmins in this country. If the Government really has the intention of annihilating caste it must have given scissors to the policemen and asked them to cut the holy thread and tuft if they spotted it.

A few days back the government put a picture. In that a Brahmin says to a Dravidian, "Don't touch. Go afar." In a picture below they show the Brahmin inside the jail because he displayed caste difference. In the picture on top the Brahmin has a sacred thread. In the picture below, it lacks the sacred thread. This is because the Brahmins in jail do not have sacred thread. According to prison rules, the sacred thread is not permitted. So, they have removed the sacred thread in accordance with the prison rules.

The government printed 10 lakh copies of this picture. One picture cost 4½ annas. On the advice of a Brahmin who said that it must not be published, and being afraid of him, the government burnt all those pictures. What does it mean if they are scared of one Brahmin and burn all the pictures?

Four crore rupees have been allocated for the development of Harijans and that amount has been given to the Brahmins by the Government. The Brahmin financially supports those Depressed Classes who betray us. For whom did we break the Pillaiyar? For whom did we burn the picture of Rama? Why do we say that caste should be annihilated in this country? For whom do we say that Brahmins should be thrown out of this country? Everyone must think over this.

[Viduthalai, 12-9-1956]

**Annihilation of Caste**

In foreign countries, if a man is a barber, his brother would be a minister, his father's younger brother would be launderer, his father's elder brother would be a judge. There is no separate caste structure for every occupation; everything is common for everyone abroad.

If caste is touched, religion will be destroyed. If caste is annihilated, harm will fall upon the Vedas, the Puranas, the Sastras, God etc. If these are harmed, the livelihood of Brahmins will be harmed. That is why when we say that caste should be annihilated, people like Rajaji say that it should be protected. In a nation that is independent, in a nation that is said to have attained Swaraj, why should there be one person called Brahmin and one person called Paraiyan?

Ask the communists why Brahmins should alone be high caste? He will say, "Don't talk about that. It is sectarian." This will be the answer if you question about caste to a man of any party. Everybody says this because they need the grace of the Brahmins and the north-Indians.
[Viduthalai, 15-9-1956]

**Depressed Classes**
When I said that only the Depressed Classes had got communal representation and the Backward Classes had not got it, the Depressed Classes feel that we are against them. I feel that I have done a lot for the Adi Dravidars. Yet, they say that the Brahmins are better, and that they face difficulties only due to the caste-Hindus. This is ingratitude!

Who has brought them to this status? Brahmins? Only due to our demonstration and agitations they have got it (communal representation). Only we have shown them the way to education. We were the first one to fight for entry into the roads around the temple. If they are educated, or occupying nice posts, or serve as ministers, it is not due to Brahmins. I do not expect gratitude from them, what I expect is no enmity.
[Viduthalai, 21-9-1956]

**Brahmin domination**
It is all right even if we lose our life or get ten years punishment, we will suffer it. But our disgrace must be eradicated. Why are we low caste in our country? Why are they high caste?

Does a man of another race rule any other country in this world? A crowd that came to beg, who only number 2 of 100, should they rule? Should they get jobs? If we ask this, they say, "Why do you separate the Tamilian and the north-Indian? Is everybody not a child of the Bharat Mata (Mother India)?" If that is so, how many husbands does Bharat Mata have? If there are Brahmins, Vaishyas, Sudras, and Panchamas, what do we say about the morality of that woman? Does a Brahmin women break stone anywhere?

Does she scavenge? Does she sweep the streets? Then how do we say that everybody is the child of Bharat Mata?

Let it come through law that there shall be no Paraiyan or Panchama or Sudra in this country. Will it come? If it comes, what work is there for our party and us? In order to get rid of the disgrace that has been imposed on us, we need to do whatever has to be done. Everyday our people should increase their feeling, "Am I a low born? Am I a coolie?"

This fifth circle is a part of the Chennai city: The people here are Paraiyars. A low caste. All those who reside in Mylapore, Tiruvelikeni, T. Nagar are Brahmins. Can this difference exist? But it is there. Why?

Who is Pundit Nehru who rules our nation today? Who is this Rajaji? Nehru's granddad was a priest. Rajaji and his father and his granddad and his race come to our homes talk of *thithi, thivasam, karumadhi, kaledupu* (rituals related with death, in order to ensure that the dead ancestor reaches heaven) and if we had not put money on their plates could they have become India's Prime Minister, Governor General, President? It is enough if we realize this. Their race will be eradicated. We have a series of plans. If all this does not work, finally there is no option but the sword.

To attain our successes, even if we sacrifice our life, we need to achieve it. Is it democracy if the majority works as a coolie, and the minority rules?
[Viduthalai, 05-10-1956]

**Brahmin domination**

I created the Self-Respect movement. I campaign that caste must be annihilated. Today, in reality, there are only two castes. One is the Brahmin caste; the other is the Sudra caste. One is the caste of the master; the other is the caste of the slave.

'If all these Sudra castes unite, our situation will become horrible'—so Brahmins divided this one caste into Chakkiliyan, Paraiyan, Chettiar, Mudaliar, but in reality there are only two castes. In the Manu dharma, Mahabharata, Ramayana, Sastras etc. only these two castes are mentioned. That is, the master caste and the slave caste. The Brahmin will not plough; he will not come and work as a constable. All the disgraceful work, that is, all the work were one must say "O' Master" is for us. The work of [authoritatively] calling us "Dey!" is for them. No risk must come for this arrangement so he is giving trouble carefully. When I see all this, what Swaraj or nonsense! Whether the white man rules, or the Russian rules, or the Chinese rules, it is okay; I feel that this north Indian and Brahmin must be annihilated. From now, we must not get afraid; we have major work to be done. I wonder why we need to fight for petty reasons. What will happen if we violate the order and hold a meeting in this Island Grounds? We will be punished for eight days. Can we not do it?

The 'tear drops' of our party get wages from the Brahmins and they say that only 'Brahminism' is not needed. The leader of the communists is a Brahmin, likewise the Socialists. We need not even ask about the Congress. They never help us.
[Viduthalai, 7-10-1956]

**Religion and caste**
Caste must be annihilated. Religion, which supports castes, must also be annihilated. The basis of this religion like puranas, Sastras must be annihilated. Brahma, Vishnu, Rudra, Yama, Vayu who keep us as Sudras are not Gods. They are Gods of the Brahmins. Ask a Turk or Christian or people of other religions. They have only one God.

"If you have the belief that there is no form, it must be in your heart. If you believe it even a little, fear god and act with honesty. Why do we need Gods with several forms? Moreover, God who says that we are Sudras, Panchamas, 4[th] or 5[th] caste, the several Gods, gods with forms, why gods with women, children, gods who demand food—all these need to be destroyed. Only then caste will be destroyed. Man's superstitions and foolishness about god will be destroyed. Who worries about all this?

I say as it is: to go for a job, education and degree are not required. "Education is not there, mark is not there, qualification and merit have gone;" shouting hoarse over this is deception. We are giving 15% representation for Adi Dravidars. How many of them have been educated? How much marks and merit and qualification do they have? What has happened because they have been given representation? No one has complained so far. Only we fought for them. When the Adi Dravidars are given 15%, must not those who are more eligible than them be given their due percentage? If the power is in our hands, there will be no jobs like pulling carts and breaking stones where the toil involved is high.

Even I have gone to foreign countries like Portugal, France and Russia. The people there don't do such work. Everything is mechanized. Only machines do even the work of transferring sacks to vehicles. Here railways, income tax and other profitable enterprises are given to north-Indians and we have to live along with the sales tax. And when they need money, they take it like alms and thus run the government. This year I have planned several important agitations. Real agitations take place often. You should support it.

[Viduthalai, 10-10-1956]

**Communal Representation**

How is it fair if Muslims, Christians and Malayalees join together with Brahmins and stomp our heads? During the British rule, Muslims and we were as intimate as brothers. By canceling proportional rights, they caused trouble. If there are three people in the Public Service Commission, should one of them be Muslim and another a Christian? We should certainly get communal representation. It is said in the Constitution that Adi Dravidars should be given 15% representation. This is a Public Service Commission order; a copy of it even came to the training school that I run. I am happy; I will give as it is. What is the way for us? What does it mean if efficiency is demanded in our matter alone? Why was it written in the Constitution that Backward classes and Depressed Classes need to be given concessions? Why was efficiency not demanded then? The Adi Dravidas are a community without 5% literacy and which industry was spoiled because it gave them jobs? What will be ruined, if all of us are given representation as per our proportion in the population?

[Viduthalai, 14-10-1956]

**Hinduism as a religion**

In foreign countries, there is no social difference that one person is high and another is low by virtue of birth. So they are qualified to talk about political and economic equality. Here they have abandoned the task that needs to be done first and instead they keep talking about this and that. There is no one to talk about caste or religion. That is why I am taking efforts for social equality and annihilation of caste, which others have abandoned, and there is nobody to talk about it. As I keep observing the world, I lose interest in other things.

In this invitation you have used the words "Reform in Hindu Religion." There is no religion called Hindu religion. In the British period, they did not accept that a religion known as Hinduism existed. While describing electoral constituencies they named it Muslim constituencies for Muslims and Christian constituencies for Christians. Otherwise, they did not call it the Hindu constituency for those who claim themselves to be Hindus. They called it Non-Muslim constituency even in the printed nomination form. The British asked, "How is someone a Hindu? What is the evidence?" Today, what is the evidence for Hindu religion? Other religions have evidence; there are historical evidences. Muslim means they can show with historical proof that their religious head Muhammad was born in this place, in this year, he did this and so on. There is history and proof to say that Jesus was born 1956 years ago, he did this, and he said this. But what proof do those who speak of Hindu religion have?

[Viduthalai, 11-11-1956]

**Eradication of Brahmins**

Even I want a good government. Only after me, the others came to the Congress; Kamaraj was then my volunteer in Coimbatore Those who were ordinary people when I was the President of the Tamil Nadu Congress Committee are ministers today. Can I not become a minister if I want to? They are not going to have a loss even if I steal money. They will give me the Minister's post.

If I come away who will do these work? Rajagopalachari carries a counter campaign that if caste is eradicated, the nation will get spoilt. No Brahmin should be there in this nation. The word Brahmin should not be there in the dictionary. The stories of Ravanan, Soorapadman etc. were created to threaten you that you shall meet their fate for having opposed the Brahmins.

[Viduthalai, 23-11-1956]

**Kings and Brahmins**

This situation came about because all the kings of this country, Cheras, Cholas, Pandiyas, Naickars were all henchmen of the Brahmins. Having a little bit of shame, if they had chased away Brahmins, the Brahmin would have run away. Instead of that, they kept the Brahmin as guru and did oblations, sacrifices and rituals and gave away their wealth to the Brahmins. Thus they remained slaves of the Brahmins.

[Viduthalai, 27-12-1956]

**Kings and Brahmins II**

What was the integrity of our kings? There was a king called Vallala Maharaja. A Brahmin came and said that he wanted a woman. He searched everywhere in the village. He did not get a woman. He wondered how he could say to a Brahmin that there was no woman? So he called his wife and asked her to be with the Brahmin. It seems this is a holy story. If they had shame or respect should they not have torn this and set it afire?

Caste is not connected with us. We don't know how we accepted it? The word '*Jati*' (caste) is not a Tamil word. There is no word in Tamil for that.

[Viduthalai, 28-12-1956]

**Buddhism**

There is no religion like Buddhism. There are only a few doctrines that the Buddha said. These are realized because of intelligence. To follow those doctrines there is no necessity to go out of Hinduism. Everything that Dr.Ambedkar said that he didn't believe in or wouldn't practice on joining Buddhism, we already don't believe it or practice it. We don't believe in *moksha* (salvation), hell, fate and 7, 8 births. We don't have faith in it, we do our work: there is no necessity to go out of Hinduism in order to tell this. No Brahmin can say that we are not Hindu just because we say that caste—which exists according to the Varnashrama religion—must be annihilated.

[Viduthalai, 28-12-1956]

**Buddhism II**

Now Buddha's doctrines have become like a religion. Just like the Muslims have an identity, likewise in countries like Burma and Ceylon the Buddhists also have an identity. They have their temples, they pray kneeling down, and with these identities it has become a religion. I said this to Ambedkar, "If we go there (to Buddhism) fearing the Brahmins in this religion, then who can exist with these Bhikkus?" If we leave this grip and talk one bad word about Hinduism tomorrow, blame will be heaped on us.

We are alone worried about this. The others are only eager to go to the Legislative Assembly. All the parties today: the communist, socialist, 'tear drops', are all caught in the hands of the Brahmins. If one goes to the Legislative Assembly itself, automatically only the love for the Brahmins comes inside them. Why? Publicity is the reason because the newspapers are only in the hands of Brahmins.

[Viduthalai, 29-12-1956]

**Representative of Depressed Classes**

I speak as a public worker; I am ready to make any sacrifice for the Depressed Classes and Sudras. I talk as a representative of all of you and not for self-esteem.

[Statement filed in the Madras High Court, 23-4-1957]

**Caste discrimination in burial grounds**

Right from the inception of the Self-Respect Movement in 1926, we have been saying, "The Brahmins, who enslave us, must be destroyed. Along with destroying Sastras, puranas and epics created by Brahmins, the Gods on whose names they fatten themselves must be broken and roads must be laid with those stones."

If we look a little earlier, I quit the Congress in 1924 and I have been writing with great clarity. Not only have I written that no man should have a title or status in the name of his caste and that if the caste name is there in the puranas and Sastras it should be destroyed but also I have spoken so in several meetings and through resolutions in conferences.

We (the Self-Respecters and the Dravidian movement) are doing this continuously and no one can consider that we have suddenly begun it today.

When I went to offer flowers and pay homage to the memorial of Comrade Pattukottai Alagiri, I found in that burial ground it was embossed stones demarcated the place where the Sudras were burnt, where Brahmins were burnt, where the Saurashtras were burnt. On seeing it, I asked them for how long this had existed. They said that is was there from the time of the Tanjore Raja. I asked them to write to the Municipality to remove these plaques and said that if they refuse to do so, we could take the necessary action. Perhaps they immediately passed a resolution and it was immediately removed. I say this because it is worrisome to see the names that disgrace us.

[Viduthalai, 5-5-1957]

**Revolution**

Statues of the British remaining in this country are itself considered harmful to self-respect and welfare of our citizens and agitations are instigated against it. What wrong can be there in the agitations against the Brahmins who are thousand times more harmful, unnecessary and troublesome to the indigenous people of the country than the British? Is the democratic rule a rule for the majority? Or is it causing trouble, shame and snatching away rights in the name of the minority? They can call it the rule of the mighty. I am afraid that if the government does not act with honesty and intelligence in this agitation, that affects the Tamil people's self-respect and feeling, people will be constrained to start a bloody revolution.

A bloody revolution might require several people to die. Did Pakistan which sacrificed and killed peoples in tens of thousands get annihilated? Or has that society disappeared?

[Viduthalai, 11-05-1957]

**Caste identity**

In those days people would see the other's ears, count the piercings and find out ones caste through that. Because the feeling of caste has reduced, they are not able to find out now.

So, piercing the ear is only a traditional ritual that protects caste. As education grows, these will be automatically destroyed.

[Viduthalai, 25-05-1957]

**Brahmin to be called Pappan**

Even if you happen to forget, the word Brahmin must not come out of your mouth. Call them "*Paappan*", "*Iyan*" and don't call them *Brahamanan* (Brahmins). Only the 'tear drops' will use the word Brahmin and not any other word.

[Viduthalai, 30-06-1957]

**Brahmin hotels**

The 1st of August (1957), must be considered as August Agitation Day. On 1st August, from four to eight in the evening, in all the villages in Tamil Nadu, in front of Brahmin hotels everywhere (whether it is called Brahmin Hotel or not), our comrades, without causing any hindrance to the traffic or to people in a public place, with peace and in a compromising manner tell them,

"Sir, this hotel is a Brahmin's hotel. This caste has ruined our life and our society. They cheat us and without effort they exploit us, they fatten and call us Sudra, son of a prostitute, Paraiyan, Chakkili, Chandala, low caste and so on. They have written so in their laws, Sastras, ithihasas and Veda puranas and are treating us likewise. To go and eat food in this Brahmin hotel is a disgrace for us! It is a great disgrace! Shameless, base disgrace! Sir, please do not go there."

You should fold your hands, bow down and beg them and see that they go back. You must also print pamphlets and distribute.

[Viduthalai, 23-07-1957]

**Annihilation of caste**

Show me one God of the Tamilians who has a Tamil name? Why does God have a marriage, a wife and children? We must look at God with intelligence; why does God with a compassionate heart require many weapons, bow, sticks and sword? God who is the embodiment of morality has been made so vulgar. Why?

The Brahmins have written a lot of puranic vulgarity and thus spoilt morality. Intelligence and honour are required for the annihilation of caste. This is very essential for people. In any way, caste has to be annihilated. For the sake of our food, for how long do we need to suffer as coolies, labourers, and physical toilers? Should we keep doing our hereditary jobs?

If caste is annihilated, we will develop, and our standard of living will also develop. So, somehow, in this period, we need to destroy caste. Even in pictures Brahmins should not exist. We did not get proportional share in employment or education. In the Government our people are there in jobs of peon and police! Only Brahmins occupy the prime jobs where one can stay comfortable and draw large salaries! Jobs should be given community-wise proportionally.

[Viduthalai, 25-07-1957]

**Caste Annihilation Struggle**

That day the opposition parties could not remove the hereditary occupation education scheme. The Dravidar Kazhagam could do it. I said this to prove that though they say that they will go to the Legislative Assembly and do something, nothing can be done by going to the Legislative Assembly. Remaining outside, only the Dravidar Kazhagam works for the people. All the others can go to the Legislative Assembly and earn salaries of Rs.150, Rs.500 and they cannot do anything else. Even if they rise up to speak, the Speaker of the Legislative Assembly will ignore them. In this country, only our party says that caste should be destroyed; the Dravidar Nadu should be partitioned; the north Indian must quit; God, Sastras, and tradition must be destroyed and works towards these. The people of other parties will speak beautifully to get votes. Also, the newspapers of this nation, with restrictions do not write about our activities. Even if they write, they will publish it in some corner where nobody can see it. Till now more than 300 people who participated in the caste eradication agitation have been arrested and punished. But, no opposition party has questioned it. Only the Dravidar Kazhagam lives without the grace of Brahmins, government, publicity or newspapers.

This year we are going to do two agitations:
1. caste annihilation struggle
2. Tamils should do pooja in the temples instead of Brahmins.

The first struggle is the Caste Annihilation struggle that is to take place in Chennai as a hotel blockade agitation. I

have planned to do this Caste Annihilation struggle for one year.
[Viduthalai, 27-7-1957]

**Destruction of Brahmins**

As far as I am concerned, only I destroyed caste. In those days, I ate in the homes of those who weaved mattresses. Somehow caste became bitter for me! In my house they would not let me inside the kitchen. My mother would not touch me till her death. They will take the tumbler I drank water in, turn it upside down and sprinkle water over it. They were so orthodox. Without any worry I ate in the homes of Muslims and Vaaniya Chettiars.

Earlier caste identity was there. I am responsible for its disappearance. After I die, people will again go and fall in the Brahmin's feet. For the last twenty years, like me who has been running Kudiarasu, Pagutharivu, Viduthalai, Self-Respect Movement? Even the Congress, has it followed one doctrine for thirty years? Only the Dravidar Kazhagam works for a single doctrine. How was our life forty years ago? Were not the Adi Dravidars unable to go to schools? We said that schools that didn't enroll Adi Dravidars would not get the grant. Only after that they enrolled…

Only if the Brahmin is destroyed, caste will be destroyed. The Brahmin is a snake entangled in our feet. He will bite. If you take off your leg, that's all. Don't leave. Brahmin is not able to dominate because power is in the hand of the Tamilian.
[Viduthalai, 30-07-1957]

**Disfiguring Gandhi's statue**

Kamaraj and I are very intimate. Can Kamaraj not destroy caste in one sentence? "Both of them spoke beforehand and they have done so," these 'tear drops' are doing a negative campaign against us.

You think over this! More than us, only for Kamaraj, caste must be destroyed. Only he should have total desire and effort. If the Minister's post goes, Kamaraj is only a Nadar. If this land comes in the hand of a Brahmin he will not even call Kamaraj as Nadar. He will only be called *saanan*! Does he not know that?

In our country's Constitution, there is no space to destroy caste. One should not consider that caste has been destroyed just because a Paraiyan is able to enter a hotel and a temple. Gandhi did this conspiracy to create such an appearance. What will be their state after the privileges given to them now are fulfilled? The world does not know of the caste system that exists in our country. But no newspaper or news agency can hide from the world Gandhi's photograph being burnt or his statue being broken.

"Caste is there in India. Brahmin is the head of the castes in India. He is also the head in politics, finance, and society. Gandhi saved the Brahmin society. He made the other communities slaves of the Brahmins. That cruelty is unbearable. The people of the south break the statues of Gandhi. Nehru is unable to rectify it."—The world will know this. Nehru cannot fly as the peace ambassador in international politics. I will certainly create this situation.
[Viduthalai, 03-10-1957]

**Gandhi's assassination**

Gandhi calculated the number of almonds that he ate. He said that he would live for 120 years. But the Brahmins shot him dead. If the same thing happens to me, what do we know? Beyond all that, who knows how many more days I will live? Doctors say that as far as the body is concerned I am fit. When I see the people's affection, I get the enthusiasm that I should live for some more time.
[Viduthalai, 7-10-1957]

**Press**

Outsiders can know information about us only if it comes in our own newspaper. Without the favor of other press people, only we, only our doctrine is gaining importance in the country. Once, in one of my demonstrations, I said that the dailies or magazines have to be slaves to the Brahmins for their existence, so they write like this. The journalists are either Brahmins or those who need a living. In foreign countries, newspapers are run for ideology. Here it is not so. If it is news about us, no one will publish it. Anyone else would have become influential at the mercy of the newspapers.

Only the press elevated Gandhi. I know that he is an ordinary man. For the welfare of the Brahmins, the journalists elevated him and made him a Mahatma. I spoke to Vinobha Bhave for more than an hour. He is nothing, just a mad man.

I say this because we are the only people who have become significant on our own because of our ideology and without favor of the press.

Comrades! You must listen carefully. Even if anyone else reports our news, he would only put it in such a way that it makes you hate us. So, you need to pay attention.

Why should caste be there? Why should there be high and low? Why are we so low to the extent that we should drink the water that was used to wash the Brahmin's feet? Even today he buys the water in which Sankaracharya's feet are washed and drinks it. What an atrocity! Must this happen in 1957? Why should the Dravidian race be Sudras?
[Viduthalai, 07-10-1957]

**Breaking Gandhi's statue**

They have made arrangements in the law to protect caste. Gandhi was the basic reason for doing so.

They are able to make arrangements to protect caste, only in the name of Gandhi; so I will say it is a great shame for Gandhi's statue to remain in this country.

We have the right to break Gandhi's statue.

Today the Congressmen say that the statue of Wellington must not be there, the statue of Queen Victoria cannot be there, Neelan's statue must not be there! Likewise, I have the right to say that Gandhi's statue must not be there in our Tamil Nadu.

Gandhi has committed the atrocities that Wellington and Neelan did not commit.

Gandhi cheated us whole-heartedly, made a strong law to protect caste, gave protection to the Brahmins, and made arrangements for us to remain slaves forever.

[Viduthalai, 08-10-1957]

**Gandhi and Untouchability**

Dr.Ambedkar has written a book. Gandhi's deception has been brought about nicely in it. That book is "What Congress and Gandhi have done for the Untouchables?"

They have an organization to eradicate untouchability called the "Untouchable League." They have also given lots of funds. They appointed a committee for this league. One of the Committee members while speaking said, "Community inter-dining must take place; caste must be eradicated; only then untouchability will go." The other members did not like this so they complained about it to Gandhi.

So Gandhi released a statement regarding that. What statement? "The eradication of untouchability and caste do not have a connection; inter-dining is different and eradication of untouchability is different." This was Gandhi's statement. After hearing this, that particular committee member submitted his resignation. After this Gandhi wrote in Young India, "Caste system is a good structure; it should be there; varnashrama should be protected."

When I was in the Congress itself, I have spoken a lot that we need equality in the social sector. For all that, Gandhi said, "If the Untouchables are not allowed to draw water, build separate wells for them. If they are not allowed into temples, build separate temples for them." He also said that he would send the money. At that time, only we opposed it. I said, "If there can be no compensation for the degradation that he cannot draw water from the well, let him die without water. His disgrace should be removed, that is important and not water."

We shouted that we need temple entry. Even Rajaji knew that my screams had some respect. We entered the temples taking the Untouchables along with us.

In Kerala it became a riot and even a murder took place. Rajaji told Gandhi, "Ramasamy's speech is respected. It will become a riot. So, they must be let into the temple." Only after that, they said that the Panchama could enter as far as the Sudra was allowed.

I said, "The Sudra and the Panchama became one and we became a little lower but the Brahmin is remaining the same. Should caste not be destroyed?"

From then onwards Gandhi was involved in the effort of saving caste, and he was hand-in-glove with the Brahmins; he committed several frauds and cheated us. Gandhi's power was used to make us slaves of the Brahmins, and to make the Brahmins remain Brahmins.

People should know of the fraud that Gandhi did. That is why we say that we will burn Gandhi's picture. Gandhi's statue should not be there in our country. It must be removed.

"We are surviving in the name of Gandhi. It looks like Gandhi's honour will be spoilt." The Central Government may get the feeling of annihilating caste.

If they do not yield even to that, we will take Nehru's effigy and drag it (along the streets). If they do not take action even then, we need to finally sharpen the plan. Let ten or hundred people die. Even now they die in Mudhukalathoor. They only say 30 people die; they do not mention which caste they belong to. Cleverly they black out the news. In an honest manner, we try to act without harm. If we cannot, what do we do? In this I had great concern that we should not cause trouble, loss of property, or harm. Now there is no other way.

If we leave it as it is, will it mean that we have crossed half the well?

They asked us to have an identity in order to let others know of our caste! Sari, jewellery, and the number of piercings in the ear: seeing these one could identify the caste! Tomorrow, how long will it take for them to make this practices return.

Aha! If we burn Gandhi's photo they say a bloody riot will take place. Let it take place, what is going to be lost?

Those who are responsible do not answer. The minister didn't say that it is wrong to burn Gandhi's photo. "If it is burnt, the people's sentiments will be wounded. We will take action so that their sentiment is not affected." They said so and they did not say that it is wrong to burn his picture. I just said the mistakes that Gandhi did. Those who say that it must not be burnt, should they not say that he did not commit any mistake? Let them at least say what they should do for caste to be destroyed! Or at least let them say that caste must not be destroyed! Or let them hold a referendum on whether caste must be destroyed or not?

Comrades! We don't have any other work now. Even we feel difficult about this. Whether we lose, or die, or go to the gallows, we need to get ready for the caste annihilation war. This is not to speak in arrogance; this is a matter relating to 100% of the people. If we leave this opportunity, there is no other. So, there is no use in worrying about 10-100 people dying.

I am 79 years old. I want to die for this issue.

They say that Hindi is going to come again. "It will not come. If it comes, we will burn the Government's flag."

Only Gandhi said that the Varnashrama Dharma must be protected. In 1927 itself I wrote, "To protect the title of Mahatma, Gandhi has become a slave of the Brahmins." Even from that time I only used to write Comrade Gandhi, not Mahatma. Later we called him Gandhiar (a respectful suffix) and not Mahatma. Gandhi has done a major harm in the social set-up.

[Viduthalai, 09-10-1957]

**Gandhi and Untouchability II**
When Gandhi and the Congress spoke of the eradication of untouchability, people were fooled into believing that it was caste annihilation! They were careful to say that untouchability was wrong, but at the same time they said that varnashrama dharma, that is the caste system must be protected as it is. Accordingly they made arrangements and found protection even through the law. Gandhi did not like caste annihilation.

The eradication of untouchability is a sham, I have written so in the Kudiarasu in 1925. In other states also, Dr.Ambedkar, Shanker Nair and others have strongly condemned the slyness of Gandhi through writing and speeches.

I also said that I know Gandhi very well. I know Gandhi better than Jinnah and Ambedkar know him. I say this mainly because people should not blindly follow Gandhi. Everyone should know that Gandhi has cheated us. Only because of Gandhi all the people of good morals have been banished and people have become cheats and dishonest persons for the sake of their living. Even after 10 years of independence, the differences like Brahmin, Sudra, Paraiyar and so on exist. The Independence Day that we are going to celebrate is a celebration only for the selfish, cunning people with a comfortable living that they achieved by cheating Sudras and Paraiyars.

[Viduthalai, 13-8-1957]

**Sepoy Mutiny and Caste**

If the Tamil kings had possessed a grain of self-respect, will the people of this country have remained tribals and low caste for 2000 to 3000 years? Did any king think that including himself, 90% of the people are low caste and degraded? Is there any evidence for this?

There is history that the Buddhists and Jains who struggled to teach our people good behaviour and knowledge were tortured, impaled, and butchered by the Tamil kings which only shows that our country was ruled by barbaric tyrants.

The revolt of 1857 (*Sepoy Mutiny*: First Indian War of Independence) did not actually take place against the domination of British in the subcontinent or the rule of foreigners is clearly evident from the present activities.

The reformation carried out by Lord Dalhousie ensured that everyone in the country was able to get the benefit of education that was the monopoly of Brahmins. The introduction of railways, steam engine, telegraph and other scientific methods made the orthodox caste people angry. They feared that such type of intellectual education would certainly shake their domination. They contemplated, "What has to be done to destroy these acts and policies?" They discovered that these poor uneducated sepoys were orthodox by nature and accordingly frightened them by saying that the British are trying to eradicate caste and religion, so they made them enemies.

The reformation and spread of scientific education would change the caste system that inculcated low caste peoples' slavish mentality and respect for Brahmins. So due to their selfishness, Brahmins instigated the Sepoys by falsely alleging that animal fat was smeared as grease on the cartridges and given to them. They informed to Hindu sepoys that cow's fat was used and to Muslim sepoys that pig's fat was smeared. Thus religious superstitions were used to bring about misunderstanding and dispute.

This nation's independence lies only in setting the foolish poor against the rich; but it has no capacity to question God, Sastras, laws or the people who made some people as 'low caste' and treated them disgracefully worse than dogs and pigs.

Can a Brahmin, Sudra, Paraiyar or Chakkiliyar exist in an independent nation? Can such a nation be an independent nation or a 'hell' nation? Do we have any right to struggle and eradicate this? Is there anyway to demand our right? The human beings of this independent nation are made slaves and fools by the Brahmins who are our ministers and leaders.

So, can this be an independent nation by writing about independence? I only feel angered, so I stop writing here.

[Viduthalai, 15-8-1957]

**Burning the Constitution**

They say that if a Brahmin needs a job, the qualification is not important, only the confidence is important. If our person is going to get that job, then they say only qualification is important. You must think about all this.

What qualification? If Adi Dravidars are given 15% without looking into their qualifications, which job that they did was spoilt? Is it not mischievous to say that if they give the jobs that have not been ruined because of the Adi Dravidars to us, as per proportion, then it will be ruined?

Our first work is to burn the Constitution; after 2-3 months we will wait and see, and burn Gandhi's picture; if they don't rectify even after that we will remove Gandhi's statue, if that doesn't work out, we will drag Nehru's effigy, later (Rajendra) Prasad's (president's) effigy, if that doesn't work out and the honest methods don't work out, then we need to go to war. Let a thousand people die! Are they not dying madly in Mudhukalathoor now? Only if it goes to that extent, something might happen.

[Viduthalai, 12-10-1957]

**Burning the Constitution**

There is not a needle size of good for us in the Constitution. There are rights in it only for the north Indians to loot. It has made way for the Brahmin to forever remain a Brahmin.

So, first the Constitution must be burnt. According to this law, the Paraiyan will exist as long as this world exists. But, he wouldn't be under the name Paraiyan, but under the name Harijan. What is the difference whether we call it a broomstick or a sweeping device? Provisions have been made in the Constitution for Brahmin, Sudra and Paraiyan to exist as long as the world exists.

The Legislative Assembly does not have the power to make a law about caste annihilation. Only the Central Government has the right to amend the Constitution. If such an amendment is going to come 251 people must vote for it. If they are going to consider it, two-thirds of the members must support it. Even for consideration, this is the state. Is it possible? Only five people made the law. To make it, five people are enough. But so much of mess to change it. What can we do in this?

Now, fights take place in Mudhukalathoor in the name of caste. Some say it is due to election. Election also takes place on the basis of caste, so numerous murders occur! Many villages have been gutted down by fire. Where is the strength for Adi Dravidars? For how long would they have fought with potency? One or two of them could have beaten others. But thousands of them would have only got beaten. Let this go anyway.

They revolt saying that five people died of police firing. Let five people have been shot dead. Because of that, several lives were saved and the carnage was contained. If they cite this and say that Kamaraj must be destroyed, what does it mean? What could have been done if Rajaji had been there? What could have been done if anyone else was there?

The opposition party is there only to destroy the ruling party. Otherwise, what work of integrity does it do? How is it justifiable if they irresponsibly speak in a manner in which the carnage will be instigated again?

Comrades! This cruel state should change. Caste must be destroyed.

Next, several agitations should take place. You must support this. Under any condition, I request our people not go to Brahmin hotels.

[Viduthalai, 13-10-1957]

**Annihilation of caste**

Only the Dravidar Kazhagam says that caste should be annihilated. Even the communist doesn't ask, "Why caste? Why so many Gods? Why does God require a wife?" Because majority of the people are poor, he will get votes only if he scolds the rich. That is their doctrine.

The Congress is a big party. Even he will not ask, "Why?" In the Constitution itself they have given the promise that they will protect caste. If they don't say, how will the 'tear-drops' say? Comrades! The Dravidar Kazhagam is doing work that no one has come forward to do and where several people have been destroyed for coming forward. If we see whether caste can be annihilated through the law, through the Parliament, one has to be afraid that it is not possible.

It is possible only if we come to the decision that something has to be done. At least thousand people should die. Only if thousand people dare to die it is possible. Only then they will conform and ask, "What should be done?" We shall tell them, if you had come yesterday, two agraharams would have survived.

[Viduthalai, 13-10-1957]

**Sacrifice in struggle**

No matter how sever the drought, only the proletariats suffer for rice and gruel. The Brahmin would eat rice and ghee in a tender plantain leaf! What labour does he do? A young Brahmin boy will see an old man (of a lower caste) and (disrespectfully) call him 'Hey!' Who thought of it? If I die, there will be no one to question, "Why the Brahmins?"

Without any reason, they are dying madly in Mudhukalathoor! They can stab a Brahmin coming in front of them and question, "Why are you a Brahmin?" Let a thousand of us go to the gallows. Yesterday, how many people died in the rail accident?

The cyclone came. How many people died? Let them perish like that. For every ten Brahmins, two of our people. What is the loss? Let our children at least remain humans.

Will the Brahmin come on his own accord and say that he will go? Only you have to question him. If it happens like Mudhukalathoor in two villages it is enough. If we announce the date it is enough. He will run away! I am not addressing this at the old people. I am not addressing this at those who go to the Legislative Assembly. I am not saying this looking at those who fatten up in the name of public service. I say this looking at the proletariat and at the youth.

[Viduthalai, 15-10-1957]

**Call to Arms**

We need to do a great revolution to change this. Not a revolution through the mouth. We must take the sword and fight. In front of our eyes itself Jinnah fought. He took the sword. He asked, "Why should we, who are 90% of the population, get stuck with the Hindus?" After the river of blood flowed, they went to his house and gave him Pakistan. Today, it is a powerful country. Nehru is shivering! The reason is that they were brave. They have only one leader. They behaved as per their leader's word. If our people are also obedient, we can very soon attain our aims.

[Viduthalai, 17-10-1957]

**Struggle to annihilate caste**

They stay behind us and say Periyar, Periyar. They get on the stage, speak four words and if two people clap hands they say that they have difference of opinion with me and run away. They go to the Legislative Assembly. We are struggling incessantly in the midst of traitors. Who knows the kind of service we are doing?

As the first step, we will burn the Constitution. Then fifteen days time. If the Government does not budge, we will burn Gandhi's photo because it was he who he created a situation for the making of a law and cheated us. Then, we will wait for fifteen more days. Then we will uproot all the Gandhi statues. Shoot ten people, put them is jail, and the world will know. The world will come to know of the cruelties of caste system. Even if that does not heed results, we will drag the effigy of Nehru. If nothing happens, we have to risk our lives. If that much emotion comes, caste will certainly be annihilated. It should come to that extent.

Let them hang 100 people. In Mudhukalathoor, without any use forty people have died because of caste. They have burnt many homes. If this takes place for a useful/necessary public work, what is the loss? I ask for all your signatures. Thousands of you must give your name. If I do anything, I will inform you and only then do it. I will not do it secretly and then cry foul.

Till now, at least 30 cases would have come against me. In any one of the cases I did not scold a witness or ask a cross-question! I have never thought that I could escape by engaging an advocate. If they question me I would say, "Yes I did so." I would not say anything else. When we do anything, we do it knowingly. We do it with the knowledge that we will be punished. "Did you do?" "Yes, I did." This punishment is correct. That's all.

So comrades! Instead of considering this as a play, consider it a great work and think of your participation. Now 400 people have signed in blood—women have also signed—that they will do whatever they are asked to do. There is no need to worry about money. Now they are going to give money equivalent to my weight, so seven thousand rupees will come. Am I going to spend it on myself? If need arises, and there is no other way, that money will be used to annihilate at least ten Brahmins.

[Viduthalai, 19-10-1957]

**Brahmin hypocrisy**

The Brahmin will instigate the Paraiyar against the Padayatchi by saying, "Look, he is calling you a Paraiyan." He will also instigate the Padayatchi against the Paraiyan. When both of them are wrestling in the ground, he will pass his time without any dirt touching his nails.

In this age, why should caste be there? What is the necessity for a community to be a low caste and lead a low life? This must be pondered upon.

[Viduthalai, 16-10-1957]

**Burning the Constitution**

Let the Government say that it will not allow the caste of Brahmin to exist in the Constitution and even if it exists, we will not allow him to lead a life of a Brahmin. That is why I burnt the Constitution.

The persons who were members of the Constitution Drafting Committee are:

1. Alladi Krishna Iyer (Brahmin)
2. T.T. Krishnamachari (Brahmin)
3. N. Gopalasamy Iyengar (Brahmin)
4. K.M. Munshi (Brahmin)
5. Dr. B. R. Ambedkar (Adi Dravida)
6. Mohammed Sadullah (Muslim)

They spent 1077 days to make the Constitution. This Constitution also applies for the Sudras who are 70% of the population, there was not even a single Sudra representative involved in making the Indian Constitution. One can see that there were 4 members for 3% population. That is, 4 out of 6 members who drafted the Constitution were Brahmins.

In the Constitution, Hindu religion has been given protection. So the castes and caste differences are completely protected. We have no opportunity to correct it or amend it. I am a person who wants to abolish caste so I carry out agitations. That is why I burn the Constitution.
[Viduthalai, 11-11-1957]

**Caste discriminatory practices**
Our society is based on religion, caste, sub-castes, customs etc. Caste-Hindus, non-caste Hindus and Panchamas are the three divisions. The Brahmins will not attend the marriages of non-caste Hindus and Panchamas. What a difference exists in our society!

Nattukkottai Chettiars went to the Privy Council to establish that they cannot be called Sudras and that they are also upper castes. Sudras cannot have marriage rituals or property or prayers for the dead. That is why if they have to get married or offer prayers to god, the officiating Brahmin priests puts a holy thread on the Sudra man. This holy thread worn by the Sudra is thrown into a river after the ritual is over.
[Vazhkai Thunainalam, 1958]

**Vaikom agitation**

I was extended an invitation to come over to the Kanyakumari District many times by your comrades. As I was busily touring in the other districts, I could not come earlier. Wherever I toured, I found a great awakening. People gathered in thousands.

Ten years ago, I addressed a meeting here in Marthandam. In those days you were citizens of a Native State. You were ruled by a Rajah whereas we were citizens of the British Government. Yet we are all "Sudras." We Dravidians were subjected to humiliation. On the outcome of the hoax played on us we continue to be "Sudras."

Today, we are citizens of one country. We are Tamilians of Tamilnadu. We are today brought together. Our unity is strengthened. We are today linked as one family because we are all, now-citizens of one country. We have to work together for achieving our ideals as we are all classified under one caste. So far as I am concerned even before 35 years ago, I led the agitation in Tamilnad to eradicate the social evils particularly the hateful 'untouchability'. For over thousand years we were not allowed to enter some of the public roads. Those who are now aged at least 50 years may recollect those days. Youngsters of this generation may not know these things of the past.

If there had not been the agitation in those days, today many of us would not have the right to pass through many of the roads. In those days conditions were very bad in this country. The Government was in the hands of the orthodox Brahmins. The Varnashrama Dharma was in its full sway. In our country, the advent of the Non-Brahmin Movement, redeemed a number of rights to the Non-Brahmins. The Non-Brahmin Movement successfully combated the Brahmin domination. That Non-Brahmin Movement was popularly known as the Justice Party; named after its journal 'Justice'.

The Brahmins too had their own organisations as Brahmana Samajam, and Brahmana Mahasabha. They worked against our interests and stood in the way of our attaining many of our legitimate rights. Brahmins were proud of calling themselves as 'High caste'. They insisted on calling themselves as 'Brahmins'. Manu law and other Sastras too termed us as 'Sudras' only. What an amount of humiliation and degradation we were subjected to!

As this state of affairs affected name for Dravidar Kazhagam or Tamilar Kazhagam, we have to choose only 'Sudra' Kazhagam, as the suitable name for our organization.

That is why, we had to change the name of the South Indian Liberal Federation and the Justice Party as the Dravidar Kazhagam, to make ourselves known to the world as to who we are. The Dravidians are a proud nation, known to the world.

On account of the efforts taken by the Non-Brahmin Movement (Justice Party) in the years 1919 and 1920 and the agitations in my Tamilnadu, the right to make use of all roads irrespective of castes, was got, not only in Tamilnadu but also in Andhra, Karnataka and Kerala.

With the powers vested in the hands of the Justice Party the right to make use of all roads by all castes was brought into practice. The Justice Party brought in a legislation even in those days permitting the so called low-caste people to make use of water from the wells, which had all along been exclusively reserved only for the use of the Brahmins.

These are all things which took place before the days of Gandhi. It is absurd and. fraudulent to say that is only Gandhi who did all these things. Not only this. The Justice Party people were the first people to pave the way for all castes to enter the Panchayats, Municipalities. Taluk Boards, District Boards and Legislative bodies even before the coming of Gandhi. They nominated the representatives of the so called low caste people as the 'Pariahs' to enter the legislatures on a par with the so called high caste 'Pariahs' 'Chakkilis', 'Pallars', were members of the legislature. I want you all to clearly understand this truth.

As a matter of fact, Gandhi's plans were different. He was not for allowing all the Sudras and Untouchables to bale out water from wells and tanks, along with the high caste Brahmins. He was not for permitting the Untouchables to enter the temples along with the high caste people. Originally, he only insisted on the continuance of certain rights exclusively for the high caste. He upheld the Manu code. He was for separate temples, tanks, wells and dwellings for the high caste Brahmins and the low-caste Sudras. That was the original plan of Gandhi. I know it. Let anyone deny. Today, false propaganda is carried on about Gandhi. Much is said about the Gandhian way and Gandhian path.

I was the Secretary of the Tamil Nadu Congress Committee. A sum of Rs.48,000 was sent to Tamilnadu as grant from the AICC, to construct separate schools and temples for the low-caste Sudras as the 'Pariahs', Chakkilis, and Pallars. It was strictly ordered that these Untouchables should not go and create trouble at the places exclusively used by the high caste Hindus.

By that time the Justice Party people had already passed orders permitting all castes to study in all schools irrespective of their castes. They made all study together. Caste restrictions were removed in the field of education long ago. This Reform was strictly enforced. There was a law compelling even private institutions to admit a percentage of 'Sudras' in their schools failing which, such schools would not be eligible for any Government grant.

At the time of inspection the officer would put the question, "How many Untouchables are studying in this institution?" If anyone should say that no untouchable approached seeking admission, the officer would say 'In that case you have got to go and get some Untouchables to study in your School'. I am telling of the conditions that prevailed in our State even before Gandhi came.

When Tamilians were so much progressive, in your Kanyakumari District things were very bad. The high caste Hindus did not even tolerate the right of low caste untouchable Hindus. Even his shadow should not fall on the so-called high castes. That was the horrible tragedy in your place. The low caste Sudras were to raise a cry from his hiding place to reveal his presence. It is on account of the laudable services of Swami Narayana Guru that the low caste Sudras were awakened. The Vaikom Agitation changed the conditions. The Untouchables here gained a lot. These things may not be known to the youngsters here.

We waged the battle at Vaikom against untouchability. We were imprisoned many times. We were severely beaten. We were put to disgrace. All these sacrifices we had to make to eradicate untouchability.

There were no classes in jail in those days. The treatment was very bad. We had to bear all these to eradicate untouchability and bring in a new change. How did this change come in? What is our present position? If you think over and search for a better position, you will in fact agree that we are very slow in eradicating casteism and its evils. We must muster strength and march on with a greater speed.

You must know the history of the Vaikom agitation. A very small incident led to the Vaikom agitation.

Comrade Madhavan was an advocate. He was to appear before the honourable judge in a case on behalf of his client. The court was in the compound of the Maharaja's palace. At that time, arrangements were made to celebrate the birthday of Rajah. The entire surrounding of the palace was thatched with Palmyrah leaves beautifully. Brahmins started chanting mantras. As Comrade Madhavan belonged to the 'Ezhava' (Nadar) community, he was not permitted to enter or pass through the place and reach the court.

It was at this time that the Justice Party was carrying on propaganda in Tamil Nadu for the abolition of castes and untouchability Inter-caste marriages were encouraged. Schools were thrown open to all. 'Samabandi Bojanam' (Inter-dining) was popular. Such intensive social reform propaganda was carried on in Tamil Nadu by the Justice Party. When Gandhi came to know of what the Justice Party was doing in Tamilnadu he started including our schemes in his constructive programme. In those days, the Justicites boldly exposed the Brahmins who were afraid to move about without company. Non-Brahmin leaders as Dr. T.M.Nair and Sir P.Theagarayar educated the masses by their incessant and extensive propaganda and secured the powers in the State. Brahmins were envious of the Justice Government. They had no platform in those days. In those days, the Brahmins cunningly took shelter under the slogans – 'we are not power-mongers. We boycott the elections.' With such false slogans they hoodwinked the people and indulged in all sorts of intrigues. Realizing the popularity of the Justice Party, Gandhi concentrated on the problem of untouchability, as the only way to bring down the Justice Party rule in Tamil Nadu.

In those days I was very familiar to the leaders of the Justice Party. They had great respect for me, because I held many posts. Mr. Rajagopalachari met me and induced me to become a follower of Gandhi. He said that Gandhi alone is capable of carrying out the much-needed social reforms. I resigned the post as chairman of Erode Municipality and joined the Congress. Before my entering the Congress, no Tamilian had the honour of becoming the Secretary or President of the Tamil Nadu Congress. I was the first Tamilian to hold these posts in the history of Tamilnadu Congress.

Com. T.V. Kalyanasundaram (Thiru Vi.Ka) was a schoolmaster. Dr. R Varadarajalu (Naidu) was the editor of 'Prapancha Mitran'. Yet Brahmins did not trust him. Com. V.0. Chidambaram (Pillai) was at the mercy of Mr. Kasturi Ranga Iyengar, after draining all his resources.

I already held big posts and hailed from a very big business community. For anything and everything Rajagopalachari believed me and reposed much confidence in me. I also believed him and reciprocated the confidence he held in me. We both worked out together and I carried on an intensive propaganda; with the result the Brahmins once again gained the platform. I was very bold in expressing the rationalist views. I openly spoke on god in all my meetings stating, "If the idol would get polluted by the touch of the people, such a god is not required and the idol has to be broken to pieces and used for constructing good roads. Otherwise it may be put near the river banks to be used for washing clothes." I was often induced to speak severely by the Brahmins. As I was not for any post or power the Brahmins remained silent.

What all I say now about God, religion and caste; I used to say even in those days. Rajagopalachari used to tell me that I have administered a very strong dose, after hearing my speeches. I used to reply that so long as people remain foolish, there is no use of our giving a light dose. On hearing me he simply smiled. In these ways, we made the Brahmins come to power in those days.

Leaders of the Ezhava community wanted to start an agitation, when Mr. Madhavan was objected from entering the court. Mr. K.P.Kesava Menon, President of the Kerala Congress Committee, T.K.Madhavan and others took the lead. They decided to launch the protest on the day of the prayers at the Raja's palace. They chose Vaikom as the fit place for the agitation. It is only in Vaikom you have a temple with four entrances on the four sides, leading to four streets around the temple. That was really convenient for the agitation. So they chose Vaikom as the place for the agitation.

There was a law that -the low caste Untouchables as 'Avarnasthans' and 'Ayithak Karans' should not enter those roads. If an untouchable has to go to the other side of the temple, he had to go two or three furlongs away from the temple and walk about a mile to reach the other side. Even the 'Asaris', Vaniars and Weavers were not allowed to enter the roads around this temple. The same conditions were prevalent in other temples also. Particularly at Suchindram it was very strictly enforced. Important government offices, courts and police station were by the side of the Vaikom temple near the entrance. Even while transferring officials, no untouchable would be transferred to these offices as they were not permitted to enter the roads around the temple. Even the coolies were debarred from entering the roads to reach the shops.

As soon as the Vaikom agitation started the Rajah ordered about 19 leaders including Advocate Madhavan, Barrister Kesava Menon, T.K.Madhavan, George Joseph, to be arrested. They were treated as special prisoners. At that time, there was a European named Mr. Pitt as the I.G. of Police, under the Rajaji's government. He ably managed the affairs of the agitators. When all the 19 agitators were put in prison, the Vaikom agitation actually fizzled out. At that time I received a letter from Mr. Kesava Menon and Barrister George Joseph.

"You must come and give life to the agitation. Otherwise we will have no other way except to tender apology to the Rajah. In that case, we will not lose anything, but, a noble cause would suffer. That is what actually worries us. So please come immediately and take up the agitation." This was what was written in the letter. They themselves chose me and wrote the letter to me because I was very vociferous in attacking the evil practice of untouchability in those days. Moreover, I have established a good name not only as a fierce propagandist but also as a successful agitator. When they sent the letter I was on tour. The letter was redirected to me from Erode and it reached me at Pannapuram in the Madurai District. On receipt of the letter I cancelled my tour and rushed up to Erode to proceed to Vaikom. I wrote a letter to Rajagopalachari requesting him to act in my place as the president of Tamilnadu Congress Committee. I pointed out the importance of the Vaikom agitation in my letter. It was a good opportunity for me. So, I did not like to miss it. I proceeded to Vaikom with two others.

Somehow the news spread that I am coming to Vaikom agitation to lead. When I reached Vaikom by boat, the police Commissioner and Tahslidar greeted us.

We were informed that the Rajah instructed them to receive us and make all arrangements for our stay, I was really very much surprised. The Rajah was so good to me because, he used to stay in our bungalow at Erode, while his officials stayed in our choultry on his way to Delhi. The Rajah and his party were always afforded a cordial treatment at Erode during their stay, before boarding the train for Delhi. That might be the real cause for the unusual treatment given to me at Vaikom. When the people of Vaikom came to know of my relations with the Rajah and the officials, they were all very happy. Even though Rajah treated me as a guest, I participated in a number of meetings supporting the Vaikom agitation. I criticized the evil practice. I said the god does not deserve to be in the temple at Vaikom, if it feels that by the touch of the Untouchables, the deity would get polluted. Such an idol should be removed to be used for washing clothes. By my propaganda more and more people were eager to join the agitation. More and more people came forward from different places. It became a problem for the Rajah. Yet he remained silent for five or six days. Many complained to him, about my speeches. Rajah could not ignore any further. So, after ten days, he permitted the police official to promulgate P.C.26 which is similar to that of Sec. 144 here.

There was no other go for me than to defy the ban. Accordingly I defied the ban and addressed a meeting and I was arrested. Mr. Ayyamuthu also defied the ban along with me. He was also arrested. We were all sentenced to undergo rigorous imprisonment for one month. I was put in Aruvikkutha jail. It was after my imprisonment, my wife Nagammai and my sister S.R.Kannammal and a few others carried on a statewide propaganda. When I was released, I once again resorted to the agitation.

When I was in prison the agitation gained momentum. Lot of people volunteered to court imprisonment. Intensive propaganda induced the people to encourage the Vaikom agitation. The enemies indulged in hooliganism. Rowdy elements tried their best to create panic and ended in failure. Even those who were in foreign countries came to know of the atrocities perpetrated in the name of caste here. They volunteered to send donations. Daily, money orders poured in. A big pandal was raised to house the volunteers. Daily, more than 300 people were provided with food. Many planters sent vegetables and coconuts daily. They were pooled together, as small hillocks. It looked like a marriage house.

At that time, Mr. Rajagopalachari wrote a letter to me. Why should you leave our country and create trouble in another country? It is wrong on your part to do so. Please leave it and come over here to take up the charge from me. This was what was contained in the letter. Mr. Srinivasa Iyengar came over from Tamilnad to meet me. He also tendered the same advice as Rajagopalachari. By that time more than 1000 people were ready .to participate in the Vaikom agitation. There were big processions and 'Bhajans' daily everywhere. The agitation gained momentum.

The news reached Punjab. There Sami Sirathananda. made an appeal. He sent about thirty Punjabis to Vaikom. They offered 2000 rupees as donation and consented to meet the catering expenses for the volunteers. Seeing this, the Brahmins here sent communication to Gandhi. They accused the Sikhs of provoking a war against Hinduism. Gandhi expressed his view. He said that Muslims, Christians, Sikhs and others who were non-Hindus should not partake in the agitation. In response to his appeal, Muslims, Christians and Sikhs withdrew from the agitation. Rajagopalachari wrote another letter to Joseph George stating that it is wrong on his part to interfere with matters pertaining to Hinduism. But, Joseph George did not pay heed to Rajagopalachari's advice. He replied that he was prepared to face expulsion from Congress. He strongly stated that he would not lose self-respect. Mr. Sen, Dr. M.E.Naidu and other leaders stood strongly in support of the agitation. But some people were afraid that Gandhi would write, condemning the agitation and stop the donations. But at that time Swamy Sirathananda came to Vaikom and assured of financial support.

The Vaikom agitation was launched in spite of opposition by Gandhi. I was once again arrested and sentenced for 6 months imprisonment. Some Namboodri Brahmins and orthodox Hindus joined together and planned to counter the Vaikom agitation through what was called as 'Satru Samhara Yagna' (Bonfire prayers for killing the enemies). They spent money like water and performed this pooja. I heard about this in prison. All on a sudden, one night I heard the sound of gunshots. I enquired the warden, whether there was any festival going on near the prison. He told me that the Raja passed away and the gunshots are fired to indicate the loss. When I learnt that the Rajah is dead, I had a melancholy feeling. Later I was glad because the prayers by the Brahmins and orthodox Hindus to destroy their enemy had resulted in the death of the Maharaja. Their prayers did not harm the Vaikom agitators. The people were also happy. Subsequently we were all released on the ceremony day of Maharaja. Our enemies also lost their tone and tenor.

Later, the Maharani wanted to settle the problem by mutual talks. She wanted to discuss the problem with me. But the Dewan of the State, who was a Brahmin stood in the way of our talks and said that the Maharani should not talk to me directly. So he wrote a letter to Rajagopalachari. Rajaji knew that I would carry the laurels and earn the credit. So he cunningly decided to make the Maharani talk to Gandhi. It is because of this trick played by Rajaji, that Gandhi's name was dragged in the history of Vaikom agitation. I did not much mind as to who personally gained the name and fame. I was not for any personal glorification. I only wanted the problem to be solved successfully.

Gandhi came and had talks with the Maharani. Maharani consented to throw open all the roads for the low-caste Sudras and untouchables. But, she expressed the fear that I would further continue the struggle seeking the right for Untouchables to enter the temple. Gandhi came to the Tourist Bungalow where I was staying and asked me to express my opinion. I said, "It is not a big thing to enable the Untouchables make use of the public roads! Even though temple entry is not at present one of the ideals of Congress, so far as I am concerned it is one of my main ideals. But, you may inform the Maharani, that for the present I have no such idea to launch the campaign for temple entry rights. Let things come to normal before I decide what to do."

Gandhi informed the Rani and she proclaimed the right for all to make use of all the roads. This is how the low caste Sudras and Untouchables got the right to use all roads, like the high caste Brahmins and orthodox Hindus.

I was for some time the Chairman of the Devasthana committee at Erode. When I was away, comrades Gurusamy and Ponnambalam and Eswaran induced two Adi-Dravidar workers in my office, to wear the sacred ashes (Vibooti) on their heads and took them inside the temple. Seeing them, the Brahmins cried loudly that they have polluted the deity. They were locked up inside and a case was instituted. They were punished In the District court. But on hearing the appeal, the High Court released them as not guilty. That was in the time of the British rule.

But, it was only at Susindram that the first agitation for temple entry right was launched publicly. A Self-Respect Conference was also conducted under my presidentship. A number of resolutions were passed urging the abolition of castes and ensuring rights for Untouchables to enter the temples.

Next, another Conference was held at Ernakulam. In that Conference a resolution was passed condemning castes and urging Hindus to become Muslims, as there are no castes in Islam. Some others recommended Christianity through amendments. At last, the option was given to join any one of the religions.

The same day about 50 Hindus joined Islam. This trend started even outside and it terrorised the orthodox Hindus and Brahmins.

One day, in Allepey one of the converts to Islam (who was a Pulayar by caste) went to purchase an article in a Nair shop. He was beaten up and it developed into a serious clash between Hindus and Muslims. Such clashes between Hindus and Muslims spread everywhere. The then Dewan, a Brahmin Sir C. R. Ramasamy Iyer put down the agitation with an iron hand. But the Rajah was later informed that most of the low caste Untouchables as, 'Ezhavars' and 'Pulayars' are becoming Muslims. He was also advised that there is no other go than to throw open all the temples for all the Untouchables to save Hinduism from the peril. At that time there was 'Yagna' by Brahmins for his long life. Moreover it was a custom to say something good to the people on his birthday. It is then that the Rajah chose it the right thing to do at the right time. He announced that on his birthday all temples would be thrown open to all, including the low caste Hindus and Untouchables. This is the history of the struggle. This is how the Untouchables were given the rights to enter the temples.

It is only after all these that Rajagopalachari and Gandhi came forward to speak of temple entry. It is quite absurd to say that these changes took place because of Gandhi. As a matter of fact Gandhi has not done even a molecule of service to the Untouchables. This you will understand well by reading 'What Congress and Gandhi have done to Untouchables', a book written by Dr. B.R.Ambedkar. When I was the Secretary of the Tamilnadu Congress Committee, out of party funds a Gurukulam was run at Cheranmadevi. As Secretary I consented to give Rs. 10,000 and paid a part sum of Rs.5,000. One Brahmin named V.V.S. Iyer took up the responsibility to run the institution. In that Gurukulam, Brahmin boys were given special preferential treatment. They were fed separately. The non-Brahmin boys in the Gurukulam were fed outside. When 'Uppuma' was given to Brahmin boys, only gruel was poured to Non-Brahmin boys. The son of Omandur Ramasamy Reddiar told the matters to me in tears. I complained to Rajagopalachari. When he contacted V.V.S. Iyer, he neither denied the charges nor repented. He stoutly refused to give equal treatment to all. He said that he cannot do anything as it was an orthodox area. Then I said that I can give the balance amount of Rs.5,000 only when the Gurukulam is reformed. He got wild. He curtly asked me 'Are you serving the nation?'' I curtly replied, ''Is what you are doing nationalism?'' This serious matter provoked me to start a party for the non-Brahmins (Tamils).

Even now you can find only 'Brahmins' preparing food at Congress congregations. Even in those days, we engaged the Virudunagar Nadars to prepare food at the Justice Party Conferences and the Self-Respect Conferences.

Why do I recollect all these past things? You should know that unless we go on agitating like this we cannot make the society progressive. Moreover you should all know that neither Congress nor Gandhi is responsible for any of the social reforms, we have been able to witness.

Even today, we are the only people who boldly ask why should the lazy Brahmins be considered as high caste, while the real tillers and toilers are considered as low caste. Why should we have a god that degrades us as Sudras?

Today they have created all safeguards to casteism in the Constitution. A Brahmin from somewhere has the boldness to come over here, and speak with temerity, tendering serious warnings. Why? Power is vested in his hands.

They tell us to meekly submit as Sudras forever. They show the prison and terrorize us.

Did any one have the guts to question?

We are the only people who are free, frank and fearless.

If we are termed as Sudras by Hinduism, what else can we do but, destroy that Hindu religion? Our D.K. is not a political organization. We do not contest elections. We do not seek votes. We are not for power. Others may hesitate to call a spade a spade. Power seekers may coax the innocent voters. They may hoodwink you for selfish gains. I am not for dragging in Gandhi's name to dupe you all for getting any post or power. I am not for that disgraceful life.

We have not made public life a profession or business to eke out our livelihood. Think over Why? We eat our own food, spend time, and bestow our energies to instill in you self-respect?

By 1938 you find all over the world wisdom having sway. But still here we are like barbarians. Our god, religion and Sastras do not elevate us from the plightful rot.

Government is also in the hands of barbarians. No one dares to question except ourselves.

We are made to be sons of prostitutes by Brahmins. Why should our sons be called sons of prostitutes. No one thinks of this disgrace. Those who survive in politics do not care about it. They implicitly obey and submit to what all Brahmins say.

When I was leading the Vaikom agitation. Mr. Sathukkutti the son of Nilamban Zamindar used to meet me often and discuss. He used to address me as 'Naicker samy'. Not only that. He used to talk high of his birth because he was born to a Namboodri Brahmin. He would often tell me that I should not consider him as one born to a Nayar. Yet he was a graduate B.A. Who is there to condemn this mentality in our people?

Consider for a moment what these 'Azhwars' have done. They attained 'Moksha' by prostituting their wives. This is revealed in the 'Purana' - Baktha Vijayam.

One 'Sudra'' an Azhwar, gained a place in heaven by allowing his wife to lead the life of a prostitute. The Nayanmars gave their wives to Brahmins. Even to this day the orthodox people propagate these things without shame or self-respect. When I point out these things, I am accused of talking damagingly of puranas. Who else talks boldly of these? These puranas ruin our morality. What else can we say?

Added to all these the Brahmins have stuck to the seats of government. Power has been transferred into the hands of Brahmins. I blame Gandhi for that? A big conspiracy was hatched to keep us eternally as 'Sudras'. Today everything is in their hands. Today the President is a Brahmin. The Vice President is a Brahmin. Premier is a Brahmin Deputy Premier is also a Brahmin. The speaker of the Parliament is also a Brahmin. Added to all these if we plead for the eradication of castes, they send the accused to be in jail for a period of three years. Who worries about all these? Most of the luminaries in public life want to safeguard the government, casteism, Sastras, puranas, religion and god. They think that there is no other go for them to live.

No one who depends on votes and bribe will question the atrocities in the name of caste, god, religion and government.

The British at least considered us as men with equal right. Today the government is in the hands of Brahmins, who call us as sons of prostitutes. That is why they have easily found safeguards in the Constitution itself. According to law those who demand castes to be abolished have to be prepared to undergo imprisonment for three years.

This casteism is a chronic disease that has eaten our society for centuries. What medicine we use for scabies and itches cannot cure cancer. We have to operate the body and remove the portion affected by cancer. Treatment will be different for different diseases. According to Hindu Law we are 'Sudras' for more than 3000 years. We are sons of prostitutes for over 3000 years. Our constitution gives full protection to this evil.

We must root out this evil. We must get rid of this ridicule. It is indeed the most difficult task. Unless you pour boiling water on the roots it will not die. Unless we take severe steps we cannot eradicate castes.

Not only in Tamil Nadu but even in the whole of India there is no force that can raise a bold voice like us. Those who are after power will never dream of raising any protest. It is only those who are dedicated to serve the people sincerely and selflessly, can dare to risk their life even for eradicating the castes. What have those who entered the legislature done so far? They cannot do anything? We can get posts by simply sending a telegram. Yet we are not prepared.

A few days ago Nehru made a sickening note about the legislatures and other elected bodies. He even threatened that he would retire and seek renunciation. What happened? He silently gulped all his remarks and stuck to power. It is all mere display of the old Gandhian tactics to gain popularity. The D.M.K. people who were with us, condemned the entry in the legislature, so long as they were here. They even wrote attacking the elected representatives and the bodies. Nay, even Nehru and Rajendra Prasad spoke against the legislature. Because they are now aware of the chances, they are quite eager to enter. They forget their past. By hook or crook they want to adorn the legislatures. They are prepared to sabotage and expose anybody. Somehow everyone wants to come up in life. No one is worried about the age long insults heaped on our nationality.

The entire country is in the clutches of three ghosts and five diseases. To believe in one thing that is not really existing is what is defined as Ghost.

God - Caste - Democracy are the three Ghosts.

Brahmins - Newspaper - Political parties - Legislatures - and Cinema are the five diseases. These diseases are preying on the human body like the diseases Cancer, Leprosy, and Malaria. If the society is to progress these things should be fought tooth and nail and destroyed completely.

[Viduthalai, 8 and 9-1-1959]

*This essay's English translation has been taken from the book Collected Works of Periyar E.V.R., published by Dravidar Kazhagam, Chennai, 2005 (Third edition).*

**Agitation against Hinduism**

To destroy and demolish the Hindu religion, one method is through law and the other is through agitation. The nine points for agitation are

1. All backward class people and Untouchables should never go to temples.
2. No one should worship Hindu gods.
3. One should not celebrate Hindu festivals.
4. No one should put any mark on the forehead.
5. No man should have a tuft.
6. No rituals should be performed after death, after birth, before birth and so on.
7. The Brahmin should not be called to perform any ritual.
8. Photos of Hindu gods should not be hung in ones home.

9. No one should go to shops or hotels run by the Brahmins.

[Viduthalai, 16-2-1959]

**Division of labour**

A labourer's son should not be a labourer. He should be a proprietor. Such stigma must vanish. No one can force a scavenger to continue his scavenging job. Such caste and professional atrocities should be banished.

[Viduthalai, 18-4-1959]

**Evils of society**

Today's Tamil society lives like barbarians. Out of the 280 crore people, only the Tamils live like fools in such a degraded status.

The 3 devils that have captured us are God-religion-Sastras, Caste and Democracy. The diseases that damage us are

1. Brahmin
2. Dailies (magazines) / Newspapers
3. Political parties
4. Elections and
5. Cinemas

Our nation will never improve; I have spoken this in several of meetings in Kerala, northern states and Karnataka.

Likewise, the three things that stand as obstructions to our development are

1. To lead our life as said by our ancestors
2. Lead our life as written by ancestors
3. Lead our life as great men

These three stumbling blocks have existed for the past 2500 years. Buddha said so. Now the Buddha of the 20$^{th}$ century, who wears black shirts, says so.

In a book "*Vazira susy pramanam*" the following is described about the birth of saints:

Kalikotu was born in the stomach of a deer, Jambukar was born from the stomach of a fox, Goutama was born of a cow, Valmiki was born to a hunter, Agastiyar was born from a pot, Vyasar was born to a fisherwoman, Vasishtar was born to a prostitute Voovasi, Koulatya was born to a widow, Narada was born to a washerwoman, Mathangar was born to a cobbler, Mandaliar born to a frog, Chasangyar was born to a untouchable, Kangayer was born to a donkey and Swanakar was born to a dog.

Thus our ancestors were born from fox, dog, donkey, pot, frog, cow etc. How can we take their words as law? Will not the foreigners who land on the moon think that we are barbarians? Though we talk of such ancestors we have not become scientific to make even a common safety pin.
[Viduthalai, 22-5-1959]

**Burning the Constitution**
To annihilate caste they burnt the Constitution. 3000 to 3500 people went to prison and 15 persons died. It is barbarous to have caste.
[Viduthalai, 5-3-1960]

**Communal representation**
To have a proper and justifiable proportion for all castes, the government should implement communal representation in education. The same should be followed in employment, postings and in ruling the country. The Brahmins who are just 3% of the population are protesting only against this.
 [Viduthalai, 13-12-1960]

**Protectors of casteism**
All kings protected caste. They demanded Swaraj to protect Manu dharma. We do not have self-respect or intellect. We sell it to satisfy our stomach. Caste will not go as long as the following five exists.

| | | |
|---|---|---|
| Our god | - | Protector of caste |
| Our religion | - | Protector casteist religion |
| Our government | - | Protector of casteist government |
| Our literature | - | Protector of casteist literature |
| Our language | - | Protector of casteist language |

What has our language done to destroy or demolish caste? What has our literature done to destroy or demolish caste? They called rationalism as atheism. So, no one had the guts to use his intellect.
[Annihilation of caste, 1961]

**Religion and casteism**
God is his (the Brahmin's) instrument to protect caste. He created religion. If we want to do any good act, first we should make a bonfire of religion, god, sastra, purana and *ithihasas* (epics). The gods killed those who acted against religion. Vishnu killed everybody who tried to destroy religion. One needs courage to destroy and annihilate religion. A Brahmin writes, "A Brahmin is equivalent to god and Sudra is equivalent to faeces."
[Annihilation of caste, 1961]

**Class and caste differences**
Can god watch a starving labourer continue in starvation and a Brahmin lead a life of luxury? If god exists, will he watch this silently? Our main duty is to educate the present generation and see that they do not take up the occupation of their father. If one asks a Chakkiliyar the cost of a slipper and requests him to lessen it, he would plead saying the cost of living is high, but on the other hand if we ask a owner or a capitalist to lessen the cost of goods he will say, "Buy it if you want, otherwise get away."
[Viduthalai, 9-5-1961]

### Hinduism

Once the Hindu religion is abolished, social equality will prevail, unity would be there among Indians, they will develop in science and technology and would also be economically well off. I am dead against all blind rituals and superstitions.

[Viduthalai, 23-11-1966]

### Rajaji, Brahmins and education

Rajaji, who came to power in 1938, was trying to put obstacles even for school education. So he declared the state a dry area (prohibition). So, there was a raise in the percentage of literates from 5% to 7%.

Infuriated by this rise in the number of literates, Rajaji closed down 2600 schools in Tamil Nadu stating that there were no funds to manage these schools because of the drastic setback due to the absence of revenue from liquor sales.

This shows the motives of Brahmins who were very much against our people getting educated even at the primary level. Thus, the only motivation of Brahmins is to see that our people are spoilt and lead a very low life. They would never show us any good path because they are least bothered about the country.

[Viduthalai, 4-3-1967]

### Caste in post-Independence India

I feel very sad that even after 16 years of independence, our nation has Brahmins and Paraiyars.

[Viduthalai, 11-3-1967]

### Communal representation

I am old at the age of 88 years and I cannot hope to live longer. I alone know the state of my health. That is why I am keen on getting communal representation.

[Viduthalai, 18-3-1967]

### Discrimination in temples

It is better if religious leaders keep their mouths shut, because if they open their mouth, it stinks. Does it imply that a believer in god should not have any sense of shame, dignity, honour or intellect? If he visits temples, they say he should not enter the temple. They say you are a Sudra or untouchable, and tell you to stand out. Why do such statements not hurt these devotees? How shameless and hard-hearted are these Brahmin people? Even if a devotee goes to temple without any sense does he not have shame?

[Viduthalai, 27-3-1967]

### Caste and rule

In a just, good rule there should not be any difference between man and man. When man is born, he is not born with caste, riches, poverty or outcasteness. If this is the case, how can low caste, Sudra, upper caste and rich exist among human beings? If the rule is for equality, how can we have high and low? If we have to find a rule that has no difference among people, it can be achieved only by killing the kings, gods and the middlemen.

[Viduthalai, 8-2-1969]

### Communal representation

I have been labouring for the Backward Classes and the Depressed Classes. I want communal representation to be given in ministerial berths and so on.

[Kudiarasu, 5-3-1969]

### Abolition of untouchability

It is said that untouchability has been legally abolished in the Constitution. Practicing untouchability is unlawful. But it has remained just a conditional implementation because it is legally punishable only as long as it does not affect the religious sentiments. That is, with the only exception of Brahmins, untouchability has been created among all class of people in order to ingrain the concept of lowness. This is the root of the philosophy of untouchability as inculcated by the Varnashrama dharma.

If you touch god, it becomes unclean i.e., untouchable (pollution). If you come to the sanctum sanctorum it is unclean (pollution). This untouchability will not be abolished; the disgrace of 97% cannot be abolished unless god is abolished. Can untouchability exist if Muslim rule comes to the nation? To think that this sort of rule is permanent is nothing but foolish, devoid of self-respect and full of superstitions. In India, Indians are divided caste-wise as Brahmins, Sudras and Panchamas. This division finds its place in law and religions. Ex-Prime Minister Nehru is a Brahmin, Gandhi is a Sudra and Dr.Ambedkar is a Panchama. Because of this division, majority of the people suffer several discriminations and injustice, which hinders their development. The terms Sudra and Panchama are there not only in the Vedas and other political Constitutions but these caste divisions are kept intact, they are made to suffer several unlawful cruelties and injustice is done to them. Thus, they should look into this structure.

Can we say with real knowledge and evidence that true democracy and independence exists in India? In independent India, 90% of the people are slaves. To be more precise, 90% of them are either coolies or salaried people who have to serve obediently to their master.

The democracy in India is that 97% of the population is low caste, degraded caste, $4^{th}$ caste (Sudras), $5^{th}$ caste (Untouchables) and Chandalas both legally and as given by the Sastras. 3% of the people who are Brahmins are upper castes equal to the caste of god. They say it is god's will and rule of the Sastras. So democracy also is only a falsehood and fraud.

[Viduthalai, 27-3-1969]

**Destruction of religion**

We Tamils have no religion. We think that Hinduism is our religion. This is a Himalayan blunder we make out of our foolishness. What is Hindu religion? What is its meaning? Does there exist any evidence for the term Hindu religion as in the case of Christianity and Islam? Does there exist any evidence for the Hindu Sastra and Veda? They say that Hindu religion is a Brahminical, Vedic religion! Brahmins call themselves Aryans and they call it the Aryan religion. In English dictionaries we see that Hinduism means it is a Brahmin religion, or a religion that is not Islam or Christianity.

All the foundations of Hinduism make us the $4^{th}$ and $5^{th}$ castes! It has disgraced us in all ways. So it is pertinent that we demolish god, religion, Sastras, puranas and ithihasas.
[Viduthalai, 17-12-1969]

**Religion and caste**
The gods propagated by these cheats, their houses, their food, their wives, their concubines, their property, their marriage, etc. and the lakhs of rupees spent on them and the people who participate on many days is only to make crores of peoples low caste and untouchable.

For instance, a man claims himself to be superior to other men by birth and so establishes a god or many gods for this act and says that only these gods have made him a high caste. This means that people who are not high caste become low caste and for this god is the basis.

So such gods should be abolished, destroyed, and those who do not come to this task of destroying and shattering gods cannot be human or possess intellect or honour. He cannot have the status to call himself a man.
[God is only imagination, 1971]

**God as basis for caste**
The upper castes enjoy life and lead a comfortable life. The men who protect god are of the low caste and undergo all problems in their life for mere sustenance. So is it not just and proper to destroy and demolish god?

Only because of god, the development of the world, social sense, morality and above all humanity have been greatly affected. Even if they realize this, no one is willing to openly express it. So, I have felt it essential to express this. The sooner man rejects and forgets god it is certain that he would become socially improved.

[God is only imagination, 1971]

**Hinduism and untouchability**

The minute you shake the concept of untouchability, the Hindu religion is shaken. That is why Brahmins use false propaganda to call us traitors. They say that religion is destroyed by our talk.

[God is only imagination, 1971]

**Partition of Tamil Nadu**

I say one thing and I say it with firmness. If Tamil Nadu had been partitioned after the British went, firstly, Tamilian social disgrace would have been wiped out. Secondly, due to his rationalism, the Tamilian would have reached great heights compared to the barbaric position, degraded status and low caste which he suffers in his own motherland today. His base attitude, of selling anything for the sake of livelihood, position, and existence, would have disappeared.

Further, the practice of low and high in birth, qualification and talent, limitless abundant rights for a very small crowd and the situation of a big or a very big crowd existing by shamelessly doing disgraceful jobs would have been wiped out. Today, the Indian government's important policy is to save caste, to ensure that upper caste people remain in high posts, save the traditions, customs; apart from that it does not have any particular good.

Why should the hereditary kings who ruled this land have to become degraded and live as sons of prostitutes?

Our people toil as builders of buildings, cultivators of lands, beautifiers of the nation; in spite of doing so much labour what is the justice in keeping us (Dravidians) as degraded people and as low castes? Who built all the big, big temples? Has any Brahmin given any land to any temple? For his comfort, temples and agraharams (exclusive Brahmin residential colonies) are built—the Sudras are ones who gave you (Brahmin) a comfortable life. Why does no one question this? A Brahmin has written that if a Brahmin ploughs it is a sin, so without working in the land he has taken a means of living. Is it virtuous to eat the yield? When one enjoys the fruits of another's labour they say this is written by god, it is in the Vedas.
[Viduthalai, 9-6-1972]

**Communal representation**

The government should forever remove the words: talent, quality, and qualification from its dictionary. Because communal representation was against the Brahmins, they removed it in 1950. They have modified it stating that in education, employment and specific professions; the backward classes may be given some concessions. This is written in the Constitution.
[Viduthalai, 18-7-1972]

**Objectives**

My only wish is that all people should be rational; caste should be abolished and the word Brahmin should not exist. This is my only principle. I joined Congress only for this.
[Viduthalai, 26-8-1972]

**Britain and India**

In Britain there is no Paraiyar, Sudra or Brahmin. But in India (the categorizations of) Brahmin, Sudra and Paraiyar exist in law, temple, lakes, even in behaviour, in marriage, wealth, work, education and god. What is the reason?
[Viduthalai, 8-10-1972]

**Caste and Class**

Wealth in society is trouble, lacks peace, can bring mental problems and it is changeable at any time. But upper caste nature in a society is extremely horrible and is a biggest crime that never changes with time. It hinders and stops development, humanity and equality. It is a criminal offence. So one can pay any cost to put an end to it.

[Unmai 14-10-1972]

### Dr.Ambedkar

In 1925 the Self Respect Movement came into being. It is not easy in this country to fight against religion and differences. Dr. Ambedkar and I were not only friends for a very long period but on many issues we both had same opinion. There is no one in India who is equivalent in status with Dr. Ambedkar.

[Viduthalai, 16-11-1972]

### Rajaji

Rajaji is the one who closed 2200 of the 6000 schools in a nation in which only 7% are literates. He was the one who said that students can read either in the morning or evening and do the profession of their father in the daytime or evening during off-school hours. Only because of this, I sent out Rajaji from politics.

He made a pact that he would not receive more than Rs.500/- per month but he took a salary of Rs.20,000/- p.m and he was the owner of several lakhs worth of property. Has he led a honest life in politics or in social life? What right has he to advice others?

[Viduthalai, 16-11-1972]

### Indian democracy

Today's independent India is ruled by 3% Brahmins who overpower the 16% Panchamas, 75% backward class people and other 7% people.

Ninety-nine percent of the social problems such as feelings of superiority, theft, worry, anger, distortion of facts, loss, cheating, cunning, rape, etc. are at their height. Instead of democracy, if there is a socialist rule then all these problems would not prevail to such a great extent. Selflessness would grow and selfishness would be reduced. The most fraudulent word in the world is democracy.

[Viduthalai, 30-4-1973]

**Rationalism**

The greatest act more important than my life is obtaining freedom for Sudras from Baniyas and Brahmins and even obtaining a separate state to live in real independence. Unlike the whites and blacks, we do not have so much of differences. And yet they lead a life of equality. On the other hand, the Brahmin claims himself to be upper caste and we are treated as lower caste. If our nation were under the rule of Muslims or Christians certainly it would have become free of these discriminations.

Dr. Ambedkar was a stalwart in rationalism. He would state every revolutionary information in a simple way.

There is no job that is very low and I am sure that there will soon be machines that will perform all these low jobs.

To develop rationality among people, I organized several feasts in which we served beef (and pork). There is no rule or law that says that one should not eat beef or pork.

It is unfortunate to state that it has become so inherent in the Hindu religion that Hindus do not eat beef and only the Christians, Muslims and the Depressed Classes feed on beef.

I am an enemy of all religions and they should be abolished. If a man seeks respect in his public life then they can be no use or real service. My only mission is to see that the Dravidians lead a life of self-respect and I will serve for this.

Whether I have the capacity to do this is not the question. As no one is interested in doing this, I continue to do it. I am fit for that duty because of my nationalism and my temper.

Not only are we ill-treated as Panchamas and Sudras in our country by the outsiders (Brahmins), but the northerners are also exploiting us.

Only human beings are there in America. Is there a Brahmin, a Sudra or an Untouchable?

Around 200 people came when I beat god with *chappals*. The awareness among people has increased.

[Last Speech of Periyar, 19-12-1973]

## 3.3  Some facts about Untouchability and its Consequences even after 57 years of Independence

Here we give some recent information about the present status of untouchability. This data is given from Dalit Right to Reservation and Employment (Advocacy Internet), Vol. 5, Issue 05, September-October, 2003.

The state has the right to make any provisions for SCs/STs (Indian Constitution, Fundamental Rights, 15-4). After 56 years of independence what are the facts? What do these Facts say?

The fact is Reservation is not applicable in

1. Judiciary
2. Trustee posts (Example: Mumbai Port trust), which determine polices in autonomous bodies.
3. Private sector
4. Defence Department, the country's largest's employer: SC share here is negligible.
5. Certain Minority Institutions: In some State Institutions (Example: Banking Industry) reservation took effect only in1972.

(Source: Suresh Mane, "Reservation Policy in the Present Scenario and Economic Development of Backward Classes", 1999, p. 5 also of National SC/ST Commission Report 1996-1997 and 1997-1998, p. 20)

**Backlog Vacancies**

| Services | in No. | in % |
|---|---|---|
| A. Government Department | | |
| Group A | 369 | 74.84 |
| Group B | 438 | 51.34 |
| Group C | 3133 | 55.87 |
| Group D | 873 | 45.70 |
| Total | 4811 | 54.30 |
| *B. Banks* | 272 | 45.10 |
| *C. Public Sector* | 2642 | 88.18 |

Note: The number and the % refer to the backlog vacancies remaining unfilled despite the Special Recruitment Drive Policy in 1996-1997.
*(Source: National SC/ST commission Report, 1996-97 and 1997-98, page, 183-184)*

Now when even in the backlog the vacancies are not filled, it only shows their least concern to do justice to SC/ST.

**SC% in Central Government Services**

| | **1965** | **1995** |
|---|---|---|
| Class I | 1.64 | 10.12 |
| Class II | 2.82 | 12.67 |
| Class III | 8.88 | 16.15 |
| Class IV | 17.75 | 21.26 |
| Total | 1317 | 17.43 |
| Sweepers | * | 44.34 |
| Grand Total | 13.17 | 18.71 |

* Figures relating to sweepers in 1965 are not available and included in the figures for class IV.

(Source: National SC/ST Commission Report, 1996-1997 and 1997-1998, p. 14)

- The lower the type, social status and salary of the class of services the greater, the number of dalits, reflecting thereby the rigidity of the hierarchically structured caste system.
- Even after a long gap of 30 years between 1963 and 1995 the increase in the intake of Dalits for all the classes of services has been very minimal.
- One striking factor is that the quantum of increase of Dalit recruits for each class of services in the same year reflects the same hierarchical pattern of the caste system and class of services.
- That 44.34% of sweepers are from the SC community is shocking! In comparison with the SC population (16.33%) in 1991 the percentage of sweepers (44.34) is almost three times as high!

## DALIT RIGHT TO EDUCATION
### National and International Standards of Measurement

- Make effective provision for securing the right to education (Indian Constitution: Directive Principles - 41).
- Provide within 10 years (1950-1959) free and compulsory education for all children up to 14 years of age (Indian Constitution: Directive Principles -45)
- Promote with special care the educational interests of SCs (Indian Constitution: Directive Principles - 46)
- Every one has the right to free and compulsory education at least in the elementary and fundamental stages. Secondary (including technical and vocational) professional and higher education shall be made generally available and equally accessible to all (UDHR-Art 26 ICESCR Art 13)

UDHR – Universal Declaration of Human Rights.
ICESCR – International Convention on Economics Socio Cultural Right.

## After 52 years of independence
## Primary Middle and Secondary Education

- In 1993 enrolment at primary level among SCs was 16.2% while among non-SCs it was 83.8% *(Annual Report 1994-1995, HRD, GOL).*
- The national dropout rate among Dalit children is 49.35% at primary 67.77% at middle and 76.65% at secondary level. *(National SC/ST Commission Report 1996-1997 and 1997-1998 p. 47)*
- While privatized education is becoming the order of the day, 99% of Dalit students come from Government schools which lack basic infrastructure adequate class rooms and teachers, teaching aid etc. *(Fraser, Sudha Sowbhagyavathy Growth of Literacy among SCs in AP (Paper) March 1999).*

### All India literacy level gaps between SC and rest of the population

| Year | Gap% |
|------|------|
| 1961 | 17.59 |
| 1971 | 19.13 |
| 1981 | 19.84 |
| 1991 | 17.20 |

*(Source: Eight Five year plan 1992-1997, Vol. 2, GOL, p, 420 and Census of India 1991, Vol. 2, page 419)*

### State-wise SC dropout rate (50% and above) in 1990-1991)

*Primary stage***:** 13 States/UTs, Andhra Pradesh, Assam, Bihar, Goa, Karnataka, Manipur, Orissa, Rajasthan Sikkim, Tripura, West Bengal, Dadra and Nagar, Haveli, Delhi.
*Middle stage***:** 19 States /UTs, Andhra Pradesh, Assam, Bihar, Goa, Gujarat, Haryana, Karnataka, Madhya Pradesh, Maharashtra, Manipur, Meghalaya, Orissa, Punjab,

Rajasthan, Sikkim, Tripura, Uttar Pradesh, West Bengal, Chandigarh.

*Secondary stage*: 23 States /UTs, Andhra Pradesh, Assam, Bihar, Goa, Gujarat Haryana, Himachal Pradesh, Jammu and Kashmir, Karnataka, Kerala, Madhya Pradesh, Maharashtra, Manipur, Orissa, Punjab, Rajasthan, Sikkim, Tamilnadu, Tripura Utter Pradesh, West Bengal, Delhi, Pondichery.

That is 11 states/UTs the drop out rate of SCs in all three stages is 50% and above is appalling and an indication and of the given attention given to the education of SCs.

**State wise % of total SC student enrolment in higher education as on 1995**

| Name of the State | % |
|---|---|
| Maharashtra | 21.74 |
| Utter Pradesh | 13.92 |
| Tamil Nadu | 9.04 |
| Andhra Pradesh | 7.68 |
| Karnataka | 7.63 |
| Gujarat | 6.80 |
| West Bengal | 5.06 |

*(Source: National SC/ST Commission Report 1996-1997 and 1997-1998 p. 77)*

When that poor percentage enrolls for higher education they are denied seats in medical colleges, engineering course and for certain bachelor degrees. They are bluntly denied entry.

**Dalit Right to life and Security**
**National and International Standards of measurement**

- Do not subject any person to any disability, liability, restriction or condition on grounds only of religion, race, caste… with regard to
    - Access to shops, public restaurants, hotels and places of public entertainment.
    - The use of wells, tanks, bathing ghats, roads and places of public resort (ICFR, 15-2)
- Untouchability is abolished, its practice in any form is forbidden. Violators will be punished in accordance with law. (ICFR, 17)
- Do not deprive life and liberty of any person (ICFR, 21)
- Everyone has the right of life, liberty and security of person. (UDHR, Art 3)
- All are equal before the law, and are entitled to equal protection of the law without and against any discrimination (UDHR, Art7)
- State parties condemn racial discrimination and undertake to eliminate it in all its forms without delay (ICERD, Art1)
- Other provisions in the Indian constitution (FR14, 15-1,25), in the Indian legislation (the protection of civil rights Act, 1955). The SC/ST (POA) ACT, 1989 and rules (1995) and in the international covenants (UDHR, Art 2)

IC = Indian Constitution,
FR = Fundamental Rights
ICERD = International convention on Elimination of Racial Discrimination
UDHR = Universal Declaration of Human Rights.

## After 52 years of India's Independence
## What are the facts? What do these facts say?

Untouchability: A challenge to Indian Nationalism and Patriotism. Continuing forms of untouchability today
- Prohibition to sit on par with dominant castes in public or private places.

- Prohibition to take out marriage procession in dominant caste locality and Dalit bridegroom prohibited from riding a horse during marriage procession.
- Prohibition to walk with footwear on roads/pathways of dominant caste villages.
- Dalit corpse prohibited from being carried through dominant caste villages, or buried in the latter's graveyards.
- Tea /coffee served in separate metal/earthenware containers in village hotels.
- Washing of such containers after use by Dalits themselves, while those used by non-Dalits are washed by someone else employed for the purpose.
- Prohibited to draw water from public village ponds/tanks/taps/wells
- Rendering forced menial services in birth, marriage and death ceremonies of dominant castes.
- Insults to and degradation of Dalit women.
- Beatings, torture, attempt to murder and actual murder by dominant castes.
- Preventing exercise of franchise at elections harassment, threats, actual murder of elected representatives, threats against, even forcible prevention from contesting democratic elections.
- Forcibly preventing occupational change /mobility.
- Rejection of demand for just wages.

## Manual Scavengers
### (the Night-soil workers) of India
### Tales of sorrow

**Rajesh**, a 25-year-old dropout has been carrying a tin of night soil in his hands every morning, every day of the year, for the past 10 years. He cleans the "vamda" (a four-walled open latrine) in his village, Rampura in Surendranagar District, carries the scooped up night soil on his head and deposits it in a place far away from the village. He earns a sum of Rs 50 (Per month) for his effort. That has been the rate for the past four years; before that the Gram Panchayat use to give us Rs25.
(Indian Express, September 22, 1998).

**Narayanamma**, 55 years, works as a safai Karmacharis (manual scavenger) in Ananthapur Municipality, Andhra Pradesh. She has been scavenging continuously for 19 years. She goes to the community dry toilet, which is one KM away from her house. She takes a bamboo basket and two small metal pieces, which she had left at the corner of the toilet…She has to clean 400 seats of dry toilet every day, having a load of 15-16 bamboo baskets of human excreta. Her health is ruined .she suffers from diarrobea and vomiting very frequently.

(R.I. Pillai former secretary General, National Human Rights commission (NHRC), Seminar paper manual scavenging a challenge to human rights 1998).

### How many manual scavengers does India have?
Conflicting Estimates and Confusing Figures

1. Says the 1995-96 Annual Report of the national commission for manual scavengers the figure for 1989 is 4.21 lakhs.
2. The Report, however expresses the commission's apprehension that the actual number of dry latrines is much more than estimated.
3. Says the 1997-98 annual report of the ministry of social justice and empowerment 8,25,572 is the estimated number of scavengers in India.

### The employment of manual scavengers and construction of Dry Latrines (Prohibition) Act, 1993
### State-wise status of implementation

*Central Act applicable (with effect from 26.01.97)*
Andhra Pradesh, Goa, Karnataka, Maharashtra, Tripura, West Bengal, All union Territories

*Central Act adopted:*
Assam, Bihar, Gujarat, Haryana, Madhya Pradesh, Orissa, Punjab

*Central Act adopted under consideration:*
Rajasthan, Tamil Nadu, Uttar Pradesh.

*No decision taken as yet:*
Arunachal Pradesh, Himachal Pradesh, Jammu and Kashmir, Kerala, Meghalaya, Manipur, Mizoram, Nagaland, Sikkim

(Human Rights Newsletter, NHRC, V5, No.11, November 1998)

1. With independence in 1947 our government system assured every citizen the protection and promotion of her/his rights. But the right of the safai karmacharis (manual scavengers) came to be officially recognized only in 1993 that is 46 years after independence! With the enactment of the employment of the manual scavengers Act. But then the act came into force only in 1997- after another gap of 4 years even then only the union territories and six states made the central Act applicable in their territories even other states have adopted the act only in principle for three others the adoption is under consideration and the rest have not taken any decision so far.

2. How to understand this state of affairs? Do not the Rights of the safai karmacharis merit any serious attention from our governance system? Have they also become untouchable to our democratic polity? Does the ruling class consider them only as dispensable and disposable entities, precisely because they are not politically organized to clamour for their rights?

3. Why has the civil society not woken up to the denial of even the basis rights of the safai karmacharis? Is it because the safai karmacharis are untouchables dealing with the "untouchable realities" of life that their dignity does not deserve any attention? Is it because they are 'useful' to do things which the non-safai karmacharis would not want to do? Is it because they are on the lowest rung of the caste ladder- so low that their rights are invisible to others, or so low that are an affront to non-dalits' rights? Have even the enlighten sections of the civil society become hapless victims of the caste system's brutal suppression of rights?

4. The last link of the colonial chain was snapped in 1947.But this was colonialism from the outside. Should not true nationalism and genuine patriotism give priority attention to the internal colonization of the safai karmacharis, through the more insidious and powerful weapon of the caste system and its practice of untouchability? Are the dominant castes willing to exhibit this brand of nationalism than any other, practiced and propagated hitherto for self-serving political ends?

5. When will Dr.Ambedkar's dream of attaining the right of social democracy dawn in the lives of the safai karmacharis? More than 50 years have passed since he sketched this dream for the dalits and for the nation. It still remains a dream. Will it remain so even in the future?

### Atrocities against Dalits Disabled Governance?

The list below refers only to those cases registered with the police. Either due to intimidation and fear, or to inaccessibility of police stations, or to loss of faith in the law enforcement machinery, a number of cases go unreported.

*Cases registered with the police under different nature of crimes and atrocities against the Dalits*

| Nature of crime | 1995 | 1996 | 1997 | Total |
|---|---|---|---|---|
| Murder | 571 | 543 | 503 | 1617 |
| Hurt | 4544 | 4585 | 3462 | 12591 |
| Rape | 873 | 949 | 1002 | 2824 |
| Kidnapping& Abduction | 276 | 281 | 242 | 799 |
| Dacoity | 70 | 90 | 57 | 217 |
| Robbery | 218 | 213 | 157 | 588 |
| Arson | 500 | 464 | 384 | 1348 |
| PCR Act 1955 | 1528 | 1417 | 1157 | 4102 |
| SC&ST(POA)ACT,1989 | 13925 | 9620 | 7831 | 31376 |
| Other Offences | 10492 | 13278 | 11693 | 35463 |
| Total | 32997 | 31440 | 26488 | 90925 |

*(Source: National SC/ST Commission Report, 1996-97 & 1997-98, Page 240.)*

1. Dalits constitute only 16.33% of the total population. Among them 62.59% are illiterates and about 50% each in rural and urban areas live below the poverty line. This being the Vulnerable position of Dalits in terms of small number, illiteracy and denial of livelihood rights, why should they be subjected to continuous perpetration of crimes and atrocities-about 3.8 crimes & atrocities per hour in 1995,3.6 in 1996 and 3 in 1997?
2. Despite the increase in the strength of the police force annually and the rapid modernization of the law enforcement machinery, why has the Indian state become so powerless to contain effectively these crimes and atrocities punish and control the perpetrators?
3. What system of governance do we have the rule of law of democracy or the reign of ''cast (e) curacy''?

## MASS MURDERS OF DALITS IN BIHAR

### Mindless Massacre in Sawaranbigha

Seven Dalits were killed in 1991 by activists of the Savarna Liberation Army (SLA), because the Dalits staked claim to 2.4 hectares of government land and cultivated it. (Front line, March 12, 1999, P.29).

### Midnight massacre in Laxmanpur-Bathe

On the night of December 01, 1997, Bhumihar, who wanted to seize 51 acres of land that was allocated to Dalits, entered the houses, shot indiscriminately, raided 14 homes, killed 67 people, injured an additional 20 people, and murdered 7 local fishermen and brutally raped and murdered 5 girls. The girls were shot in the chest and vagina (Smita Narula, "broken people-caste violence against India's untouchables", Human rights watch, 1998, pp.61)

### Cold-blooded murders in Shankarbigha

Twenty-three scheduled caste landless agricultural workers were murdered in cold blood on 25 January 1999 by the private army of upper caste bhumihar landlords. Five women, and seven children, including a 10-month-old, were among those killed in order to terrorize the residents who were getting attracted to the ideology of two prominent Naxalite groups- the communist party of India (Marxist-Leninist) party unity and to seek to establish the supremacy of landlords (Frontline, February 26,1999 p.37)

### Heartless carnage in Narayanpur

On Feb.10, 1999 barely a fortnight after the republic day-eve massacre of 23 Dalits at Shankarbigha village the private army of upper caste Bhumihar landlords struck again killing 12 Dalits at Narayanpur village in the district (Frontline, March 12.1999, p.29)

### The bloody rivers of Bihar

Much blood has been split in Bihar in caste violence over the past three decades. Between the first reported caste based massacre, at Rupaspur Chandwa in Purnca district in 1971, and the latest bloodbath, at Narayanpur village in Jehanabad district on February 10, 1999, there were 59 recorded instances of mass murders in which about 600 people were killed. The majority of these massacres was directed at dalits and was carried out by the private armies of the upper castes, such as the Ranvir Sena, the Bhoomi Sena, the Brahmarshi Sena, the Sunlight Sena and the Savarna Liberation Army.

The period between 1990 and 1999 witnessed 35 instances of caste-based massacres, the total number of victims being about 400. More than 350 of those killed were from among the lower castes. (Frontline, March 12, 1999, p.30)

### The Killing Fields of India

In the recent history of the central Bihar districts, killings have occurred with frightening regularity. A few such ghastly massacres

| | |
|---|---|
| Rupaspur Chandwa | 1971 |
| Arwal&Kansara | 1986 |
| Golakpur | 1987 |
| Lalibigha | 1988 |
| Lakhawar | 1990 |
| Sawanbigha | 1992 |
| Aiara | 1994 |
| Khadasin | 1997 |
| Lakshamanpur-bathe | 1997 |
| Chou ram and Rampur | 1998 |
| Shankarbigha | 1999 |
| Narayanpur | 1999 |

*(Source: Frontline, Feb 26,1999)*

## Is there a future for Dalits in Bihar?

1. Sociologists have pointed out that resort to measures that merely address violence, as a law and –order problem will not be enough to smash these Senas.
2. Such steps have to be coupled with bold and far reaching measures such as land reforms which address the fundamental problem of economic exploitation and social discrimination of landless agricultural labourers from among the scheduled castes by upper-caste feudal landowners

(Source: Frontline, March 12, 1999, p. 31)

*State-wise incidence of crimes against Dalits*

| STATES | Incidence of total cognizable crimes | % contribution to all-India total |
| --- | --- | --- |
| Uttar Pradesh | 10963 | 34.9 |
| Rajasthan | 6623 | 21.1 |
| Madhya Pradesh | 4075 | 13.0 |
| Tamil Nadu | 1812 | 5.8 |
| Gujarat | 1764 | 5.6 |
| Andhra Pradesh | 1629 | 5.2 |
| Maharashtra | 1352 | 4.3 |
| Karnataka | 1089 | 3.5 |
| Bihar | 810 | 2.6 |
| Kerala | 640 | 2.0 |
| Orissa | 486 | 1.5 |

(This list contains only states having incidence above 100. Source: National SC/ST Commission Report, 1996-97 & 1997-98, p. 309)

**The Ranvir Sena in Bihar: Symbol of Jungle law?**

The Narayanpur killings is the 19$^{th}$ massacre perpetrated by the Sena since it was founded in August 1994 by Brahmeshwar Singh and Dharicharan Choudhry prosperous landlords of Belier village in bhojpur district. So far the Sena has killed 277 persons, almost all of them poor, landless and oppressed Dalits.

1. Naxalite groups in recent years have mobilized agricultural workers against their social persecution and economic exploitation by the landlords. The Sena was formed to protect the economic and feudal interests of the upper caste landlords and has systematically targeted Naxalite and their suspected sympathizers.
2. The Sena has 300 well-trained Bhumihar youth as its members and has sophisticated arms in its possession. It has insured the lives of its activists and provides those monthly allowances and other benefits.
3. It depends on the Bhumihar community for financial support
4. That political patronage, cutting across party lines, is available to the Sena is evident from the fact that despite an official ban no major crackdown has been launched against it.
5. No important member of the Sena has been tried in court till date. Brahmeshwar Singh who owns 97 bighas of land was arrested on two occasions earlier but was released.
6. The Sena leaders boast that at least 125 of the killings were carried out after July 1995 when the Bihar government banned the group .As the communist party of India (Marxist-Leninist) liberation has pointed out repeatedly, the ban existed only on paper

(Source: Frontline, February 26 and March 12, 1999)

## Dalit right To Livelihood
## National and International standards of Measurement

1. Ensure the right to adequate means of livelihood (IC DP,39-a)

2. Provide, within 10 years (1950-59) free and compulsory education for all children up to 14 years of age (ICDP, 45).
3. Promote with special care the educational and economic interests of SCs/STs.(ICDP,46)
4. Everyone has the right to a standard of living adequate for the health and well-being of herself/himself and of her/him family, including food, clothing, housing and medical care and necessary social services. (UDHR art 25)
5. Everyone has the right to free and compulsory education, at least in the elementary and fundamental stages (UDHR Art 26, ICESCR Art 13)
6. Other provisions in the Indian constitution (DP, 39-c, 39-f, 41, 47) and in the international conventions (UDHR, Art 25,ICESCR Art 12, ICERD, Art 5)

IC = Indian Constitution, DP = Directive Principles
ICERD = International convention on Elimination of Racial Discrimination
UDHR = Universal Declaration of Human Rights.
ICESCR=International Convention on Economic, Socio-cultural Rights

### After 52 Years of India's Independence: what are the facts? Basic Amenities: Ever Attainable?

Even after 52 years of independence, a majority of Dalits have been deprived of electricity & sanitation facilities. And compared to non-Dalits their plight is indeed shocking!

*Household in % having amenities*

|  | Electricity | Sanitation |
| --- | --- | --- |
| SCs | 30.91 | 9.34 |
| Non-SC's | 61.31 | 26.76 |

*(Source: "Demographic and Occupational Charactistics in India, ''M.Thangaraj (working paper 134 of MIDS), p. 6)*

In many of the rural areas, the scheduled caste continue to have separate wells or source of drinking water, the quality of which is considerably poorer than the general source of drinking water in the village.

*Poverty Line a Blot on Indian Nationalism*

*SCs and all population in % below poverty line*

| Year | Total Population | SCs |
|---|---|---|
| 1977-78 | 51.2 | 64.6 |
| 1983-84 | 40.4 | 53.1 |
| 1987-88 | 33.4 | 44.7 |

*(Source: Eighth Five year plan, GOI, 1992-97, Vol. 2, p. 420)*

1. Despite overall declining rates, incidence of poverty among SCs remains comparatively higher than in the total population
2. That about 50% of the Dalits in rural as well as urban areas live below the poverty line is an appalling situation!
3. That the number of Dalits below the poverty line is significantly greater than that of non Dalits is more shocking

**% of SCs and all population below poverty line**

| | SCs | | Total Population | |
|---|---|---|---|---|
| Year | Urban | Rural | Urban | Rural |
| 1993-94 | 49.48 | 48.11 | 37.27 | 32.36 |

*(Source: National SC/ST Commission, 96-97 & 97-98, p.82)*

1. That the % of Dalits below the poverty line in 1993-94(49.48% in Urban and 48.11% in rural areas) is an increase over that of 1987-88(44.7%) is very alarming!

2. That this sudden increase is in sharp contrast to the declining rate in the preceding years is quite significant.
3. That this increase has taken place during the privatization liberalization- globalization regime indicates how vulnerable dalits could become to the market forces.

### Basic Literacy A Receding Dream?
### Basic literacy rate of dalits (%)

| All India (%) | | Urban (%) | | Rural (%) | |
|---|---|---|---|---|---|
| 37.41 | | 55.11 | | 33.25 | |
| Male | Female | Male | Female | Male | Female |
| 49.91 | 23.75 | 66.60 | 42.29 | 45.94 | 19.46 |

*(Source: 1991 Census)*

1. Only 1/3 of the SC population is literate.
2. Female SC literacy falls short by more than half of male SC literacy
3. In both the rural and urban areas female SC literacy is far behind the male SC literacy rate
4. So also the rural SC literacy rate lags far behind the urban SC literacy rate.

### Government-Recognized Atrocities against the Dalits

What constitutes an offence under Indian law are certain acts and practices of non Dalits treading Dalits as unequal which contravene the fundamental right under the Indian constitution that all citizens are equal having equal opportunity and legal protection under the law. The following acts and practices are considered as acts of "Atrocity" under Chapter 2, No 3 of the "Atrocities on Scheduled Castes and Scheduled Tribes (prevention) act, 1988". What is important to note is that the acts mentioned here continue to take place every day in villages and cities throughout India?

*Punishment for offences of atrocities-*

**Whoever, not being a member of a Scheduled Caste or a Scheduled Tribe, -**

1. Forces a member of a Scheduled Caste or a Scheduled Tribe to drink or eat any inedible or obnoxious substance

2. Acts with intent to cause injury, insult or annoyance to any member of a Scheduled Caste or a Scheduled Tribe by dumping excreta, waste matter, carcasses or any other obnoxious substance in his premises or neighborhood

3. Forcibly removes clothes from the person of a member of a Scheduled Caste or a Scheduled Tribe or parades him naked or with painted face or body or commits any similar act which is derogatory to human dignity.

4. Wrongfully occupies or cultivates any land owned by or allotted to or notified by any competent authority to be allotted to a, member of a Scheduled Caste or a Scheduled Tribe or gets the land allotted to him transferred.

5. Wrongfully dispossess a member of a Scheduled Caste or a Scheduled Tribe from his land or premises or interferes with the enjoyment of his rights over any land, premises or water.

6. Compels or entices a member of a Scheduled Caste or a Scheduled Tribe to do begging or other similar forms of forces or bonded labor other than any compulsory service for public purposes imposed by government.

7. Forces or intimidates a member of a Scheduled Caste or a Scheduled Tribe not to vote or to vote to a particular

candidate or to vote in a manner other than that provided by law.

8. Institutes false, malicious or veracious suit or criminal or other legal proceedings against a member of a Scheduled Caste or a Scheduled Tribe.

9. Gives any false or frivolous information to any public servant and thereby causes such public servant to use his lawful power to the injury or annoyance of a member of a Scheduled Caste or a Scheduled Tribe.

10. Intentionally insults or intimidates with intent to humility a member of a scheduled caste or scheduled tribe in any place within public view.

11. Assaults or uses force to any woman belonging to a Scheduled Caste or Scheduled Tribe with intent to dishonor or outrage her modesty.

12. Being in a position to dominate the will of a woman belonging to a Scheduled Caste or a Scheduled Tribe and uses that position to exploit her sexually to which she would not have otherwise agreed.

13. Corrupts or fouls the water of any spring, reservoir or any other source ordinarily used by members of the Scheduled Castes or Scheduled Tribes so as to render it less fit for the purpose for which it is ordinarily used.

14. Denies a member of a Scheduled Caste or a Scheduled Tribe any customary right of passage to a place of public resort or obstructs such member so as to prevent him from using or having access to a place of public resort to which other members of public or any section thereof have a right to use or access to.

15. Forces or causes a member of a Scheduled Caste or a Scheduled Tribe to leave his house, village or other place of residence

**Shall be punishable with imprisonment for a term which shall not be less than six months but which may extend to five years with fine.**

### Human Rights Violations of the Dalits

The core of a Dalit life in India means the following.

1. *Untouchable Status*: While in the urban areas, the form of the untouchability practice has changed in rural areas it continues to be openly practiced. Every Indian village has a segregated and isolated Dalit housing locality known as 'Vas'. A Dalit cannot dream of hiring building or purchasing a house amongst the populace of non Dalits. Even a large scale natural disaster like the earth quake in Gujarat in 2001 could not ensure that the new upcoming housing colonies built during the rehabilitation effort from public funds could be mixed Dalits and non-Dalit housing colonies. Two other basic parameters are used to decide an individual caste status inter dining and inter caste marriage. If both these become a reality between those considered untouchables and those who are not, than the traditional meaning of untouchables caste status does not hold ground. But in the case of Dalits even when a touchable girl may marry an untouchable boy, she retains her touchable caste status but her children carry the untouchable caste status of their father.
2. *Untouchability:* The touch of a Dalits whether by accident or on purpose brings impurity for the non-Dalit? In some places Non-dalits will immediately rush home and take a bath to purify their body and soul. In other places, it is sufficient for them to dip their thump even in water, even if it be dirty water, or have their wives or relatives sprinkle water on them before they enter their house. In July 1998 in the state of Uttar Pradesh, a Brahmin district judge had his entire court chamber washed with the holy water of the Ganges River because the previous occupant of that position (as district judge) happened to be a Dalit. It is common knowledge that after the former home minister of India Mr.Jagjivanram, also a Dalit, visited the famous Hindu temple in bananas; the priest had the entire temple washed upon his departure. It is not unusual to witness a non-Dalit lighting a match and touching the flame to his skins an act of purification after having come in contact with a Dalit. The touch of a Dalit, even to food, drinking water, a smoking pipe, mattress, towel, cooking pan, water pot, is considered defiling its non-Dalit user.

3. *Discrimination*: Even today most of the villages have restriction on using neither public roads. Dalits can neither take out their marriage procession nor carry their dead on the main village roads .Any act of a Dalit that can symbolize them as non Dalits is prohibited or becomes the cause of violence on Dalits. This includes playing a music band in the wedding ceremony. Tucking ones shirt in, wearing sunglasses, riding a horse, riding in the vehicle from the village main square, even having big moustache like the Kshatriya (Feudal).

4. *Forced impure occupations:* Disposal of a dead animal is thrust upon Dalits. So are all the occupations, which are considered to be unclean and filthy like manual scavenging. In India, most manual scavengers are employed by the state agencies in spite of a law that bans the same. The state under the rules of the village level social justice committees that is to be chaired by person belonging to a scheduled caste or a scheduled tribe has fixed the responsibility of ensuring disposal of the dead animals on the committee. This classification of caste-based occupations is rigid. In India all kinds of work which involve some handling of dirt are allotted to some caste or other among Harijan (the word used by Gandhi for the untouchable, it means children of god), Examples (1) Removal of carcasses and skinning them, (2) Tanning the hides, (3) Manufacture of leather goods, (4) Sweeping the streets and (5) Scavenging. These are all occupations, which are quite essential to the well being of society, but as they involve the physical handling of dirty or quickly putrefying matter, no non-Harijan, as a rule, will do any of the aforesaid jobs as a profession. That would be below his dignity.

*Government's Failure to meet CERD Obligations*

Pre and post independence, the state promised to remove untouchability and caste- based discrimination. It is true that there are adequate laws and legislations in place to redress the same. The question pertains to implementation — merely filing a complaint does not ensure justice. There is no legal system that can keep vigil on the social behaviors of the implementing authorities, each of which emerges with his or her social prejudices.

The Indian government has passed a special law to prevent commit atrocities on Dalits but over the years its implementation has been diluted by controversial judgments. The conviction rate in India in criminal cases in any case does not exceed four per cent. Dalit activists have come across several cases of atrocities in which the court has acquitted the accused merely on the grounds that the investigation of the case was carried out by an officer of lower rank than the one ought to have investigated!

The awareness on the parts of Dalits in fact has caused more problems for them. In a case where there was a social boycott of Dalits for nearly three years, the National

*Human Rights Commission made these observations*

"Due to education and marginal cultural development, when some youths, either assert their right to equal treatment or attempt to protect to dignity of their person or of their women, or resist the perpetration of the practice of untouchability or atrocities being committed on Dalits, They are often branded as "Naxalite" or "Extremists", they are implicated in false crimes and killed in false encounters. When they resist as a group, mass killings, arson of their hamlets, mass rape of their women and stripping them naked and parading them in the village are regular features." The lack of political-will is more apparent on the part of the state when it comes to the question of implementing economic programs, such as land reforms.

## A Profile of the Dalits of India

1. *Poverty*: People below the poverty line among SCs (49.48% urban areas and 48.11% in rural areas) is much higher than that of the average Indian population (37.27% in urban areas and 32.36% in rural areas).
2. *Income*: The income levels disaggregated by social group suggest that both the total household and the per-capita income levels are least for the SC followed by the ST (when compared with the national average). SCs have a total household income of Rs. 17,465, a mere 68% of the national average, and a per capita income of Rs. 3, 237, 72% of the national average.
3. *Wage Economy:* The share of income derived from wage labour (both agriculture and non- agriculture) is highest among SCs (at about 33% followed by approximately 20% among STs).
4. *Land Ownership:* STs reported a higher level land ownership(69%) with an average holding of 4.3 acres compared with SCs who owned the least land among the groups surveyed(only 47%), reporting an average holding of only 2.8 acres. The land ownership pattern seems to confirm the historical and domiciliary or residential patterns that affect specific caste groups in India.
5. *Housing:* The village development index is also associated with the percentage of kutcha houses( low cost house often made from a mixture of mud and tin) over 70% of landless labourers live in kutcha houses, as do a majority of both SCs and STs(74% and 67% respectively).
6. *Amenities:* STs and SCs are considerably more disadvantaged when compared with all other social groups in regards to ownership of and accessibility to amenities such as an electricity connection, piped water, and toilets.
7. *Dependency ratio and poverty:* The dependency ratios are very low among landless wage earners, STs and SCs. The decline in dependency may be attributed to the higher participation of females from these groups in employment and income earning activities. This apparent paradox supports the hypothesis that many among the marginal groups are at the risk of economic stress resulting in a higher level participation in the workforce as a coping mechanism. The evidence also suggests that low economic dependency among the low income and vulnerable population groups is the result of poverty.
8. *Participation in wage earning activities:* is higher among SCs and STs (58% and 55% respectively). The SCs are mainly landless (69.6%) with little control over resources such as land, forest and water. There has been a marked rise in the number of agricultural labourers (49.1%), casual labourers (72%), industrial labourers (17.3%), plantation labour (6.1%) and fishing labour (92.5%).
9. *Child labour:* Child labour exists in 58.75 of the SCs communities.
10. *Social groups and literacy:* STs and SCs recorded a literacy level of about 40% in comparison to the national Indian average of 54%. Gender disparity in terms of literacy is high among both SCs and STs. The school drop out rates are substantially higher among the lower income groups, landless wage earners, females, STs and SCs. STs and SCs have lower levels of literacy, especially at the level of matriculation and above for example only about 5% of girls among these communities complete matriculation.

11. *Undernourishment and infant mortality:* 57.5% SC children under the age of 4 were undernourished in 1992 while the infant mortality rate among SCs was 91 per 1000 live births.

To redistribute resources currently under the domain of upper castes. Even after 54 years of independence there is little progress on the question of land reforms. Furthermore, the planning commission of India has stated. The Programme of ceiling set out in the plan had been diluted in implementation. There were deficiencies in the law and delays in its enactment and implementation resulting in large scale evasions. Several states had made provisions for disregarding transfers made after a certain date, but often these provisions proved to be ineffective and not much surplus land has been available for redistribution. The main object of ceiling which is to redistribute land to the landless at a reasonable price on a planned basis has thus been largely defeated. The core questions, in context of discriminations faced by Dalits on the basis of descent are

1. Why are Dalits treated as lesser humans?
2. Why do not Dalits have the same de facto rights as other Indian citizens have?
3. Why is strict implementation of the laws viewed as only way to protect Dalits from being discriminated?
4. Why are affirmative action program seen as the only way to ensure Dalit representation in higher education jobs and political offices?
5. Why, in spite of constitutional provisions abolishing untouchability, do government institutions themselves fail to abolish discrimination in institutions set up and run out of public funds?

These questions are central to the more than 250 million people in south Asia whose life, dignity and rights are affected by caste based discrimination. In granting them the status of scheduled castes the Indian constitution has given special status to them as scheduled castes taking the sole criteria into consideration that they form the disabilities arising out of the caste system and for having been born into a particular caste community. The constitution recognized the caste system as an integral part of Hinduism and therefore, the Dalits who converted to another religion were not granted this same special status following their conversion. However, as late as 1990 the Indian government granted Dalits who had converted to Buddhism the same special status as Hindu Dalits, on the grounds that their conversion to a non Hindu religion did not change their status within the caste system and Indian society and the resultant discriminatory treatment meted out to them by the non-Dalits.

## Conclusion

Dalits have been the victims of caste-based discrimination based on descent. The major problem they face is of denial and the power to deny rests largely with those who discriminate. The progress in the country has been lopsided and vulnerable sections of Dalits continue to suffer inequalities.

What more remains to be proved when the official data suggest that 50 per cent of Dalit children can not complete primary education and in the case of the Dalit girls the drop out rate is 64%? What more does one need to prove on descent based discrimination when over one million people working as manual scavengers in the country belong only to the Dalits? The best proof of descent-based discrimination is that only Dalits are untouchables.

Chapter Four

# OBSERVATIONS

This chapter gives observations and conclusions made on the mathematical analysis of the social evil untouchability. Almost all the directed graphs given by the experts had many edges. They were in fact very dense graphs. Further most of the hidden patterns resulted in the ON state of all the nodes. Moreover, all the hidden patterns resulted only in fixed points. No hidden pattern in all these models was a limit cycle, which shows the permanent impact of the social evils rather than the recurrence after a time period.

It is very important to note that the attributes related with each and every type of untouchability is so strongly interlinked, it is evident from the fact that in many cases the ON state of a single node converts to ON state all the other nodes. This shows that many of the nodes are highly powerful and play a vital role.

It is pertinent to mention here that in most of the dynamical systems given by the experts (or modeled using the experts' opinions) all the nodes connected with any form of untouchability are so potent that the resultant of any state vector makes ON all the nodes of the resultant state vector. This factor is very unusual in the study of any another social problem.

As we had our own limitations, we did not ask the expert to modify or give lesser number of attributes or links or interrelations. In several cases these experts gave over 20 nodes. On request and after re-discussion we coupled a few of the nodes (which had the same shades of meaning). Thus at almost every stage the ON state of any node invariably made all the other nodes ON. This mathematically proves that the inter-relation between the set of nodes under investigation are densely interrelated and intertwined. As the number of nodes is many, the reader is advised to use C program given in [69. 76] for verification.

Another important observation was the experts accepted that certain nodes or relations couldn't be determined. Those concepts remained indeterminate, so most of the experts welcomed the idea of using indeterminacy. They felt at home to use this neutrosophic notion i.e., the concept of indeterminacy. We found in our analysis that only few of the relations or interrelations remained indeterminate. Hence we made the Neutrosophic Cognitive Map Model that followed the Fuzzy Cognitive Map Model. Our study confirms the reality of the situation of caste system and untouchability in India: they are extremely central to the problems faced by the people. We could conclude that any manifestation of untouchability was as strong as the practice of untouchability.

The heinous practice untouchability did not stop from the days of Manu and it is continuing even 57 years after Indian independence. In Tamil Nadu, it was Periyar's social revolution that led to a great deal of change. He was highly responsible for the development in the state of Dalits, Sudras and Women. For more than half a century, he voiced against the cruelties done in the name of the caste system.

**Observations about Periyar and untouchability practiced in religion**

It was Periyar who led the first historic and successful struggle at Vaikom for the rights of Backward Classes and Dalits to enter into temple streets. The success of this agitation led to the temples itself being thrown open to everybody.

He acutely understood that the source of untouchability was religion. That is why he countered a stigmatized caste identity with an identity of Self-Respect.

The Manu Smriti is a very inhuman law that differentiates people on the basis of their birth. This is why Periyar also advocated the concept of Samadharma: equality for one and all. It is clear that the concept of caste is actually ingrained in the minds of people because of religious brainwashing for centuries. There was really no basis for the caste system except the various epics, puranas, Vedas and other Hindu texts brought in by the Aryans. Therefore Periyar sought to annihilate the basis of casteism by annihilating Hinduism. He did not accept god, because all the cruelties against people were done in the name of god. His anti-theistic quote is renowned: "There is no god/ There is no god/ There is no god at all. He who created god was a fool/ He who propagates god is a scoundrel/ He who worships god is a barbarian."

Superstitions were highly dominant among the people in Tamil Nadu. Even for small acts like starting a new venture, travels, building a new house etc. people would seek the astrologers (who are invariably Brahmins) for auspicious time and date. By squandering money on rituals, rites, superstitions these non-Brahmins lost the little money they earned. Only Periyar's Self-Respect movement managed to eradicate the baseless superstitions prevailing in the society to a large extent.

Several experts observed the terrible impact of the caste system and Varnashrama Dharma that has infiltrated into all other major religions practiced in India. For instance in Christianity, there are Nadar Christians, Dalit Christians, Syrian Christians, Vellala Christians, Kallar Christians and so on. Untouchability is practiced against the Dalit Christians in some rural areas, who are allocated separate seating arrangements in churches, separate entry and exit doors and separate burial grounds. Categories like Dalit and Backward Class Muslims have also come into existence. This proves that caste system has seeped into every other religious structure in India. Our experts said that Periyar had mentioned about the existence of the caste system in Christianity in many speeches. It was with this motive that Periyar founded the Self Respect Movement for people to renounce their castes. Also, he advocated religious conversions to Islam to the Dalits as a means of liberation from caste.

We have concluded from our study that religious untouchability is the root cause of all other manifestations of untouchability. Even today people always fear to question anything that is done in the name of religion or god. It is because of their blind faith on god and religion is so great. Dalits and Sudras who were subjugated because of the caste system failed to question it because it was a matter of Hindu religion. This was the reason why Periyar first called for the annihilation of Hindu religion in order to annihilate caste.

Hindu religion attributed the caste system to the karma (deeds of the previous birth). The deeds of the previous birth decided the caste status in the next rebirth. Thus people were made to complacently accept the caste system because each of them were held personally responsible for the 'high' or 'low' status given to them in society. People could not challenge the social setup because it was religiously constructed. Periyar exposed the falsity of the notion of karma. He rationally pointed out that if indeed the cycle of birth and rebirth were in practice, it was entirely futile to offer prayers for the dead and perform rituals for dead ancestors through the Brahmins.

The Vedas, puranas, Sastras, Manu Smriti and epics like the Ramayana were all full of contradictions. Periyar held that they were written for the sole purpose of subjugating the Dravidian people and it made them into untouchables and Sudras, that is prostitute sons. He burnt the Manu Smriti and the Ramayana; and broke statues of Vinayaka. Periyar understood that it was impossible to improve the conditions of the Sudras without first annihilating untouchability. He said, "If you believe that without the label of Paraiyan being removed, the Sudra label will go away, then you are downright idiots."

Periyar condemned the use of Sanskrit as a language of worship. The Brahmins had claimed that it was a divine language, whereas Tamil was a 'low' language. He demanded the use of Tamil as a language of worship, because he wanted the Sudra and Dalit to understand what the officiating Brahmin priest prays on their behalf.

His last statewide agitation was to obtain the rights of priesthood for Dalits and Sudras. Periyar viewed it as an essential prerequisite to reorganize the social order on the basis of equality. It was more than twenty years after his death did the Supreme Court of India open priesthood to all the castes. Though the law officially recognizes this right, it is yet to be put into practice. Compared to Brahmin priests, the number of non-Brahmins who are priests at Hindu temples in negligible. Even now, the Dravidar Kazhagam organized demonstrations to enter the sanctum sanctorum of temples.

The neutrosophic model was used in our study and the node Samadharma was included as one of the attributes. When the only node rule of Manu Dharma was in the ON state, the hidden pattern of the dynamical system was a fixed point that showed all other nodes to be in the ON state. However the node Samadharma remained as an indeterminate node. Further when we worked with the ON state of the only node Samadharma we saw the resultant was a fixed point with all nodes in the OFF state and Samadharma alone was in the ON state. Thus it clearly proves none of these evils will prevail if Samadharma is present, a complete contrast with role or rule of Manu dharma.

**Observations about Periyar and untouchability practiced in education**

Dalits and Sudras were forbidden to learn anything or recite the sacred texts: this was the rule of the Manu Smriti. If they heard the Vedas being recited, molten lead was poured into their years.

If they dared to recite any of the Vedas, their tongue would be cut off. Only the coming of the British to India and the advent of missionary run schools provide basic education to people who had been unlettered because of the strictures of the caste system. Even today, denial of education continues to be the most sophisticated manifestation of untouchability.

We can conclude from our mathematical study that unless policy and decision makers on education are drawn from the Dalit and Sudra communities, the educational status of these downtrodden people cannot be improved. It must also be noted that although reservation in education is followed properly in Tamil Nadu state-run educational institutions, it is not implemented in the Central Government educational institutions.

Even though reservation might be followed, the need of the hour is the presence of Dalit and Backward class teachers in these institutions. Only such teachers will nurture the students from the downtrodden sections and give them adequate encouragement. A vast majority of students from the marginalized sections of society complain of the harassment that they are made to undergo in the hands of the teachers from the Brahmin community.

When non-Brahmin students at the Cheranamadevi Gurukulam where discriminated, he fought against it and succeeded in putting an end to it. Periyar vehemently opposed the introduction of hereditary education system. He rightly felt that it will perpetuate the occupation and birth-based caste system. The campaign he built was so powerful that the then Chief Minister Rajagopalachari was forced to resign his post and the scheme was withdrawn immediately.

He was also against the imposition of Hindi as a compulsory language of study in schools. He felt that learning Hindi would automatically lead to an imposition of Aryan values and domination, and would further oppress the Dravidians. Periyar vigorously fought against the monopoly of Brahmins in education, administration and employment, who occupied more than 90% of the coveted jobs.

He seldom spoke without quoting the exact statistics to prove his point. As an alternative, he campaigned for communal representation in employment and education to be put into practice, so that it allowed the Dalits and Sudras to achieve equality. In fact, it was his sustained and powerful crusade that brought about the first amendment to the Indian Constitution that allowed the State to make provisions for providing reservation to the Backward Classes and Dalits in education.

**Observations about Periyar and untouchability practiced in society**

India's development remains stagnated mainly because of the caste system. He said untouchability was crueler than all the cruelties of the world. Periyar, as a social revolutionary, perfectly understood the graveness and ugliness of the caste system. He believed that it was concocted by the Brahmins to establish their supremacy. As a rationalist, he fought for the equality of one and all. He felt that the minority Brahmins, who were just 3% of the population, was oppressing the majority of the people who were from the Backward and Dalit communities. In the Manu Smriti, it is described that the untouchable people had to wear the clothes of the dead. They could not own any property and they had to only use broken vessels. That is why, as a part of his struggle against Brahmin supremacy, Periyar burnt the Manu Smriti.

He led the historical struggle in Vaikom that obtained Dalits and Backward communities the right of entry into temple streets. He argued that if dogs and pigs could use the streets, why should human beings be denied their rights?

Similarly, at the Cheranmadevi school hostel, he intervened and condemned in strong words the discriminatory treatment meted out to non-Brahmins. Both these incidents proved remarkably his commitment to social justice. He was always against all kinds of social discrimination and fought against segregation including separate wells, separate burial grounds, etc. Because inter-dining and inter-marriage were forbidden between the castes, Periyar promoted common dining and inter-caste marriages.

He asked the Self Respecters to drop the caste tags from their names and to stop wearing religious symbols on their foreheads. Likewise, he was responsible for the Self Respect marriage that was devoid of any rituals, and didn't involve any priests or Brahmins and was merely a declaration between the man and the woman who wanted to be life-partners.

He launched agitations against the practice of serving Brahmins and 'others' separately in railway station restaurants and succeeded in putting an end to this practice. Also, he demonstrated against the use of the caste-name 'Brahmin' in the name boards of Brahmin hotels. He and his followers blackened such name boards and the Brahmins were forced to stop using their caste-name on hotel name boards.

Cadres of the Dravidar Kazhagam burnt sections of the Constitution that protected the caste system. All this proves that Periyar was vehement in eradicating the caste system and social inequalities prevalent in the Hindu social set-up.

**Observations about Periyar and untouchability practiced in politics**

Though Periyar entered political life by joining the Congress, in later days he condemned Gandhi who supported the Varnashrama Dharma. Periyar viewed the struggle against the caste system to be more important than the struggle against British imperialism. He said that the Indians, who practice untouchability don't have any right to condemn the practice of apartheid in South Africa. He attacked the Congress party that he felt existed in order to uphold Brahmin supremacy.

Since he felt that the oppression of the Dravidians (Sudras and Dalits) was because of Aryan (Brahmin) domination, he demanded the creation of a separate Dravida Nadu, and he put forth the slogan: "Tamil Nadu is for Tamils alone." In fact his separatist demand preceded even

Jinnah's demand for the creation of an independent Muslim nation, Pakistan.

Periyar declared, "I accept Dr. Ambedkar as my leader and all of us accept him as our leader." He supported Dr.Ambedkar's conversion to Buddhism. Periyar rued the failure of the Round Table Conference where the Depressed Classes demanded separate electorates.

He decried the absence of any Sudra member in the Constitution Drafting Committee. Three Brahmins, one Dalit, one Muslim and one Christian were members of the committee. He said that the Constitution, which was drafted without taking into consideration the Sudra community, could not be accepted.

Even now, political untouchability is rampant against the Dalits. Speaking about the state of Dalits in politics, Periyar said: How do we maintain the Untouchables? Even after 200 years of British rule, can an untouchable stand in an election against a caste-Hindu and win? Does any Indian think of making the Depressed Classes, who are one-third of the population, stand in an election, support and elect them? Have anyone of the leaders shown this way? So, the only state of mind in which everyone lives is, "I will not give my support for your equality and if you wish to attain equality, we will destroy it." Is this human?

In Tamil Nadu alone, ghastly and shocking anti-Dalit violence has taken place. In September 1996, the village of Melavalavu in Madurai district was declared a reserved constituency under Article 243D of the Indian Constitution. The Thevars, a higher caste of the area were incited because of this and resolved that they wouldn't allow any Dalit to occupy the posts in the panchayat and they threatened dire consequences. Murugesan bravely contested the election and won. After a few months, on 30 June 1997, the Thevars attacked and they beheaded him and six other Dalits. In a society which celebrated Dalit servitude, self-determination was not acceptable. Eight years after the Melavalavu massacre, in the reserved villages of Keeripatti, Paapapatti, Nattamangalam, the upper castes see to it that Dalits are not even allowed to contest the seats.

## Observations about Periyar and untouchability practiced in the form of economic oppression

It is mentioned in the Manu Smriti that "No collection of wealth must be made by a Sudra, even though he be able (to do it); for a Sudra who has acquired wealth, gives pain to Brahmins. (10:129)" This clearly proves that Sudras and Dalits were denied even the rights to own land or property.

The landlessness of the Dalit communities forced them to be economically dependent on the caste-Hindus for their survival. They continue to be agricultural labourers who are daily wagers who work in the fields of landlords who invariably are caste-Hindus. That is why Periyar had asked for land to be distributed among the Dalits.

Due to caste atrocities, and the discrimination they suffer in the villages, Dalits migrate to cities and other towns in search of jobs. This makes them more vulnerable to caste-Hindus. They inhabit city slums because of their poverty.

Periyar decried the practice of manual scavenging. He said that it was very unjust that people were given these menial jobs on account of their birth in a particular caste. He was well aware that economic improvement of the Dalits and Sudras was integral to their empowerment.
The entire project of division of labour among the various caste groups was intrinsic to the protection of caste system. This was only meant to economically oppress them.

Thus we see all these manifestations of untouchability mentioned in chapter two of this book are closely inter-related. When we studied the node 'religious untouchability' alone in the ON state and all other nodes were in the OFF state we found that the resulting hidden pattern was a fixed point. When we get a hidden pattern to be fixed point it implies that there is no periodical recurrence of the situation, and the resultant remains a standardized value.

# FURTHER READING

1. **Ambedkar B.R.,** *Writings and Speeches,* Vols. 1 to 16, Government of Maharashtra, Mumbai, (1993).

2. **Anaimuthu. V.,** *EVR Sinthanaikal (Thoughts of Periyar),* Vols. 1 to 3, Sinthanaiyalar Pathippagam, Trichy, (1974).

3. **Arignar Anna**, *My Leader Periyar*, Translated from Tamil by Venu, Periyar Self Respect Propaganda Institution, Madras, (1981).

4. **Arunachalam, M.**, *An Introduction to the History of Tamil Literature*, Gandhi Vidyalayam, (1974).

5. **Balasubramaniam, K.M.**, *Periyar E.V.Ramasamy*, Periyar Self Respect Propaganda Institution, Trichy, (1973).

6. **Balu, M.S.,** *Application of Fuzzy Theory to Indian Politics*, Masters Dissertation, Guide: Dr. W. B. Vasantha Kandasamy, Department of Mathematics, Indian Institute of Technology, (2001).

7. **Bechtel, J.H.,** *An Innovative Knowledge Based System using Fuzzy Cognitive Maps for Command and Control*, Storming Media, (1997). http://www.stormingmedia.us/cgi-bin/32/3271/A327183.php

8. **Brannback, M., L. Alback, T. Finne and R. Rantanen,** Cognitive Maps: An Attempt to Trace Mind and Attention in Decision Making, *in* C. Carlsson ed. *Cognitive Maps and Strategic Thinking,* Meddelanden

92. **Viduthalai Rajendran**, *Periyarai Kochaip Paduthum Kuzhappavadhikal* (The Confused who Vulgarise Periyar), Dravidar Kazhagam Publications, Madras, (1995).

93. **Viswanathan, E. Sa.**, *The Political Career of E. V. Ramasamy Naicker—A Study in the Politics of Tamil Nadu 1920-1949*, A thesis submitted in partial fulfillment of requirements for the degree of the doctor of Philosophy in the Australian National University, Canberra, (1973).

94. **Vysoký, P.** *Fuzzy Cognitive Maps and their Applications in Medical Diagnostics*. http://www.cbmi.cvut.cz/lab/publikace/30/Vys98_11.doc

95. **William Taylor**, *Oriental Historical Manuscripts in The Tamil Language*, Chennai, (1835).

96. http://www.countercurrents.org/dalit-periyar280603.htm

97. http://www.periyar.org

# INDEX

# About the Authors

**Dr.W.B.Vasantha Kandasamy** is an Associate Professor in the Department of Mathematics, Indian Institute of Technology Madras, Chennai, where she lives with her husband Dr.K.Kandasamy and daughters Meena and Kama. Her current interests include Smarandache algebraic structures, fuzzy theory, coding/ communication theory. In the past decade she has guided nine Ph.D. scholars in the different fields of non-associative algebras, algebraic coding theory, transportation theory, fuzzy groups, and applications of fuzzy theory of the problems faced in chemical industries and cement industries. Currently, six Ph.D. scholars are working under her guidance. She has to her credit 612 research papers of which 209 are individually authored. Apart from this, she and her students have presented around 329 papers in national and international conferences. She teaches both undergraduate and post-graduate students and has guided over 45 M.Sc. and M.Tech. projects. She has worked in collaboration projects with the Indian Space Research Organization and with the Tamil Nadu State AIDS Control Society. This is her 21$^{st}$ book.

She can be contacted at vasantha@iitm.ac.in
You can visit her work on the web at: http://mat.iitm.ac.in/~wbv

**Dr.Florentin Smarandache** is an Associate Professor of Mathematics at the University of New Mexico in USA. He published over 75 books and 100 articles and notes in mathematics, physics, philosophy, psychology, literature, rebus. In mathematics his research is in number theory, non-Euclidean geometry, synthetic geometry, algebraic structures, statistics, neutrosophic logic and set (generalizations of fuzzy logic and set respectively), neutrosophic probability (generalization of classical and imprecise probability). Also, small contributions to nuclear and particle physics, information fusion, neutrosophy (a generalization of dialectics), law of sensations and stimuli, etc.).

He can be contacted at smarand@unm.edu

**Dr. K. Kandasamy** is a staunch Periyarist. He worked as a guest professor in the Tamil Department of the University of Madras till recently. He holds postgraduate M.A. degrees in Tamil Literature, Political Science, Saiva Siddhanta, Defence Studies, Education and History.

He can be contacted at dr.k.kandasamy@gmail.com